# Scaling limit of a weakly asymmetric simple exclusion process in the framework of regularity structures

May 1, 2025

**Ruojun Huang**[1], **Konstantin Matetski**[2] **and Hendrik Weber**[3]

[1] University of Münster, Email: ruojun.huang@uni-muenster.de
[2] Michigan State University, Email: matetski@msu.edu
[3] University of Münster, Email: hendrik.weber@uni-muenster.de

### Abstract

We prove that a parabolically rescaled and suitably renormalised height function of a weakly asymmetric simple exclusion process on a circle converges to the Cole-Hopf solution of the KPZ equation. This is an analogue of the celebrated result by Bertini and Giacomin from 1997 for the exclusion process on a circle with any particle density. The main goal of this article is to analyse the interacting particle system using the framework of regularity structures without applying the Gärtner transform, a discrete version of the Cole-Hopf transform which linearises the KPZ equation.

Our analysis relies on discretisation framework for regularity structures developed by Erhard and Hairer [EH19] as well as estimates for iterated integrals with respect to càdlàg martingales derived by Grazieschi, Matetski and Weber [GMW24]. The main technical challenge addressed in this work is the renormalisation procedure which requires a subtle analysis of regularity preserving discrete convolution operators.

## Contents

## 1 Introduction

The KPZ equation is a non-linear Stochastic PDE in $1 + 1$ dimensions for the function $h(t, x)$, depending on the time and space variables $t$ and $x$,

$$\partial_t h = \alpha\, \partial_x^2 h + \lambda(\partial_x h)^2 + \sigma\xi, \qquad h(0, \cdot) = h_0(\cdot), \tag{1.1}$$

where $\alpha, \sigma > 0$, $\lambda \neq 0$ and $\xi$ is a space-time white noise. It was introduced in [KPZ86] as a phenomenological model of random growing interfaces with lateral growth. It is now well-understood that the KPZ equation is a member of the KPZ universality class named after it. In particular, after the so-called 1:2:3 rescaling and suitable recentering, its solution converges to the *KPZ fixed point* [QS23, ACH24]. The latter is a conjectural universal scaling limit of a large class of random growing interfaces, which was first completely characterised in [MQR21]. We refer to the review papers [Qua12, Cor12] for a detailed exposition on the KPZ equation and its scaling limit.

Among the most important mathematically rigorous results about the KPZ equation are its derivation as the scaling limit of a whole class of microscopic surface growth models in suitable parameter regimes. These microscopic models are often described by discrete interacting particle systems. On the one hand, establishing such convergence results justifies the central importance of the KPZ equation—the emergence of the equation as a scaling limit for many *weakly asymmetric* growth models is sometimes dubbed *weak universality*. On the other hand, particle approximations can be used to rigorously derive properties of the solutions in their own right. For example, the convergences of such exactly solvable particle systems as ASEP and $q$-TASEP were used in [ACQ11, OQR16] and [Cor18] to derive exact formulas for the KPZ equation with special initial conditions.

The first convergence result of an interacting particle system to the KPZ equation was obtained in the celebrated work by Bertini and Giacomin [BG97] who showed that the large scale fluctuations of a *weakly asymmetric simple exclusion process (WASEP)* (which can be also described as a *solid-on-solid* model) converge to the KPZ equation. This article made a tremendous impact and analogous convergence results have been proved for many other discrete models. We mention only some of them because the list is very long: variants of weakly asymmetric exclusion processes [GJ14, GJS15, DT16, GJ17, CST18, CS18, Par19, GPS20, CT20, BFS21, Yan22, Par23, Yan23a, Yan23b, Yan24], weakly asymmetric random interfaces [EL15, Lab17], directed random polymers [AKQ14, Nak19, JMF20, SSSX21], higher spin exclusion processes [CT17], stochastic six vertex models [Lin20, CGST20], weakly asymmetric interacting Brownian motions [DGP17].

A central mathematical challenge in proving such scaling limits is the solution theory for the equation itself: the KPZ equation is classically ill-posed, in the sense that the solution of the linearised equation (the stochastic heat equation) has the spatial Hölder regularity arbitrarily close but strictly less than $\frac{1}{2}$ and the product $(\partial_x h)^2$ cannot be defined classically using theory of distributions. A renormalisation procedure which amounts to removing an infinite term has to be performed to define and construct solutions.

In Bertini-Giacomin's work [BG97], this renormalisation was implemented rigorously by introducing the so-called *Cole-Hopf solution*, which is defined by using the non-linear Cole-Hopf transform for parabolic PDEs with quadratic non-linearities [Hop50, Col51]. This transform allows to convert the KPZ equation into a stochastic heat equation with a multiplicative noise, which can be solved using the standard SPDEs techniques [Wal86].

Proving the convergence of particle approximations to the KPZ equation through the Cole-Hopf transform is challenging in general. Each individual approximation requires a particular discrete analogue of the Cole-Hopf transform, which converts an approximate KPZ equation into a simpler equation. The article [BG97] uses the Gärtner transform [Gär88] of WASEP. Finding such a transform for each approximation is a non-trivial problem, and moreover the converted equation usually appears to be non-linear and gets linearised only in the limit; see e.g. the works [DT16, Yan23b] which deal with generalisations of [BG97] to exclusion processes with non-nearest-neighbour jumps.

A natural way to avoid these complications is to work with the KPZ equation "directly", i.e. without using the Cole-Hopf transform. Defining a direct notion of solution for the KPZ equation was an open problem for a long time, and it was first solved by Hairer [Hai13] using the *theory of controlled rough paths*. This approach gave a path-wise solution to the KPZ equation, and it was moreover shown that for a suitable choice of the renormalisation constant the new solution almost surely coincides with the Cole-Hopf solution. Starting with [Hai13], the theory of such *singular stochastic PDEs* has seen massive developments in the last decade and today there are essentially two well-developed approaches to deal with the KPZ equation directly: on the one hand, there are several path-wise theories extending Hairer's work [Hai13]; these include the *theory of regularity structures* [Hai14, FH20], the theory of *paracontrolled distributions* [GIP15, GP17] and the *renormalisation group approach* [Kup16, KM17]. At a high level the idea of these path-wise solution theories is to develop a local approximation of the solution in terms of polynomials in the underlying noise, perform the renormalisation at the level of these approximate solutions and then use PDE regularity estimates to control a "remainder". On the other hand, a notion of probabilistically weak solution, *energy solutions*, for the KPZ and stochastic Burgers equations was introduced in [GJ14] and [GJ13] and uniqueness was shown in [GP18]; see also [GP20] for a characterisation of solutions in terms of their Markov generator. At a high level these theories use the explicit knowledge of invariant measures (e.g. white noise for the stochastic Burgers equation) to get control on singular non-linearity.

Going back to particle approximations, there are essentially two strategies that have been implemented successfully to date. All of the convergence results mentioned above have been derived based on one of these approaches:

- Going through the Cole-Hopf transform on the level of particle system and showing convergence at this level, following and extending the original idea in [BG97] as discussed above.
- If an explicit invariant measure is known on the level of particle dynamics (such as Bernoulli product measures for WASEP), the notion of weak solutions is a powerful tool to prove convergence.

The aim of this article is to open up third avenue to proving the large-scale convergence of interacting particle systems to the KPZ equation using a path-wise solution theory, namely the theory of regularity structures. We re-visit the classical setting of WASEP and prove the convergence to KPZ without using the Gärtner transform[1] and without making use of the explicit invariant measure. In principle, the theory of regularity structures is ideally suited to prove such scaling limits because the solution theory is itself very stable and does not rely on specific symmetries to make sense of an equation. In particular, no knowledge of invariant measures or of good behaviour under a non-linear transform is required of the particle model. However, a few technical obstacles have to be overcome in order to implement such a convergence analysis, as both PDE regularity arguments and the construction of models have to be implemented on a discrete model. These techniques have both been developed recently, most notably in [HM18, EH19] where a discretisation framework for regularity structures is given that permits to set up the PDE arguments on a lattice and in [GMW24] where a framework, which allows to prove moment bounds on stochastic integrals driven by càdlàg martingales in a similar way as it is done for multiple Wiener integrals in [HQ18], is developed. The combination of these techniques has already been successfully applied to show the convergence of suitably scaled Glauber dynamics for a Kac-Ising model to the dynamic $\Phi^4$ equation in three dimensions [GMW25]. The present article builds on these technical developments and extends them to prove convergence to the KPZ equation. The main challenge addressed is the renormalisation procedure, which is much more involved than in the Kac-Ising case, as we explain below in Section 1.2. This difficulty arises, because, in contrast to the Kac-Ising model, the exclusion process is not of a (local) mean-field type. We expect that the method developed in this article is highly flexible and that it can be transferred to a large class of particle models.

### 1.1 Model and main convergence result

We fix $N \in \mathbf{N}$ and set $\varepsilon = \frac{2}{2N+1}$ and the discrete circle $\mathbb{T}_N := \mathbb{Z}/(2N+1)\mathbb{Z}$. Let $\sigma^\varepsilon(t)$ be a WASEP on $\mathbb{T}_N$, i.e. $\sigma^\varepsilon(t)$ is a configuration of particles $(\sigma^\varepsilon(t,x) \in \{-1,1\} : x \in \mathbb{T}_N)$, where $\sigma^\varepsilon(t,x) = 1$ means that there is a particle at $x$ and $\sigma^\varepsilon(t,x) = -1$ means that the site $x$ is unoccupied. The evolution of particles is defined as follows: particles try to make unit jumps independently to the left with rate $r^\varepsilon_{x\to x-1} := \frac{1}{2} + \sqrt{\varepsilon}$ and to the right with rate $r^\varepsilon_{x\to x+1} := \frac{1}{2}$. If the destination site

[1] with the only difference that due to a limitation of the theory our WASEP is on a circle rather than whole line

is unoccupied then the particle makes a jump, and it stays put otherwise. Hence, at any time each site contains at most one particle. Here, the value $\sqrt{\varepsilon} > 0$ is a "weak asymmetry". Let $\mathcal{L}^\varepsilon$ be the generator of the Markov process $\sigma^\varepsilon(t)$. Then it acts on functions $f$ of configurations $\sigma \in \{-1, 1\}^{\mathbb{T}_N}$ as

$$\mathcal{L}^\varepsilon_{\text{WASEP}} f(\sigma) = \sum_{x \sim y \in \mathbb{T}_N} r^\varepsilon_{x\to y} \mathcal{L}^\varepsilon_{x \to y} f(\sigma), \tag{1.2}$$

where $x \sim y$ means that $x$ and $y$ are neighboring points on $\mathbb{T}_N$ and the operator $\mathcal{L}^\varepsilon_{x\to y}$ describes a jump of a particle from $x$ to $y$ as

$$\mathcal{L}^\varepsilon_{x\to y} f(\sigma) = \frac{1+\sigma(x)}{2}\frac{1-\sigma(y)}{2}\big(f(\sigma^{x,y}) - f(\sigma)\big).$$

Here, the configuration $\sigma^{x,y}$ is obtained from $\sigma$ by swapping the two values at $x$ and $y$. We note that the role of the multiplier $\frac{1+\sigma(x)}{2}\frac{1-\sigma(y)}{2}$ is to make sure that particles do not jump on top of each other; namely, this multiplier is the indicator that there is a particle at $x$ and the destination site $y$ is empty.

We are going to consider initial configurations of particles $\sigma^\varepsilon(0)$, parametrised by $\varepsilon$. Let us then define the densities of particles

$$\varrho_\varepsilon := \frac{1}{2N+1} \sum_{x \in \mathbb{T}_N} \sigma^\varepsilon(t,x) \;\in\; [-1,1], \tag{1.3}$$

which are preserved under the dynamics.

**Assumption 1** *We assume that the particle densities* (1.3) *are such that the limit* $\varrho = \lim_{\varepsilon\to 0}\varrho_\varepsilon$ *exists and is in* $(-1,1)$.

It will be convenient in what follows to identify $\mathbb{T}_N$ with $\{-N, \ldots, N\}$. We extend $\sigma^\varepsilon(t)$ periodically to $\mathbb{Z}$ and denote by $\zeta^\varepsilon(t)$ the respective $(1+1)$-dimensional periodic *solid-on-solid growth process* [BG97] on $\mathbb{Z}$. We recall that $\zeta^\varepsilon(t,x) - \zeta^\varepsilon(t, x-1) = \sigma^\varepsilon(t,x)$ for each $x \in \mathbb{Z}$, and we identify $\zeta^\varepsilon(t,x)$ with its piece-wise linear extension to all $x \in \mathbb{R}$. We note that the function $\zeta^\varepsilon(t)$ is defined in this way up to a vertical shift, and to define the process $\zeta^\varepsilon(t)$ without ambiguity we assume that it satisfies $\zeta^\varepsilon(0,0) = 0$. Its dynamics is described as follows: each local minimum on $\{-N,\ldots,N\}$ becomes a local maximum with rate $\frac{1}{2} + \sqrt{\varepsilon}$, and each local maximum becomes a local minimum with rate $\frac{1}{2}$, independently. When such an update happens at $x \in \{-N, \ldots, N\}$, it gets repeated periodically to all $(2N+1)\mathbb{Z} \setminus \{x\}$.

One can go back from $\sigma^\varepsilon$ to $\zeta^\varepsilon$ using the identity

$$\zeta^\varepsilon(t,x) = 2J^\varepsilon_t + \begin{cases} \sum_{0<y\le x} \sigma^\varepsilon(t,y) & \text{for } x > 0, \\ 0 & \text{for } x = 0, \\ -\sum_{x<y\le 0}\sigma^\varepsilon(t,y) & \text{for } x < 0, \end{cases} \tag{1.4}$$

where $J^\varepsilon_t$ is the *current of particles* through the origin on the time interval $[0,t]$, i.e. $J^\varepsilon$ increases by $1$ at time $s$ if a particle jumps from $1$ to $0$ at time $s$, and decreases by $1$ if a particle jumps from $0$ to $1$ at time $s$. We note that $\zeta^\varepsilon(t,x)$ changes by $2$ only if the spins $\sigma^\varepsilon$ at sites $x$ and $x+1$ get swapped. Hence, we get

$$\begin{aligned}\mathcal{L}^\varepsilon \zeta^\varepsilon(t,x) &= r^\varepsilon_{x\to x-1}\frac{1-\sigma^\varepsilon(t,x)}{2}\frac{1+\sigma^\varepsilon(t,x+1)}{2}2 - r^\varepsilon_{x\to x+1}\frac{1+\sigma^\varepsilon(t,x)}{2}\frac{1-\sigma^\varepsilon(t,x+1)}{2}2 \\ &= \frac{1}{2}(1+\sqrt{\varepsilon})(\sigma^\varepsilon(t,x+1) - \sigma^\varepsilon(t,x)) - \frac{1}{2}\sqrt{\varepsilon}\sigma^\varepsilon(t,x)\sigma^\varepsilon(t,x+1) + \frac{1}{2}\sqrt{\varepsilon}. \end{aligned} \tag{1.5}$$

Our goal is to study the scaling limit of the height function

$$h^\varepsilon(t,x) := \sqrt{\varepsilon}\zeta^\varepsilon(\varepsilon^{-2}t, \varepsilon^{-1}x) \tag{1.6}$$

as $\varepsilon \to 0$, parabolically rescaled by the same parameter $\varepsilon$ used in the weak asymmetry. Here, $t \geq 0$ and $x \in \Lambda_\varepsilon := \varepsilon\mathbb{Z}$. However, due to the asymmetry one can see that this function tends to $+\infty$ in the limit, and one has to subtract a deterministic drift to get a finite limit. We note furthermore that the function $\zeta^\varepsilon(t)$ is quasi-periodic, i.e. it satisfies

$$\zeta^\varepsilon(t, x + 2N + 1) = \zeta^\varepsilon(t, x) + \varrho_\varepsilon(2N + 1)$$

for all $x \in \mathbb{Z}$. Hence, when studying a scaling limit of this function, we will consider its tilted version $\zeta^\varepsilon(t, x) - \varrho_\varepsilon x$ which is $(2N + 1)$-periodic. More precisely, we are going to study the following rescaled and renormalised function:

$$\hat{h}^\varepsilon(t, x) := h^\varepsilon(t, x) - \varrho_\varepsilon \varepsilon^{-\frac{1}{2}} x + \frac{1}{2}(\varrho_\varepsilon^2 - 1)\varepsilon^{-1} t, \tag{1.7}$$

where $t \geq 0$ and $x \in \Lambda_\varepsilon$. We note that the tilting by $\varrho_\varepsilon \varepsilon^{-\frac{1}{2}} x$ makes $\hat{h}^\varepsilon(t, x)$ periodic with period 2. However, this tilting forces us to make a transformation of the spatial variable. For this, it will be convenient to denote by $\lfloor x \rfloor_\varepsilon$ the point on $\Lambda_\varepsilon$ which is the closest to $x \in \mathbb{R}$. Then we will study the limit of the function

$$\tilde{h}^\varepsilon(t, x) := \hat{h}^\varepsilon(t, x + \lfloor \varrho_\varepsilon \varepsilon^{-\frac{1}{2}} t \rfloor_\varepsilon). \tag{1.8}$$

This transformation leaves the function 2-periodic in the spatial variable.

The following is the main result of this article.

**Theorem 1.1** *Let Assumption 1 be satisfied. Let us identify the function $\tilde{h}^\varepsilon$, defined in (1.8), with its piece-wise linear extension to $\mathbb{R}$. Let there be values $\frac{2}{5} < \alpha < \bar{\alpha} < \frac{1}{2}$ and $\varepsilon_\star > 0$, and a 2-periodic function $h_0 \in \mathcal{C}^{\bar{\alpha}}(\mathbb{R})$ such that the initial state $\tilde{h}_0^\varepsilon$ of the model (1.8) at time $t = 0$ satisfies*

$$\sup_{\varepsilon \in (0, \varepsilon_\star)} \|\tilde{h}_0^\varepsilon\|_{\mathcal{C}^{\bar{\alpha}}} < \infty, \qquad \lim_{\varepsilon \to 0} \|\tilde{h}_0^\varepsilon - h_0\|_{\mathcal{C}^\alpha} = 0. \tag{1.9}$$

*Then there is a constant $c_0^\varepsilon$, bounded uniformly in $\varepsilon$, such that*

$$t \mapsto \tilde{h}^\varepsilon(t, \bullet) - c_0^\varepsilon t \tag{1.10}$$

*converge in law as $\varepsilon \to 0$ to $t \mapsto h_{\mathrm{CH}}(t, \bullet)$ with respect to the topology of the Skorokhod space $\mathcal{D}(\mathbb{R}_+, \mathcal{C}(\mathbb{R}))$ of càdlàg functions, where $h_{\mathrm{CH}}$ is the Cole-Hopf solution of the KPZ equation*

$$\partial_t h = \frac{1}{2}\Delta h - \frac{1}{2}(\partial_x h)^2 + \sqrt{1 - \varrho^2}\, \xi, \qquad h(0, \bullet) = h_0(\bullet), \tag{1.11}$$

*on $\mathbb{R}_+ \times \mathbb{R}$ driven by a spatially 2-periodic Gaussian space-time white noise $\xi$, and where $\varrho$ is from Assumption 1.*

**Remark 1.2** A special case is $\varrho_\varepsilon = \mathcal{O}(\varepsilon)$, when we consider configurations in which approximately half of sites are occupied. Then the function (1.6) is periodic and we do not need to do the tilting as in (1.8). We have $\varrho = 0$ in this case and the constant multiplier of $\xi$ in (1.11) equals 1. This is the limit proved in [BG97] but for the exclusion process on $\mathbb{Z}$.

**Remark 1.3** The constants, which appear as parameters in the definition of the height-functions (1.7)-(1.8), and the coefficients in the limiting equation (1.11) can be written as derivatives of the current of particles at the equilibrium. More precisely, the Bernoulli product measure on $\mathbb{T}_N$ with any parameter $\varrho \in [-1, 1]$ is invariant for this process [Lig05, Ch. VIII], and we denote by $\mathbf{P}_\varrho$ and $\mathbf{E}_\varrho$ this measure and the expectation with respect to it. Then

$$\mathbf{P}_\varrho(\sigma(x) = 1) = \frac{1 + \varrho}{2}, \qquad \mathbf{P}_\varrho(\sigma(x) = -1) = \frac{1 - \varrho}{2}.$$

We define the *average steady state current* of particles as in [Spo14]

$$\mathrm{Curr}_\varepsilon(\varrho) := \mathbf{E}_\varrho[r^\varepsilon_{1\to0}\mathcal{L}^\varepsilon_{1\to0}\sigma(1) - r^\varepsilon_{0\to1}\mathcal{L}^\varepsilon_{0\to1}\sigma(0)], \tag{1.12}$$

where we write $\mathcal{L}^\varepsilon_{x\to y}\sigma(x)$ as a shorthand of $\mathcal{L}^\varepsilon_{x\to y}f(\sigma)$ with $f(\sigma) = \sigma(x)$. Then we compute the constants

$$\lambda^\varepsilon_\ell := \varepsilon^{-\frac{1}{2}} \frac{2}{\ell!} \frac{d^\ell}{d\varrho^\ell} \mathrm{Curr}_\varepsilon(\varrho)\Big|_{\varrho=\varrho_\varepsilon}, \tag{1.13}$$

for integer $\ell \geq 0$. We clearly have $\mathrm{Curr}_\varepsilon(\varrho) = \sqrt{\varepsilon}\frac{\varrho^2-1}{4}$ and we can readily see that $\lambda^\varepsilon_0$ and $\lambda^\varepsilon_1$ coincide with the two constants in the divergent terms in (1.7). Moreover, we have $\lambda^\varepsilon_2 = \frac{1}{2}$ and $\lambda^\varepsilon_\ell = 0$ for all $\ell \geq 3$. The constant $\lambda^\varepsilon_2$ coincides with the coefficient $\frac{1}{2}$ in front of the non-linearity in (1.11). The diffusion coefficient $\frac{1}{2}$ in front of the Laplacian in (1.11) is usually computed by the Green-Kubo formula and coincides in this case with the symmetric jump rates of the model. One expects to be able to compute these constants from the equilibrium current also for more complicated exclusion processes with non-nearest neighbors jumps and inhomogeneous jump rates. Formal computations support this conjecture but we leave a rigorous proof for future work.

**Remark 1.4** The constant $c^\varepsilon_0$ in (1.10) does not seem to have a simple interpretation in the hydrodynamic theory and is needed to get the Cole-Hopf solution in the limit. In the case $\varrho = 0$ we have $\lim_{\varepsilon\to0} c^\varepsilon_0 = \frac{1}{24}$ where the latter is the same constant as in [BG97, Thm. 2.3]. Identifying this constant for arbitrary limiting density and giving it a physical interpretation is an interesting problem.

**Remark 1.5** The regularity on the initial conditions in [BG97, Def. 2.2] can be written in terms of the rescaled height function as follows: for any $n \geq 1$ there are constants $C$ and $a$, depending on $n$, such that

$$\mathbb{E}[e^{-n\tilde{h}^\varepsilon_0(x)}] \leq Ce^{a|x|}, \qquad \mathbb{E}[|\tilde{h}^\varepsilon_0(x) - \tilde{h}^\varepsilon_0(y)|^{2n}] \leq Ce^{a(|x|+|y|)}|x-y|^n$$

for all $x, y \in \mathbb{R}$. There are the exponential weights in these bounds because the model is studied on the whole line. By the Kolmogorov continuity test these bounds yield the Hölder regularity which we require in (1.9).

**Remark 1.6** The restriction $\alpha < \frac{1}{2}$ on the regularity of the initial condition is natural because it coincides with the regularity of the solution to the KPZ equation. The bound $\alpha > \frac{2}{5}$ is dictated by the method which we are using the prove the convergence. More precisely, the proof of moment bounds in Section 6.4.3 requires to have $1 - 2\alpha < 2\kappa$, where $\kappa$ is the value used in (5.2). As we explain in Section 5, we need to impose $\kappa < \frac{1}{10}$ in order to have a minimal number of basis elements of the regularity structure. This yields $1 - 2\alpha < 2\kappa < \frac{1}{5}$ and hence $\alpha > \frac{2}{5}$.

## 1.2 Methodology

We show in Section 2 that the martingale problem yields a discrete equation for the rescaled hight function (1.8) in the form

$$d\tilde{h}^\varepsilon(t,x) = \frac{1}{2}((1+\sqrt{\varepsilon})\Delta_\varepsilon\tilde{h}^\varepsilon - \nabla^-_\varepsilon\tilde{h}^\varepsilon\nabla^+_\varepsilon\tilde{h}^\varepsilon)(t,x)dt + E^\varepsilon(\tilde{h}^\varepsilon)(t,x)dt + d\widetilde{M}^\varepsilon(t,x), \tag{1.14}$$

where the discrete operators approximate the respective continuous differential operators. The "error" term $E^\varepsilon(\tilde{h}^\varepsilon)$ appears because of the shift of the spatial coordinate in (1.8) and it vanishes in the limit as $\varepsilon \to 0$; proof of this fact is however somewhat involved and it forces us to work with a larger regularity structure than we would otherwise. The driving random processes $\widetilde{M}^\varepsilon$ are càdlàg martingales, whose predictable covariations can be written explicitly in terms of the function $\tilde{h}^\varepsilon$ (see (2.13)). From this equation one can clearly see the structure of the limiting KPZ equation (1.11). To prove this limit we use the theory of regularity structures [Hai14] and its discretisation framework [EH19].

To control stochastic objects, which appear in the definition of a discrete "model" (it is one of the key objects in the theory of regularity structures), we use the result of [GMW24] which we slightly improve in Appendix C. This result allows to bound multiple stochastic integrals with respect to càdlàg martingales similarly to how it is done for Wiener stochastic integrals. Application of these frameworks forces us to deal with several difficulties, and we are going to describe only two of them which are the most complicated ones:

1. *Renormalisation of the non-linearity:* The product $\nabla_\varepsilon^- \tilde h^\varepsilon \nabla_\varepsilon^+ \tilde h^\varepsilon$ in the discrete equation needs to be renormalised. This is done similarly to how renormalisation is performed for continuous equations [Hai14] as

$$\nabla_\varepsilon^- \tilde h^\varepsilon(t,x) \nabla_\varepsilon^+ \tilde h^\varepsilon(t,x) - \int_{-\infty}^t \varepsilon^2 \sum_{y \in \Lambda_\varepsilon} \nabla_\varepsilon^- G^\varepsilon_{t-s}(x-y) \nabla_\varepsilon^+ G^\varepsilon_{t-s}(x-y)\, d\langle \widetilde M^\varepsilon(y)\rangle_s, \quad (1.15)$$

where $G^\varepsilon$ is the Green function of the linear operator and where $\widetilde M^\varepsilon$ is suitably extended to the time interval $(-\infty, t]$. If $\widetilde M^\varepsilon$ was a spatially mollified cylindrical Wiener process, the bracket $d\langle \widetilde M^\varepsilon(y)\rangle_s$ would be $\varepsilon^{-1} ds$, and the subtracted term would be a renormalisation constant for the non-linearity. Let us define the kernel

$$Q^\varepsilon_t(x) = \nabla_\varepsilon^- G^\varepsilon_t(x) \nabla_\varepsilon^+ G^\varepsilon_t(x).$$

One of the key properties of this discrete kernel is $\int_0^\infty \sum_{x\in\Lambda_\varepsilon} Q^\varepsilon_t(x) dt = 0$, as already observed in [BG97, Lem. A.1] (see also Lemma 4.4 below) . Then the renormalisation constant in (1.15) vanishes and, surprisingly, the non-linearity does not need to be renormalised.

In our case, $\widetilde M^\varepsilon$ are martingales such that

$$d\langle \widetilde M^\varepsilon(y)\rangle_s := \varepsilon^{-1}((1+\sqrt{\varepsilon})(1-\varrho_\varepsilon^2 - \varepsilon \nabla_\varepsilon^- \tilde h^\varepsilon(t,x) \nabla_\varepsilon^+ \tilde h^\varepsilon(t,x) + \ldots)) ds \quad (1.16)$$

for $s > 0$, where "…" contains several terms depending on $\tilde h^\varepsilon$, which are not so important at this moment. The constant term $(1+\sqrt{\varepsilon})(1-\varrho_\varepsilon^2)$ in (1.16) vanishes after the convolution with the kernel in (1.15), and (1.15) thus suggests to define the renormalised product

$$\nabla_\varepsilon^- \tilde h^\varepsilon \diamond \nabla_\varepsilon^+ \tilde h^\varepsilon = \nabla_\varepsilon^- \tilde h^\varepsilon \nabla_\varepsilon^+ \tilde h^\varepsilon + (1+\sqrt{\varepsilon})\varepsilon Q^\varepsilon \star_\varepsilon^+ (\nabla_\varepsilon^- \tilde h^\varepsilon \nabla_\varepsilon^+ \tilde h^\varepsilon) + \ldots,$$

where $\star_\varepsilon^+$ is a space-time convolution on $\mathbb{R}_+ \times \Lambda_\varepsilon$ (see (1.18)) and where "…" again contains several terms depending on $\tilde h^\varepsilon$. The bracket process of the martingale contains the same term $\nabla_\varepsilon^- \tilde h^\varepsilon \nabla_\varepsilon^+ \tilde h^\varepsilon$ which we are renormalising in (1.15). Hence, this product has to be renormalised again, which yields an infinite series of renormalisations. More precisely, we can write

$$\begin{aligned} \nabla_\varepsilon^- \tilde h^\varepsilon \nabla_\varepsilon^+ \tilde h^\varepsilon &= \nabla_\varepsilon^- \tilde h^\varepsilon \diamond \nabla_\varepsilon^+ \tilde h^\varepsilon - (1+\sqrt{\varepsilon})\varepsilon Q^\varepsilon \star_\varepsilon^+ (\nabla_\varepsilon^- \tilde h^\varepsilon \nabla_\varepsilon^+ \tilde h^\varepsilon) + \ldots \\ &= (\mathscr{Q}^\varepsilon)^{-1} (\nabla_\varepsilon^- \tilde h^\varepsilon \diamond \nabla_\varepsilon^+ \tilde h^\varepsilon) + \ldots, \end{aligned} \quad (1.17)$$

where

$$(\mathscr{Q}^\varepsilon)^{-1} f(z) = f(z) + \sum_{n=1}^\infty (-\varepsilon)^n (1+\sqrt{\varepsilon})^n (Q^\varepsilon)^{(\star_\varepsilon^+)^n} \star_\varepsilon^+ f(z)$$

is a convolution operator acting on suitable functions $f$. Another remarkable property of the function $Q^\varepsilon$ is the bound

$$\int_0^\infty \varepsilon \sum_{x \in \Lambda_\varepsilon} |Q^\varepsilon_t(x)| dt \le \varepsilon^{-1} \Theta$$

for some $\Theta \in (0,1)$, which was observed in [BG97, Lem. A.2] (see also Lemma 4.5 below). This bound guarantees that the series in the definition of $(\mathscr{Q}^\varepsilon)^{-1}$ converges absolutely in $L^1$ when $\varepsilon$ is taken sufficiently small. Plugging (1.17) into (1.14) and writing the equation in the mild form, we get

$$\tilde{h}^\varepsilon(t,x) = (G^\varepsilon_t *_\varepsilon \tilde{h}^\varepsilon_0)(x) - \frac{1}{2}\widetilde{G}^\varepsilon \star^+_\varepsilon (\nabla^-_\varepsilon \tilde{h}^\varepsilon \diamond \nabla^+_\varepsilon \tilde{h}^\varepsilon)(t,x) + (G^\varepsilon \star^+_\varepsilon d\widetilde{M}^\varepsilon)(t,x) + \ldots,$$

with a new kernel $\widetilde{G}^\varepsilon = (\mathscr{Q}^\varepsilon)^{-1} G^\varepsilon$ and where $*_\varepsilon$ is a convolution in the spatial variable. This is the actual equation whose convergence to the KPZ equation (1.11) we prove. In particular, we show that both $G^\varepsilon$ and $\widetilde{G}^\varepsilon$ converge to the same limit in a suitable sense, which is the heat kernel (see (4.37) and (8.21)), and that the error terms hidden in "..." vanish in the limit (see (8.5c)).

2. *Extensions of the discrete kernels:* A technical difficulty we encounter is in defining extensions of the discrete kernels off the grid and bounding them. In order to have Schauder-type estimates on convolutions with the discrete kernels, we need to show that their regular extensions have bounds similar to the continuous heat kernel. We define such extensions in Appendix A using the idea of [HM18, Sec. 5.1]. While proving the necessary bounds on the extension of the kernel $G^\varepsilon$ is pretty straightforward (see Section 4.2), the kernel $\widetilde{G}^\varepsilon$ is much more complicated to analyse. The necessary bounds are proved in Section 4.5 exploiting a recursive relation for the kernel.

## 1.3 Notation

We denote by $\mathbf{N}$ the set of integers $\{0,1,2,\ldots\}$ and we define $\mathbb{R}_+ := [0,\infty)$. We are going to work on the space-time domain $\mathbb{R}\times\mathbb{R}$ with the parabolic scaling $\mathfrak{s} = (2,1)$. Then for $(t,x)\in\mathbb{R}\times\mathbb{R}$ we define the parabolic distance

$$\|(t,x)\|_\mathfrak{s} := \sqrt{|t|} + |x|$$

and for a multiindex $j = (j_0, j_1) \in \mathbf{N}^2$ we set

$$|j|_\mathfrak{s} := 2j_0 + j_1$$

and the mixed derivative $D^j := \partial_t^{j_0}\partial_x^{j_1}$.

For $\varepsilon > 0$ we will work on the grid $\Lambda_\varepsilon := \varepsilon\mathbb{Z}$ of mesh size $\varepsilon$, and for two functions $f, g : \Lambda_\varepsilon \to \mathbb{R}$ we define their convolution on this grid as

$$(f *_\varepsilon g)(x) := \varepsilon \sum_{y\in\Lambda_\varepsilon} f(x-y)g(y).$$

For the domain $D_\varepsilon := \mathbb{R}\times\Lambda_\varepsilon$ and functions $f, g : D_\varepsilon \to \mathbb{R}$ we define the space-time convolution

$$(f \star_\varepsilon g)(t,x) := \int_\mathbb{R} \varepsilon \sum_{y\in\Lambda_\varepsilon} f(x-y, t-s) g(y,s) ds.$$

We denote the space-time domain $D^+_\varepsilon := \mathbb{R}_+ \times \Lambda_\varepsilon$. When working with functions $f, g : D^+_\varepsilon \to \mathbb{R}$ on a space-time domain, we define

$$(f \star^+_\varepsilon g)(t,x) := \int_0^t \varepsilon \sum_{y\in\Lambda_\varepsilon} f(x-y, t-s) g(y,s) ds. \tag{1.18}$$

It will be convenient to use the shorthands

$$\int_{D^+_\varepsilon} f(z)dz := \int_0^\infty \varepsilon \sum_{y\in\Lambda_\varepsilon} f(t,x)dt, \qquad \int_{D_\varepsilon} f(z)dz := \int_\mathbb{R} \varepsilon \sum_{y\in\Lambda_\varepsilon} f(t,x)dt. \tag{1.19}$$

A semidiscrete analogue of the mixed derivative $D^j$ will be $D^j_\varepsilon := \partial_t^{j_0}(\nabla^-_\varepsilon)^{(j_1)}$, where the discrete operator $\nabla^-_\varepsilon$ acts on the spatial variable as $\nabla^-_\varepsilon f(t,x) = \varepsilon^{-1}(f(t,x) - f(t,x-\varepsilon))$.

Sometimes we will use an operator $A^\varepsilon$ acting on functions $f : D_\varepsilon^+ \to \mathbb{R}$, which will be defined as a convolution operator with a kernel, which we denote by $A^\varepsilon : D_\varepsilon^+ \times D_\varepsilon^+ \to \mathbb{R}$, so that

$$(A^\varepsilon f)(t, x) := (A^\varepsilon \star_\varepsilon^+ f)(t, x),$$

provided the convolution converges absolutely. The only exception will be the identity operator Id, mapping $f$ to itself.

For integer $n \geq 0$, we denote by $\mathcal{C}_0^n$ the set of compactly supported $\mathcal{C}^n$ functions $\varphi : \mathbb{R} \to \mathbb{R}$, and the set $\mathcal{B}^n$ contains all functions $\varphi \in \mathcal{C}_0^n$, which are supported on $[-1, 1]$, and which satisfy $\|\varphi\|_{\mathcal{C}^n} \leq 1$. For a function $\varphi \in \mathcal{B}^n$, for $x \in \mathbb{R}$ and for $\lambda \in (0, 1]$, we define its rescaled and recentered version

$$\varphi_x^\lambda(y) := \frac{1}{\lambda}\varphi\Big(\frac{y - x}{\lambda}\Big). \tag{1.20}$$

We denote by $\mathscr{D}'(\mathbb{R})$ the space of distributions on $\mathbb{R}$. For $\alpha < 0$ we define the Besov space $\mathcal{C}^\alpha$ as a completion of smooth compactly supported functions $f : \mathbb{R} \to \mathbb{R}$, under the seminorm

$$\|f\|_{\mathcal{C}^\alpha} := \sup_{\varphi \in \mathcal{B}^r} \sup_{x \in \mathbb{R}} \sup_{\lambda \in (0,1]} \lambda^{-\alpha} |f(\varphi_x^\lambda)| < \infty,$$

for $r$ being the smallest integer such that $r > -\alpha$, where we write $f(\varphi_x^\lambda) = \langle f, \varphi_x^\lambda \rangle$ for the duality pairing. It is important to define these spaces as completions of smooth functions, because this makes the spaces separable.

We need to define discrete norms and seminorms for discrete functions $f^\varepsilon : \Lambda_\varepsilon \to \mathbb{R}$. For $\alpha \in (0, 1)$, a discrete analogue of the $\alpha$-Hölder norm is

$$\|f^\varepsilon\|_{\mathcal{C}^\alpha}^{(\varepsilon)} := \sup_{x \in \Lambda_\varepsilon} |f^\varepsilon(x)| + \sup_{\substack{x, y \in \Lambda_\varepsilon \\ 0 < |x - y| \leq 1}} \frac{|f^\varepsilon(x) - f^\varepsilon(y)|}{|x - y|^\alpha}. \tag{1.21}$$

We map a discrete function $f^\varepsilon : \Lambda_\varepsilon \to \mathbb{R}$ to a distribution as

$$(\iota_\varepsilon f^\varepsilon)(\varphi) := \varepsilon \sum_{x \in \Lambda_\varepsilon} f^\varepsilon(x) \varphi(x), \tag{1.22}$$

for any continuous and compactly supported function $\varphi$. In the case $\varepsilon = 1$ we will write $\iota = \iota_1$. Then for $\alpha < 0$ and for the smallest integer $r$ such that $r > -\alpha$, we define the seminorm

$$\|f^\varepsilon\|_{\mathcal{C}^\alpha}^{(\varepsilon)} := \sup_{\varphi \in \mathcal{B}^r} \sup_{x \in \Lambda_\varepsilon} \sup_{\lambda \in (0,1]} (\lambda \vee \varepsilon)^{-\alpha} |(\iota_\varepsilon f^\varepsilon)(\varphi_x^\lambda)|. \tag{1.23}$$

We note that this definition implies that on scales larger than $\varepsilon$ the distribution $\iota_\varepsilon f^\varepsilon$ behaves as a Besov distribution from $\mathcal{C}^\alpha$. On the other hand, on the scale $\varepsilon$ we have $|f^\varepsilon(x)| \lesssim \varepsilon^\alpha$.

When we write $f^\varepsilon \in \mathcal{C}_\varepsilon^\alpha$, we mean that there is $\varepsilon_0 > 0$ and we consider a sequence of functions $f^\varepsilon$ parametrised by $\varepsilon \in (0, \varepsilon_0)$, such that the norms $\|f^\varepsilon\|_{\mathcal{C}^\alpha}^{(\varepsilon)}$ are bounded uniformly in $\varepsilon \in (0, \varepsilon_0)$.

Depending on context, we will use both notations $f(t, x)$ and $f_t(x)$ for functions depending on space and time variables.

The notation $a \lesssim b$ means that there is a constant $C$, independent of relevant quantities, such that $A \leq Cb$. These "relevant quantities" will be usually $\varepsilon$ and space-time points and they will be always clear from the context. When we want to specify dependence of $C$ on parameters $j, k, \ldots$, we will write $C = C(j, k, \ldots)$.

## 1.4 Structure of the article

A discrete equation for the rescaled and renormalised height function is derived in Section 2. We recall the solution theory for the KPZ equation and essential properties of the solution in Section 3. Discrete kernels and operators involved in the discrete equation are analysed in Section 4. In particular, we define extensions of the discrete heat kernel and its modifications off the grid. For this we use a quite general construction of regular extensions described in Appendix A. A regularity structure for the discrete equation is define in Section 5, and Section 6 is devoted to constructing of a renormalised discrete model and proving its moment bounds. The moment bounds follow from an estimate on generalised convolutions from Appendix C. Properties of the solution to the discrete equation are studied in Section 7 and the main convergence result is proved in Section 8. We prove bounds on some function, defined via the symmetric exclusion process, in Appendix B.

**Acknowledgments**

HW is funded by the European Union (ERC, GE4SPDE, No. 101045082). HW and RH are funded by the Deutsche Forschungsgemeinschaft (DFG, German Research Foundation) under Germany's Excellence Strategy EXC 2044-390685587, Mathematics Münster: Dynamics-Geometry-Structure. KM is partially supported by an MSU startup grant.

The authors thank Jeremy Quastel for numerous discussions on this and related problems.

## 2 Semidiscrete equations from the particle system

We can write a semidiscrete equation for the function $h^\varepsilon(t,x)$ defined in (1.6). We first compute the drift by applying the generator $\mathcal{L}^\varepsilon$. Namely, applying rescaling (1.6) to (1.5) we get

$$\begin{aligned}\mathcal{L}^\varepsilon h^\varepsilon(t,x) = \frac{1}{2}\sqrt{\varepsilon}(1+\sqrt{\varepsilon})(\sigma^\varepsilon(\varepsilon^{-2}t,\varepsilon^{-1}x+1) - \sigma^\varepsilon(\varepsilon^{-2}t,\varepsilon^{-1}x))\\ - \frac{1}{2}\varepsilon\sigma^\varepsilon(\varepsilon^{-2}t,\varepsilon^{-1}x)\sigma^\varepsilon(\varepsilon^{-2}t,\varepsilon^{-1}x+1) + \frac{1}{2}\varepsilon.\end{aligned}$$

Defining the discrete forward/backward derivative and Laplacian

$$\nabla_\varepsilon^\pm g(x) = \pm\,\varepsilon^{-1}(g(x\pm\varepsilon) - g(x)), \qquad \Delta_\varepsilon g(x) = \varepsilon^{-2}(g(x+\varepsilon) - 2g(x) + g(x-\varepsilon)), \tag{2.1}$$

and using (1.4) together with (1.6), we can write the preceding expression as

$$\mathcal{L}^\varepsilon h^\varepsilon(t,x) = \frac{1}{2}\varepsilon^2(1+\sqrt{\varepsilon})\Delta_\varepsilon h^\varepsilon(t,x) - \frac{1}{2}\varepsilon^2\nabla_\varepsilon^- h^\varepsilon(t,x)\nabla_\varepsilon^+ h^\varepsilon(t,x) + \frac{1}{2}\varepsilon.$$

The martingale problem allows to write the semi-discrete equation for $h^\varepsilon$:

$$\begin{aligned}dh^\varepsilon(t,x) &= \varepsilon^{-2}\mathcal{L}^\varepsilon h^\varepsilon(t,x)dt + dM^\varepsilon(t,x)\\ &= \frac{1}{2}((1+\sqrt{\varepsilon})\Delta_\varepsilon h^\varepsilon(t,x) - \nabla_\varepsilon^- h^\varepsilon(t,x)\nabla_\varepsilon^+ h^\varepsilon(t,x) + \varepsilon^{-1})dt + dM^\varepsilon(t,x),\end{aligned} \tag{2.2}$$

on $t \geq 0$ and $x \in \Lambda_\varepsilon$, and with a fixed initial state $h^\varepsilon(0,x) = h_0^\varepsilon(x)$ at time $0$. Here, the processes $t \mapsto M^\varepsilon(t,x)$ are càdlàg martingales whose predictable quadratic covariations can be written in terms of the carré-du-champ operator. We note that $M^\varepsilon(t,x)$ is 2-periodic in $x$, and hence it is sufficient to characterize these martingales only restricted to $x \in \mathbb{T}_\varepsilon$. In particular, we have

$$\langle M^\varepsilon(\bullet, x+2m), M^\varepsilon(\bullet, x'+2m')\rangle_t = \langle M^\varepsilon(\bullet,x), M^\varepsilon(\bullet,x')\rangle_t$$

for any $m, m' \in \mathbb{Z}$ and $x, x' \in \mathbb{T}_\varepsilon$. Two martingales $M^\varepsilon(t,x)$ and $M^\varepsilon(t,x')$ with $x \neq x' \in \mathbb{T}_\varepsilon$ a.s. do not change at the same time, because $h^\varepsilon(t,x)$ and $h^\varepsilon(t,x')$ a.s. do not change at the same time. This implies that $\langle M^\varepsilon(\bullet,x), M^\varepsilon(\bullet,x')\rangle_t = 0$. Furthermore,

$$\begin{aligned}d\langle M^\varepsilon(\bullet,x)\rangle_t = \varepsilon^{-2}\bigg( r^\varepsilon_{x\to x-1}\frac{1-\sigma^\varepsilon(\varepsilon^{-2}t,\varepsilon^{-1}x)}{2}\frac{1+\sigma^\varepsilon(\varepsilon^{-2}t,\varepsilon^{-1}x+1)}{2}(2\sqrt{\varepsilon})^2\\ + r^\varepsilon_{x\to x+1}\frac{1+\sigma^\varepsilon(\varepsilon^{-2}t,\varepsilon^{-1}x)}{2}\frac{1-\sigma^\varepsilon(\varepsilon^{-2}t,\varepsilon^{-1}x+1)}{2}(2\sqrt{\varepsilon})^2\bigg)dt\end{aligned}$$

$$= \varepsilon^{-1}\mathbf{C}_\varepsilon(t,x)dt, \tag{2.3}$$

with

$$\mathbf{C}_\varepsilon(t,x) := (1+\sqrt{\varepsilon})(1-\varepsilon\nabla_\varepsilon^- h^\varepsilon(t,x)\nabla_\varepsilon^+ h^\varepsilon(t,x)) + \varepsilon^2\Delta_\varepsilon h^\varepsilon(t,x). \tag{2.4}$$

We can get from (2.2) as similar equation for the function (1.7). We have $\nabla_\varepsilon^\pm \hat{h}^\varepsilon(t,x) = \nabla_\varepsilon^\pm h^\varepsilon(t,x) - \varrho_\varepsilon\varepsilon^{-\frac{1}{2}}$ which yields

$$\nabla_\varepsilon^- h^\varepsilon(t,x)\nabla_\varepsilon^+ h^\varepsilon(t,x) = \nabla_\varepsilon^- \hat{h}^\varepsilon\nabla_\varepsilon^+\hat{h}^\varepsilon + 2\varrho_\varepsilon\varepsilon^{-\frac{1}{2}}\nabla_\varepsilon\hat{h}^\varepsilon + \varrho_\varepsilon^2\varepsilon^{-1},$$

where $\nabla_\varepsilon := \frac{1}{2}(\nabla_\varepsilon^- + \nabla_\varepsilon^+)$ is the symmetric finite difference derivative, which can be written explicitly as

$$\nabla_\varepsilon g(x) = \frac{1}{2}\varepsilon^{-1}(g(x+\varepsilon) - g(x-\varepsilon)). \tag{2.5}$$

Using furthermore the identity $\Delta_\varepsilon\hat{h}^\varepsilon(t,x) = \Delta_\varepsilon h^\varepsilon(t,x)$, we get

$$d\hat{h}^\varepsilon(t,x) = \frac{1}{2}((1+\sqrt{\varepsilon})\Delta_\varepsilon\hat{h}^\varepsilon - \nabla_\varepsilon^-\hat{h}^\varepsilon\nabla_\varepsilon^+\hat{h}^\varepsilon - 2\varrho_\varepsilon\varepsilon^{-\frac{1}{2}}\nabla_\varepsilon\hat{h}^\varepsilon)(t,x)dt + dM^\varepsilon(t,x) \tag{2.6}$$

with the initial state $\hat{h}^\varepsilon(0,x)$ at time $0$. The function (2.4) can be written as

$$\mathbf{C}_\varepsilon(t,x) = \nu^\varepsilon - (1+\sqrt{\varepsilon})(\varepsilon\nabla_\varepsilon^-\hat{h}^\varepsilon(t,x)\nabla_\varepsilon^+\hat{h}^\varepsilon(t,x) + 2\varrho_\varepsilon\varepsilon^{\frac{1}{2}}\nabla_\varepsilon\hat{h}^\varepsilon(t,x)) + \varepsilon^2\Delta_\varepsilon\hat{h}^\varepsilon(t,x), \tag{2.7}$$

where

$$\nu^\varepsilon := (1+\sqrt{\varepsilon})(1-\varrho_\varepsilon^2). \tag{2.8}$$

Assumption 1 yields existence of the limit

$$\nu := \lim_{\varepsilon\to 0}\nu^\varepsilon = 1-\varrho^2, \tag{2.9}$$

which will be used throughout the paper.

Now, we will derive an equation for the function (1.8). The goal of the change of the spatial variable in (1.8) is to remove the $\varepsilon^{-\frac{1}{2}}$-term on the right-hand side of (2.6). We note that the shift of the spatial variable in (1.8) has effect only at the time points $t_n^\varepsilon = \varrho_\varepsilon^{-1}\varepsilon^{\frac{3}{2}}n$ with $n \in \mathbb{Z}\setminus\{0\}$, at which the spatial variable gets changed by $\varepsilon$. We note that $\varrho_\varepsilon \neq 0$ because there is an odd number of terms in the sum in (1.3). Let us derive the equation in the case $\varrho_\varepsilon > 0$ and for $\varrho_\varepsilon < 0$ it can be derived similarly. For $\varrho_\varepsilon > 0$ and a fixed $t > 0$, let $n_\star$ be the maximal $n \in \mathbf{N}$ such that $t_n^\varepsilon \le t$. Let furthermore $x_n^\varepsilon = x + \varepsilon n$. Then we have

$$\begin{aligned}\tilde{h}^\varepsilon(t,x) - \tilde{h}_0^\varepsilon(x) &= (\hat{h}^\varepsilon(t,x_{n_\star}^\varepsilon) - \hat{h}^\varepsilon(t_{n_\star}^\varepsilon, x_{n_\star}^\varepsilon)) \\ &\quad + \sum_{n=1}^{n_\star}(\hat{h}^\varepsilon(t_n^\varepsilon, x_{n-1}^\varepsilon) - \hat{h}^\varepsilon(t_{n-1}^\varepsilon, x_{n-1}^\varepsilon)) + \sum_{n=1}^{n_\star}(\hat{h}^\varepsilon(t_n^\varepsilon, x_n^\varepsilon) - \hat{h}^\varepsilon(t_n^\varepsilon, x_{n-1}^\varepsilon)).\end{aligned}$$

Using (2.6) to write each time increment of the function $\hat{h}^\varepsilon$, we get

$$\begin{aligned}\tilde{h}^\varepsilon(t,x) - \tilde{h}_0^\varepsilon(x) &= \int_{t_{n_\star}^\varepsilon}^t \frac{1}{2}((1+\sqrt{\varepsilon})\Delta_\varepsilon\hat{h}^\varepsilon - \nabla_\varepsilon^-\hat{h}^\varepsilon\nabla_\varepsilon^+\hat{h}^\varepsilon - 2\varrho_\varepsilon\varepsilon^{-\frac{1}{2}}\nabla_\varepsilon\hat{h}^\varepsilon)(s,x_{n_\star}^\varepsilon)ds \\ &\quad + \sum_{n=1}^{n_\star}\int_{t_{n-1}^\varepsilon}^{t_n^\varepsilon} \frac{1}{2}((1+\sqrt{\varepsilon})\Delta_\varepsilon\hat{h}^\varepsilon - \nabla_\varepsilon^-\hat{h}^\varepsilon\nabla_\varepsilon^+\hat{h}^\varepsilon - 2\varrho_\varepsilon\varepsilon^{-\frac{1}{2}}\nabla_\varepsilon\hat{h}^\varepsilon)(s,x_{n-1}^\varepsilon)ds \\ &\quad + \sum_{n=1}^{n_\star}(\hat{h}^\varepsilon(t_n^\varepsilon,x_n^\varepsilon) - \hat{h}^\varepsilon(t_n^\varepsilon,x_{n-1}^\varepsilon)) + \widetilde{M}^\varepsilon(t,x)\end{aligned}$$

with the martingale

$$\widetilde{M}^\varepsilon(t,x) = M^\varepsilon(t, x^\varepsilon_{n_\star}) - M^\varepsilon(t^\varepsilon_{n_\star}, x^\varepsilon_{n_\star}) + \sum_{n=1}^{n_\star}(M^\varepsilon(t^\varepsilon_n, x^\varepsilon_{n-1}) - M^\varepsilon(t^\varepsilon_{n-1}, x^\varepsilon_{n-1}))$$
$$= \int_0^t dM^\varepsilon(s, x + \lfloor \varrho_\varepsilon \varepsilon^{-\frac{1}{2}} s \rfloor_\varepsilon).$$

Hence, we get

$$d\tilde{h}^\varepsilon(t,x) = \frac{1}{2}((1+\sqrt{\varepsilon})\Delta_\varepsilon \tilde{h}^\varepsilon - \nabla^-_\varepsilon \tilde{h}^\varepsilon \nabla^+_\varepsilon \tilde{h}^\varepsilon)(t,x)dt + E^\varepsilon(\tilde{h}^\varepsilon)(t,x)dt + d\widetilde{M}^\varepsilon(t,x), \tag{2.10}$$

where the error term $E^\varepsilon(\tilde{h}^\varepsilon)$ equals

$$E^\varepsilon(\tilde{h}^\varepsilon)(t,x) = \sum_{s \in \varrho_\varepsilon^{-1}\varepsilon^{3/2}\mathbf{N}} (\tilde{h}^\varepsilon(s,x) - \tilde{h}^\varepsilon(s, x-\varepsilon))\delta(t-s) - \varrho_\varepsilon \varepsilon^{-\frac{1}{2}} \nabla_\varepsilon \tilde{h}^\varepsilon(t,x) \tag{2.11}$$

if $\varrho_\varepsilon > 0$. In the case $\varrho_\varepsilon < 0$ we get the same equation with the error term

$$E^\varepsilon(\tilde{h}^\varepsilon)(t,x) = \sum_{s \in |\varrho_\varepsilon|^{-1}\varepsilon^{3/2}\mathbf{N}} (\tilde{h}^\varepsilon(s,x) - \tilde{h}^\varepsilon(s, x+\varepsilon))\delta(t-s) - \varrho_\varepsilon \varepsilon^{-\frac{1}{2}} \nabla_\varepsilon \tilde{h}^\varepsilon(t,x). \tag{2.12}$$

One can now see that (2.10) is a discretized KPZ equation driven by the processes $\widetilde{M}^\varepsilon$. The processes $t \mapsto \widetilde{M}^\varepsilon(t,x)$ are martingales which have the same properties as $M^\varepsilon$ but with the function

$$\widetilde{\mathbf{C}}_\varepsilon(t,x) := \nu^\varepsilon - (1+\sqrt{\varepsilon})(\varepsilon \nabla^-_\varepsilon \tilde{h}^\varepsilon(t,x) \nabla^+_\varepsilon \tilde{h}^\varepsilon(t,x) + 2\varrho_\varepsilon \varepsilon^{\frac{1}{2}} \nabla_\varepsilon \tilde{h}^\varepsilon(t,x)) + \varepsilon^2 \Delta_\varepsilon \tilde{h}^\varepsilon(t,x) \tag{2.13}$$

instead of (2.7).

### 2.1 Renormalised semidiscrete equation

In order to prove convergence of the solution of equation (2.10), one should renormalise the non-linearity to compensate its divergence in the limit. Renormalisation of the non-linearity in (2.10) is however quite involved because the bracket process of the driving martingales, which is used in the renormalisation, produces a term containing exactly the non-linearity which we want to renormalise. We explain it in more details below.

We utilize the Green function $G^\varepsilon$ of the linear operator $\frac{d}{dt} - \frac{1}{2}(1+\sqrt{\varepsilon})\Delta_\varepsilon$ to write equation (2.10) in a mild form:

$$\tilde{h}^\varepsilon(t,x) = (G^\varepsilon_t *_\varepsilon \tilde{h}^\varepsilon_0)(x) - \frac{1}{2} G^\varepsilon \star^+_\varepsilon (\nabla^-_\varepsilon \tilde{h}^\varepsilon \nabla^+_\varepsilon \tilde{h}^\varepsilon)(t,x) + (G^\varepsilon \star^+_\varepsilon E^\varepsilon(\tilde{h}^\varepsilon))(t,x) + (G^\varepsilon \star^+_\varepsilon d\widetilde{M}^\varepsilon)(t,x), \tag{2.14}$$

where we use the convolutions $*_\varepsilon$ and $\star^+_\varepsilon$ defined in Section 1.3.

We can write $\tilde{h}^\varepsilon = \widetilde{X}^\varepsilon + \widetilde{Y}^\varepsilon$, where $\widetilde{X}^\varepsilon$ is a linear approximation of the solution to (2.10) given by

$$d\widetilde{X}^\varepsilon(t,x) = \frac{1}{2}(1+\sqrt{\varepsilon})\Delta_\varepsilon \widetilde{X}^\varepsilon(t,x)dt + d\widetilde{M}^\varepsilon(t,x)$$

with the initial condition $\widetilde{X}^\varepsilon(0,x) \equiv 0$. The process $\widetilde{Y}^\varepsilon$ is usually expected to have a higher regularity than $\widetilde{X}^\varepsilon$, which we do not specify precisely at this point (such linear approximation is an analogue of the Da Prato-Debussche ansatz used first in [DPD03] and then for several other equations). We can write explicitly

$$\widetilde{X}^\varepsilon(t,x) = \int_0^t \varepsilon \sum_{y \in \Lambda_\varepsilon} G^\varepsilon_{t-s}(x-y)\, d\widetilde{M}^\varepsilon(s,y),$$

where $G^\varepsilon$ is the Green function of the linear operator on $D_\varepsilon^+$ and the stochastic integral is with respect to the time variable $s$. Then the non-linearity in (2.2) is approximated by $\nabla_\varepsilon^- \widetilde{X}^\varepsilon \nabla_\varepsilon^+ \widetilde{X}^\varepsilon$ and the latter needs to be renormalised in order to have a non-trivial limit as $\varepsilon \to 0$. More precisely, we need to consider the renormalised non-linearity

$$\nabla_\varepsilon^- \tilde{h}^\varepsilon(t,x) \nabla_\varepsilon^+ \tilde{h}^\varepsilon(t,x) - \int_0^t \varepsilon^2 \sum_{y \in \Lambda_\varepsilon} \nabla_\varepsilon^- G_{t-s}^\varepsilon(x-y) \nabla_\varepsilon^+ G_{t-s}^\varepsilon(x-y)\, d\langle \widetilde{M}^\varepsilon(y)\rangle_s. \tag{2.15}$$

We will refer to the last term as a *renormalisation function*, in contrast to the renormalisation constants which are typically used for stationary Gaussian noises.

**Remark 2.1** It is important to note that we do not set the last term in (2.15) to be

$$\int_0^t \varepsilon^2 \sum_{y_1, y_2 \in \Lambda_\varepsilon} \nabla_\varepsilon^- G_{t-s}^\varepsilon(x-y_1) \nabla_\varepsilon^+ G_{t-s}^\varepsilon(x-y_2)\, d\langle \widetilde{M}^\varepsilon(y_1), \widetilde{M}^\varepsilon(y_2)\rangle_s$$

which may look more natural. This is because we show in Lemma 4.3 that a singular part of the preceding expression is exactly the diagonal sum as in (2.15), while the other terms in the sum are bounded uniformly in $\varepsilon$. (There will be non-vanishing terms for the points $y_2 - y_1 \in 2\mathbb{Z}$ due to periodicity of the martingales.)

We recall that $d\langle \widetilde{M}^\varepsilon(\bullet, x)\rangle_t = \varepsilon^{-1} \widetilde{\mathbf{C}}_\varepsilon(t,x) dt$ with the function $\widetilde{\mathbf{C}}_\varepsilon$ defined in (2.13). Then the renormalisation function (2.15) contains the term

$$-(1+\sqrt{\varepsilon}) \int_0^t \varepsilon^2 \sum_{y \in \Lambda_\varepsilon} \nabla_\varepsilon^- G_{t-s}^\varepsilon(x-y) \nabla_\varepsilon^+ G_{t-s}^\varepsilon(x-y) \nabla_\varepsilon^- \tilde{h}^\varepsilon(s,x) \nabla_\varepsilon^+ \tilde{h}^\varepsilon(s,x) ds$$

and hence we renormalise the non-linearity in (2.15) by a function written in terms of the same non-linearity. This forces us to define the renormalisation of (2.15) in a more complicated way.

Let us define the kernel

$$Q_t^\varepsilon(x) := \nabla_\varepsilon^- G_t^\varepsilon(x) \nabla_\varepsilon^+ G_t^\varepsilon(x). \tag{2.16}$$

Then, using (2.13) we define the renormalised product

$$\nabla_\varepsilon^- \tilde{h}^\varepsilon \diamond \nabla_\varepsilon^+ \tilde{h}^\varepsilon := \nabla_\varepsilon^- \tilde{h}^\varepsilon \nabla_\varepsilon^+ \tilde{h}^\varepsilon - \widetilde{C}_\varepsilon, \tag{2.17}$$

where $\widetilde{C}_\varepsilon$ is a random function defined as

$$\widetilde{C}_\varepsilon(t,x) = Q^\varepsilon \star_\varepsilon^+ \Big(-(1+\sqrt{\varepsilon})(\varepsilon \nabla_\varepsilon^- \tilde{h}^\varepsilon \nabla_\varepsilon^+ \tilde{h}^\varepsilon + 2\varrho_\varepsilon \varepsilon^{\frac{1}{2}} \nabla_\varepsilon \tilde{h}^\varepsilon) + \varepsilon^2 \Delta_\varepsilon \tilde{h}^\varepsilon\Big)(t,x). \tag{2.18}$$

We note that we have not included the constant $\nu^\varepsilon$ from (2.13) into this renormalisation function. This is because we prove in Lemma 4.4 that the convolution of $Q^\varepsilon$ with a constant vanishes. Hence, we obtain the equation for the non-linearity

$$\begin{aligned}\nabla_\varepsilon^- \tilde{h}^\varepsilon \nabla_\varepsilon^+ \tilde{h}^\varepsilon = \nabla_\varepsilon^- \tilde{h}^\varepsilon \diamond \nabla_\varepsilon^+ \tilde{h}^\varepsilon + \widetilde{C}_\varepsilon &= \nabla_\varepsilon^- \tilde{h}^\varepsilon \diamond \nabla_\varepsilon^+ \tilde{h}^\varepsilon - (1+\sqrt{\varepsilon})\varepsilon Q^\varepsilon \star_\varepsilon^+ (\nabla_\varepsilon^- \tilde{h}^\varepsilon \nabla_\varepsilon^+ \tilde{h}^\varepsilon) \\ &\quad + Q^\varepsilon \star_\varepsilon^+ \Big(-2\varrho_\varepsilon \varepsilon^{\frac{1}{2}} (1+\sqrt{\varepsilon}) \nabla_\varepsilon \tilde{h}^\varepsilon + \varepsilon^2 \Delta_\varepsilon \tilde{h}^\varepsilon\Big).\end{aligned} \tag{2.19}$$

Let us moreover define the convolution operator

$$\mathscr{Q}^\varepsilon f(z) := f(z) + (1+\sqrt{\varepsilon})\varepsilon Q^\varepsilon \star_\varepsilon^+ f(z)$$

on $L^\infty(D_\varepsilon^+)$. We can conclude from Lemma 4.5 below that for any $\varepsilon \in (0, \varepsilon_0]$ the operator norm of $\mathscr{Q}^\varepsilon - \mathrm{Id}$ is bounded by $(1+\sqrt{\varepsilon})\Theta$, and there is $0 < \varepsilon_1 \leq \varepsilon_0$ such that $(1+\sqrt{\varepsilon})\Theta < 1$ for all $\varepsilon \in (0, \varepsilon_1]$. Then for such values of $\varepsilon$ the operator has the inverse which can be written as

$$(\mathscr{Q}^\varepsilon)^{-1} f(z) = f(z) + \sum_{n=1}^{\infty} (-\varepsilon)^n (1+\sqrt{\varepsilon})^n (Q^\varepsilon)^{(\star_\varepsilon^+)^n} \star_\varepsilon^+ f(z). \tag{2.20}$$

This is the Neumann series of the operator $(I+A)^{-1} = I + \sum_{n=1}^\infty (-A)^n$ provided it converges in the operator norm. Applying the inverse to (2.19) yields

$$\begin{aligned} \nabla_\varepsilon^- \tilde{h}^\varepsilon \nabla_\varepsilon^+ \tilde{h}^\varepsilon &= (\mathscr{Q}^\varepsilon)^{-1} (\nabla_\varepsilon^- \tilde{h}^\varepsilon \diamond \nabla_\varepsilon^+ \tilde{h}^\varepsilon) \\ &\quad + (\mathscr{Q}^\varepsilon)^{-1} Q^\varepsilon \star_\varepsilon^+ \Big( -2\varrho_\varepsilon \varepsilon^{\frac{1}{2}} (1+\sqrt{\varepsilon}) \nabla_\varepsilon \tilde{h}^\varepsilon + \varepsilon^2 \Delta_\varepsilon \tilde{h}^\varepsilon \Big). \end{aligned}$$

Plugging this into equation (2.10) gives

$$\begin{aligned} d\tilde{h}^\varepsilon(t,x) &= \frac{1}{2} \Big( (1+\sqrt{\varepsilon}) \Delta_\varepsilon \tilde{h}^\varepsilon - (\mathscr{Q}^\varepsilon)^{-1} (\nabla_\varepsilon^- h^\varepsilon \diamond \nabla_\varepsilon^+ h^\varepsilon) \Big)(t,x) dt + E^\varepsilon(\tilde{h}^\varepsilon)(t,x) dt \\ &\quad - \frac{1}{2} (\mathscr{Q}^\varepsilon)^{-1} Q^\varepsilon \star_\varepsilon^+ \Big( -2\varrho_\varepsilon \varepsilon^{\frac{1}{2}} (1+\sqrt{\varepsilon}) \nabla_\varepsilon \tilde{h}^\varepsilon + \varepsilon^2 \Delta_\varepsilon \tilde{h}^\varepsilon \Big)(t,x) dt + d\widetilde{M}^\varepsilon(t,x). \end{aligned}$$

It will be convenient to define new kernels

$$\widetilde{G}^\varepsilon := (\mathscr{Q}^\varepsilon)^{-1} G^\varepsilon, \qquad \widetilde{\widetilde{G}}^\varepsilon := \widetilde{G}^\varepsilon \star_\varepsilon^+ \varepsilon Q^\varepsilon, \tag{2.21}$$

to write this equation in a mild form as

$$\begin{aligned} \tilde{h}^\varepsilon(t,x) &= (G_t^\varepsilon *_\varepsilon \tilde{h}_0^\varepsilon)(x) - \frac{1}{2} \widetilde{G}^\varepsilon \star_\varepsilon^+ (\nabla_\varepsilon^- \tilde{h}^\varepsilon \diamond \nabla_\varepsilon^+ \tilde{h}^\varepsilon)(t,x) + (G^\varepsilon \star_\varepsilon^+ E^\varepsilon(\tilde{h}^\varepsilon))(t,x) \\ &\quad - \frac{1}{2} \widetilde{\widetilde{G}}^\varepsilon \star_\varepsilon^+ \Big( -2\varrho_\varepsilon \varepsilon^{-\frac{1}{2}} (1+\sqrt{\varepsilon}) \nabla_\varepsilon \tilde{h}^\varepsilon + \varepsilon \Delta_\varepsilon \tilde{h}^\varepsilon \Big)(t,x) + (G^\varepsilon \star_\varepsilon^+ d\widetilde{M}^\varepsilon)(t,x). \end{aligned}$$

We write furthermore $2\varrho_\varepsilon (1+\sqrt{\varepsilon}) \varepsilon^{-\frac{1}{2}} \nabla_\varepsilon \tilde{h}^\varepsilon = 2\varrho_\varepsilon \nabla_\varepsilon \tilde{h}^\varepsilon + 2\varrho_\varepsilon \varepsilon^{-\frac{1}{2}} \nabla_\varepsilon \tilde{h}^\varepsilon$ and group the error terms

$$\begin{aligned} \tilde{h}^\varepsilon(t,x) &= (G_t^\varepsilon *_\varepsilon \tilde{h}_0^\varepsilon)(x) - \frac{1}{2} \widetilde{G}^\varepsilon \star_\varepsilon^+ (\nabla_\varepsilon^- \tilde{h}^\varepsilon \diamond \nabla_\varepsilon^+ \tilde{h}^\varepsilon)(t,x) \\ &\quad + (G^\varepsilon \star_\varepsilon^+ d\widetilde{M}^\varepsilon)(t,x) + \widetilde{E}^\varepsilon(\tilde{h}^\varepsilon)(t,x), \end{aligned} \tag{2.22}$$

where

$$\widetilde{E}^\varepsilon(\tilde{h}^\varepsilon) = \widetilde{E}_1^\varepsilon(\tilde{h}^\varepsilon) + \widetilde{E}_2^\varepsilon(\tilde{h}^\varepsilon) \tag{2.23}$$

with

$$\widetilde{E}_1^\varepsilon(\tilde{h}^\varepsilon)(t,x) := (G^\varepsilon \star_\varepsilon^+ E^\varepsilon(\tilde{h}^\varepsilon))(t,x), \tag{2.24a}$$

$$\widetilde{E}_2^\varepsilon(\tilde{h}^\varepsilon)(t,x) := \varrho_\varepsilon (1+\varepsilon^{-\frac{1}{2}}) (\widetilde{\widetilde{G}}^\varepsilon \star_\varepsilon^+ \nabla_\varepsilon \tilde{h}^\varepsilon)(t,x) - \frac{\varepsilon}{2} (\nabla_\varepsilon^+ \widetilde{\widetilde{G}}^\varepsilon \star_\varepsilon^+ \nabla_\varepsilon^- \tilde{h}^\varepsilon)(t,x). \tag{2.24b}$$

Here, we used summation by parts to move a discrete derivative on the kernel $\widetilde{\widetilde{G}}^\varepsilon$. We could move the discrete derivative $\nabla_\varepsilon^-$ from $\tilde{h}^\varepsilon$ to $\widetilde{\widetilde{G}}^\varepsilon$ too, but as we will see in Section 7.3 this way of writing will be more convenient.

Let us discuss the error terms appearing in equation (2.22) in more detail. The term $E^\varepsilon$ in (2.24a) is the same as in (2.11) and it appears because of the transformation of the spatial variable performed in (1.8). The error term (2.24b) contains the terms from the bracket process (2.13) which are not handled by renormalisation.

We can furthermore write the error term $\widetilde{E}_1^\varepsilon$ in a more explicit form. For this, we define the lattice (note that $\varrho_\varepsilon \neq 0$)

$$\bar{\Lambda}_\varepsilon^+ := |\varrho_\varepsilon|^{-1} \varepsilon^{\frac{3}{2}} \mathbf{N},$$

and a convolution on this lattice

$$(f \bar{\star}_\varepsilon g)(z) := |\varrho_\varepsilon|^{-1} \varepsilon^{5/2} \sum_{z' \in \bar{\Lambda}_\varepsilon^+ \times \Lambda_\varepsilon} f(z - z') g(z'). \tag{2.25}$$

Then the next result follows from the definitions (2.11)-(2.12) and (2.24a).

**Lemma 2.2** *One has*

$$\widetilde{E}_1^\varepsilon(\tilde{h}^\varepsilon)(t,x) = \varepsilon^{-\frac{1}{2}} \varrho_\varepsilon (\nabla_\varepsilon^- G^\varepsilon \bar{\star}_\varepsilon \tilde{h}^\varepsilon - \nabla_\varepsilon G^\varepsilon \star_\varepsilon^+ \tilde{h}^\varepsilon + \nabla_\varepsilon \widetilde{\widetilde{G}}^\varepsilon \star_\varepsilon^+ \tilde{h}^\varepsilon)(t,x) \tag{2.26}$$

*if* $\varrho_\varepsilon > 0$, *and*

$$\widetilde{E}_1^\varepsilon(\tilde{h}^\varepsilon)(t,x) = \varepsilon^{-\frac{1}{2}} \varrho_\varepsilon (\nabla_\varepsilon^+ G^\varepsilon \bar{\star}_\varepsilon \tilde{h}^\varepsilon - \nabla_\varepsilon G^\varepsilon \star_\varepsilon^+ \tilde{h}^\varepsilon + \nabla_\varepsilon \widetilde{\widetilde{G}}^\varepsilon \star_\varepsilon^+ \tilde{h}^\varepsilon)(t,x) \tag{2.27}$$

*if* $\varrho_\varepsilon < 0$.

We are going to prove convergence of $\tilde{h}^\varepsilon$ to the solution of the KPZ equation (1.11) using the semidiscrete equation (2.22). For this, we will show that the kernel $\widetilde{G}^\varepsilon$ behaves in the limit similarly to the discrete heat kernel $G^\varepsilon$, and converges in a suitable sense as $\varepsilon \to 0$ to the heat kernel. We will show furthermore, that each error term $\widetilde{E}_i^\varepsilon(\tilde{h}^\varepsilon)$ vanishes in the limit in a suitable topology. After proving that $\widetilde{M}^\varepsilon$ converges to a cylindrical Wiener process, we can formally see that (2.22) should converge to the mild form of the KPZ equation (1.11). Proving all these convergences rigorously is quite involved.

## 2.2 A priori bounds

The limiting KPZ equation (1.11) can be guessed from the discrete equation (2.22). However, to show this limit we need to control several objects, for example those which appear in the bracket process (2.3) of the martingales. In order to introduce a priori bounds, we are going to work with a "stopped" process (1.8).

In the proof of Theorem 1.1, we are going to show convergence of $\tilde{h}^\varepsilon(t)$ in a stronger topology than $\mathcal{C}(\mathbb{R})$. For this we need to control this process using the seminorm (1.21). More precisely, for a fixed constant $\mathfrak{m} \geq 1$ and the value $\alpha$ as in the statement of Theorem 1.1 we define the stopping time

$$\tau_{\varepsilon,\mathfrak{m}}^{(1)} := \inf\Big\{ t \geq 0 : \|\tilde{h}^\varepsilon(t)\|_{\mathcal{C}^\alpha}^{(\varepsilon)} \geq \mathfrak{m} \Big\}. \tag{2.28}$$

On the random time interval $[0, \tau_{\varepsilon,\mathfrak{m}}^{(1)})$ we have the bound $\|\tilde{h}^\varepsilon(t)\|_{\mathcal{C}^\alpha}^{(\varepsilon)} \leq \mathfrak{m}$, while on the closed interval $[0, \tau_{\varepsilon,\mathfrak{m}}^{(1)}]$ the bound is $\|\tilde{h}^\varepsilon(t)\|_{\mathcal{C}^\alpha}^{(\varepsilon)} \leq \mathfrak{m} + 2$ almost surely because the height function $\tilde{h}^\varepsilon(t)$ may make a jump of size $2\sqrt{\varepsilon}$ at time $t = \tau_{\varepsilon,\mathfrak{m}}^{(1)}$.

To control the bracket process (2.13) we need to control the product $\nabla_\varepsilon^- \tilde{h}^\varepsilon(t,x) \nabla_\varepsilon^+ \tilde{h}^\varepsilon(t,x)$. Since we are going to prove convergence of $\tilde{h}^\varepsilon$ to a function from $\mathcal{C}^\alpha$, one can expect that the renormalised product (2.17) behaves as a function of regularity $2(\alpha - 1)$. Then for fixed $\underline{\kappa} \in (0, 2(\alpha - \frac{2}{5}))$ we define the stopping time

$$\tau_{\varepsilon,\mathfrak{m}}^{(2)} := \inf\Big\{ t \geq 0 : \|(\nabla_\varepsilon^- \tilde{h}^\varepsilon \diamond \nabla_\varepsilon^+ \tilde{h}^\varepsilon)(t)\|_{\mathcal{C}^{2(\alpha-1)}}^{(\varepsilon)} \geq \mathfrak{m} \varepsilon^{-\underline{\kappa}} \Big\}. \tag{2.29}$$

Then on the random time interval $[0, \tau_{\varepsilon,\mathfrak{m}}^{(1)})$ we have the bound $\|(\nabla_\varepsilon^- \tilde{h}^\varepsilon \diamond \nabla_\varepsilon^+ \tilde{h}^\varepsilon)(t)\|_{\mathcal{C}^{2(\alpha-1)}}^{(\varepsilon)} \leq \mathfrak{m}\varepsilon^{-\underline{\kappa}}$. We expect to have a blow-up of this norm as $\varepsilon \to 0$ because $(\nabla_\varepsilon^- \tilde{h}^\varepsilon \diamond \nabla_\varepsilon^+ \tilde{h}^\varepsilon)(t, \bullet)$ does not converge to a function in the time variable and hence it cannot be bounded uniformly in time and $\varepsilon > 0$. The restriction $\underline{\kappa} < 2(\alpha - \frac{2}{5})$ is dictated by our analysis in Section 6.4.3, where we need $1 - 2\alpha + \underline{\kappa} < 2\kappa < \frac{1}{5}$.

To combine these two stopping times, we set

$$\tau_{\varepsilon,\mathfrak{m}} := \tau^{(1)}_{\varepsilon,\mathfrak{m}} \wedge \tau^{(2)}_{\varepsilon,\mathfrak{m}}, \tag{2.30}$$

and we restrict the time variable to the interval $[0, \tau_{\varepsilon,\mathfrak{m}}]$. For this we consider a stopped process $\sigma^\varepsilon(t)$, extended beyond the random time $\tau_{\varepsilon,\mathfrak{m}}$. To define such extension, we introduce a new process $\sigma^{\varepsilon,\mathfrak{m}}_{\text{sym}}$ which starts from the configuration $\sigma^{\varepsilon,\mathfrak{m}}_{\text{sym}}(\varepsilon^{-2}\tau_{\varepsilon,\mathfrak{m}}) = \sigma^\varepsilon(\varepsilon^{-2}\tau_{\varepsilon,\mathfrak{m}}-)$ and which for the times $t > \varepsilon^{-2}\tau_{\varepsilon,\mathfrak{m}}$ it evolves as a symmetric simple exclusion process, i.e., its infinitesimal generator $\mathcal{L}^\varepsilon_{\text{sym}}$ is given by (1.2) with constant jump rates $r^\varepsilon_{\text{sym}} = \frac{1}{2}\nu^\varepsilon$ to the neighboring sites, where we use the constant (2.8). We choose these values for the jump rates to have the same constant term in the bracket process of the martingales $\widetilde{M}^{\varepsilon,\mathfrak{m}}$ defined below. Then we set

$$\sigma^{\varepsilon,\mathfrak{m}}(t) := \begin{cases} \sigma^\varepsilon(t) & \text{for } t < \varepsilon^{-2}\tau_{\varepsilon,\mathfrak{m}}, \\ \sigma^{\varepsilon,\mathfrak{m}}_{\text{sym}}(t) & \text{for } t \geq \varepsilon^{-2}\tau_{\varepsilon,\mathfrak{m}}. \end{cases}$$

We note that the constants (1.13), defined via the average steady state current of particles, change their values after the stopping time $\varepsilon^{-2}\tau_{\varepsilon,\mathfrak{m}}$. Namely, the current (1.12) of the symmetric simple exclusion process equals 0 and the constants (1.13) vanish.

We define the function $h^{\varepsilon,\mathfrak{m}}_{\text{sym}}$ by (1.6) via the spin system $\sigma^{\varepsilon,\mathfrak{m}}_{\text{sym}}$, and we define the martingales $M^{\varepsilon,\mathfrak{m}}$ via the process $\sigma^{\varepsilon,\mathfrak{m}}$ in the same way as we defined $M^\varepsilon$ in (2.2) via the process $\sigma^\varepsilon$. For $t < \tau_{\varepsilon,\mathfrak{m}}$ the martingale $M^{\varepsilon,\mathfrak{m}}(t)$ coincides with $M^\varepsilon(t)$, while for $t \geq \tau_{\varepsilon,\mathfrak{m}}$ we denote $M^{\varepsilon,\mathfrak{m}}(t) = M^{\varepsilon,\mathfrak{m}}_{\text{sym}}$, where the latter has the predictable quadratic covariations (which are derived by analogy with (2.3))

$$\begin{aligned}
&\langle M^{\varepsilon,\mathfrak{m}}_{\text{sym}}(\bullet, x), M^{\varepsilon,\mathfrak{m}}_{\text{sym}}(\bullet, x')\rangle_t \\
&\quad = \varepsilon^{-2}\mathbf{1}_{x=x'} \int_{\tau_{\varepsilon,\mathfrak{m}}}^t \Bigg( r^\varepsilon_{\text{sym}} \frac{1-\sigma^{\varepsilon,\mathfrak{m}}_{\text{sym}}(\varepsilon^{-2}s, \varepsilon^{-1}x)}{2} \frac{1+\sigma^{\varepsilon,\mathfrak{m}}_{\text{sym}}(\varepsilon^{-2}s, \varepsilon^{-1}x+1)}{2} (2\sqrt{\varepsilon})^2 \\
&\qquad\qquad + r^\varepsilon_{\text{sym}} \frac{1+\sigma^{\varepsilon,\mathfrak{m}}_{\text{sym}}(\varepsilon^{-2}s, \varepsilon^{-1}x)}{2} \frac{1-\sigma^{\varepsilon,\mathfrak{m}}_{\text{sym}}(\varepsilon^{-2}s, \varepsilon^{-1}x+1)}{2} (2\sqrt{\varepsilon})^2 \Bigg) ds \\
&\quad = \varepsilon^{-1}\mathbf{1}_{x=x'} \nu^\varepsilon \int_{\tau_{\varepsilon,\mathfrak{m}}}^t (1 - \varepsilon \nabla^+_\varepsilon h^{\varepsilon,\mathfrak{m}}_{\text{sym}}(s,x) \nabla^+_\varepsilon h^{\varepsilon,\mathfrak{m}}_{\text{sym}}(s,x)) ds.
\end{aligned}$$

Then for $t \geq \tau_{\varepsilon,\mathfrak{m}}$ we have

$$\langle M^{\varepsilon,\mathfrak{m}}(\bullet, x), M^{\varepsilon,\mathfrak{m}}(\bullet, x')\rangle_t = \langle M^\varepsilon(\bullet, x), M^\varepsilon(\bullet, x')\rangle_{\tau_{\varepsilon,\mathfrak{m}}} + \langle M^{\varepsilon,\mathfrak{m}}_{\text{sym}}(\bullet, x), M^{\varepsilon,\mathfrak{m}}_{\text{sym}}(\bullet, x')\rangle_t.$$

We define as in (2.10) the martingales

$$\widetilde{M}^{\varepsilon,\mathfrak{m}}(t,x) = \int_0^t dM^{\varepsilon,\mathfrak{m}}(s, x + \lfloor \varrho_\varepsilon \varepsilon^{-\frac{1}{2}} s \rfloor_\varepsilon),$$

whose predictable quadratic variations are $\langle \widetilde{M}^{\varepsilon,\mathfrak{m}}(\bullet, x), \widetilde{M}^{\varepsilon,\mathfrak{m}}(\bullet, x')\rangle_t = 0$ for $x \neq x'$ and

$$d\langle \widetilde{M}^{\varepsilon,\mathfrak{m}}(\bullet, x)\rangle_t = \varepsilon^{-1}\widetilde{\mathbf{C}}_{\varepsilon,\mathfrak{m}}(t,x)dt, \tag{2.31}$$

where $\widetilde{\mathbf{C}}_{\varepsilon,\mathfrak{m}}(t,x) = \widetilde{\mathbf{C}}_\varepsilon(t,x)$ for $t < \tau_{\varepsilon,\mathfrak{m}}$ and

$$\widetilde{\mathbf{C}}_{\varepsilon,\mathfrak{m}}(t,x) = \nu^\varepsilon(1 - \varepsilon \nabla^-_\varepsilon \tilde{h}^{\varepsilon,\mathfrak{m}}_{\text{sym}}(t,x) \nabla^+_\varepsilon \tilde{h}^{\varepsilon,\mathfrak{m}}_{\text{sym}}(t,x)) \tag{2.32}$$

for $t \geq \tau_{\varepsilon,\mathfrak{m}}$, where $\tilde{h}^{\varepsilon,\mathfrak{m}}_{\text{sym}}(t,x) := h^{\varepsilon,\mathfrak{m}}_{\text{sym}}(t, x + \lfloor \varrho_\varepsilon \varepsilon^{-\frac{1}{2}} \tau_{\varepsilon,\mathfrak{m}} \rfloor_\varepsilon)$. We note that the constant term in the bracket processes (2.32) and (2.13) are the same and equal $\nu^\varepsilon$, which was the reason to introduce the symmetric exclusion process with such particular jump rates.

We define a version of the renormalised equation (2.22) but driven by the new martingale

$$\tilde h^{\varepsilon,\mathfrak{m}}(t,x) = (G^\varepsilon_t *_\varepsilon \tilde h^\varepsilon_0)(x) - \frac{1}{2}\widetilde{G}^\varepsilon \star^+_\varepsilon (\nabla^-_\varepsilon \tilde h^{\varepsilon,\mathfrak{m}} \diamond \nabla^+_\varepsilon \tilde h^{\varepsilon,\mathfrak{m}})(t,x) + \varrho_\varepsilon(\widetilde{\widetilde{G}}^\varepsilon \star^+_\varepsilon \nabla_\varepsilon \tilde h^{\varepsilon,\mathfrak{m}})(t,x) + \widetilde{E}^{\varepsilon,\mathfrak{m}}(\tilde h^{\varepsilon,\mathfrak{m}})(t,x) + (G^\varepsilon \star^+_\varepsilon d\widetilde{M}^{\varepsilon,\mathfrak{m}})(t,x), \tag{2.33}$$

where the error term is the same as in (2.23) and the renormalised product $\nabla^-_\varepsilon \tilde h^{\varepsilon,\mathfrak{m}} \diamond \nabla^+_\varepsilon \tilde h^{\varepsilon,\mathfrak{m}}$ is defined as in (2.17) (see also (2.15)) but via the renormalisation function

$$\widetilde{C}_{\varepsilon,\mathfrak{m}}(t,x) = \int_0^t \varepsilon \sum_{y\in\Lambda_\varepsilon} Q^\varepsilon_{t-s}(x-y)(\widetilde{\mathbf{C}}_{\varepsilon,\mathfrak{m}}(s,y) - \nu^\varepsilon)\, ds. \tag{2.34}$$

We clearly have $\tilde h^{\varepsilon,\mathfrak{m}}(t,x) = \tilde h^\varepsilon(t,x)$ for $t \in [0, \tau_{\varepsilon,\mathfrak{m}})$.

We are going to prove convergence of the function $\tilde h^{\varepsilon,\mathfrak{m}}$ to the solution of the KPZ equation (1.11). The desired limit for the functions (1.8) will be then obtained by taking $\mathfrak{m} \to \infty$ so that $\tau_{\varepsilon,\mathfrak{m}} \to \infty$ in probability.

## 3 Solution of the KPZ equation

A solution theory for the KPZ equation using regularity structures is developed in [FH20, Ch. 15]. In this section we recall some key concepts and notations of this solution theory in order to motivate the definition of the regularity structure for the discrete equation in Section 5, the renormalisation of the discrete model in Section 6 and the solution maps in Section 7.

It will be convenient to denote

$$\zeta := \sqrt{\nu}\xi \tag{3.1}$$

as the driving noise in (1.11), which will prevent us from writing the constant factor everywhere. We use here the constant $\nu$ defined in (2.9). The driving noise used in [FH20, Ch. 15] is the space-time noise $\xi$. The constant multiplier of the noise does not create any difficulties and yields respective constant multipliers of the renormalisation constants in (3.20).

### 3.1 Model space

Let us recall the construction of the model space $\mathcal{T}$ because we will use a more general construction for the discrete equation. Let $\Xi$ will represent the driving noise $\zeta$, let the integration map $\mathcal{I}$ represent the space-time convolution with the heat kernel, and let $\partial$ represent a spatial derivative. We denote $\mathcal{I}' := \partial\mathcal{I}$. Let the symbols $X_i$, $i = 0, 1$, will represent the time and space variables, and for $\ell = (\ell_0, \ell_1) \in \mathbf{N}^2$ we will use the shorthand $X^\ell = X_0^{\ell_0} X_1^{\ell_1}$, with the special unit symbol $\mathbf{1} := X^0$. We define the set of all monomials $\mathcal{W}_{\text{poly}} := \{X^\ell : \ell \in \mathbf{N}^2\}$. Then we define the minimal sets $\mathcal{V}, \mathcal{U}$ and $\mathcal{U}'$ of formal expressions such that $\Xi \in \mathcal{V}$, $\mathcal{W}_{\text{poly}} \subset \mathcal{V} \cap \mathcal{U} \cap \mathcal{U}'$ and the following implications hold:

$$\tau \in \mathcal{V} \quad \Rightarrow \quad \mathcal{I}(\tau) \in \mathcal{U}, \quad \mathcal{I}'(\tau) \in \mathcal{U}', \tag{3.2a}$$

$$\tau_1, \tau_2 \in \mathcal{U}' \quad \Rightarrow \quad \tau_1\tau_2 \in \mathcal{V}, \tag{3.2b}$$

where the product of symbols is commutative with the convention $\mathbf{1}\tau = \tau$. We postulate $\mathcal{I}(X^\ell) = \mathcal{I}'(X^\ell) = 0$ and do not include such zero elements into the sets.

We set $\mathscr{W} := \mathcal{U} \cup \mathcal{V} \cup \mathcal{U}'$, and for a fixed $\kappa \in (0, \frac{1}{10})$ we define the homogeneity $|\bullet| : \mathscr{W} \to \mathbb{R}$ of each element of $\mathscr{W}$ by the recurrent relations

$$|X^\ell| = 2\ell_0 + \ell_1 \quad \text{for} \quad \ell = (\ell_0, \ell_1) \in \mathbf{N}^2, \tag{3.3a}$$

$$|\Xi| = -\frac{3}{2} - \kappa, \tag{3.3b}$$

$$|\tau_1\tau_2| = |\tau_1| + |\tau_2|, \tag{3.3c}$$

$$|\mathcal{I}(\tau)| = |\tau| + 2, \quad |\mathcal{I}'(\tau)| = |\tau| + 1, \quad \tau \notin \mathcal{W}_{\text{poly}}. \tag{3.3d}$$

The restriction $\kappa < \frac{1}{10}$ allows us to use the minimal number of elements from $\mathscr{W}$ to solve the equation. We define $\mathscr{T}$ to contain all finite linear combinations of the elements in $\mathscr{W}$, and we view $\mathcal{I}$ and $\mathcal{I}'$ as linear maps $\tau \mapsto \mathcal{I}(\tau), \mathcal{I}'(\tau)$ on the respective subspaces of $\mathscr{T}$. We define the action of the map $\partial$ on monomials naturally as $\partial X^\ell = 0$ if $\ell_1 = 0$ and $\partial X^\ell = \ell_1 X^{\ell-(0,1)}$ if $\ell_1 \geq 1$. Then we extend $\partial$ linearly to the respective subspace of $\mathscr{T}$. The set $\mathfrak{A}$ contains the homogeneities $|\tau|$ for all $\tau \in \mathscr{W}$.

In order to solve equation (1.11), it is enough to consider the elements in $\mathscr{W}$ with negative homogeneities to describe the right-hand side, while the solution of this equation is described by the elements of homogeneities smaller than $\frac{3}{2}$. Hence, we define

$$\mathcal{W} := \{\tau \in \mathcal{V} : |\tau| \leq 0\} \cup \Big\{\tau \in \mathcal{U}' : |\tau| < \frac{1}{2}\Big\} \cup \Big\{\tau \in \mathcal{U} : |\tau| < \frac{3}{2}\Big\}. \tag{3.4}$$

We define $\mathcal{T}$ to be the linear span of the elements in $\mathcal{W}$, and the set $\mathcal{A}$ contains the homogeneities $|\tau|$ for all elements $\tau \in \mathcal{W}$.

It is convenient to represent the elements of $\mathcal{W}$ as trees. Namely, we denote $\Xi$ by a node . When a map $\mathcal{I}$ is applied to a symbol $\tau$, we draw a red edge from the root of the tree representing this symbol $\tau$. Respectively, the map $\mathcal{I}'$ is represented by a blue edge. For example, the symbols $\mathcal{I}(\Xi)$ and $\mathcal{I}'(\Xi)$ are represented by the diagrams  and  respectively. The product of symbols $\tau_1$ and $\tau_2$ is represented by the tree, obtained from the trees of these symbols by drawing them from the same root. For example, 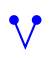 is the diagram for $\mathcal{I}'(\Xi)^2$. We use the symbols for the polynomials as before. Then $\mathcal{W}$ consist of 16 elements, which can be represented as symbols in increasing order of homogeneities

$\mathcal{W} :=$ 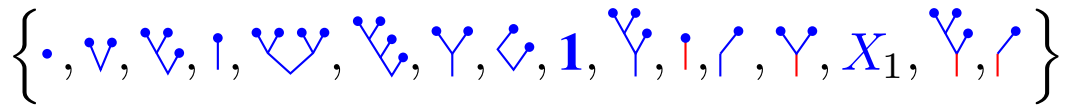 (3.5)

The homogeneities of these elements are provided in Table 1.

| Element | Homogeneity | Set |
|---|---|---|
|  | $-\frac{3}{2} - \kappa$ | $\mathcal{V}$ |
| 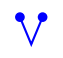 | $-1 - 2\kappa$ | $\mathcal{V}$ |
| 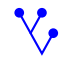 | $-\frac{1}{2} - 3\kappa$ | $\mathcal{V}$ |
|  | $-\frac{1}{2} - \kappa$ | $\mathcal{U}'$ |
| 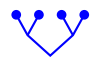 | $-4\kappa$ | $\mathcal{V}$ |
| 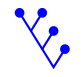 | $-4\kappa$ | $\mathcal{V}$ |
| 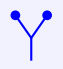 | $-2\kappa$ | $\mathcal{U}'$ |
| 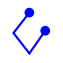 | $-2\kappa$ | $\mathcal{V}$ |

| Element | Homogeneity | Set |
|---|---|---|
| **1** | $0$ | $\mathcal{V} \cap \mathcal{U} \cap \mathcal{U}'$ |
| 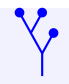 | $\frac{1}{2} - 3\kappa$ | $\mathcal{U}'$ |
|  | $\frac{1}{2} - \kappa$ | $\mathcal{U}$ |
| 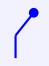 | $\frac{1}{2} - \kappa$ | $\mathcal{U}'$ |
| 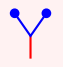 | $1 - 2\kappa$ | $\mathcal{U}$ |
| $X_1$ | $1$ | $\mathcal{U}$ |
| 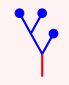 | $\frac{3}{2} - 3\kappa$ | $\mathcal{U}$ |
| 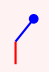 | $\frac{3}{2} - \kappa$ | $\mathcal{U}$ |

Table 1: The elements of $\mathcal{W}$ and their homogeneities. The third column indicates in which set the element is: $\mathcal{V}, \mathcal{U}$ or $\mathcal{U}'$. The backgrounds of the rows corresponding to the elements from the last two sets (except the element **1**) are red and blue respectively.

Every element $f \in \mathscr{T}$ can be uniquely written as $f = \sum_{\tau \in \mathscr{W}} f_\tau \tau$ for $f_\tau \in \mathbb{R}$, and we define

$$|f|_\alpha := \sum_{\tau \in \mathscr{W} : |\tau| = \alpha} |f_\tau|, \tag{3.6}$$

postulating $|f|_\alpha = 0$ if the sum runs over the empty set. We also introduce the projections

$$\mathcal{Q}_{<\alpha} f := \sum_{\tau \in \mathscr{W} : |\tau| < \alpha} f_\tau \tau, \qquad \mathcal{Q}_{\leq \alpha} f := \sum_{\tau \in \mathscr{W} : |\tau| \leq \alpha} f_\tau \tau \tag{3.7}$$

and

$$\mathrm{Proj}_\tau f := f_\tau. \tag{3.8}$$

Let the model space $\mathscr{T}_{<\alpha}$ contain all the elements $f \in \mathscr{T}$ satisfying $f = \mathcal{Q}_{<\alpha} f$. All these definitions can be immediately transferred to $\mathcal{T}$.

## 3.2 A structure group

In order to use the results of [Hai14], we need to define a structure group $\mathcal{G}$. For this, we need to introduce another set of basis elements

$$\mathcal{W}_+ := \left\{\mathbf{1}, \mathsf{Y}, \uparrow, \wr, \mathsf{Y}, X_1, \mathsf{Y}, \wr \right\}, \tag{3.9}$$

which is a subset of $\mathcal{W}$ containing the elements with non-negative homogeneities. Then we define $\mathcal{T}_+$ to be the free commutative algebra generated by the elements of $\mathcal{W}_+$.

We define a linear map $\Delta : \mathcal{T} \to \mathcal{T} \otimes \mathcal{T}_+$ by the identities

$$\Delta \mathbf{1} = \mathbf{1} \otimes \mathbf{1}, \qquad \Delta X_1 = X_1 \otimes \mathbf{1} + \mathbf{1} \otimes X_1, \qquad \Delta \Xi = \Xi \otimes \mathbf{1}, \tag{3.10a}$$

and then recursively by (we denote by $I$ the identity operator on $\mathcal{T}_+$)

$$\Delta(\tau_1\tau_2) = (\Delta\tau_1)(\Delta\tau_2), \tag{3.10b}$$
$$\Delta\mathcal{I}(\tau) = (\mathcal{I} \otimes I)\Delta\tau + \mathbf{1} \otimes \mathcal{I}(\tau), \qquad \text{if } |\tau| < -1, \tag{3.10c}$$
$$\Delta\mathcal{I}'(\tau) = (\mathcal{I}' \otimes I)\Delta\tau, \qquad \text{if } |\tau| < -1, \tag{3.10d}$$
$$\Delta\mathcal{I}(\tau) = (\mathcal{I} \otimes I)\Delta\tau + \mathbf{1} \otimes \mathcal{I}(\tau) + (\Delta X_1)(\mathbf{1} \otimes \mathcal{I}'(\tau)), \quad \text{if } |\tau| > -1, \tag{3.10e}$$
$$\Delta\mathcal{I}'(\tau) = (\mathcal{I}' \otimes I)\Delta\tau + \mathbf{1} \otimes \mathcal{I}'(\tau), \qquad \text{if } |\tau| > -1, \tag{3.10f}$$

for respective elements $\tau_i, \tau \in \mathcal{W}$. We provide $\Delta\tau$ for all $\tau \in \mathcal{W}$ in Table 2.

| Element $\tau$ | $\Delta\tau$ |
|---|---|
| $\cdot$ | $\cdot \otimes \mathbf{1}$ |
| $\vee$ | $\vee \otimes \mathbf{1}$ |
| $\curlyvee$ | $\curlyvee \otimes \mathbf{1}$ |
| $\uparrow$ | $\uparrow \otimes \mathbf{1}$ |
| $\vee\vee$ | $\vee\vee \otimes \mathbf{1}$ |
| $\curlyvee$ | $\curlyvee \otimes \mathbf{1} + \uparrow \otimes \mathsf{Y}$ |
| $\mathsf{Y}$ | $\mathsf{Y} \otimes \mathbf{1}$ |
| $\diamond$ | $\diamond \otimes \mathbf{1} + \uparrow \otimes \wr$ |

| Element $\tau$ | $\Delta\tau$ |
|---|---|
| $\mathbf{1}$ | $\mathbf{1} \otimes \mathbf{1}$ |
| $\mathsf{Y}$ | $\mathsf{Y} \otimes \mathbf{1} + \mathbf{1} \otimes \mathsf{Y}$ |
| $\uparrow$ | $\uparrow \otimes \mathbf{1} + \mathbf{1} \otimes \uparrow$ |
| $\wr$ | $\wr \otimes \mathbf{1} + \mathbf{1} \otimes \wr$ |
| $\mathsf{Y}$ | $\mathsf{Y} \otimes \mathbf{1} + \mathbf{1} \otimes \mathsf{Y}$ |
| $X_1$ | $X_1 \otimes \mathbf{1} + \mathbf{1} \otimes X_1$ |
| $\mathsf{Y}$ | $\mathsf{Y} \otimes \mathbf{1} + \mathbf{1} \otimes \mathsf{Y} + X_1 \otimes \mathsf{Y} + \mathbf{1} \otimes X_1 \mathsf{Y}$ |
| $\wr$ | $\wr \otimes \mathbf{1} + \mathbf{1} \otimes \wr + X_1 \otimes \wr + \mathbf{1} \otimes X_1 \wr$ |

Table 2: The image of the operator $\Delta$.

For any linear functional $f : \mathcal{T}_+ \to \mathbb{R}$ we define the map $\Gamma_f : \mathcal{T} \to \mathcal{T}$ as

$$\Gamma_f \tau := (I \otimes f)\Delta\tau. \tag{3.11}$$

Then the structure group $\mathcal{G}$ is defined as $\mathcal{G} := \{\Gamma_f : f \in \mathcal{G}_+\}$, where $\mathcal{G}_+$ contains all linear multiplicative functionals $f : \mathcal{T}_+ \to \mathbb{R}$ satisfying $f(\mathbf{1}) = 1$. Multiplicativity means $f(\tau\bar\tau) = f(\tau)f(\bar\tau)$ for $\tau, \bar\tau \in \mathcal{T}_+$.

Each functional $f \in \mathcal{G}_+$ is defined by its action on the elements of $\mathcal{W}_+$. Let us denote $a_i = f(\tau_i)$ for $i = 0, \ldots, 7$, where $\tau_i$ is the $i$-th element in (3.9) in the written order. In particular, $\tau_0 = \mathbf{1}$ and $a_0 = 1$. Then the linear map $\Gamma_f \in \mathcal{G}$ has the matrix in Figure 1 with respect to the basis (3.5).

In the rest of this section we use the framework of [Hai14] to work with the regularity structure $\mathscr{T} = (\mathcal{A}, \mathcal{T}, \mathcal{G})$ just introduced.

$$
\begin{bmatrix}
1&0&0&0&0&0&0&0&0&0&0&0&0&0&0&0\\
0&1&0&0&0&0&0&0&0&0&0&0&0&0&0&0\\
0&0&1&0&0&0&0&0&0&0&0&0&0&0&0&0\\
0&0&0&1&0&a_1&0&a_3&0&0&0&0&0&0&0&0\\
0&0&0&0&1&0&0&0&0&0&0&0&0&0&0&0\\
0&0&0&0&0&1&0&0&0&0&0&0&0&0&0&0\\
0&0&0&0&0&0&1&0&0&0&0&0&0&0&0&0\\
0&0&0&0&0&0&0&1&0&0&0&0&0&0&0&0\\
0&0&0&0&0&0&0&0&1&a_1&a_2&a_3&a_4&a_5&a_6+a_1a_5&a_7+a_3a_5\\
0&0&0&0&0&0&0&0&0&1&0&0&0&0&0&0\\
0&0&0&0&0&0&0&0&0&0&1&0&0&0&0&0\\
0&0&0&0&0&0&0&0&0&0&0&1&0&0&0&0\\
0&0&0&0&0&0&0&0&0&0&0&0&1&0&0&0\\
0&0&0&0&0&0&0&0&0&0&0&0&0&1&a_1&a_3\\
0&0&0&0&0&0&0&0&0&0&0&0&0&0&1&0\\
0&0&0&0&0&0&0&0&0&0&0&0&0&0&0&1
\end{bmatrix}
$$

Figure 1: The matrix of the linear map $\Gamma_f \in \mathcal{G}$ with respect to the basis (3.5).

### 3.3 Renormalisation group

We recall that the renormalisation group $\mathfrak{R}$ for the regularity structure $\mathscr{T}$ is defined as a character group of $\mathcal{T}_-$, which is a free unital algebra generated by

$$\mathcal{W}_- := \{ \dots \}. \tag{3.12}$$

We note that our set $\mathcal{W}_-$ is different from the more general set of generators of $\mathcal{T}_-$ in [FH20, Sec. 15.5.1], because the noise • and the planted trees and (which have only one edge incident to the root) do not get renormalised, and the tree with three leaves should be renormalised by a third cumulant of the noise (see Section 3.2 in [HS17]), which is zero due to its Gaussianity. Then we define the map $\Delta_- : \mathcal{T} \to \mathcal{T}_- \otimes \mathcal{T}$, as in [FH20, Eq. 15.25], which for every $\tau \in \mathcal{W}$ roughly speaking removes from the tree representing $\tau$ all subgraphs of negative homogeneities. More precisely, for $\tau = \mathbf{1}$ or $\tau = X_1$ we set $\Delta_- \tau = \mathbf{1} \otimes \tau$. For other elements $\tau \in \mathcal{W}$ let $\mathcal{N}(\tau)$ be the set containing $\mathbf{1}$ and all elements $A = \tau_1 \cdots \tau_n$ for $n \geq 1$, where $\tau_1, \ldots, \tau_n \in \mathcal{W}_-$, such that there is an injective map, mapping all vertices in $A$ to the vertices in $\tau$ and respecting connectivity. For such $A$, if $A$ contains the whole tree $\tau$ then we denote $\tau' = \mathbf{1}$. Otherwise, we denote by $\tau'$ the graph obtained from $\tau$ by contracting the subtrees $\tau_1, \ldots, \tau_n$ into non-bold vertices (recall that all trees representing the elements of $\mathcal{W}$ have two types of vertices: bold leaves and non-bold inner vertices). Then we postulate $\mathrm{Cut}_A(\tau) = \tau'$ if $\tau' \in \mathcal{W}$ and $\mathrm{Cut}_A(\tau) = 0$ otherwise. Moreover, we set $\mathrm{Cut}_{\mathbf{1}}(\tau) = \tau$ and define

$$\Delta_- \tau := \sum_{A \in \mathcal{N}(\tau)} A \otimes \mathrm{Cut}_A(\tau). \tag{3.13}$$

The action of this map on the elements of the set (3.5) is provided in Table 3. For every $g \in \mathfrak{R}$ we define a linear map $M_g : \mathcal{T} \to \mathcal{T}$, giving renormalisation of each element $\tau \in \mathcal{W}$ as

$$M_g \tau := (g \otimes I) \Delta_- \tau. \tag{3.14}$$

| Element $\tau$ | $\Delta_- \tau$ | Element $\tau$ | $\Delta_- \tau$ |
|---|---|---|---|
| $\bullet$ | $\mathbf{1} \otimes \bullet$ | $\mathbf{1}$ | $\mathbf{1} \otimes \mathbf{1}$ |
| V | $\mathbf{1} \otimes$ V $+$ V $\otimes \mathbf{1}$ | Ẏ | $\mathbf{1} \otimes$ Ẏ $+ 2$ ◇ $\otimes$ ↗ |
| V̇ | $\mathbf{1} \otimes$ V̇ $+ 2$ ◇ $\otimes$ ↑ | ↑ | $\mathbf{1} \otimes$ ↑ |
| ↑ | $\mathbf{1} \otimes$ ↑ | ↗ | $\mathbf{1} \otimes$ ↗ |
| VV | $\mathbf{1} \otimes$ VV $+$ VV $\otimes \mathbf{1}$ | Y | $\mathbf{1} \otimes$ Y |
| V̇̇ | $\mathbf{1} \otimes$ V̇̇ $+$ V̇̇ $\otimes \mathbf{1} + 2$ ◇ $\otimes$ ◇ $+$ ◇ $\otimes$ Y | $X_1$ | $\mathbf{1} \otimes X_1$ |
| Y | $\mathbf{1} \otimes$ Y | Ẏ | $\mathbf{1} \otimes$ Ẏ $+ 2$ ◇ $\otimes$ ↗ |
| ◇ | $\mathbf{1} \otimes$ ◇ $+$ ◇ $\otimes \mathbf{1}$ | ↗ | $\mathbf{1} \otimes$ ↗ |

Table 3: The image of the operator $\Delta_-$.

For $g \in \mathfrak{R}$, we have $g(\mathbf{1}) = 1$ and we denote

$$g(◇) = -C_0, \qquad g(V) = -C_1, \qquad g(VV) = -C_2, \qquad g(\dot{V}\dot{}) = -C_3. \tag{3.15}$$

Moreover, when we consider spatially symmetric approximations of the driving noise, the element ◇ does not need to be renormalised, i.e. $C_0 = 0$, and the renormalisation group becomes 3-dimensional. The map $M_g$ then satisfies

$$M_g V = V - C_1 \mathbf{1}, \qquad M_g VV = VV - C_2 \mathbf{1}, \qquad M_g \dot{V}\dot{} = \dot{V}\dot{} - C_3 \mathbf{1}$$

and it maps all other elements of $\mathcal{W}$ to themselves.

Given a smooth model $Z = (\Pi, \Gamma)$, the map $M_g$ allows to define a renormalised model $Z^g = (\Pi^g, \Gamma^g)$ on the regularity structure $\mathscr{T}$ (see [FH20, Thm. 15.21]) as

$$\Pi^g_z := \Pi_z M_g, \qquad \Gamma^g_{z\bar{z}} := M_g^{-1} \Gamma_{z\bar{z}} M_g. \tag{3.16}$$

### 3.4 Solution map

Let $G$ be the heat kernel on $\mathbb{R}_+ \times \mathbb{R}$, which is the solution of the PDE $\partial_t G = \frac{1}{2}\Delta G$ with the initial condition $G(0, \bullet) = \delta(\bullet)$. Then [Hai14, Lem. 7.7] allows to write $G = K + R$, where $K$ is a singular compactly supported function and $R$ is smooth. For a constant $\mu > 0$, we consider the equation

$$H = \mathcal{Q}_{<\gamma}\Big(Gh_0 + \mathscr{P}\mathbf{1}_+ \mathcal{Q}_{\le 0}\Big(-\frac{1}{2}(\partial H)^2 + \Xi\Big)\Big) \tag{3.17}$$

on a suitable space $H \in \mathcal{D}^{\gamma,\eta}$ of modelled distributions with respect to the model $Z_{\mathrm{KPZ}}$ defined in [Hai14, Thm. 3.10], where $Gh_0$ is a polynomial lift of the spatial convolution $Gh_0(t,x) := G_t * h_0(x)$ to the regularity structure (defined in [Hai14, Lem. 7.5]), $\mathscr{P}$ is the abstract integration operator

$$\mathscr{P} := \mathscr{K}_\kappa + \mathscr{R}_{\frac{3}{2}+2\kappa} \mathcal{R}_{\mathrm{KPZ}}, \tag{3.18}$$

where the operator $\mathscr{K}_\kappa$ is defined in [Hai14, Eq. 5.15] via the kernel $K$ for the values $\beta = 2$ and $\gamma = \kappa$, the operator $\mathscr{R}_{\frac{3}{2}+2\kappa}$ is defined in [Hai14, Eq. 7.7] as a Taylor's expansion of the function $R$ up to the order $\frac{3}{2} + 2\kappa$, and $\mathcal{R}_{\mathrm{KPZ}}$ is the reconstruction map for the model $Z_{\mathrm{KPZ}}$. We take the value of $\kappa$ to be the same as in (1.11). Furthermore, $\mathbf{1}_+$ is the projection of modelled distributions to positive times.

We also consider a mollified noise $\zeta^\delta = \varrho_\delta \star \zeta$, where $\zeta$ is defined in (3.1) and where the mollifier $\varrho_\delta$ is defined for $\delta > 0$ as

$$\varrho_\delta(t,x) := \frac{1}{\delta^3} \varrho\Big(\frac{t}{\delta^2}, \frac{x}{\delta}\Big), \tag{3.19}$$

with $\varrho : \mathbb{R}^2 \to \mathbb{R}$ being a symmetric, smooth, supported on the ball of radius 1 (with respect to the parabolic distance $\| \bullet \|_{\mathfrak{s}}$) and satisfying $\int_{\mathbb{R}^2} \varrho(z)dz = 1$. This noise can be canonically lifted to a smooth model $Z^\delta$ (see [FH20]). In this case the renormalisation constants (3.15) are given by (see Section 15.5.3 in [FH20] and the definition of the noise in (3.1)) $C_0^\delta = 0$,

$$\begin{aligned} C_1^\delta &= \nu \int_{\mathbb{R}^2} (\partial_x K^\delta)^2 dz \quad \sim \quad \delta^{-1}, \\ C_2^\delta &= 2\nu^2 \int_{\mathbb{R}^2} \partial_x K(z) \partial_x K(\bar z) Q_\delta(z - \bar z)^2 dz d\bar z \quad \sim \quad \log \delta, \\ C_3^\delta &= 2\nu^2 \int_{\mathbb{R}^2} \partial_x K(z) \partial_x K(\bar z) Q_\delta(\bar z) Q_\delta(z + \bar z) dz d\bar z \quad \sim \quad \log \delta, \end{aligned} \tag{3.20}$$

where $K^\delta := K \star \varrho_\delta$ and

$$Q_\delta(z) = \int_{\mathbb{R}^2} \partial_x K^\delta(\tilde z) \partial_x K^\delta(\tilde z - z) d\tilde z.$$

The derivative $\partial_x$ in these formulas is with respect to the spatial variable. An important observation is

$$C_2^\delta + 4 C_3^\delta \quad \sim \quad 1, \tag{3.21}$$

i.e. the logarithmic divergences cancel out. Due to this cancelation there is no logarithmically divergent part in the renormalisation constant of the KPZ equation (3.22) below. These renormalisation constants define a character $g_\delta$ by (3.15) and this character allows to renormalise the smooth model $Z^\delta$ as in (3.16), which yields a renormalised model $Z^\delta_{\mathrm{KPZ}} = (\Pi^\delta, \Gamma^\delta)$. Moreover, $\lim_{\delta \to 0} Z^\delta_{\mathrm{KPZ}} = Z_{\mathrm{KPZ}}$ in probability in the metric of models.

Let $H^\delta$ be the solution of equation (3.17), defined with respect to a 2-periodic initial condition $h_0^\delta$ and the renormalised model $Z^\delta_{\mathrm{KPZ}}$. Then one can see from [FH20, Prop. 15.12] that the process $h^\delta = \mathcal{R}^\delta H^\delta$, where $\mathcal{R}^\delta$ is the reconstruction map for the model $Z^\delta_{\mathrm{KPZ}}$ from [Hai14, Thm. 3.10], is the classical solution of the SPDE

$$\partial_t h^\delta = \frac{1}{2} \Delta h^\delta - \frac{1}{2}((\partial_x h^\delta)^2 - C^\delta) + \zeta^\delta, \qquad h^\delta(0, \bullet) = h_0^\delta(\bullet). \tag{3.22}$$

The renormalisation constant $C^\delta \sim \delta^{-1}$ is such that the solution of (3.22) converges as $\delta \to 0$ in a suitable space of functions. This constant may be written explicitly as

$$C^\delta = C_1^\delta + \frac{1}{4}(C_2^\delta + 4 C_3^\delta). \tag{3.23}$$

The relation (3.21) implies that there is no logarithmic sub-divergence in $C^\delta$. More precisely, we have the following result:

**Theorem 3.1** *Equation* (3.17) *has a local in time solution* $H \in \mathcal{D}^{\gamma,\eta}(Z_{\mathrm{KPZ}})$ *with* $\gamma = \frac{3}{2} + 2\kappa$ *and* $\eta = \alpha$, *where* $\alpha$ *is as in the statement of Theorem 1.1. The solution map* $H = \mathcal{S}(h_0, Z_{\mathrm{KPZ}})$ *is locally Lipschitz continuous with respect to the initial state* $h_0 \in \mathcal{C}^\eta$ *and the model* $Z_{\mathrm{KPZ}}$.

*Then the solution of the KPZ equation* (1.11) *is defined as the reconstruction* $h = \mathcal{R}_{\mathrm{KPZ}} H$. *There is a constant* $c_0$ *such that* $h_{\mathrm{CH}}(t,x) = h(t,x) - c_0 t$ *coincides with the Cole-Hopf solution of the KPZ equation [BG97], and hence it almost surely exists globally in time.*

*Let* $h^\delta$ *be the solution of* (3.22) *with an initial state satisfying* $\lim_{\delta \to 0} \|h_0 - h_0^\delta\|_{\mathcal{C}^\alpha} = 0$. *Then for any* $T > 0$ *and* $\eta > 0$ *one has*

$$\lim_{\delta \to 0} \mathbb{P}\left( \sup_{t \in [0,T]} \|h(t) - h^\delta(t)\|_{\mathcal{C}^\alpha} \geq \eta \right) = 0. \tag{3.24}$$

**Remark 3.2** The renormalisation of the non-linearity in the KPZ equation depends on the decomposition of the heat kernel $G = K + R$ (see the definition of the renormalisation constants (3.20)), and for a different decomposition $G = \bar{K} + \bar{R}$ the theory of regularity structures would yield a different solution $\bar{h}$. The two solutions are different from each other by a linear shift, i.e. $h(t,x) - \bar{h}(t,x) = \bar{c}_0 t$ for a constant $\bar{c}_0$ depending on the decompositions of the kernel. Hence, if we want to get a Cole-Hopf solution, we need to perform a respective shift $h(t,x) - c_0 t$. The precise value of $c_0$ is usually hard to determine and we are not going to do it.

*Proof of Theorem 3.1.* Existence of a local solution $H \in \mathcal{D}^{\gamma,\eta}$ and the local Lipschitz continuity of the solution map follow from Theorem 7.8, Proposition 7.11 and Corollary 7.12 in [Hai14]. The fact that we get the Cole-Hopf solution after a linear shift is explained in Remark 3.2. The convergence in probability (3.24) follows from continuity of the solution map and the convergences of the initial conditions and the models. □

# 4 Properties of discrete operators and kernels

The goal of this section is to prove some properties of the discrete operators and kernels involved in equation (2.22). Throughout this paper we are going to use the Fourier transform

$$(\mathcal{F}g)(\omega) := \int_{-\infty}^{\infty} g(x) e^{-2\pi i \omega x} dx \tag{4.1}$$

of an $L^1$ function or a tempered distribution $g$ on $\mathbb{R}$. The inverse Fourier transform is

$$g(x) = \mathcal{F}^{-1}(\mathcal{F}g)(x) = \int_{-\infty}^{\infty} (\mathcal{F}g)(\omega) e^{2\pi i \omega x} d\omega.$$

## 4.1 Properties of the finite difference derivatives

The finite difference derivatives (2.1) have the Fourier multipliers

$$\mathcal{F}(\nabla_\varepsilon^\pm g)(\omega) = i\omega\, m^\pm(\varepsilon\omega)(\mathcal{F}g)(\omega), \qquad \mathcal{F}(\Delta_\varepsilon g)(k) = -\omega^2 f(\varepsilon\omega)(\mathcal{F}g)(\omega), \tag{4.2}$$

with the Fourier multipliers

$$m^+(x) := \frac{e^{2\pi i x} - 1}{ix}, \qquad m^-(x) := \frac{1 - e^{-2\pi i x}}{ix}, \qquad f(x) := \frac{2(1 - \cos(2\pi x))}{x^2}, \tag{4.3}$$

with the definition at $x = 0$ by the L'Hopital's rule. The following properties of these Fourier multipliers can be used throughout the article.

**Lemma 4.1** *1. The function $f$ is smooth and satisfies $16 \leq f(x) \leq 4\pi^2$ for all $x \in [-\frac{1}{2}, \frac{1}{2}]$. Moreover, for any integer $n \geq 1$ and all $x \in \mathbb{R}$ we have*

$$|f^{(n)}(x)| \leq \begin{cases} \frac{2^{n+3}\pi^{n+2}}{(n+1)(n+2)} \cosh(2\pi x) & \text{if } n \text{ is even}, \\ \frac{2^{n+3}\pi^{n+2}}{(n+1)(n+2)} \sinh(2\pi |x|) & \text{if } n \text{ is odd}. \end{cases} \tag{4.4}$$

2. *The functions $m^\pm$ are smooth and satisfy $|m^\pm(x)| \leq 2\pi$ for all $x \in [-\frac{1}{2}, \frac{1}{2}]$. Moreover, for any integer $n \geq 1$ and all $x \in \mathbb{R}$ we have*

$$|(m^\pm)^{(n)}(x)| \leq \frac{(2\pi)^{n+1}}{n+1} \sqrt{\cosh(4\pi x)}. \tag{4.5}$$

3. *The following identities hold: $m^\pm(-x) = m^\mp(x)$ and*

$$\overline{m^\pm(x)} = m^\mp(x), \qquad m^\pm(x) m^\mp(x) = f(x), \qquad m^\pm(x)^2 = e^{\pm 2\pi i x} f(x). \tag{4.6}$$

*Proof.* All properties, except the bounds (4.4) and (4.5), are straightforward. The Maclaurin series of $\cos(2\pi x)$ yields

$$f^{(n)}(x) = 8\pi^2(2\pi)^n \sum_{m\geq n/2} \frac{(-1)^m(2\pi x)^{2m-n}}{(2m+2)!}\frac{(2m)!}{(2m-n)!}.$$

Estimating this series absolutely and using $\frac{(2m)!}{(2m+2)!} \leq \frac{1}{(n+1)(n+2)}$ we get (4.4). Similarly, we can write

$$(m^+)^{(n)}(x) = 2\pi(2\pi i)^n \sum_{m\geq 0} \frac{(2\pi i x)^m}{m!}\frac{1}{m+n+1}.$$

The real and imaginary parts of the series are absolutely bounded by $\frac{\cosh(4\pi|x|)}{n+1}$ and $\frac{\sinh(4\pi|x|)}{n+1}$ respectively. This yields the bound (4.5) for $(m^+)^{(n)}$, and the bound for $(m^-)^{(n)}$ follows from the relation $m^-(x) = m^+(-x)$. □

## 4.2 The discrete heat kernel $G^\varepsilon$

The Green function $G^\varepsilon$ of the linear operator $\partial_t - \frac{1}{2}\Delta_\varepsilon$ is the solution of

$$\frac{d}{dt}G^\varepsilon(t,x) = \frac{1}{2}\Delta_\varepsilon G^\varepsilon(t,x), \qquad G^\varepsilon(0,x) = \varepsilon^{-1}\mathbf{1}_{x=0}, \tag{4.7}$$

on $\mathbb{R}_+ \times \Lambda_\varepsilon$. One way to solve this equation is by computing its Fourier transform, which yields

$$G^\varepsilon(t,x) = \int_{-\frac{1}{2\varepsilon}}^{\frac{1}{2\varepsilon}} e^{-f_\varepsilon(\omega)t}e^{2\pi i\omega x}d\omega, \tag{4.8}$$

where

$$f_\varepsilon(\omega) := \frac{1}{2}\omega^2 f(\varepsilon\omega) = \varepsilon^{-2}(1-\cos(2\pi\varepsilon\omega)), \tag{4.9}$$

and we used the function $f$ defined in (4.3). We are going to use this formula later when proving bounds on the kernel. Another useful formula is

$$G^\varepsilon(t,x) = \varepsilon^{-1}q^\circ(\varepsilon^{-2}t, \varepsilon^{-1}x), \tag{4.10}$$

where $q^\circ$ is the solution of (4.7) with $\varepsilon = 1$, i.e., $q^\circ$ is a transition kernel of a continuous-time random walk on $\mathbb{Z}$ with symmetric and independent nearest-neighbors jumps. The latter is a discrete heat kernel whose definition and properties can be found in [Bar17]. In particular, $q^\circ(t,x)$ has an exponential decay in the spatial variable $x$ and it vanishes as $t\to\infty$.

Sometimes we will use a periodised discrete kernel, defined on $\mathbb{T}_\varepsilon := \{\varepsilon k : k \in \mathbb{T}_N\}$ as

$$P^\varepsilon(t,x) := \sum_{m\in\mathbb{Z}} G^\varepsilon(t, x+2m). \tag{4.11}$$

This kernel $P^\varepsilon(t,x)$ can be expressed as the Fourier series

$$P^\varepsilon(t,x) = \frac{1}{2}\sum_{k\in\mathbb{T}_N} e^{-f_\varepsilon(k)t}e^{\pi i k x}. \tag{4.12}$$

It is important for the theory of regularity structures to extend the discrete kernel $G^\varepsilon$ off the lattice, so that the extension, which we denote by $G^\varepsilon_{\text{ext}}$, is sufficiently regular and its derivatives satisfy bounds similar to the bounds on the continuous heat kernel. We need this extension because we test the discrete models in (5.15a) with respect to sufficiently regular functions, and having such extension of $G^\varepsilon$, the kernel can be treated as a test function. Then we can get an analogue of the Schauder estimates. We define in Appendix A such extensions for quite general discrete kernels, and as an application we define $G^\varepsilon_{\text{ext}}$ in Proposition A.9.

The following result provides the main properties of the kernel $G^\varepsilon_{\text{ext}}$ which are necessary in the framework of regularity structures.

**Proposition 4.2** *For any integer $r \geq 2$ we can write $G^\varepsilon_{\text{ext}} = K^\varepsilon + R^\varepsilon$, where*

1. *The function $R^\varepsilon$ satisfies $R^\varepsilon(t,x) = 0$ if $t < 0$. Moreover, for every $j \in \mathbf{N}^2$ such that $|j|_\mathfrak{s} \leq r$ there is a constant $C = C(r) > 0$ such that*

$$|D^j R^\varepsilon(z)| \leq C, \qquad \int_\mathbb{R} |D^j R^\varepsilon(t,x)|\,dx \leq C, \tag{4.13}$$

*uniformly over $z \in \mathbb{R}^2$ and $t \in \mathbb{R}$.*

2. *$K^\varepsilon$ satisfies $K^\varepsilon(t,x) = 0$ if $t < 0$, is even in the spatial variable, and can be written as $K^\varepsilon = \sum_{n=0}^M K^{\varepsilon,n}$ with $M = -\lfloor \log_2 \varepsilon \rfloor$, where the functions $\{K^{\varepsilon,n}\}_{n=0}^M$ are defined on $\mathbb{R}^2$ and have the following properties:*

   *(a) There is a constant $c \geq 1$ such that for every $0 \leq n \leq M$, the function $K^{\varepsilon,n}$ is compactly supported on $\{z \in \mathbb{R}^2 : \|z\|_\mathfrak{s} \leq c2^{-n}\}$.*

   *(b) There is a constant $C$ independent of $\varepsilon$ such that*

$$|D^j K^{\varepsilon,n}(z)| \leq C 2^{n(1+|j|_\mathfrak{s})}, \tag{4.14}$$

   *uniformly in $z \in \mathbb{R}^2$, $j \in \mathbf{N}^2$ with $|j|_\mathfrak{s} \leq r$ and $0 \leq n < M$. For $n = M$ and $j \neq 0$ the bound (4.14) holds only for $z \in \mathbb{R}_+ \times \mathbb{R}$.*

   *(c) For all $0 \leq n < M$ and $j \in \mathbf{N}^2$ with $|j|_\mathfrak{s} \leq r$, we have that*

$$\int_{D_\varepsilon} D^j K^{\varepsilon,n}(z)dz = 0$$

   *For $n = M$ the identity only holds for $j = 0$.*

*Proof.* This decomposition follows from the proof of [HM18, Lem. 5.4] and the bounds provided in Proposition A.9. The only difference in the argument is that we do not restrict the domain of the function $R^\varepsilon$. □

Having the derived bounds on the discrete heat kernel, we can show that the terms discussed in Remark 2.1 and which are not included into the renormalisation in (2.15) have a very good bound, in conjunction with the uniform bound (2.4) on the quadratic variation.

**Lemma 4.3** *There exists a constant $C > 0$ such that for any $T > 0$ one has*

$$\int_0^T \varepsilon \sum_{x \in \Lambda_\varepsilon} \sum_{m \in \mathbb{Z}\setminus\{0\}} |\nabla_\varepsilon^- G^\varepsilon(t,x) \nabla_\varepsilon^+ G^\varepsilon(t,x+2m)|dt \leq C(1+\sqrt{T}) \tag{4.15}$$

*uniformly in $\varepsilon \in (0,1)$.*

*Proof.* We write the left-hand side of (4.15) as

$$\begin{aligned}
&\sum_{m \in \mathbb{Z}\setminus\{0\}} \int_{[0,T]\times\Lambda_\varepsilon : \|z\|_\mathfrak{s} \leq |m|} |\nabla_\varepsilon^- G^\varepsilon(z) \nabla_\varepsilon^+ G^\varepsilon(z+(0,2m))|dz \\
&\quad + \sum_{m \in \mathbb{Z}\setminus\{0\}} \int_{[0,T]\times\Lambda_\varepsilon : |m| < \|z\|_\mathfrak{s} \leq 3|m|} |\nabla_\varepsilon^- G^\varepsilon(z) \nabla_\varepsilon^+ G^\varepsilon(z+(0,2m))|dz \\
&\quad + \sum_{m \in \mathbb{Z}\setminus\{0\}} \int_{[0,T]\times\Lambda_\varepsilon : \|z\|_\mathfrak{s} > 3|m|} |\nabla_\varepsilon^- G^\varepsilon(z) \nabla_\varepsilon^+ G^\varepsilon(z+(0,2m))|dz.
\end{aligned} \tag{4.16}$$

We bound the first term by

$$\sum_{m\in\mathbb{Z}\setminus\{0\}}\left(\int_{[0,T]\times\Lambda_\varepsilon:\|z\|_\mathfrak{s}\le|m|}|\nabla_\varepsilon^- G^\varepsilon(z)|dz\right)\left(\sup_{[0,T]\times\Lambda_\varepsilon:\|z\|_\mathfrak{s}\le|m|}|\nabla_\varepsilon^+ G^\varepsilon(z+(0,2m))|\right). \tag{4.17}$$

By Lemma A.8, the integral is bounded by a constant times $\int_{\|z\|_\mathfrak{s}\le|m|}(\|z\|_\mathfrak{s}+\varepsilon)^{-2}dz \lesssim |m|$. Using Lemma A.8 again, we bound

$$|\nabla_\varepsilon^+ G^\varepsilon(z+(0,2m))| \lesssim \frac{(\sqrt{t}+\varepsilon)^{n-2}}{(\|(t,x+2m)\|_\mathfrak{s}+\varepsilon)^n} \lesssim \frac{(\sqrt{t}+1)^{n-2}}{|m|^n}, \tag{4.18}$$

where $z=(t,x)$, and where the last bound holds for $z$ and $m$ as in the first term in (4.16). Then the first term in (4.16) is estimated by a constant multiple of

$$\sum_{m\in\mathbb{Z}\setminus\{0\}}\frac{(\sqrt{T}+1)^{n-2}}{|m|^{n-1}} \lesssim (\sqrt{T}+1)^{n-2},$$

for any $n\ge 3$. Taking $n=3$ we get a bound of order $\sqrt{T}+1$.

Now, we will bound the second term in (4.16). We use Lemma A.8 to bound

$$|\nabla_\varepsilon^- G^\varepsilon(z)| \lesssim \frac{(\sqrt{t}+\varepsilon)^{\bar{n}-2}}{(\|z\|_\mathfrak{s}+\varepsilon)^{\bar{n}}} \lesssim \frac{(\sqrt{t}+\varepsilon)^{\bar{n}-2}}{|m|^{\bar{n}}}, \tag{4.19}$$

for $z$ and $m$ as in the second term in (4.16). Applying Lemma A.8 again, we bound $|\nabla_\varepsilon^+ G^\varepsilon(z+(0,2m))| \lesssim (\sqrt{t}+\varepsilon)^{-2}$. Then the second term in (4.16) is bounded by a constant multiple of

$$\sum_{m\in\mathbb{Z}\setminus\{0\}}\frac{1}{|m|^{\bar{n}}}\int_{[0,T]\times\Lambda_\varepsilon:|m|<\|z\|_\mathfrak{s}\le 3|m|}(\sqrt{t}+\varepsilon)^{\bar{n}-4}dz.$$

Taking $\bar{n}=5$, this expression is of the order $\sum_{m\in\mathbb{Z}\setminus\{0\}}\frac{1}{|m|^2}(\sqrt{T}+1) \lesssim \sqrt{T}+1$.

Finally, we will bound the last term in (4.16). The estimate (4.18) holds in this case. Then, using the intermediate estimate in (4.19), we bound the last term in (4.16) by a constant times

$$\sum_{m\in\mathbb{Z}\setminus\{0\}}\frac{1}{|m|^n}\int_{[0,T]\times\Lambda_\varepsilon:\|z\|_\mathfrak{s}>3|m|}\frac{(\sqrt{t}+1)^{\bar{n}+n-4}}{(\|z\|_\mathfrak{s}+\varepsilon)^{\bar{n}}}dz \lesssim \sum_{m\in\mathbb{Z}\setminus\{0\}}\frac{(\sqrt{T}+1)^{\bar{n}+n-4}}{|m|^{\bar{n}+n-3}},$$

where the last bound holds for $\bar{n}\ge 4$. Taking $\bar{n}=4$ and $n=1$, the last expression becomes of the order $\sqrt{T}+1$. □

### 4.3 Properties of the kernel $Q^\varepsilon$

The following two lemmas state some useful properties of the kernel $Q^\varepsilon$ defined in (2.16), which follow readily from [BG97, Lems. A.1 & A.2]. We prefer however to provide a complete proof of Lemma 4.4 to demonstrate how the formula follows only from the discrete equation (4.7), without using the Fourier transform.

As explained in Section 1.2, these properties of the function $Q^\varepsilon$ play a crucial role in our analysis of the discrete equation (2.33). More precisely, the identity from Lemma 4.4 guarantees that the renormalisation constant $C_\varepsilon(\vee)$, defined in (6.16), does not diverge (see Lemma 6.5). If that was not the case, we would have to subtract another $\varepsilon^{-1}$-divergent constant in (1.7) to renormalise the non-linearity in equation (2.33). Moreover, identity (4.20) is used in the proof of the bound on the kernel in Lemma 4.13 which guarantees that the error terms in (2.33) containing $\widetilde{\widetilde{G}}^\varepsilon$ vanish in the limit.

Lemma 4.5 guarantees that the convolution operator (2.20) is defined. We use this lemma in Section 4.4 to prove bounds on this operator.

**Lemma 4.4** *The function* (2.16) *satisfies*

$$\int_{\mathbb{R}_+ \times \Lambda_\varepsilon} Q^\varepsilon(z) dz = 0. \tag{4.20}$$

*Proof.* Using the definition (2.16) and applying summation by parts, we get

$$\sum_{x \in \Lambda_\varepsilon} Q^\varepsilon(t,x) = -\sum_{x \in \Lambda_\varepsilon} \Delta_\varepsilon G^\varepsilon_t(x) G^\varepsilon_t(x+\varepsilon) = -\sum_{x \in \Lambda_\varepsilon} G^\varepsilon_t(x) \Delta_\varepsilon G^\varepsilon_t(x+\varepsilon).$$

Adding the last two expressions and using (4.7) and the chain rule yields

$$\sum_{x \in \Lambda_\varepsilon} Q^\varepsilon(t,x) = -\sum_{x \in \Lambda_\varepsilon} \frac{d}{dt}(G^\varepsilon_t(x) G^\varepsilon_t(x+\varepsilon)).$$

Using the Gaussian decay of $G^\varepsilon$ and its time derivative, we can swap the sum over $x$ and the integral over $t$. Then the left-hand side of (4.20) equals

$$-\sum_{x \in \Lambda_\varepsilon} \int_0^\infty \frac{d}{dt}(G^\varepsilon_t(x) G^\varepsilon_t(x+\varepsilon)) dt = \sum_{x \in \Lambda_\varepsilon} G^\varepsilon_0(x) G^\varepsilon_0(x+\varepsilon),$$

where we used $\lim_{t\to\infty} G^\varepsilon_t(x) = 0$. The initial condition in (4.7) implies that the last sum vanishes. □

**Lemma 4.5** *There exists some* $\Theta \in (0,1)$ *such that for any* $\varepsilon > 0$ *the kernel* (2.16) *satisfies*

$$\int_{\mathbb{R}_+ \times \Lambda_\varepsilon} |Q^\varepsilon(z)| dz \leq \varepsilon^{-1} \Theta. \tag{4.21}$$

*Proof.* Let us denote by $\nabla^\pm$ the operators $\nabla^\pm_\varepsilon$ with $\varepsilon = 1$, and recall the identity (4.10). Then

$$\int_0^\infty \varepsilon \sum_{x \in \Lambda_\varepsilon} |Q^\varepsilon(t,x)| dt = \varepsilon^{-1} \int_0^\infty \sum_{x \in \mathbb{Z}} |\nabla^- q^\circ(t,x) \nabla^+ q^\circ(t,x)| dt$$

and the claim then follow from [BG97, Lem. A.2]. □

Lemma A.8 implies the decay bounds for the kernel $Q^\varepsilon$ defined in (2.16).

**Lemma 4.6** *For any* $j \in \mathbf{N}^2$ *there is a constant* $C = C(j) > 0$ *such that*

$$|D^j_\varepsilon Q^\varepsilon(z)| \leq C(\|z\|_\mathfrak{s} + \varepsilon)^{-|j|_\mathfrak{s} - 4}, \tag{4.22}$$

*uniformly over* $z \in \mathbb{R}_+ \times \Lambda_\varepsilon$.

*Proof.* It is straightforward to check the following identity for the discrete derivative:

$$(\nabla^-_\varepsilon)^n (fg)(x) = \sum_{m=0}^n \binom{n}{m} (\nabla^-_\varepsilon)^{n-m} f(x - \varepsilon m) (\nabla^-_\varepsilon)^m g(x). \tag{4.23}$$

Combining it with the Leibniz rule, we write the derivative of (2.16) as

$$D^j_\varepsilon Q^\varepsilon(t,x) = \sum_{0 \leq k \leq j} \binom{j}{k} D^{j-k}_\varepsilon (\nabla^-_\varepsilon G^\varepsilon(t, x - \varepsilon k_1)) D^k_\varepsilon (\nabla^+_\varepsilon G^\varepsilon(t,x)), \tag{4.24}$$

where $k = (k_0, k_1) \in \mathbf{N}^2$. We have $\nabla_\varepsilon^+ G^\varepsilon(t,x) = \nabla_\varepsilon^- G^\varepsilon(t,x+\varepsilon)$ and by Lemma A.8 we get

$$\begin{aligned}|D_\varepsilon^{j-k}(\nabla_\varepsilon^- G^\varepsilon(t,x-\varepsilon k_1))| &= |D_\varepsilon^{j-k+(0,1)} G^\varepsilon(t,x-\varepsilon k_1)| \\ &\leq C(\|(t,x-\varepsilon k_1)\|_{\mathfrak{s}} + \varepsilon)^{-|j-k|_{\mathfrak{s}}-2}, \\ |D_\varepsilon^k(\nabla_\varepsilon^+ G^\varepsilon(t,x))| = |D_\varepsilon^{k+(0,1)} G^\varepsilon(t,x+\varepsilon)| &\leq C(\|(t,x+\varepsilon)\|_{\mathfrak{s}} + \varepsilon)^{-|k|_{\mathfrak{s}}-2}.\end{aligned}$$

Furthermore, for $0 \leq k \leq j$ we have $(\|(t,x-\varepsilon k_1)\|_{\mathfrak{s}} + \varepsilon)^{-|j-k|_{\mathfrak{s}}-2} \lesssim (\|(t,x)\|_{\mathfrak{s}} + \varepsilon)^{-|j-k|_{\mathfrak{s}}-2}$ with a proportionality constant depending on $j$. Then the bound (4.22) follows from (4.24). □

Lemma 4.5 yields the bound $\|Q^\varepsilon \star_\varepsilon^+ f^\varepsilon\|_{L^\infty(\mathbb{R}_+\times\Lambda_\varepsilon)} \lesssim \varepsilon^{-1}\|f^\varepsilon\|_{L^\infty(\mathbb{R}_+\times\Lambda_\varepsilon)}$ for a function $f^\varepsilon$ on the space-time domain. If however this function has some Hölder regularity, this bound cam be improved, what is demonstrated in Lemma 4.7 below. More precisely, we use functions $f^\varepsilon : \mathbb{R}\times\Lambda_\varepsilon \to \mathbb{R}$, such that the following quantity is bounded uniformly in $\varepsilon \in (0,1)$:

$$\|f^\varepsilon\|_{\mathcal{C}^\alpha_{\mathfrak{s},\varepsilon}} := \sup_{z\in\mathbb{R}\times\Lambda_\varepsilon} |f^\varepsilon(z)| + \sup_{\substack{z,\bar{z}\in\mathbb{R}\times\Lambda_\varepsilon \\ \|z-\bar{z}\|_{\mathfrak{s}}>0}} \frac{|f^\varepsilon(z) - f^\varepsilon(\bar{z})|}{\|z-\bar{z}\|_{\mathfrak{s}}^\alpha}. \tag{4.25}$$

In this case we write $f^\varepsilon \in \mathcal{C}^\alpha_{\mathfrak{s},\varepsilon}$.

**Lemma 4.7** *Let $f^\varepsilon \in \mathcal{C}^\alpha_{\mathfrak{s},\varepsilon}$ for some $\alpha \in [0,1)$. Then*

$$\|Q^\varepsilon \star_\varepsilon^+ f^\varepsilon\|_{L^\infty(\mathbb{R}_+\times\Lambda_\varepsilon)} \leq C\varepsilon^{\alpha-1}\|f^\varepsilon\|_{\mathcal{C}^\alpha_{\mathfrak{s},\varepsilon}}, \tag{4.26}$$

*with a constant $C > 0$ independent of $\varepsilon$ and $f^\varepsilon$.*

*Proof.* Using Lemma 4.4, we write

$$(Q^\varepsilon \star_\varepsilon^+ f^\varepsilon)(z) = \int_{\mathbb{R}_+\times\Lambda_\varepsilon} Q^\varepsilon(\bar{z})(f^\varepsilon(z-\bar{z}) - f^\varepsilon(z))d\bar{z}.$$

Using the regularity assumption for $f^\varepsilon$ and Lemma 4.6, we bound

$$|(Q^\varepsilon \star_\varepsilon^+ f^\varepsilon)(z)| \lesssim \|f^\varepsilon\|_{\mathcal{C}^\alpha_{\mathfrak{s},\varepsilon}} \int_{D_\varepsilon^+} (\|\bar{z}\|_{\mathfrak{s}} + \varepsilon)^{-4}\|\bar{z}\|_{\mathfrak{s}}^\alpha d\bar{z}.$$

We split the integration domain into $\|\bar{z}\|_{\mathfrak{s}} \leq \varepsilon$ and $\|\bar{z}\|_{\mathfrak{s}} > \varepsilon$. In the first case the integral is simply bounded by a constant multiple of $\varepsilon^{\alpha-1}$, while in the second case the integral is bounded by

$$\int_{\|\bar{z}\|_{\mathfrak{s}}>\varepsilon} \|\bar{z}\|_{\mathfrak{s}}^{\alpha-4} d\bar{z} \lesssim \varepsilon^{\alpha-1}$$

as desired. □

### 4.4 Properties of the kernel of $(\mathscr{Q}^\varepsilon)^{-1}$

We can show that the operator $(\mathscr{Q}^\varepsilon)^{-1}$, defined in (2.20), converges to the identity.

**Lemma 4.8** *Let $f^\varepsilon \in \mathcal{C}^\alpha_{\mathfrak{s},\varepsilon}$ for some $\alpha \in (0,1)$, where we use the space defined below* (4.25). *Then there is $\varepsilon_0 \in (0,1)$ such that*

$$\|(\mathscr{Q}^\varepsilon)^{-1} f^\varepsilon - f^\varepsilon\|_{L^\infty(\mathbb{R}_+\times\Lambda_\varepsilon)} \leq C\varepsilon^\alpha\|f^\varepsilon\|_{\mathcal{C}^\alpha_{\mathfrak{s},\varepsilon}}, \tag{4.27}$$

*for all $\varepsilon \in (0,\varepsilon_0)$ with a constant $C > 0$ independent of $\varepsilon$ and $f^\varepsilon$.*

*Proof.* The definition (2.20) yields

$$\|(\mathscr{Q}^\varepsilon)^{-1} f - f\|_{L^\infty(\mathbb{R}_+\times\Lambda_\varepsilon)} \leq \sum_{n=1}^\infty \varepsilon^n (1+\sqrt{\varepsilon})^n \|Q^\varepsilon\|^{n-1}_{L^1(\mathbb{R}_+\times\Lambda_\varepsilon)} \|Q^\varepsilon \star^+_\varepsilon f^\varepsilon\|_{L^\infty(\mathbb{R}_+\times\Lambda_\varepsilon)}.$$

Lemmas 4.5 and 4.7 allow to bound this expression by a constant times

$$\varepsilon^\alpha \|f^\varepsilon\|_{\mathcal{C}^\alpha_{\mathfrak{s},\varepsilon}} \sum_{n=1}^\infty (1+\sqrt{\varepsilon})^n \Theta^{n-1}.$$

Recalling that $\Theta < 1$, we choose $\varepsilon_0 \in (0,1)$ such that $(1+\sqrt{\varepsilon_0})\Theta < 1$. Then for any $\varepsilon \in (0,\varepsilon_0)$ the preceding geometric series converges and we get the desired bound. □

The definition (2.20) implies that the operator $(\mathscr{Q}^\varepsilon)^{-1} - \mathrm{Id}$ has an integral kernel on $\mathbb{R}_+ \times \Lambda_\varepsilon$, which is a function. We will need to use a decay of this kernel.

**Lemma 4.9** *For any $j \in \mathbf{N}^2$ there is a constant $C = C(j) > 0$ and $\varepsilon_0 \in (0,1)$ such that*

$$|D^j_\varepsilon((\mathscr{Q}^\varepsilon)^{-1} - \mathrm{Id})(z)| \leq C\varepsilon(\|z\|_{\mathfrak{s}} + \varepsilon)^{-|j|_{\mathfrak{s}}-4}, \tag{4.28}$$

*uniformly over $z \in \mathbb{R}_+ \times \Lambda_\varepsilon$ and $\varepsilon \in (0,\varepsilon_0)$.*

*Proof.* The definition (2.20) yields the recursion on $\mathbb{R}_+ \times \Lambda_\varepsilon$:

$$(\mathscr{Q}^\varepsilon)^{-1} = \mathrm{Id} - \varepsilon(1+\sqrt{\varepsilon}) Q^\varepsilon \star^+_\varepsilon (\mathscr{Q}^\varepsilon)^{-1}.$$

Denoting by $S^\varepsilon$ the kernel of the operator $(\mathscr{Q}^\varepsilon)^{-1} - \mathrm{Id}$, we get

$$S^\varepsilon(z) = -\varepsilon(1+\sqrt{\varepsilon})(Q^\varepsilon \star^+_\varepsilon S^\varepsilon)(z) - \varepsilon(1+\sqrt{\varepsilon}) Q^\varepsilon(z). \tag{4.29}$$

Then

$$\begin{aligned}
&\sup_{z\in\mathbb{R}_+\times\Lambda_\varepsilon} \Big\{ |D^j_\varepsilon S^\varepsilon(z)| \left(\|z\|_{\mathfrak{s}} + \varepsilon\right)^{4+|j|_{\mathfrak{s}}} \Big\} \\
&\quad \leq \varepsilon(1+\sqrt{\varepsilon}) \sup_{z\in\mathbb{R}_+\times\Lambda_\varepsilon} \Big\{ |D^j_\varepsilon (Q^\varepsilon \star^+_\varepsilon S^\varepsilon)(z)| \left(\|z\|_{\mathfrak{s}} + \varepsilon\right)^{4+|j|_{\mathfrak{s}}} \Big\} \\
&\qquad + \varepsilon(1+\sqrt{\varepsilon}) \sup_{z\in\mathbb{R}_+\times\Lambda_\varepsilon} \Big\{ |D^j_\varepsilon Q^\varepsilon(z)| \left(\|z\|_{\mathfrak{s}} + \varepsilon\right)^{4+|j|_{\mathfrak{s}}} \Big\},
\end{aligned} \tag{4.30}$$

where the last supremum is bounded by a constant $C$ by applying Lemma 4.6.

Now, we are going to bound the first supremum on the right-hand side of (4.30). We clearly have $D^j_\varepsilon(Q^\varepsilon \star^+_\varepsilon S^\varepsilon)(0) = 0$, and we need to bound this convolution evaluated at $z \neq 0$. For this we fix $z \neq 0$. Let $\varphi : \mathbb{R}^2 \to [0,1]$ be a smooth function such that $\varphi(y) = 0$ for $\|y\|_{\mathfrak{s}} \geq 1$ and $\varphi(y) = 1$ for $\|y\|_{\mathfrak{s}} \leq \frac{1}{2}$. We set $r = \delta\|z\|_{\mathfrak{s}} + \eta\varepsilon$ for some $\delta \in (0,\frac{1}{2})$ and $\eta \in [0,\delta]$ whose values will be specified below. We define $\varphi_r(t',x') = \varphi(r^{-2}t', r^{-1}x')$ and write

$$\begin{aligned}
(Q^\varepsilon \star^+_\varepsilon S^\varepsilon)(z) &= \int_{\mathbb{R}_+\times\Lambda_\varepsilon} \varphi_r(\bar z) Q^\varepsilon(\bar z) S^\varepsilon(z-\bar z)d\bar z + \int_{\mathbb{R}_+\times\Lambda_\varepsilon} \varphi_r(\bar z) Q^\varepsilon(z-\bar z) S^\varepsilon(\bar z)d\bar z \\
&\quad + \int_{\mathbb{R}_+\times\Lambda_\varepsilon} (1 - \varphi_r(\bar z) - \varphi_r(z-\bar z)) Q^\varepsilon(z-\bar z) S^\varepsilon(\bar z)d\bar z,
\end{aligned}$$

where we changed the integration variable $\bar z$ in the first integral to $z - \bar z$. For any $j \in \mathbf{N}^2$, we have

$$D^j_\varepsilon(Q^\varepsilon \star^+_\varepsilon S^\varepsilon)(z) = \int_{\mathbb{R}_+\times\Lambda_\varepsilon} \varphi_r(\bar z) Q^\varepsilon(\bar z) D^j_\varepsilon S^\varepsilon(z-\bar z)d\bar z + \int_{\mathbb{R}_+\times\Lambda_\varepsilon} \varphi_r(\bar z) D^j_\varepsilon Q^\varepsilon(z-\bar z) S^\varepsilon(\bar z)d\bar z$$

$$+\int_{\mathbb{R}_+\times\Lambda_\varepsilon}(1-\varphi_r(\bar z)-\varphi_r(z-\bar z))D^j_\varepsilon Q^\varepsilon(z-\bar z)S^\varepsilon(\bar z)d\bar z$$
$$-\sum_{0<k\le j}\binom{j}{k}\int_{\mathbb{R}_+\times\Lambda_\varepsilon}D^k_\varepsilon\varphi_r(z-\bar z)D^{j-k}_\varepsilon Q^\varepsilon(z-\bar z-(0,\varepsilon k_1))S^\varepsilon(\bar z)d\bar z,$$

where in the last term we used the Leibniz rule in the time variable and (4.23) in the spatial variable. The sum runs over $k=(k_0,k_1)\in\mathbf{N}^2$. Then we have

$$\begin{aligned}(\|z\|_{\mathfrak{s}}+\varepsilon)^{4+|j|_{\mathfrak{s}}}|D^j_\varepsilon(Q^\varepsilon\star^+_\varepsilon S^\varepsilon)(z)|&\le\int_{\mathbb{R}_+\times\Lambda_\varepsilon}\varphi_r(\bar z)|Q^\varepsilon(\bar z)|\,|D^j_\varepsilon S^\varepsilon(z-\bar z)|\,(\|z\|_{\mathfrak{s}}+\varepsilon)^{4+|j|_{\mathfrak{s}}}d\bar z\\&+\int_{\mathbb{R}_+\times\Lambda_\varepsilon}\varphi_r(\bar z)|D^j_\varepsilon Q^\varepsilon(z-\bar z)|\,(\|z\|_{\mathfrak{s}}+\varepsilon)^{4+|j|_{\mathfrak{s}}}|S^\varepsilon(\bar z)|d\bar z\\&+\int_{\mathbb{R}_+\times\Lambda_\varepsilon}|1-\varphi_r(\bar z)-\varphi_r(z-\bar z)|\,|D^j_\varepsilon Q^\varepsilon(z-\bar z)|\,(\|z\|_{\mathfrak{s}}+\varepsilon)^{4+|j|_{\mathfrak{s}}}|S^\varepsilon(\bar z)|d\bar z\\&+\sum_{0<k\le j}\binom{j}{k}\int_{\mathbb{R}_+\times\Lambda_\varepsilon}|D^k_\varepsilon\varphi_r(z-\bar z)|\,|D^{j-k}_\varepsilon Q^\varepsilon(z-\bar z-(0,\varepsilon k_1))|\,(\|z\|_{\mathfrak{s}}+\varepsilon)^{4+|j|_{\mathfrak{s}}}\,|S^\varepsilon(\bar z)|d\bar z,\end{aligned}\tag{4.31}$$

and we will analyse each of these terms. The definition of the function $\varphi_r$ implies that $\varphi_r(\bar z)\neq 0$ only for $\|\bar z\|_{\mathfrak{s}}\le r=\delta\|z\|_{\mathfrak{s}}+\eta\varepsilon$. Then in the first term in (4.31) we have $\|z-\bar z\|_{\mathfrak{s}}\ge(1-\delta)\|z\|_{\mathfrak{s}}-\eta\varepsilon$, which implies $\|z-\bar z\|_{\mathfrak{s}}+\varepsilon\ge(1-\delta)(\|z\|_{\mathfrak{s}}+\varepsilon)+(\delta-\eta)\varepsilon\ge(1-\delta)(\|z\|_{\mathfrak{s}}+\varepsilon)$ because $\delta\ge\eta$. Then the whole term is bounded by

$$\begin{aligned}&\sup_{\substack{\bar z\in\mathbb{R}_+\times\Lambda_\varepsilon\\ \|z-\bar z\|_{\mathfrak{s}}+\varepsilon\ge(1-\delta)(\|z\|_{\mathfrak{s}}+\varepsilon)}}\Big\{|D^j_\varepsilon S^\varepsilon(z-\bar z)|\,(\|z\|_{\mathfrak{s}}+\varepsilon)^{4+|j|_{\mathfrak{s}}}\Big\}\int_{\mathbb{R}_+\times\Lambda_\varepsilon}|Q^\varepsilon(\bar z)|d\bar z\\&\qquad\le(1-\delta)^{-4-|j|_{\mathfrak{s}}}\sup_{z'\in\mathbb{R}_+\times\Lambda_\varepsilon}\Big\{|D^j_\varepsilon S^\varepsilon(z')|\,(\|z'\|_{\mathfrak{s}}+\varepsilon)^{4+|j|_{\mathfrak{s}}}\Big\}\int_{\mathbb{R}_+\times\Lambda_\varepsilon}|Q^\varepsilon(\bar z)|d\bar z,\end{aligned}$$

where we changed the variable in the supremum to $z'=z-\bar z$ and estimated $\|z\|_{\mathfrak{s}}+\varepsilon\le\frac{1}{1-\delta}(\|z'\|_{\mathfrak{s}}+\varepsilon)$. Using Lemma 4.5, this expression is bounded by

$$\varepsilon^{-1}\Theta(1-\delta)^{-4-|j|_{\mathfrak{s}}}\sup_{z'\in\mathbb{R}_+\times\Lambda_\varepsilon}\Big\{|D^j_\varepsilon S^\varepsilon(z')|\,(\|z'\|_{\mathfrak{s}}+\varepsilon)^{4+|j|_{\mathfrak{s}}}\Big\}.\tag{4.32}$$

The second term in (4.31) is bounded in exactly the same way by

$$(1-\delta)^{-4-|j|_{\mathfrak{s}}}\sup_{z'\in\mathbb{R}_+\times\Lambda_\varepsilon}\Big\{|D^j_\varepsilon Q^\varepsilon(z')|\,(\|z'\|_{\mathfrak{s}}+\varepsilon)^{4+|j|_{\mathfrak{s}}}\Big\}\int_{\mathbb{R}_+\times\Lambda_\varepsilon}|S^\varepsilon(\bar z)|d\bar z.\tag{4.33}$$

By Lemma 4.6, the supremum is bounded by a constant $C>0$. Moreover, using the definition (2.20) and Lemma 4.5, we bound

$$\int_{\mathbb{R}_+\times\Lambda_\varepsilon}|S^\varepsilon(\bar z)|d\bar z\le\sum_{n=1}^\infty\varepsilon^n(1+\sqrt{\varepsilon})^n\|Q^\varepsilon\|^n_{L^1(\mathbb{R}_+\times\Lambda_\varepsilon)}\le\sum_{n=1}^\infty(1+\sqrt{\varepsilon})^n\Theta^n.$$

Taking as usual $\varepsilon<\varepsilon_0$ so that $(1+\sqrt{\varepsilon_0})\Theta<1$, the preceding expression is bounded by $\frac{(1+\sqrt{\varepsilon_0})\Theta}{1-\Theta}$. Then (4.33) is bounded by

$$C(1-\delta)^{-4-|j|_{\mathfrak{s}}}\frac{(1+\sqrt{\varepsilon_0})\Theta}{1-\Theta}.\tag{4.34}$$

The third term in (4.31) can be bounded similarly to the second one, due to the cutoff. Namely, if we have $\|z-\bar z\|_{\mathfrak{s}}<\frac{r}{2}$ in the integral, then $\|\bar z\|_{\mathfrak{s}}\ge\|z\|_{\mathfrak{s}}-\|z-\bar z\|_{\mathfrak{s}}\ge\frac{r}{\delta}-\frac{r}{2}$, where we

used the definition of $r$. Hence, $\|\bar z\|_{\mathfrak s} > r$ for $\delta < \frac{2}{3}$ which implies $1 - \varphi_r(\bar z) - \varphi_r(z-\bar z) = 0$. It means that we have $\|z-\bar z\|_{\mathfrak s} \geq \frac{r}{2} = \frac{\delta}{2}\|z\|_{\mathfrak s} + \frac{\eta}{2}\varepsilon$ in the third term in (4.31), which implies $\|z-\bar z\|_{\mathfrak s} + \varepsilon \geq \frac{\delta}{2}(\|z\|_{\mathfrak s}+\varepsilon)$. Then we bound this term in the same way as we bounded the second term by

$$
\begin{aligned}
2 &\sup_{\substack{\bar z \in \mathbb{R}_+\times\Lambda_\varepsilon \\ \|z-\bar z\|_{\mathfrak s}+\varepsilon \geq \frac{\delta}{2}(\|z\|_{\mathfrak s}+\varepsilon)}} \Big\{ |D_\varepsilon^j Q^\varepsilon(z-\bar z)| \, (\|z\|_{\mathfrak s}+\varepsilon)^{4+|j|_{\mathfrak s}} \Big\} \int_{\mathbb{R}_+\times\Lambda_\varepsilon} |S^\varepsilon(\bar z)| d\bar z \\
&\leq 2(2/\delta)^{4+|j|_{\mathfrak s}} \sup_{z' \in \mathbb{R}_+\times\Lambda_\varepsilon} \Big\{ |D_\varepsilon^j Q^\varepsilon(z')| \, (\|z'\|_{\mathfrak s}+\varepsilon)^{4+|j|_{\mathfrak s}} \Big\} \int_{\mathbb{R}_+\times\Lambda_\varepsilon} |S^\varepsilon(\bar z)| d\bar z.
\end{aligned}
$$

We bound this expression as we did for (4.33) by

$$
2C(2/\delta)^{4+|j|_{\mathfrak s}} \frac{(1+\sqrt{\varepsilon_0})\Theta}{1-\Theta}. \tag{4.35}
$$

Now, we turn to the last term in (4.31). From the definition of the function $\varphi_r$ we have $\varphi_r(z-\bar z) = 1$ for $\|z-\bar z\|_{\mathfrak s} \leq \frac{r}{2}$. Hence, $D_\varepsilon^k \varphi_r(z-\bar z) = 0$ for $\|z-\bar z\|_{\mathfrak s} \leq \frac{r}{2} - \varepsilon j_1$ because $k \leq j$. Thus, we have $\|z-\bar z\|_{\mathfrak s} \geq \frac{r}{2} - \varepsilon j_1$ in the last integral in (4.31). This implies

$$
\|z-\bar z-(0,\varepsilon k_1)\|_{\mathfrak s} + (1+2j_1)\varepsilon \geq \|z-\bar z\|_{\mathfrak s} + (1+j_1)\varepsilon \geq \frac{\delta}{2}(\|z\|_{\mathfrak s}+\varepsilon). \tag{4.36}
$$

Note also that $|D_\varepsilon^k \varphi_r(z-\bar z)| \lesssim r^{-|k|_s}$ uniformly over $z-\bar z$. Thus, the last integral in (4.31) is bounded by a constant times

$$
\begin{aligned}
&\sup_{\substack{\bar z \in \mathbb{R}_+\times\Lambda_\varepsilon \\ \|z-\bar z-(0,\varepsilon k_1)\|_{\mathfrak s}+(1+2j_1)\varepsilon \geq \frac{\delta}{2}(\|z\|_{\mathfrak s}+\varepsilon)}} \Big\{ |D_\varepsilon^{j-k} Q^\varepsilon(z-\bar z-(0,\varepsilon k_1))| \, (\|z\|_{\mathfrak s}+\varepsilon)^{4+|j-k|_{\mathfrak s}} \Big\} \int_{\mathbb{R}_+\times\Lambda_\varepsilon} |S^\varepsilon(\bar z)| d\bar z \\
&\leq (2/\delta)^{4+|j|_{\mathfrak s}} \sup_{z' \in \mathbb{R}_+\times\Lambda_\varepsilon} \Big\{ |D_\varepsilon^{j-k} Q^\varepsilon(z')| \, (\|z'\|_{\mathfrak s}+(1+2j_1)\varepsilon)^{4+|j-k|_{\mathfrak s}} \Big\} \int_{\mathbb{R}_+\times\Lambda_\varepsilon} |S^\varepsilon(\bar z)| d\bar z,
\end{aligned}
$$

which is bounded by (4.35) but probably with a different constant multiplier.

Combining (4.30) with the derived bounds on the terms in (4.31), we obtain

$$
\begin{aligned}
&\sup_{z \in \mathbb{R}_+\times\Lambda_\varepsilon} \Big\{ |D_\varepsilon^j S^\varepsilon(z)| \, (\|z\|_{\mathfrak s}+\varepsilon)^{4+|j|_{\mathfrak s}} \Big\} \\
&\qquad \leq (1+\sqrt{\varepsilon})\Theta(1-\delta)^{-4-|j|_{\mathfrak s}} \sup_{z' \in \mathbb{R}_+\times\Lambda_\varepsilon} \Big\{ |D_\varepsilon^j S^\varepsilon(z')| \, (\|z'\|_{\mathfrak s}+\varepsilon)^{4+|j|_{\mathfrak s}} \Big\} \\
&\qquad\quad + \varepsilon(1+\sqrt{\varepsilon})C(1-\delta)^{-4-|j|_{\mathfrak s}} \frac{(1+\sqrt{\varepsilon_0})\Theta}{1-\Theta} \\
&\qquad\quad + \varepsilon(1+\sqrt{\varepsilon})\tilde C(2/\delta)^{4+|j|_{\mathfrak s}} \frac{(1+\sqrt{\varepsilon_0})\Theta}{1-\Theta} + \varepsilon C,
\end{aligned}
$$

for some constant $\tilde C > 0$. Since $\Theta \in (0,1)$, we can take $\varepsilon$ and $\delta$ sufficiently small such that $(1+\sqrt{\varepsilon})\Theta(1-\delta)^{-4-|j|_{\mathfrak s}} < 1$. Then the supremum on the right-hand side can be moved to the left, which yields the bound (4.28). □

## 4.5 Extensions of the other discrete kernels

Using the definition (2.21), we aim to deduce the decay of $\widetilde G^\varepsilon$ from Lemmas A.8 and 4.9.

**Lemma 4.10** *For any $j \in \mathbf{N}^2$ and $\alpha \in [0,1]$ there is $\varepsilon_0 \in (0,1)$ and a constant $C = C(j,\alpha) > 0$ such that*

$$
|D_\varepsilon^j(\widetilde G^\varepsilon - G^\varepsilon)(z)| \leq C\varepsilon^\alpha(\|z\|_{\mathfrak s}+\varepsilon)^{-|j|_{\mathfrak s}-1-\alpha}, \tag{4.37}
$$

*uniformly over $z \in \mathbb{R}_+ \times \Lambda_\varepsilon$ and $\varepsilon \in (0, \varepsilon_0)$. Moreover,*

$$|D_\varepsilon^j \widetilde{G}^\varepsilon(z)| \leq C(\|z\|_\mathfrak{s} + \varepsilon)^{-|j|_\mathfrak{s} - 1}, \tag{4.38}$$

*uniformly over $z \in \mathbb{R}_+ \times \Lambda_\varepsilon$ and $\varepsilon \in (0, \varepsilon_0)$.*

*Proof.* The bound (4.38) follows from (4.37) and Lemma A.8. Hence, we need to prove (4.37).

Recall from the proof of Lemma 4.9 that we denoted by $S^\varepsilon$ the kernel of $(\mathscr{Q}^\varepsilon)^{-1} - \mathrm{Id}$, and hence the definition (2.21) yields

$$\widetilde{G}^\varepsilon(z) - G^\varepsilon(z) = (G^\varepsilon \star_\varepsilon^+ S^\varepsilon)(z).$$

Let $\varphi : \mathbb{R}^2 \to [0, 1]$ be a smooth function such that $\varphi(\bar z) = 0$ for $\|\bar z\|_\mathfrak{s} \geq 1$ and $\varphi(\bar z) = 1$ for $\|\bar z\|_\mathfrak{s} \leq \frac{1}{2}$. For any fixed $z \in \mathbb{R}_+ \times \Lambda_\varepsilon$ we let $r = \frac{1}{2}(\|z\|_\mathfrak{s} + \varepsilon)$ and $\varphi_r(s, y) = \varphi(r^{-2}s, r^{-1}y)$. Then we can write

$$\begin{aligned}(G^\varepsilon \star_\varepsilon^+ S^\varepsilon)(z) = \int_{\mathbb{R}_+ \times \Lambda_\varepsilon} \varphi_r(\bar z)\, G^\varepsilon(\bar z)\, S^\varepsilon(z - \bar z) d\bar z + \int_{\mathbb{R}_+ \times \Lambda_\varepsilon} \varphi_r(\bar z)\, G^\varepsilon(z - \bar z)\, S^\varepsilon(\bar z) d\bar z \\ + \int_{\mathbb{R}_+ \times \Lambda_\varepsilon} (1 - \varphi_r(\bar z) - \varphi_r(z - \bar z))\, G^\varepsilon(z - \bar z)\, S^\varepsilon(\bar z) d\bar z.\end{aligned}$$

For any $j \in \mathbf{N}^2$, we have

$$\begin{aligned}D_\varepsilon^j (G^\varepsilon \star_\varepsilon^+ S^\varepsilon)(z) = \int_{\mathbb{R}_+ \times \Lambda_\varepsilon} \varphi_r(\bar z)\, G^\varepsilon(\bar z)\, D_\varepsilon^j S^\varepsilon(z - \bar z) d\bar z + \int_{\mathbb{R}_+ \times \Lambda_\varepsilon} \varphi_r(\bar z)\, D_\varepsilon^j G^\varepsilon(z - \bar z)\, S^\varepsilon(\bar z) d\bar z \\ + \int_{\mathbb{R}_+ \times \Lambda_\varepsilon} (1 - \varphi_r(\bar z) - \varphi_r(z - \bar z))\, D_\varepsilon^j G^\varepsilon(z - \bar z)\, S^\varepsilon(\bar z) d\bar z \\ - \sum_{0 < k \leq j} \binom{j}{k} \int_{\mathbb{R}_+ \times \Lambda_\varepsilon} D_\varepsilon^k \varphi_r(z - \bar z)\, D_\varepsilon^{j-k} G^\varepsilon(z - \bar z - (0, \varepsilon k_1))\, S^\varepsilon(\bar z) d\bar z,\end{aligned} \tag{4.39}$$

where in the last term we used the Leibniz rule and (4.23). Thus we need to absolutely bound the integrals on the right-hans side of (4.39).

In the first term in (4.39), due to the cutoff we have $\|\bar z\|_\mathfrak{s} \leq r$ which implies $\|z - \bar z\|_\mathfrak{s} + \varepsilon \geq \|z\|_\mathfrak{s} + \varepsilon - r = \frac{1}{2}(\|z\|_\mathfrak{s} + \varepsilon)$. Then by Lemmas 4.9 and A.8 this terms is absolutely bounded by a constant multiple of

$$\begin{aligned}\varepsilon(\|z\|_\mathfrak{s} + \varepsilon)^{-|j|_\mathfrak{s} - 4} \int_{\|\bar z\|_\mathfrak{s} \leq r} |G^\varepsilon(\bar z)| d\bar z \lesssim \varepsilon(\|z\|_\mathfrak{s} + \varepsilon)^{-|j|_\mathfrak{s} - 4} \int_{\|\bar z\|_\mathfrak{s} \leq r} (\|\bar z\|_\mathfrak{s} + \varepsilon)^{-1} d\bar z \\ \lesssim \varepsilon(\|z\|_\mathfrak{s} + \varepsilon)^{-|j|_\mathfrak{s} - 4}(r + \varepsilon)^2 \lesssim \varepsilon(\|z\|_\mathfrak{s} + \varepsilon)^{-|j|_\mathfrak{s} - 2} \lesssim \varepsilon^\alpha(\|z\|_\mathfrak{s} + \varepsilon)^{-|j|_\mathfrak{s} - 1 - \alpha},\end{aligned}$$

for any $\alpha \in [0, 1]$, where the bound $\int_{\|\bar z\|_\mathfrak{s} \leq r} (\|\bar z\|_\mathfrak{s} + \varepsilon)^{-1} d\bar z \lesssim (r + \varepsilon)^2$ follows from the definitions of the norm $\|\bar z\|_\mathfrak{s}$ and where in the second to the last bound we used our definition of $r$.

The second term in (4.39) is more complicated to bound. For this, we are going to use the expansion of $G^\varepsilon$ provided in Proposition 4.2. Namely, we write the second term as

$$\sum_{n=0}^M \int_{\mathbb{R}_+ \times \Lambda_\varepsilon} \varphi_r(\bar z)\, D_\varepsilon^j K^{\varepsilon, n}(z - \bar z)\, S^\varepsilon(\bar z) d\bar z + \int_{\mathbb{R}_+ \times \Lambda_\varepsilon} \varphi_r(\bar z)\, D_\varepsilon^j R^\varepsilon(z - \bar z)\, S^\varepsilon(\bar z) d\bar z. \tag{4.40}$$

We note that the kernel $R^\varepsilon$ is defined on $\mathbb{R}^2$, which allows us to use Lemma 4.8 and bound

$$\Big| \int_{\mathbb{R}_+ \times \Lambda_\varepsilon} \varphi_r(\bar z)\, D_\varepsilon^j R^\varepsilon(z - \bar z)\, S^\varepsilon(\bar z) d\bar z \Big| \lesssim \varepsilon^\alpha \|\varphi_r(\bullet)\, D_\varepsilon^j R^\varepsilon(z - \bullet)\|_{\mathcal{C}^\alpha_{\mathfrak{s}, \varepsilon}}.$$

For $k \in \mathbf{N}^2$ we have $\|D_\varepsilon^k \varphi_r\|_{L^\infty} \lesssim r^{-|k|_\mathfrak{s}} \lesssim (\|z\|_\mathfrak{s} + \varepsilon)^{-|k|_\mathfrak{s}}$. Moreover, we have $\|\bar z\|_\mathfrak{s} \leq r$ inside the integral, which implies $\|z - \bar z\|_\mathfrak{s} + \varepsilon \geq \frac{1}{2}(\|z\|_\mathfrak{s} + \varepsilon)$ as above. Then Proposition 4.2 yields $|D_\varepsilon^k R^\varepsilon(z - \bar z)| \lesssim (\|z\|_\mathfrak{s} + \varepsilon)^{-|k|_\mathfrak{s} - 1}$. Hence,

$$|D_\varepsilon^k (\varphi_r(\bar z)\, D_\varepsilon^j R^\varepsilon(z - \bar z))| \lesssim (\|z\|_\mathfrak{s} + \varepsilon)^{-|j|_\mathfrak{s} - |k|_\mathfrak{s} - 1},$$

where the operator $D_\varepsilon^k$ is applied to the variable $\bar z$. Using interpolation, we estimate the norm (4.25) as $|\varphi_r(\bullet)\, D_\varepsilon^j R^\varepsilon(z - \bullet)\|_{\mathcal{C}^\alpha_{\mathfrak{s},\varepsilon}} \lesssim r^{-\alpha} \lesssim (\|z\|_\mathfrak{s} + \varepsilon)^{-|j|_\mathfrak{s} - 1 - \alpha}$. Hence, the last term in (4.40) is absolutely bounded by a constant multiple of $\varepsilon^\alpha(\|z\|_\mathfrak{s} + \varepsilon)^{-|j|_\mathfrak{s} - 1 - \alpha}$ for any $\alpha \in [0, 1]$.

We will now bound the first term in (4.40). Let us look at the term for $0 \leq n < M$. Since $K^{\varepsilon,n}(z - \bar z)$ is supported on $\|z - \bar z\|_s \leq c2^{-n}$ and $\|z - \bar z\|_\mathfrak{s} + \varepsilon \geq \frac{1}{2}(\|z\|_\mathfrak{s} + \varepsilon)$ in the integral, we conclude that the integral vanishes if $\|z\|_\mathfrak{s} + \varepsilon > (1 + c)2^{1-n}$. Hence, we should consider $\|z\|_\mathfrak{s} + \varepsilon \leq (1 + c)2^{1-n}$. For $k \in \mathbf{N}^2$, Proposition 4.2 yields $|D_\varepsilon^k K^{\varepsilon,n}(z - \bar z)| \lesssim 2^{n(1+|k|_\mathfrak{s})} \lesssim (\|z\|_\mathfrak{s} + \varepsilon)^{-|k|_\mathfrak{s} - 1}$ and hence

$$|D_\varepsilon^k (\varphi_r(\bar z)\, D_\varepsilon^j K^{\varepsilon,n}(z - \bar z))| \lesssim (\|z\|_\mathfrak{s} + \varepsilon)^{-|j|_\mathfrak{s} - |k|_\mathfrak{s} - 1}.$$

Then the $n$-th term in (4.40) is bounded, in the same way as we bounded the last term, by a constant multiple of $\varepsilon^\alpha(\|z\|_\mathfrak{s} + \varepsilon)^{-|j|_\mathfrak{s} - 1 - \alpha}$ for any $\alpha \in [0, 1]$.

Let as loot at the term in (4.40) with $n = M$. This term is different from the others because the kernel $K^{\varepsilon,M}$ is not extended smoothly to $\mathbb{R}^2$ (see Proposition 4.2). In particular, we cannot apply Lemma 4.8, similarly to how we used it for the other terms. We absolutely bound this term by

$$\int_{\|\bar z\|_\mathfrak{s} \leq r} |D_\varepsilon^j K^{\varepsilon,M}(z - \bar z)|\, |S^\varepsilon(\bar z)|\, d\bar z \lesssim \int_{\|\bar z\|_\mathfrak{s} \leq r} \varepsilon^{-|j|_\mathfrak{s}} (\|\bar z\|_\mathfrak{s} + \varepsilon)^{-4}\, d\bar z \lesssim \varepsilon^{-|j|_\mathfrak{s} - 1},$$

where we used (4.14) and Lemma 4.9. As we explained in the previous paragraph, this term may be non-vanishing only for $\|z\|_\mathfrak{s} + \varepsilon \leq 2(1 + c)\varepsilon$. Hence, the preceding expression is bounded by a constant multiple of $\varepsilon^\alpha(\|z\|_\mathfrak{s} + \varepsilon)^{-|j|_\mathfrak{s} - 1 - \alpha}$ for any $\alpha \in [0, 1]$.

In the third term in (4.39), we have $1 - \varphi_r(\bar z) - \varphi_r(z - \bar z) = 0$ if $\|\bar z\|_\mathfrak{s} \leq \frac{r}{2}$ or $\|z - \bar z\|_\mathfrak{s} \leq \frac{r}{2}$, and thus we can restrict the integral to those $\bar z \in \mathbb{R}_+ \times \Lambda_\varepsilon$ such that $\|\bar z\|_\mathfrak{s} \geq \frac{1}{4}(\|z\|_\mathfrak{s} + \varepsilon)$ and $\|z - \bar z\|_\mathfrak{s} \geq \frac{1}{4}(\|z\|_\mathfrak{s} + \varepsilon)$, hence by Lemmas A.8 and 4.9 the third term in (4.39) is absolutely bounded by a constant times

$$\begin{aligned}(\|z\|_\mathfrak{s} + \varepsilon)^{-|j|_\mathfrak{s} - 1} &\int_{\|\bar z\|_\mathfrak{s} \geq \frac{1}{4}(\|z\|_\mathfrak{s} + \varepsilon)} |S^\varepsilon(\bar z)| d\bar z \\ &\lesssim \varepsilon(\|z\|_\mathfrak{s} + \varepsilon)^{-|j|_\mathfrak{s} - 1} \int_{\|\bar z\|_\mathfrak{s} \geq \frac{1}{4}(\|z\|_\mathfrak{s} + \varepsilon)} (\|\bar z\|_\mathfrak{s} + \varepsilon)^{-4} d\bar z \\ &\lesssim \varepsilon(\|z\|_\mathfrak{s} + \varepsilon)^{-|j|_\mathfrak{s} - 2} \lesssim \varepsilon^\alpha(\|z\|_\mathfrak{s} + \varepsilon)^{-|j|_\mathfrak{s} - 1 - \alpha},\end{aligned}$$

for any $\alpha \in [0, 1]$.

Now we turn to the last term in (4.39). In the same way as we did in (4.36), we conclude that $\|z - \bar z\|_\mathfrak{s} \geq \frac{r}{2} - \varepsilon j_1$ in the last integral in (4.39), which implies

$$\|z - \bar z - (0, \varepsilon k_1)\|_\mathfrak{s} + (1 + 2j_1)\varepsilon \geq \|z - \bar z\|_\mathfrak{s} + (1 + j_1)\varepsilon \geq \frac{1}{4}(\|z\|_\mathfrak{s} + \varepsilon).$$

Furthermore, we have $|D_\varepsilon^k \varphi_r(z - \bar z)| \lesssim r^{-|k|_\mathfrak{s}} \lesssim (\|z\|_\mathfrak{s} + \varepsilon)^{-|k|_\mathfrak{s}}$ uniformly over $z - \bar z$. Then the last term in (4.39) is bounded, in the same way as we bounded the second term, by a constant multiple of $\varepsilon^\alpha(\|z\|_\mathfrak{s} + \varepsilon)^{-|j|_\mathfrak{s} - 1 - \alpha}$ for any $\alpha \in [0, 1]$. Combining the preceding bounds on the terms in (4.39) we get (4.37). □

Combining Lemma 4.10 and Corollary A.3 we get the following bound on the extension (A.5) of the first kernel in (2.21), which we denote by $\widetilde{G}^\varepsilon_{\text{ext}} = \mathsf{Ext}(\widetilde{G}^\varepsilon)$.

**Lemma 4.11** *For any $j \in \mathbf{N}^2$ and $\alpha \in [0,1]$ there is $\varepsilon_0 \in (0,1)$ and a constant $C = C(j,\alpha) > 0$ such that*

$$|D^j(\widetilde{G}^\varepsilon_{\text{ext}} - G^\varepsilon_{\text{ext}})(z)| \leq C\varepsilon^\alpha(\|z\|_\mathfrak{s} + \varepsilon)^{-|j|_\mathfrak{s}-1-\alpha},$$

*uniformly over $z \in \mathbb{R}_+ \times \Lambda_\varepsilon$ and $\varepsilon \in (0,\varepsilon_0)$.*

Using these bounds, we get an analogue of Proposition 4.2 for this kernel.

**Proposition 4.12** *For any integer $r \geq 2$ we can write $\widetilde{G}^\varepsilon_{\text{ext}} = \widetilde{K}^\varepsilon + \widetilde{R}^\varepsilon$, where the functions have the same properties as listed in Proposition 4.2. Moreover, these functions are related to the respective functions of the kernel $G^\varepsilon$, provided by Proposition 4.2, as follows: for any $j \in \mathbf{N}^2$ such that $|j|_\mathfrak{s} \leq r$ and any $\alpha \in [0,1]$ there is $\varepsilon_0 \in (0,1)$ and a constant $C = C(r,\alpha) > 0$ such that*

$$|D^j(\widetilde{R}^\varepsilon - R^\varepsilon)(z)| \leq C\varepsilon^\alpha, \qquad \int_\mathbb{R} |D^j(\widetilde{R}^\varepsilon - R^\varepsilon)(z)|\, dx \leq C\varepsilon^\alpha$$

*and*

$$|D^j(\widetilde{K}^{\varepsilon,n} - K^{\varepsilon,n})(z)| \leq C\varepsilon^\alpha 2^{n(1+\alpha+|j|_\mathfrak{s})}$$

*uniformly over $z \in \mathbb{R}^2$, $0 \leq n < M$ and $\varepsilon \in (0,\varepsilon_0)$.*

*Proof.* The proof is analogous to that of Proposition 4.2, where we only need to show that the function $\widetilde{R}^\varepsilon$ satisfies (4.13). Indeed, these bounds for the function $R^\varepsilon$ follow from the estimate (A.35) on the discrete heat kernel. We have not proved such estimates for the kernel $\widetilde{G}^\varepsilon$. These bounds will follow if we show that the derivatives of $\widetilde{G}^\varepsilon$ are absolutely integrable in the spatial variable.

Let us fix $t \geq 0$ and use the definition to write

$$\begin{aligned}\widetilde{G}^\varepsilon(t,x) &= \int_0^t \varepsilon \sum_{y\in\Lambda_\varepsilon} (\mathscr{Q}^\varepsilon)^{-1}(s,y) G^\varepsilon(t-s,x-y)ds \\ &= G^\varepsilon(t,x) + \int_0^t \varepsilon \sum_{y\in\Lambda_\varepsilon} ((\mathscr{Q}^\varepsilon)^{-1} - \text{Id})(s,y) G^\varepsilon(t-s,x-y)ds.\end{aligned}$$

This yields

$$\varepsilon \sum_{x\in\Lambda_\varepsilon} |D^j_\varepsilon \widetilde{G}^\varepsilon(t,x)| \leq \varepsilon \sum_{x\in\Lambda_\varepsilon} |D^j_\varepsilon G^\varepsilon(t,x)| + \int_0^t \varepsilon^2 \sum_{x,y\in\Lambda_\varepsilon} |((\mathscr{Q}^\varepsilon)^{-1} - \text{Id})(s,y)||D^j_\varepsilon G^\varepsilon(t-s,x-y)|\, ds.$$

By (A.35), the first term is bounded uniformly in $\varepsilon$. Using Lemma 4.9, the second term is bounded too. □

We can now define an extension of the kernel $\widetilde{\widetilde{G}}^\varepsilon$, similarly to how we did it for $\widetilde{G}^\varepsilon$.

**Lemma 4.13** *For any $j \in \mathbf{N}^2$ and $\alpha \in [0,1]$ there is $\varepsilon_0 \in (0,1)$ and a constant $C = C(j,\alpha) > 0$ such that*

$$|D^j_\varepsilon \widetilde{\widetilde{G}}^\varepsilon(z)| \leq C\varepsilon^\alpha(\|z\|_\mathfrak{s} + \varepsilon)^{-|j|_\mathfrak{s}-1-\alpha}, \tag{4.41}$$

*uniformly over $z \in \mathbb{R}_+ \times \Lambda_\varepsilon$ and $\varepsilon \in (0,\varepsilon_0)$.*

*Proof.* Using the definition (2.21) we write

$$\widetilde{\widetilde{G}}^\varepsilon(z) = (\widetilde{G}^\varepsilon \star^+_\varepsilon \varepsilon Q^\varepsilon)(z).$$

We note that the kernel $\varepsilon Q^\varepsilon$ satisfies the same bounds (4.22) and (4.26) as the kernel of $(\mathscr{Q}^\varepsilon)^{-1} - \text{Id}$ (cf. (4.27) and (4.28)). Hence, we can repeat the proof of Lemma 4.10, but using the kernels $\widetilde{G}^\varepsilon$ and $\varepsilon Q^\varepsilon$ in place of $G^\varepsilon$ and $S^\varepsilon$ respectively. This yields (4.41). □

Combining this lemma with Corollary A.3, we get the following bound on the extension $\widetilde{\widetilde{G}}^\varepsilon_{\text{ext}} = \mathsf{Ext}(\widetilde{\widetilde{G}}^\varepsilon)$.

**Lemma 4.14** *For any $j \in \mathbf{N}^2$ and $\alpha \in [0, 1]$ there is $\varepsilon_0 \in (0, 1)$ and a constant $C = C(j, \alpha) > 0$ such that*

$$|D^j \widetilde{\widetilde{G}}^\varepsilon_{\text{ext}}(z)| \leq C\varepsilon^\alpha (\|z\|_\mathfrak{s} + \varepsilon)^{-|j|_\mathfrak{s}-1-\alpha},$$

*uniformly over $z \in \mathbb{R}_+ \times \Lambda_\varepsilon$ and $\varepsilon \in (0, \varepsilon_0)$.*

Using these bounds, the following result is proved similarly to Proposition 4.12.

**Proposition 4.15** *For any integer $r \geq 2$ we can write $\widetilde{\widetilde{G}}^\varepsilon_{\text{ext}} = \widetilde{\widetilde{K}}^\varepsilon + \widetilde{\widetilde{R}}^\varepsilon$, where the functions have the same properties as listed in Proposition 4.2. Moreover, for any $j \in \mathbf{N}^2$ such that $|j|_\mathfrak{s} \leq r$ and any $\alpha \in [0, 1]$ there is $\varepsilon_0 \in (0, 1)$ and a constant $C = C(r, \alpha) > 0$ such that*

$$|D^j \widetilde{\widetilde{R}}^\varepsilon(z)| \leq C\varepsilon^\alpha, \qquad |D^j \widetilde{\widetilde{K}}^{\varepsilon,n}(z)| \leq C\varepsilon^\alpha 2^{n(1+\alpha+|j|_\mathfrak{s})}$$

*uniformly over $z \in \mathbb{R}^2$, $0 \leq n < M$ and $\varepsilon \in (0, \varepsilon_0)$.*

### 4.6 Periodisation of the kernels

As in (4.11), we will also need the periodised kernels (2.21) defined as

$$\widetilde{P}^\varepsilon(t, x) := \sum_{m\in\mathbb{Z}} \widetilde{G}^\varepsilon(t, x + 2m), \qquad \widetilde{\widetilde{P}}^\varepsilon(t, x) := \sum_{m\in\mathbb{Z}} \widetilde{\widetilde{G}}^\varepsilon(t, x + 2m). \tag{4.42}$$

They can be written, similarly to (4.12), as the Fourier series

$$\widetilde{P}^\varepsilon(t, x) = \frac{1}{2} \sum_{k\in\mathbb{T}_N} (\mathcal{F}\widetilde{G}^\varepsilon_t)(k)\, e^{\pi i k x}, \qquad \widetilde{\widetilde{P}}^\varepsilon(t, x) = \frac{1}{2} \sum_{k\in\mathbb{T}_N} (\mathcal{F}\widetilde{\widetilde{G}}^\varepsilon_t)(k)\, e^{\pi i k x}, \tag{4.43}$$

where $\mathcal{F}$ is the Fourier transform (4.1).

## 5 Regularity structure for the discrete equation

The goal of this section is to write the discrete equation (2.33) in the framework of regularity structures.

The regularity structure for the equation (2.33) has three main differences from the regularity structure for the KPZ equation (1.11): (1) we need to introduce two "abstract" derivatives $\partial_-$ and $\partial_+$ which correspond to the two discrete derivatives $\nabla^-_\varepsilon$ and $\nabla^+_\varepsilon$ respectively, (2) we need to introduce three "abstract integration maps" $\mathcal{I}$, $\widetilde{\mathcal{I}}$ and $\widetilde{\widetilde{\mathcal{I}}}$ corresponding to $G^\varepsilon$, $\widetilde{G}^\varepsilon$ and $\widetilde{\widetilde{G}}^\varepsilon$ respectively, (3) we need to describe the multiplication operator by $\varepsilon$ on the level of regularity structure, which will allow to describe the error (2.24b). We assume that the maps $\partial_-$ and $\partial_+$ act on monomials $X^n$ in the same way as $\partial$ (see Section 3.1). Moreover, the symbol $\mathcal{E}$ will correspond to multiplication by $\varepsilon$.

Using the notation of Section 3.1, we define the minimal sets $\widetilde{\mathcal{V}}$, $\widetilde{\mathcal{U}}$, $\widetilde{\mathcal{U}}'_-$ and $\widetilde{\mathcal{U}}'_+$ of formal expressions such that $\mathcal{W}_{\text{poly}} \subset \widetilde{\mathcal{V}} \cap \widetilde{\mathcal{U}} \cap \widetilde{\mathcal{U}}'_- \cap \widetilde{\mathcal{U}}'_+$,

$$\Xi \in \widetilde{\mathcal{V}}, \qquad \mathcal{I}(\Xi) \in \widetilde{\mathcal{U}}, \quad \partial_-\mathcal{I}(\Xi) \in \widetilde{\mathcal{U}}'_-, \quad \partial_+\mathcal{I}(\Xi) \in \widetilde{\mathcal{U}}'_+,$$

and the following implications hold:

$$\tau \in \widetilde{\mathcal{V}} \setminus \{\Xi\} \quad \Rightarrow \quad \widetilde{\mathcal{I}}(\tau) \in \widetilde{\mathcal{U}}, \quad \partial_-\widetilde{\mathcal{I}}(\tau) \in \widetilde{\mathcal{U}}'_-, \quad \partial_+\widetilde{\mathcal{I}}(\tau) \in \widetilde{\mathcal{U}}'_+, \tag{5.1a}$$

$$\tau_1 \in \widetilde{\mathcal{U}}'_-,\ \tau_2 \in \widetilde{\mathcal{U}}'_+ \quad \Rightarrow \quad \tau_1\tau_2 \in \widetilde{\mathcal{V}}, \tag{5.1b}$$

$$\tau \in \widetilde{\mathcal{U}}'_- \cup \widetilde{\mathcal{U}}'_+ \quad \Rightarrow \quad \widetilde{\widetilde{\mathcal{I}}}(\tau) \in \widetilde{\mathcal{U}}, \quad \partial_-\widetilde{\widetilde{\mathcal{I}}}(\tau) \in \widetilde{\mathcal{U}}'_-, \quad \partial_+\widetilde{\widetilde{\mathcal{I}}}(\tau) \in \widetilde{\mathcal{U}}'_+, \tag{5.1c}$$

$$\tau \in \widetilde{\mathcal{U}}'_- \setminus \mathcal{W}_{\text{poly}} \quad \Rightarrow \quad \mathcal{E}(\partial_+ \widetilde{\widetilde{\mathcal{I}}}(\tau)) \in \widetilde{\mathcal{U}}, \quad \partial_- \mathcal{E}(\partial_+ \widetilde{\widetilde{\mathcal{I}}}(\tau)) \in \widetilde{\mathcal{U}}'_-, \quad \partial_+ \mathcal{E}(\partial_+ \widetilde{\widetilde{\mathcal{I}}}(\tau)) \in \widetilde{\mathcal{U}}'_+, \tag{5.1d}$$

where as before the product of symbols is commutative with the convention $\mathbf{1}\tau = \tau$. We postulate that $\mathcal{I}$, $\widetilde{\mathcal{I}}$ and $\widetilde{\widetilde{\mathcal{I}}}$ map monomials to zero and do not include such zero elements into the sets.

We set $\widetilde{\mathscr{W}} := \widetilde{\mathcal{U}} \cup \widetilde{\mathcal{V}} \cup \widetilde{\mathcal{U}}'_- \cup \widetilde{\mathcal{U}}'_+$, and similarly to (3.3) we define the homogeneities

$$|X^\ell| = 2\ell_0 + \ell_1 \quad \text{for} \quad \ell = (\ell_0, \ell_1) \in \mathbf{N}^2, \tag{5.2a}$$
$$|\Xi| = -\frac{3}{2} - \kappa, \tag{5.2b}$$
$$|\tau_1 \tau_2| = |\tau_1| + |\tau_2|, \tag{5.2c}$$
$$|\mathcal{I}(\tau)| = |\tau| + 2, \tag{5.2d}$$
$$|\partial_- \bar\tau| = |\partial_+ \bar\tau| = |\bar\tau| - 1, \tag{5.2e}$$

where $\kappa \in (\frac{1}{2} - \alpha, \frac{1}{10})$. The increase of homogeneities by the maps $\widetilde{\mathcal{I}}$ and $\widetilde{\widetilde{\mathcal{I}}}$ is the same as by $\mathcal{I}$, and we also have the identity

$$|\mathcal{E}(\tau)| = |\tau| + 1. \tag{5.3}$$

The restriction on $\kappa < \frac{1}{10}$ is explained below and is necessary to have a minimal number of basis elements in the regularity structure. The bound $\kappa > \frac{1}{2} - \alpha$ will be used in Section 6.4.3 when proving moment bounds on the discrete model.

The set $\widetilde{\mathscr{T}}$ is defined to contain all finite linear combinations of the elements in $\widetilde{\mathscr{W}}$, and we view the maps $\mathcal{E}$, $\mathcal{I}$, $\widetilde{\mathcal{I}}$ and $\widetilde{\widetilde{\mathcal{I}}}$ and all their defined derivatives as linear maps on the respective subspaces of $\widetilde{\mathscr{T}}$. We also extend $\partial_-$ and $\partial_+$ linearly to the respective subspaces of $\widetilde{\mathscr{T}}$. The set $\widetilde{\mathfrak{A}}$ contains all homogeneities $|\tau|$ for all $\tau \in \widetilde{\mathscr{T}}$.

Similarly to (3.4) we define

$$\widetilde{\mathcal{W}} := \{\tau \in \widetilde{\mathcal{V}} : \mathcal{E}(\tau) \notin \widetilde{\mathcal{V}},\ |\tau| \le 0\} \cup \{\tau \in \widetilde{\mathcal{V}} : \mathcal{E}(\tau) \in \widetilde{\mathcal{V}},\ |\tau| \le -1\} \cup \Big\{\tau \in \widetilde{\mathcal{U}}'_- \cup \widetilde{\mathcal{U}}'_+ : |\tau| < \frac{1}{2}\Big\} \cup \Big\{\tau \in \widetilde{\mathcal{U}} : |\tau| < \frac{3}{2}\Big\}. \tag{5.4}$$

The set $\widetilde{\mathcal{W}}$ contains $48$ elements which are the symbols from Table 1 with different combinations of integration maps and their derivatives and $5$ new symbols containing $\mathcal{E}$. These symbols are

$$\partial_+ \mathcal{I}(\Xi) \partial_- \mathcal{E}(\partial_+ \widetilde{\widetilde{\mathcal{I}}}(\partial_- \mathcal{I}(\Xi))), \qquad \partial_- \mathcal{I}(\Xi) \partial_+ \mathcal{E}(\partial_+ \widetilde{\widetilde{\mathcal{I}}}(\partial_- \mathcal{I}(\Xi))),$$
$$\partial_- \mathcal{E}(\partial_+ \widetilde{\widetilde{\mathcal{I}}}(\partial_- \mathcal{I}(\Xi))), \qquad \partial_+ \mathcal{E}(\partial_+ \widetilde{\widetilde{\mathcal{I}}}(\partial_- \mathcal{I}(\Xi))), \qquad \mathcal{E}(\partial_+ \widetilde{\widetilde{\mathcal{I}}}(\partial_- \mathcal{I}(\Xi))). \tag{5.5}$$

The three symbols in the first line of (5.5) have homogeneity $-2\kappa$, and the two symbols in the second line have homogeneity $\frac{1}{2} - \kappa$. The element which is not included into $\widetilde{\mathcal{W}}$ and which has the minimal homogeneity, is an element of $\widetilde{\mathcal{V}}$ with homogeneity $\frac{1}{2} - 5\kappa$. The reason of having $\kappa < \frac{1}{10}$ is to have $\frac{1}{2} - 5\kappa > 0$ and hence not to include it into the set $\widetilde{\mathcal{W}}$.

We define $\widetilde{\mathcal{T}}$ to be the linear span of the elements in $\widetilde{\mathcal{W}}$, and the set $\widetilde{\mathcal{A}}$ contains all homogeneities $|\tau|$ for elements $\tau \in \widetilde{\mathcal{W}}$. In order to represent symbols as trees, as we did in (3.5), an edge will correspond to the symbol $\mathcal{I}$, an edge will correspond to the symbol $\widetilde{\mathcal{I}}$ and an edge will correspond to the symbol $\widetilde{\widetilde{\mathcal{I}}}$. Integration against a discrete derivative of one of the kernels will be represented by a respective blue edge with a label either "+" or "−", corresponding to one of the two discrete derivatives. We use the quantity (3.6) and the projections (3.7) also for these sets.

To define a structure group $\widetilde{\mathcal{G}}$, we introduce another set of elements

$$\widetilde{\mathcal{W}}_+ := \{\tau \in \widetilde{\mathcal{W}} : |\tau| \ge 0\}. \tag{5.6}$$

This set contains in total 25 elements, which are the elements of (3.9) with different combinations of integration maps and their derivatives and the 3 symbols from the second line of (5.5). Then we define $\widetilde{\mathcal{T}}_+$ is the free commutative algebra generated by the elements of $\widetilde{\mathcal{W}}_+$.

By analogy with (3.10), we define a linear map $\widetilde{\Delta} : \widetilde{\mathcal{T}} \to \widetilde{\mathcal{T}} \otimes \widetilde{\mathcal{T}}_+$ by the identities

$$\widetilde{\Delta}\mathbf{1} = \mathbf{1} \otimes \mathbf{1}, \qquad \widetilde{\Delta} X_1 = X_1 \otimes \mathbf{1} + \mathbf{1} \otimes X_1, \qquad \widetilde{\Delta}\Xi = \Xi \otimes \mathbf{1}, \tag{5.7a}$$

and then recursively by

$$\widetilde{\Delta}(\tau_1\tau_2) = (\widetilde{\Delta}\tau_1)(\widetilde{\Delta}\tau_2), \tag{5.7b}$$
$$\widetilde{\Delta}\widetilde{\mathcal{I}}(\tau) = (\widetilde{\mathcal{I}} \otimes I)\widetilde{\Delta}\tau + \mathbf{1} \otimes \widetilde{\mathcal{I}}(\tau), \qquad \text{if } |\tau| < -1, \tag{5.7c}$$
$$\widetilde{\Delta}\partial_\pm\widetilde{\mathcal{I}}(\tau) = (\partial_\pm\widetilde{\mathcal{I}} \otimes I)\widetilde{\Delta}\tau, \qquad \text{if } |\tau| < -1, \tag{5.7d}$$
$$\widetilde{\Delta}\widetilde{\mathcal{I}}(\tau) = (\widetilde{\mathcal{I}} \otimes I)\widetilde{\Delta}\tau + \mathbf{1} \otimes \widetilde{\mathcal{I}}(\tau) + (\widetilde{\Delta}X_1)(\mathbf{1} \otimes \partial\widetilde{\mathcal{I}}(\tau)), \quad \text{if } |\tau| > -1, \tag{5.7e}$$
$$\widetilde{\Delta}\partial_\pm\widetilde{\mathcal{I}}(\tau) = (\partial_\pm\widetilde{\mathcal{I}} \otimes I)\widetilde{\Delta}\tau + \mathbf{1} \otimes \partial\widetilde{\mathcal{I}}(\tau), \qquad \text{if } |\tau| > -1, \tag{5.7f}$$

for respective elements $\tau_i, \tau \in \widetilde{\mathcal{W}}$. We use the symbol $\partial = \frac{1}{2}(\partial_- + \partial_+)$ in (5.7e), which means

$$\partial\widetilde{\mathcal{I}}(\tau) = \frac{1}{2}(\partial_-\widetilde{\mathcal{I}}(\tau) + \partial_+\widetilde{\mathcal{I}}(\tau)) \tag{5.7g}$$

by analogy with the symmetric discrete derivative (2.5). We now have three symbols corresponding to derivatives: $\partial$, $\partial_-$ and $\partial_+$, and we could choose any of them to use on the right-hand side of (5.7e). We think that it is slightly more natural to use the symmetric version $\partial$ which corresponds to the symmetric discrete derivative in our later definition of the model in (6.4d). Having the derivative $\partial$ in (5.7e) will be also more convenient because the reconstruction map will satisfy (6.5) which will simplify computations in the proof of Lemma 7.3. We use respectively $\partial$ in (5.7f), which corresponds to the last two identities in (6.4d). The action of $\widetilde{\Delta}$ on the map $\mathcal{I}$ and $\widetilde{\widetilde{\mathcal{I}}}$ is defined by analogy with its action on $\widetilde{\mathcal{I}}$. Since the map $\mathcal{I}$ appears only as $\mathcal{I}(\Xi)$, the action of $\widetilde{\Delta}$ on is defined by (5.7c) and (5.7d). Finally, we define

$$\widetilde{\Delta}\partial_\pm\mathcal{E}(\tau) = (\partial_\pm\mathcal{E} \otimes I)\widetilde{\Delta}\tau + \mathbf{1} \otimes \partial_\pm\mathcal{E}(\tau), \tag{5.7h}$$
$$\widetilde{\Delta}\mathcal{E}(\tau) = (\mathcal{E} \otimes I)\widetilde{\Delta}\tau + \mathbf{1} \otimes \mathcal{E}(\tau) + (\widetilde{\Delta}X_1)(\mathbf{1} \otimes \partial_-\mathcal{E}(\tau)). \tag{5.7i}$$

The action of $\widetilde{\Delta}$ on $\mathcal{E}$ is defined by analogy with $\widetilde{\mathcal{I}}$, although here we use the derivative $\partial_-$. This choice of the derivative will be convenient because the discrete model $Z^{\varepsilon,\mathrm{m}}$, defined in Section 6, will satisfy the assumptions of Lemma 7.2.

The maps $\Gamma_f : \widetilde{\mathcal{T}} \to \widetilde{\mathcal{T}}$ are defined as in (3.11) via $\widetilde{\Delta}$, i.e.

$$\Gamma_f\tau := (I \otimes f)\widetilde{\Delta}\tau, \tag{5.8}$$

for any linear functional $f : \widetilde{\mathcal{T}}_+ \to \mathbb{R}$. The structure group $\widetilde{\mathcal{G}}$ is defined respectively as $\{\Gamma_f : f \in \widetilde{\mathcal{G}}_+\}$, where $\widetilde{\mathcal{G}}_+$ contains all linear multiplicative functionals $f : \widetilde{\mathcal{T}}_+ \to \mathbb{R}$ satisfying $f(\mathbf{1}) = 1$. This yields a new regularity structure $\widetilde{\mathscr{T}} = (\widetilde{\mathcal{A}}, \widetilde{\mathcal{T}}, \widetilde{\mathcal{G}})$.

### 5.1 Discrete models

Let $\mathcal{B}^2_\mathfrak{s}$ be the set of all test functions $\varphi \in \mathcal{C}^2(\mathbb{R}^2)$, compactly supported in the ball of radius 1 around the origin (with respect to the parabolic distance $\|\bullet\|_\mathfrak{s}$ defined in Section 1.3), and satisfying $\|\varphi\|_{\mathcal{C}^2} \leq 1$. By analogy with (1.20), for $\varphi \in \mathcal{B}^2_\mathfrak{s}$, $\lambda \in (0,1]$ and $(s,y) \in \mathbb{R}^2$ we define a rescaled and recentered function

$$\varphi^\lambda_{(s,y)}(t,x) := \frac{1}{\lambda^3}\varphi\Big(\frac{t-s}{\lambda^2}, \frac{x-y}{\lambda}\Big). \tag{5.9}$$

In order to use the results of [EH19], we need to define a discretisation for the regularity structure $\widetilde{\mathscr{T}}$ according to [EH19, Def. 2.1]. The following definition corresponds to the "semidiscrete" case in [EH19, Sec. 2]. We use the space-time domain $D_\varepsilon := \mathbb{R} \times \Lambda_\varepsilon$ in these definitions.

**Definition 5.1** 1. We define the space $\mathcal{X}_\varepsilon := L^\infty(D_\varepsilon)$, and we extend the operator (1.22) to $\iota_\varepsilon : \mathcal{X}_\varepsilon \hookrightarrow L^\infty(\mathbb{R}, \mathscr{D}'(\mathbb{R}))$ as

$$(\iota_\varepsilon f_\varepsilon)(t, \bullet) := (\iota_\varepsilon f_\varepsilon(t))(\bullet) \tag{5.10}$$

for $f \in \mathcal{X}_\varepsilon$. For any smooth compactly supported function $\varphi : \mathbb{R}^2 \to \mathbb{R}$ it will be convenient to write

$$(\iota_\varepsilon f)(\varphi) := \varepsilon \sum_{x \in \Lambda_\varepsilon} \int_{\mathbb{R}} f(t,x)\varphi(t,x)dt. \tag{5.11}$$

2. For any $\gamma \in \mathbb{R}$, $z \in D_\varepsilon$ and a compact set $K_\varepsilon \subset \mathbb{R}^2$ of diameter at most $2\varepsilon$, we define the following seminorm $f_\varepsilon \in \mathcal{X}_\varepsilon$:

$$\|f_\varepsilon\|_{\gamma;K_\varepsilon;z;\varepsilon} := \varepsilon^{-\gamma} \sup_{z \in K_\varepsilon \cap D_\varepsilon} |f_\varepsilon(z)|. \tag{5.12}$$

This seminorm is local in the sense that if $f_\varepsilon, g_\varepsilon \in \mathcal{X}_\varepsilon$ and $(\iota_\varepsilon f_\varepsilon)(\varphi) = (\iota_\varepsilon g_\varepsilon)(\varphi)$ for every $\varphi \in \mathcal{C}^2$ supported in $K_\varepsilon$, then $\|f_\varepsilon - g_\varepsilon\|_{\gamma;K_\varepsilon;z;\varepsilon} = 0$.

3. For any function $\Gamma : D_\varepsilon^2 \to \widetilde{\mathcal{G}}$, any compact set $K \subset \mathbb{R}^2$ and any $\gamma \in \mathbb{R}$, we define the following seminorm on the functions $f_\varepsilon : D_\varepsilon \to \widetilde{\mathcal{T}}_{<\gamma}$:

$$|\!|\!| f_\varepsilon |\!|\!|_{\gamma;K;\varepsilon} := \sup_{\substack{z,\bar{z} \in K \cap D_\varepsilon \\ \|z - \bar{z}\|_{\mathfrak{s}} \leq \varepsilon}} \sup_{m < \gamma} \varepsilon^{m-\gamma} |f_\varepsilon(z) - \Gamma_{z\bar{z}} f_\varepsilon(\bar{z})|_m, \tag{5.13}$$

where we use the quantity (3.6). For a second function $\bar{\Gamma} : D_\varepsilon^2 \to \widetilde{\mathcal{G}}$ and for $\bar{f}_\varepsilon : D_\varepsilon \to \widetilde{\mathcal{T}}_{<\gamma}$ we also define

$$|\!|\!| f_\varepsilon; \bar{f}_\varepsilon |\!|\!|_{\gamma;K;\varepsilon} := \sup_{\substack{z,\bar{z} \in K \cap D_\varepsilon \\ \|z - \bar{z}\|_{\mathfrak{s}} \leq \varepsilon}} \sup_{m < \gamma} \varepsilon^{m-\gamma} \left|f_\varepsilon(z) - \Gamma_{z\bar{z}} f_\varepsilon(\bar{z}) - \bar{f}_\varepsilon(z) + \bar{\Gamma}_{z\bar{z}} \bar{f}_\varepsilon(\bar{z})\right|_m. \tag{5.14}$$

Both seminorms depend only on the values of $f_\varepsilon$ and $\bar{f}_\varepsilon$ in a neighbourhood of size $c\varepsilon$ around $K$, for a fixed constant $c > 0$. The seminorms (5.13) and (5.14) depends on the functions $\Gamma$ and $\bar{\Gamma}$. However, we prefer not to indicate it to have a lighter notation. The choice of these functions will be always clear from the context.

**Remark 5.2** Let the function $\varphi_z^\varepsilon$ be defined by (5.9) with $\lambda = \varepsilon$, and let $[\varphi_z^\varepsilon] \subset \mathbb{R}^2$ denote its support. Then from the definition (5.12) we get the bound

$$|(\iota_\varepsilon f_\varepsilon)(\varphi_z^\varepsilon)| \leq \Big(\sup_{\bar{z} \in [\varphi_z^\varepsilon] \cap D_\varepsilon} |f_\varepsilon(\bar{z})|\Big) \varepsilon \sum_{x \in \Lambda_\varepsilon} \int_{\mathbb{R}} |\varphi_z^\varepsilon(t,x)| dt \lesssim \varepsilon^\gamma \|f_\varepsilon\|_{\gamma;[\varphi_z^\varepsilon];z;\varepsilon},$$

uniformly over $f_\varepsilon \in \mathcal{X}_\varepsilon$, $z \in D_\varepsilon$, $\gamma \in \mathbb{R}$ and $\varphi \in \mathcal{B}^2_{\mathfrak{s}}$.

Following [EH19, Def. 2.5], we can define a discrete model on the regularity structure $\widetilde{\mathscr{T}}$. However, our definition will be slightly different because we distinguish the element $\Xi$. This is because $\Xi$ corresponds to a time derivative of the driving martingale in (2.33), which is a distribution in the time variable. On the other hand, all the other elements of $\widetilde{\mathcal{W}}$ correspond to functions in time.

**Definition 5.3** A *discrete model* $(\Pi^\varepsilon, \Gamma^\varepsilon)$ on the regularity structure $\widetilde{\mathscr{T}}$ consists of a collection of maps $D_\varepsilon \ni z \mapsto \Pi_z^\varepsilon \in \mathcal{L}(\mathrm{span}\{\widetilde{\mathcal{W}} \setminus \{\Xi\}\}, \mathcal{X}_\varepsilon)$ and $D_\varepsilon \times D_\varepsilon \ni (z, \bar{z}) \mapsto \Gamma_{z\bar{z}}^\varepsilon \in \widetilde{\mathcal{G}}$ with the following algebraic properties:

1. $\Gamma^\varepsilon_{zz} = \mathrm{Id}$ and $\Gamma^\varepsilon_{z\bar z}\Gamma^\varepsilon_{\bar z\tilde z} = \Gamma^\varepsilon_{z\tilde z}$ for all $z, \bar z, \tilde z \in D_\varepsilon$.
2. $\Pi^\varepsilon_{\bar z} = \Pi^\varepsilon_z \Gamma^\varepsilon_{z\bar z}$ for all $z, \bar z \in D_\varepsilon$. (We note that these identities make sense on $\mathrm{span}\{\widetilde{\mathcal{W}} \setminus \{\Xi\}\}$ because, as follows from the definitions (3.11) and (5.7), $\Gamma^\varepsilon_{z\bar z}\tau$ contains a term proportional to $\Xi$ only for $\tau = \Xi$.)
3. $(\Pi^\varepsilon_z \partial_\pm \tau)(z) = \nabla^\pm_\varepsilon (\Pi^\varepsilon_z \tau)(\bar z)|_{\bar z = z}$ for any $\tau$ from the domain of $\partial_\pm$, where the discrete derivatives act on the spatial component of $\bar z$.
4. $(\Pi^\varepsilon_z \mathbf{1})(\bar z) = 1$, $(\Pi^\varepsilon_z X_1)(\bar z) = \bar x - x$ for all $z, \bar z \in D_\varepsilon$ with the spatial components $x, \bar x$.

Furthermore, for any compact set $K \subset \mathbb{R}^2$ the following bounds hold

$$\sup_{\varphi \in \mathcal{B}^2_\mathfrak{s}} \sup_{z \in K \cap D_\varepsilon} |(\iota_\varepsilon \Pi^\varepsilon_z \tau)(\varphi^\lambda_z)| \lesssim \lambda^{|\tau|}, \qquad \sup_{K_\varepsilon \subset K} \sup_{z \in K \cap D_\varepsilon} \|\Pi^\varepsilon_z \tau\|_{|\tau|; K_\varepsilon; z; \varepsilon} \lesssim 1, \tag{5.15a}$$

uniformly over $\lambda \in [\varepsilon, 1]$ and $\tau \in \widetilde{\mathcal{W}} \setminus \{\Xi\}$, where the supremum in the second bound is over compact sets $K_\varepsilon \subset K$ with the diameter not exceeding $2\varepsilon$. For the function $f^{\tau, \Gamma^\varepsilon}_{\bar z}(z) := \Gamma^\varepsilon_{z\bar z}\tau - \tau$ one has

$$|\Gamma^\varepsilon_{z\bar z}\tau|_m \lesssim \|z - \bar z\|_\mathfrak{s}^{|\tau| - m}, \qquad \sup_{\bar z \in D_\varepsilon} \|f^{\tau, \Gamma^\varepsilon}_{\bar z}\|_{|\tau|; K; \varepsilon} \lesssim 1, \tag{5.15b}$$

uniformly over $\tau \in \widetilde{\mathcal{W}}$, $m < |\tau|$ and $z, \bar z \in K \cap D_\varepsilon$ such that $\|z - \bar z\|_\mathfrak{s} \in [\varepsilon, 1]$. In the second bound in (5.15b) we consider the seminorm (5.13) with respect to the map $\Gamma^\varepsilon$.

Finally, the maps $\Pi^\varepsilon_z$ are defined also on $\Xi$ such that $(\Pi^\varepsilon_z \Xi)(s, y)$ is a distribution in $s \in \mathbb{R}$ and a function in $y \in \Lambda_\varepsilon$, and the first bound in (5.15a) holds with $\tau = \Xi$.

**Remark 5.4** The last two algebraic properties are not included to the definition of a discrete model [EH19, Def. 2.5]. Since these properties will be crucial in our analysis in Section 7, we prefer to make them a part of the definition. We note that the last property implies

$$\Gamma^\varepsilon_{z\bar z}\mathbf{1} = \mathbf{1}, \qquad \Gamma^\varepsilon_{z\bar z} X_1 = X_1 - (\bar x - x)\mathbf{1}.$$

**Remark 5.5** The first bounds in (5.15) control the model on the scale above $\varepsilon$ similarly to continuous models in [Hai14], and the second bounds in (5.15) control the model on the scale below $\varepsilon$.

We denote by $\|\Pi^\varepsilon\|^{(\varepsilon)}_K$ and $\|\Gamma^\varepsilon\|^{(\varepsilon)}_K$ the smallest proportionality constants such that the bounds (5.15a) and (5.15b) hold respectively. Then for the model $Z^\varepsilon = (\Pi^\varepsilon, \Gamma^\varepsilon)$ we set

$$|\!|\!| Z^\varepsilon |\!|\!|^{(\varepsilon)}_K := \|\Pi^\varepsilon\|^{(\varepsilon)}_K + \|\Gamma^\varepsilon\|^{(\varepsilon)}_K.$$

For a second model $\bar Z^\varepsilon = (\bar\Pi^\varepsilon, \bar\Gamma^\varepsilon)$ we define the "distance"

$$|\!|\!| Z^\varepsilon; \bar Z^\varepsilon |\!|\!|^{(\varepsilon)}_K := \|\Pi^\varepsilon - \bar\Pi^\varepsilon\|^{(\varepsilon)}_K + \|\Gamma^\varepsilon; \bar\Gamma^\varepsilon\|^{(\varepsilon)}_K,$$

where $\|\Gamma^\varepsilon; \bar\Gamma^\varepsilon\|^{(\varepsilon)}_K$ is the smallest proportionality constant such that the following bounds hold

$$\|(\Gamma^\varepsilon_{z\bar z} - \bar\Gamma^\varepsilon_{z\bar z})\tau\|_m \lesssim \|z - \bar z\|_\mathfrak{s}^{|\tau| - m}, \qquad \sup_{\bar z \in D_\varepsilon} |\!|\!| f^{\tau, \Gamma^\varepsilon}_{\bar z}; f^{\tau, \bar\Gamma^\varepsilon}_{\bar z} |\!|\!|_{|\tau|; K; \varepsilon} \lesssim 1,$$

uniformly over the same quantities as in (5.15b), where in the second bound we consider the distance (5.14) with respect to $\Gamma^\varepsilon$ and $\bar\Gamma^\varepsilon$.

**Remark 5.6** We will often work with models on the set $K = [-T, T] \times [-2, 2]$. In this case we prefer to remove the set $K$ from the notation and write $\|\Pi^\varepsilon\|^{(\varepsilon)}_T$, $\|\Gamma^\varepsilon\|^{(\varepsilon)}_T$, etc.

### 5.2 Modelled distributions and a reconstruction theorem

By analogy with [Hai14, Sec. 6], we are going to define a weighted norm for $\widetilde{\mathcal{T}}$-valued functions with a weight at time $0$. For this we define the following quantities for $z, \bar{z} \in \mathbb{R}^2$:

$$\|z\|_0 := |t|^{\frac{1}{2}} \wedge 1, \qquad \|z, \bar{z}\|_0 := \|z\|_0 \wedge \|\bar{z}\|_0,$$

where $z = (t, x)$ with $t \in \mathbb{R}$.

For $\gamma, \eta \in \mathbb{R}$ and for a compact set $K \subset \mathbb{R}^2$, we define in the context of Definition 5.1(3) the following quantities, as in [EH19, Eqs. 3.21, 3.22]:

$$\|\!| f_\varepsilon \|\!|_{\gamma,\eta;K;\varepsilon} := \sup_{\substack{z \in K \cap D_\varepsilon \\ \|z\|_{\mathfrak{s}} \leq \varepsilon}} \sup_{m < \gamma} \frac{|f_\varepsilon(z)|_m}{\varepsilon^{(\eta - m)\wedge 0}} + \sup_{\substack{z, \bar{z} \in K \cap D_\varepsilon \\ \|z - \bar{z}\|_{\mathfrak{s}} \leq \varepsilon}} \sup_{m < \gamma} \frac{|f_\varepsilon(z) - \Gamma^\varepsilon_{z\bar{z}} f_\varepsilon(\bar{z})|_m}{\varepsilon^{\gamma - m}(\|z, \bar{z}\|_0 + \varepsilon)^{\eta - \gamma}}, \tag{5.16}$$

and

$$\begin{aligned} \|\!| f_\varepsilon; \bar{f}_\varepsilon \|\!|_{\gamma,\eta;K;\varepsilon} := &\sup_{\substack{z \in K \cap D_\varepsilon \\ \|z\|_{\mathfrak{s}} \leq \varepsilon}} \sup_{m < \gamma} \frac{|f_\varepsilon(z) - \bar{f}_\varepsilon(z)|_m}{\varepsilon^{(\eta - m)\wedge 0}} \\ &+ \sup_{\substack{z, \bar{z} \in K \cap D_\varepsilon \\ \|z - \bar{z}\|_{\mathfrak{s}} \leq \varepsilon}} \sup_{m < \gamma} \frac{|f_\varepsilon(z) - \Gamma^\varepsilon_{z\bar{z}} f_\varepsilon(\bar{z}) - \bar{f}_\varepsilon(z) + \bar{\Gamma}^\varepsilon_{z\bar{z}} \bar{f}_\varepsilon(\bar{z})|_m}{\varepsilon^{\gamma - m}(\|z, \bar{z}\|_0 + \varepsilon)^{\eta - \gamma}}. \end{aligned} \tag{5.17}$$

Let us now take a discrete model $Z^\varepsilon = (\Pi^\varepsilon, \Gamma^\varepsilon)$. A *discrete modelled distribution* is an element of the space $\mathcal{D}^{\gamma,\eta}_\varepsilon(\Gamma^\varepsilon)$, containing the maps $f : D_\varepsilon \to \widetilde{\mathcal{T}}_{<\gamma}$ such that, for any compact set $K \subseteq \mathbb{R}^2$,

$$\begin{aligned} \|\!| f_\varepsilon \|\!|^{(\varepsilon)}_{\gamma,\eta;K} := &\sup_{\substack{z \in K \cap D_\varepsilon \\ \|z\|_{\mathfrak{s}} > \varepsilon}} \sup_{m < \gamma} \frac{|f_\varepsilon(z)|_m}{\|z\|_0^{(\eta - m)\wedge 0}} \\ &+ \sup_{\substack{z, \bar{z} \in K \cap D_\varepsilon \\ \|z - \bar{z}\|_{\mathfrak{s}} > \varepsilon}} \sup_{m < \gamma} \frac{|f_\varepsilon(z) - \Gamma^\varepsilon_{z\bar{z}} f_\varepsilon(\bar{z})|_m}{\|z - \bar{z}\|_{\mathfrak{s}}^{\gamma - m}(\|z, \bar{z}\|_0 + \varepsilon)^{\eta - \gamma}} + \|\!| f_\varepsilon \|\!|_{\gamma,\eta;K;\varepsilon} < \infty, \end{aligned} \tag{5.18}$$

where the last term is define by (5.16) via $\Gamma^\varepsilon$. Sometimes it will be convenient to write $\mathcal{D}^{\gamma,\eta}_\varepsilon(Z^\varepsilon)$ for $\mathcal{D}^{\gamma,\eta}_\varepsilon(\Gamma^\varepsilon)$, and when the model is clear from the context we will omit it from the notation and will simply write $\mathcal{D}^{\gamma,\eta}_\varepsilon$. Observe that the first two terms in (5.18) are the same as in the definition of the modelled distributions in [Hai14, Def. 6.2], except that we look at the scale above $\varepsilon$. The last term measures regularity of $f$ on scale below $\varepsilon$.

For another discrete model $\bar{Z}^\varepsilon = (\bar{\Pi}^\varepsilon, \bar{\Gamma}^\varepsilon)$ and a modelled distribution $\bar{f}_\varepsilon \in \mathcal{D}^\gamma_\varepsilon(\bar{Z}^\varepsilon)$, we set

$$\begin{aligned} \|\!| f_\varepsilon; \bar{f}_\varepsilon \|\!|^{(\varepsilon)}_{\gamma,\eta;K} := &\sup_{\substack{z \in K \cap D_\varepsilon \\ \|z\|_{\mathfrak{s}} > \varepsilon}} \sup_{m < \gamma} \frac{|f_\varepsilon(z) - \bar{f}_\varepsilon(z)|_m}{\|z\|_0^{(\eta - m)\wedge 0}} \\ &+ \sup_{\substack{z, \bar{z} \in K \cap D_\varepsilon \\ \|z - \bar{z}\|_{\mathfrak{s}} > \varepsilon}} \sup_{m < \gamma} \frac{|f_\varepsilon(z) - \Gamma^\varepsilon_{z\bar{z}} f_\varepsilon(\bar{z}) - \bar{f}_\varepsilon(z) + \bar{\Gamma}^\varepsilon_{z\bar{z}} \bar{f}_\varepsilon(\bar{z})|_m}{\|z - \bar{z}\|_{\mathfrak{s}}^{\gamma - m}(\|z, \bar{z}\|_0 + \varepsilon)^{\eta - \gamma}} + \|\!| f_\varepsilon; \bar{f}_\varepsilon \|\!|_{\gamma;K;\varepsilon}, \end{aligned}$$

where the last term is defined by (5.17) via $\Gamma^\varepsilon$ and $\bar{\Gamma}^\varepsilon$.

**Remark 5.7** When we work on the compact set $K = [-T, T] \times [-2, 2]$, we simply write $\|\!| f_\varepsilon \|\!|^{(\varepsilon)}_{\gamma,\eta;T}$ and $\|\!| f_\varepsilon; \bar{f}_\varepsilon \|\!|^{(\varepsilon)}_{\gamma,\eta;T}$. The space of modelled distributions, restricted to this set $K$ we denote by $\mathcal{D}^{\gamma,\eta}_{\varepsilon,T}$.

For a discrete model $(\Pi^\varepsilon, \Gamma^\varepsilon)$ and for a modelled distribution $f_\varepsilon \in \mathcal{D}^{\gamma,\eta}_\varepsilon$ we follow [HM18, Defn. 4.5] to define a *reconstruction map* $\mathcal{R}^\varepsilon : \mathcal{D}^{\gamma,\eta}_\varepsilon \to \mathcal{X}_\varepsilon$ as

$$(\mathcal{R}^\varepsilon f_\varepsilon)(z) := (\Pi^\varepsilon_z f_\varepsilon(z))(z). \tag{5.19}$$

Assumpsions 3.6 and 3.12 in [EH19] follow readily from our definitions, and the following reconstruction theorem [EH19, Thm. 3.13] holds.

**Proposition 5.8** *Let $(\Pi^\varepsilon, \Gamma^\varepsilon)$ be a discrete model and let $f_\varepsilon \in \mathcal{D}_\varepsilon^{\gamma,\eta}(\Gamma^\varepsilon)$ be a modelled distribution, taking values in a sector of regularity $\alpha \leq 0$ and such that $\gamma > 0$, $\eta \leq \gamma$ and $\alpha \wedge \eta > -2$. Then for any compact set $K \subset \mathbb{R}^2$ one has*

$$|\iota_\varepsilon(\mathcal{R}^\varepsilon f_\varepsilon)(\varphi_z^\lambda)| \lesssim (\lambda + \varepsilon)^{\alpha\wedge\eta} \|\Pi^\varepsilon\|_{\bar K}^{(\varepsilon)} \|\!|f_\varepsilon|\!\|_{\gamma;[\varphi_z^\lambda]}^{(\varepsilon)},$$

*uniformly over $\varphi \in \mathcal{B}_\mathfrak{s}^2$, $\lambda \in (0,1]$, $z \in D_\varepsilon$, and $\varepsilon \in (0,1]$. Here, $\bar K$ is the 1-fattening of $K$, $[\varphi_z^\lambda]$ is the support of $\varphi_z^\lambda$, and we used the map (5.11). For a second discrete model $(\bar\Pi^\varepsilon, \bar\Gamma^\varepsilon)$ and for $\bar f_\varepsilon \in \mathcal{D}_\varepsilon^\gamma(\bar\Gamma^\varepsilon)$ one has*

$$|\iota_\varepsilon(\mathcal{R}^\varepsilon f_\varepsilon - \bar{\mathcal{R}}^\varepsilon \bar f_\varepsilon)(\varphi_z^\lambda)| \lesssim (\lambda + \varepsilon)^{\alpha\wedge\eta} \Big( \|\bar\Pi^\varepsilon\|_{\bar K}^{(\varepsilon)} \|\!|f_\varepsilon; \bar f_\varepsilon|\!\|_{\gamma;[\varphi_z^\lambda]}^{(\varepsilon)} + \|\Pi^\varepsilon - \bar\Pi^\varepsilon\|_{\bar K}^{(\varepsilon)} \|\!|f_\varepsilon|\!\|_{\gamma;[\varphi_z^\lambda]}^{(\varepsilon)} \Big),$$

*uniformly over the same quantities.*

## 6 Lift of martingales to a discrete model

Our next goal is to construct a lift of the martingales $\widetilde M^{\varepsilon,\mathfrak{m}}$ to a discrete model and define a solution map for the discrete equation (2.33) by analogy with the KPZ equation. We are going first to construct a discrete model $Z^{\varepsilon,\mathfrak{m}} = (\Pi^{\varepsilon,\mathfrak{m}}, \Gamma^{\varepsilon,\mathfrak{m}})$ and then its renormalisation.

By Propositions 4.2, 4.12 and 4.15 we can decompose $G^\varepsilon = K^\varepsilon + R^\varepsilon$, $\widetilde G^\varepsilon = \widetilde K^\varepsilon + \widetilde R^\varepsilon$ and $\widetilde{\widetilde G}{}^\varepsilon = \widetilde{\widetilde K}{}^\varepsilon + \widetilde{\widetilde R}{}^\varepsilon$, where $K^\varepsilon$, $\widetilde K^\varepsilon$ and $\widetilde{\widetilde K}{}^\varepsilon$ are compactly supported and contain the singularities of the kernels, while $R^\varepsilon$, $\widetilde R^\varepsilon$ and $\widetilde{\widetilde R}{}^\varepsilon$ are regular and integrable in the spatial variable.

We proceed to extend the definition of the martingale $\widetilde M^{\varepsilon,\mathfrak{m}}(t,x)$ from $\mathbb{R}_+ \times \Lambda_\varepsilon$ to $\mathbb{R} \times \Lambda_\varepsilon$. Let $\underline{\tilde h}^{\varepsilon,\mathfrak{m}}(t,x)$ be an independent copy of $\tilde h^{\varepsilon,\mathfrak{m}}(t,x)$, satisfying the discrete equation (2.33) driven by a martingale $\underline{\widetilde M}^{\varepsilon,\mathfrak{m}}(t,x)$. We denote by $\underline{\widetilde{\mathbf{C}}}_{\varepsilon,\mathfrak{m}}$ and $\underline{\widetilde C}_{\varepsilon,\mathfrak{m}}$ the respective functions (2.31) and (2.34) defined via $\underline{\tilde h}^{\varepsilon,\mathfrak{m}}$. We extend the martingale associated with $\tilde h^{\varepsilon,\mathfrak{m}}(t,x)$ to all $t \in \mathbb{R}$ by setting $\widetilde M^{\varepsilon,\mathfrak{m}}(t,x) = \underline{\widetilde M}^{\varepsilon,\mathfrak{m}}(-t,x)$ for $t < 0$. Of course, $t \mapsto \widetilde M^{\varepsilon,\mathfrak{m}}(t,x)$ is not a martingale on $\mathbb{R}$, but its projections to $\mathbb{R}_+$ and $\mathbb{R}_-$ are martingales, and in order to control moments of $\widetilde M^{\varepsilon,\mathfrak{m}}$ we can apply the Burkholder-Davis-Gundy inequality to the two projections separately. The reason to extend the martingale to all $t \in \mathbb{R}$ is technical and we have to do it because the continuous model $Z_{\mathrm{KPZ}}$ is defined via a noise on the whole time line. A particular form of extension is not very important as the noise $\widetilde M^{\varepsilon,\mathfrak{m}}(t,x)$ with $t < 0$ will not affect the solution to the discrete equation. We have chosen such a form of extension of the martingales because it was used in [GMW24] and [GMW25] and proved to be quite easy to work with. In what follows we write for $t < 0$

$$\widetilde{\mathbf{C}}_{\varepsilon,\mathfrak{m}}(t,x) = \underline{\widetilde{\mathbf{C}}}_{\varepsilon,\mathfrak{m}}(-t,x), \qquad \widetilde C_{\varepsilon,\mathfrak{m}}(t,x) = \underline{\widetilde C}_{\varepsilon,\mathfrak{m}}(-t,x). \tag{6.1}$$

Before defining the two maps $\Pi^{\varepsilon,\mathfrak{m}}$ and $\Gamma^{\varepsilon,\mathfrak{m}}$, it will be convenient to define an auxiliary linear map $\mathbf{\Pi}^{\varepsilon,\mathfrak{m}}$ which maps elements of $(\widetilde{\mathcal{W}} \cup \widetilde{\mathcal{W}}_+) \setminus \{\Xi\}$ to functions on $D_\varepsilon$. We define $\mathbf{\Pi}^{\varepsilon,\mathfrak{m}}\Xi$ by its action on any continuous compactly supported $\varphi : \mathbb{R}^2 \to \mathbb{R}$,

$$\iota_\varepsilon(\mathbf{\Pi}^{\varepsilon,\mathfrak{m}}\Xi)(\varphi) := \int_{D_\varepsilon} \varphi(z)\, d\widetilde M^{\varepsilon,\mathfrak{m}}(z),$$

where we use the integral (1.19) but over $\mathbb{R}$ in the time variable. We have to define this object in this way, because $\Xi$ represents $\partial_t \widetilde M^{\varepsilon,\mathfrak{m}}(t,x)$ which is a distribution in the time variable. All the other basis elements in $\widetilde{\mathcal{W}} \cup \widetilde{\mathcal{W}}_+$ correspond to semidiscrete functions. Let furthermore

$$(\mathbf{\Pi}^{\varepsilon,\mathfrak{m}}\mathbf{1})(z) := 1, \qquad (\mathbf{\Pi}^{\varepsilon,\mathfrak{m}} X_1)(z) := x.$$

The action of $\mathbf{\Pi}^{\varepsilon,\mathrm{m}}$ on all other elements in $\widetilde{\mathcal{W}} \cup \widetilde{\mathcal{W}}_+$ is defined recursively, using the following rules according to the definitions (5.1):

1) For the elements of the form $\widetilde{\mathcal{I}}(\tau)$ and $\widetilde{\widetilde{\mathcal{I}}}(\tau)$, and the element $\mathcal{I}(\Xi)$ we set

$$(\mathbf{\Pi}^{\varepsilon,\mathrm{m}}\mathcal{I}(\Xi))(z) = \int_{D_\varepsilon} K^\varepsilon(z-\tilde{z})(\mathbf{\Pi}^{\varepsilon,\mathrm{m}}\Xi)(\tilde{z})d\tilde{z},$$

$$(\mathbf{\Pi}^{\varepsilon,\mathrm{m}}\widetilde{\mathcal{I}}(\tau))(z) = \int_{D_\varepsilon} \widetilde{K}^\varepsilon(z-\tilde{z})(\mathbf{\Pi}^{\varepsilon,\mathrm{m}}\tau)(\tilde{z})d\tilde{z},$$

$$(\mathbf{\Pi}^{\varepsilon,\mathrm{m}}\widetilde{\widetilde{\mathcal{I}}}(\tau))(z) = \int_{D_\varepsilon} \widetilde{\widetilde{K}}^\varepsilon(z-\tilde{z})(\mathbf{\Pi}^{\varepsilon,\mathrm{m}}\tau)(\tilde{z})d\tilde{z}.$$

Since the kernel $K^\varepsilon$ is regular for $\varepsilon > 0$, the first integral on the right-hand side is well-defined in the sense of (6.4a).

2) For the maps $\partial_\pm$ we set

$$(\mathbf{\Pi}^{\varepsilon,\mathrm{m}}\partial_\pm\tau)(z) = \nabla_\varepsilon^\pm(\mathbf{\Pi}^{\varepsilon,\mathrm{m}}\tau)(z).$$

Using (5.7g) we immediately get $(\mathbf{\Pi}^{\varepsilon,\mathrm{m}}\partial\tau)(z) = \nabla_\varepsilon(\mathbf{\Pi}^{\varepsilon,\mathrm{m}}\tau)(z)$, where we use the symmetric discrete derivative (2.5).

3) For elements of the form $\tau_1\tau_2$ we set

$$(\mathbf{\Pi}^{\varepsilon,\mathrm{m}}\tau_1\tau_2)(z) = (\mathbf{\Pi}^{\varepsilon,\mathrm{m}}\tau_1)(z)(\mathbf{\Pi}^{\varepsilon,\mathrm{m}}\tau_2)(z).$$

4) For the only element $\mathcal{E}(\tau)$ we set

$$(\mathbf{\Pi}^{\varepsilon,\mathrm{m}}\mathcal{E}(\tau))(z) = \varepsilon(\mathbf{\Pi}^{\varepsilon,\mathrm{m}}\tau)(z), \qquad (\mathbf{\Pi}^{\varepsilon,\mathrm{m}}\partial_\pm\mathcal{E}(\tau))(z) = \varepsilon\nabla_\varepsilon^\pm(\mathbf{\Pi}^{\varepsilon,\mathrm{m}}\tau)(z). \tag{6.2}$$

Then for every $z \in D_\varepsilon$ we can define the map $F_z^{\varepsilon,\mathrm{m}} \in \widetilde{\mathcal{G}}$ as $F_z^{\varepsilon,\mathrm{m}} = \Gamma_{f_z^{\varepsilon,\mathrm{m}}}$, where the latter is (5.8) with the linear multiplicative functionals $f_z^{\varepsilon,\mathrm{m}} : \widetilde{\mathcal{T}}_+ \to \mathbb{R}$ given by

$$f_z^{\varepsilon,\mathrm{m}}(\tau) = -(\mathbf{\Pi}^{\varepsilon,\mathrm{m}}\tau)(z)$$

for $\tau \in \widetilde{\mathcal{W}}_+$ not of the form $\tau_1\tau_2$, and

$$f_z^{\varepsilon,\mathrm{m}}(\tau_1\tau_2) = f_z^{\varepsilon,\mathrm{m}}(\tau_1)f_z^{\varepsilon,\mathrm{m}}(\tau_2).$$

for $\tau_1\tau_2 \in \widetilde{\mathcal{W}}_+$. The functional $f_z^{\varepsilon,\mathrm{m}}$ is the same as in [Hai14, Eqs. 8.20 & 8.21], except its action on the symbols $\mathcal{E}(\tau)$ and $\partial_\pm\mathcal{E}(\tau)$. By analogy with [Hai14, Eqs. 8.22 & 8.26] we define the discrete model $Z^{\varepsilon,\mathrm{m}} = (\Pi^{\varepsilon,\mathrm{m}}, \Gamma^{\varepsilon,\mathrm{m}})$ as

$$\Gamma_{z\bar{z}}^{\varepsilon,\mathrm{m}} = (F_z^{\varepsilon,\mathrm{m}})^{-1} \circ F_{\bar{z}}^{\varepsilon,\mathrm{m}}, \qquad \Pi_z^{\varepsilon,\mathrm{m}} = \mathbf{\Pi}^{\varepsilon,\mathrm{m}}F_z^{\varepsilon,\mathrm{m}}. \tag{6.3}$$

The algebraic identities $\Gamma_{zz}^{\varepsilon,\mathrm{m}} = \mathrm{Id}$, $\Gamma_{z\bar{z}}^{\varepsilon,\mathrm{m}}\Gamma_{\bar{z}\tilde{z}}^{\varepsilon} = \Gamma_{z\tilde{z}}^{\varepsilon,\mathrm{m}}$ and $\Pi_{\bar{z}}^{\varepsilon,\mathrm{m}} = \Pi_z^{\varepsilon,\mathrm{m}}\Gamma_{z\bar{z}}^{\varepsilon,\mathrm{m}}$ follow immediately from these definitions. However, these maps do not satisfy the required analytic bounds uniformly in $\varepsilon$ and they have to be renormalised.

We do not write explicitly the action of the map $\Gamma_{z\bar{z}}^{\varepsilon,\mathrm{m}}$ on each element, because there are 48 of them. One can readily check that the map $\Pi^{\varepsilon,\mathrm{m}}$ is given by

$$\iota_\varepsilon(\Pi_z^{\varepsilon,\mathrm{m}}\Xi)(\varphi) := \int_{D_\varepsilon} \varphi(\bar{z})\, d\widetilde{M}^{\varepsilon,\mathrm{m}}(\bar{z}), \tag{6.4a}$$

and

$$(\Pi_z^{\varepsilon,\mathrm{m}}\mathbf{1})(\bar{z}) := 1, \qquad (\Pi_z^{\varepsilon,\mathrm{m}}X_1)(\bar{z}) := \bar{x} - x.$$

The action of $\Pi_z^{\varepsilon,\mathfrak{m}}$ on all other elements in $\widetilde{\mathcal{W}}$ satisfy the recursions:

$$\begin{aligned}(\Pi_z^{\varepsilon,\mathfrak{m}}\mathcal{I}(\Xi))(\bar{z}) &= \int_{D_\varepsilon}(K^\varepsilon(\bar{z}-\tilde{z})-K^\varepsilon(z-\tilde{z}))(\Pi_z^{\varepsilon,\mathfrak{m}}\Xi)(\tilde{z})d\tilde{z},\\ (\Pi_z^{\varepsilon,\mathfrak{m}}\partial_-\mathcal{I}(\Xi))(\bar{z}) &= \int_{D_\varepsilon}\nabla_\varepsilon^- K^\varepsilon(\bar{z}-\tilde{z})(\Pi_z^{\varepsilon,\mathfrak{m}}\Xi)(\tilde{z})d\tilde{z},\\ (\Pi_z^{\varepsilon,\mathfrak{m}}\partial_+\mathcal{I}(\Xi))(\bar{z}) &= \int_{D_\varepsilon}\nabla_\varepsilon^+ K^\varepsilon(\bar{z}-\tilde{z})(\Pi_z^{\varepsilon,\mathfrak{m}}\Xi)(\tilde{z})d\tilde{z}.\end{aligned} \tag{6.4b}$$

For $\widetilde{\mathcal{I}}(\tau)\in\widetilde{\mathcal{W}}$, if $|\widetilde{\mathcal{I}}(\tau)|<1$ then

$$\begin{aligned}(\Pi_z^{\varepsilon,\mathfrak{m}}\widetilde{\mathcal{I}}(\tau))(\bar{z}) &= \int_{D_\varepsilon}(\widetilde{K}^\varepsilon(\bar{z}-\tilde{z})-\widetilde{K}^\varepsilon(z-\tilde{z}))(\Pi_z^{\varepsilon,\mathfrak{m}}\tau)(\tilde{z})d\tilde{z},\\ (\Pi_z^{\varepsilon,\mathfrak{m}}\partial_-\widetilde{\mathcal{I}}(\tau))(\bar{z}) &= \int_{D_\varepsilon}\nabla_\varepsilon^-\widetilde{K}^\varepsilon(\bar{z}-\tilde{z})(\Pi_z^{\varepsilon,\mathfrak{m}}\tau)(\tilde{z})d\tilde{z},\\ (\Pi_z^{\varepsilon,\mathfrak{m}}\partial_+\widetilde{\mathcal{I}}(\tau))(\bar{z}) &= \int_{D_\varepsilon}\nabla_\varepsilon^+\widetilde{K}^\varepsilon(\bar{z}-\tilde{z})(\Pi_z^{\varepsilon,\mathfrak{m}}\tau)(\tilde{z})d\tilde{z},\end{aligned} \tag{6.4c}$$

and if $1<|\widetilde{\mathcal{I}}(\tau)|<\frac{3}{2}$ then

$$\begin{aligned}(\Pi_z^{\varepsilon,\mathfrak{m}}\widetilde{\mathcal{I}}(\tau))(\bar{z}) &= \int_{D_\varepsilon}(\widetilde{K}^\varepsilon(\bar{z}-\tilde{z})-\widetilde{K}^\varepsilon(z-\tilde{z})-(\bar{x}-x)\nabla_\varepsilon\widetilde{K}^\varepsilon(z-\tilde{z}))(\Pi_z^{\varepsilon,\mathfrak{m}}\tau)(\tilde{z})d\tilde{z},\\ (\Pi_z^{\varepsilon,\mathfrak{m}}\partial_-\widetilde{\mathcal{I}}(\tau))(\bar{z}) &= \int_{D_\varepsilon}(\nabla_\varepsilon^-\widetilde{K}^\varepsilon(\bar{z}-\tilde{z})-\nabla_\varepsilon\widetilde{K}^\varepsilon(z-\tilde{z}))(\Pi_z^{\varepsilon,\mathfrak{m}}\tau)(\tilde{z})d\tilde{z},\\ (\Pi_z^{\varepsilon,\mathfrak{m}}\partial_+\widetilde{\mathcal{I}}(\tau))(\bar{z}) &= \int_{D_\varepsilon}(\nabla_\varepsilon^+\widetilde{K}^\varepsilon(\bar{z}-\tilde{z})-\nabla_\varepsilon\widetilde{K}^\varepsilon(z-\tilde{z}))(\Pi_z^{\varepsilon,\mathfrak{m}}\tau)(\tilde{z})d\tilde{z},\end{aligned} \tag{6.4d}$$

where we use the symmetric discrete derivative (2.5). We could have used any of the discrete derivatives in the first identity in (6.4d), and we have chosen the symmetric one because it looks more natural. Because of this choice we use the symmetric map $\partial$ in the definition (5.7e). Since we want the discrete model to have property 3 in Definition 5.3, the first identity in (6.4d) yields the last two.

The action on the elements $\widetilde{\widetilde{\mathcal{I}}}(\tau)\in\widetilde{\mathcal{W}}$ are given in the same way but via the kernel $\widetilde{\widetilde{K}}^\varepsilon$. For elements $\tau_1\tau_2\in\widetilde{\mathcal{W}}$ we have

$$(\Pi_z^{\varepsilon,\mathfrak{m}}\tau_1\tau_2)(\bar{z}) = (\Pi_z^{\varepsilon,\mathfrak{m}}\tau_1)(\bar{z})(\Pi_z^{\varepsilon,\mathfrak{m}}\tau_2)(\bar{z}).$$

For the only element $\mathcal{E}(\tau)\in\widetilde{\mathcal{W}}$ we have

$$(\Pi_z^{\varepsilon,\mathfrak{m}}\mathcal{E}(\tau))(\bar{z}) := \varepsilon\Big((\Pi_z^{\varepsilon,\mathfrak{m}}\tau)(\bar{z})-(\Pi_z^{\varepsilon,\mathfrak{m}}\tau)(z)-(\bar{x}-x)\nabla_\varepsilon^-(\Pi_z^{\varepsilon,\mathfrak{m}}\tau)(\bar{z})\big|_{\bar{z}=z}\Big), \tag{6.4e}$$

where the discrete operator $\nabla_\varepsilon^-$ acts on the spatial component of the variable $\bar{z}$.

The reconstruction map $\mathcal{R}^{\varepsilon,\mathfrak{m}}$ for the discrete model is defined by (5.19).

**Remark 6.1** From (6.4d) we clearly have $(\Pi_z^{\varepsilon,\mathfrak{m}}\partial_\pm\widetilde{\mathcal{I}}(\tau))(z)\neq 0$ and $(\Pi_z^{\varepsilon,\mathfrak{m}}\partial_\pm\widetilde{\widetilde{\mathcal{I}}}(\tau))(z)\neq 0$ in the case $1<|\widetilde{\mathcal{I}}(\tau)|, |\widetilde{\widetilde{\mathcal{I}}}(\tau)|<\frac{3}{2}$. This implies that the reconstruction map does not vanish when applied to elements of positive homogeneities (the exact values of the reconstruction map applied to these elements are provided in Lemma 6.2). This will create a difficulty when we reconstruct the solution of the discrete equation (7.14).

The following lemma explains the action of the reconstruction map on the elements with positive homogeneities.

**Lemma 6.2** *For the elements of the form $\partial_\pm\widetilde{\mathcal{I}}(\bar\tau)$ and $\partial_\pm\widetilde{\widetilde{\mathcal{I}}}(\bar\tau)$ with $1 < |\widetilde{\mathcal{I}}(\bar\tau)|, |\widetilde{\widetilde{\mathcal{I}}}(\bar\tau)| < \frac{3}{2}$ one has*

$$(\mathcal{R}^{\varepsilon,\mathfrak{m}}\partial_\pm\widetilde{\mathcal{I}}(\bar\tau))(z) = \pm\frac{\varepsilon}{2}\int_{D_\varepsilon}\Delta_\varepsilon\widetilde{K}^\varepsilon(z-\tilde z)(\mathcal{R}^{\varepsilon,\mathfrak{m}}\bar\tau)(\tilde z)d\tilde z, \tag{6.5}$$

$$(\mathcal{R}^{\varepsilon,\mathfrak{m}}\partial_\pm\widetilde{\widetilde{\mathcal{I}}}(\bar\tau))(z) = \pm\frac{\varepsilon}{2}\int_{D_\varepsilon}\Delta_\varepsilon\widetilde{\widetilde{K}}^\varepsilon(z-\tilde z)(\mathcal{R}^{\varepsilon,\mathfrak{m}}\bar\tau)(\tilde z)d\tilde z, \tag{6.6}$$

*which, using* (5.7g)*, yields*

$$(\mathcal{R}^{\varepsilon,\mathfrak{m}}\partial\widetilde{\mathcal{I}}(\bar\tau))(z) = (\mathcal{R}^{\varepsilon,\mathfrak{m}}\partial\widetilde{\widetilde{\mathcal{I}}}(\bar\tau))(z) = 0. \tag{6.7}$$

*For other elements $\tau \in \widetilde{\mathcal{W}}$, such that $|\tau| > 0$, one has $(\mathcal{R}^{\varepsilon,\mathfrak{m}}\tau)(z) = 0$.*

*Proof.* The recursive definition of the discrete model implies that $(\mathcal{R}^{\varepsilon,\mathfrak{m}}\tau)(z) = 0$ for all $\tau \in \widetilde{\mathcal{W}}$, such that $|\tau| > 0$ and which are not of the form (6.5)-(6.6). For $\tau = \partial_\pm\widetilde{\mathcal{I}}(\bar\tau)$ with $1 < |\widetilde{\mathcal{I}}(\bar\tau)| < \frac{3}{2}$ we get from (6.4d)

$$(\mathcal{R}^{\varepsilon,\mathfrak{m}}\partial_\pm\widetilde{\mathcal{I}}(\bar\tau))(z) = (\Pi_z^{\varepsilon,\mathfrak{m}}\partial_\pm\widetilde{\mathcal{I}}(\bar\tau))(z) = \int_{D_\varepsilon}(\nabla_\varepsilon^\pm\widetilde{K}^\varepsilon(z-\tilde z) - \nabla_\varepsilon\widetilde{K}^\varepsilon(z-\tilde z))(\Pi_z^{\varepsilon,\mathfrak{m}}\bar\tau)(\tilde z)d\tilde z,$$

which equals (6.5) because we have $(\Pi_z^{\varepsilon,\mathfrak{m}}\bar\tau)(\tilde z) = (\mathcal{R}^{\varepsilon,\mathfrak{m}}\bar\tau)(\tilde z)$ for $\bar\tau$ in these expressions. Identity (6.6) can be proved in the same way. □

## 6.1 Renormalised discrete model

As we will see later, some of the products in the definition of the discrete model need to be renormalised if we want to obtain a non-trivial limit as $\varepsilon \to 0$. The definition of the renormalisation group, described in Section 3.3, naturally extends to the regularity structure $\widetilde{\mathscr{T}}$. However, the actual renormalisation of the discrete model is more complicated because of the two new elements 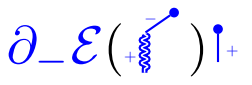 and 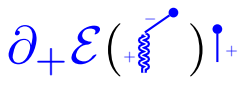, and because the product (2.17) should be renormalised by a function rather than a constant.

Let us define the linear map $L : \widetilde{\mathcal{T}} \to \mathcal{T}$ by the following rules: (1) for $\bar\tau \in \widetilde{\mathcal{W}}$ we set $L(\bar\tau) = 0$ if $\bar\tau$ contains at least one instance of $\mathcal{E}$, (2) if $\bar\tau$ does not contain $\mathcal{E}$, then $L(\bar\tau)$ is an element of $\mathcal{W}$ obtained by replacing all $\widetilde{\mathcal{I}}$ and $\widetilde{\widetilde{\mathcal{I}}}$ by $\mathcal{I}$, and all $\partial_-$ and $\partial_+$ by $\partial$. For example, using the graphical representation of the elements, we get

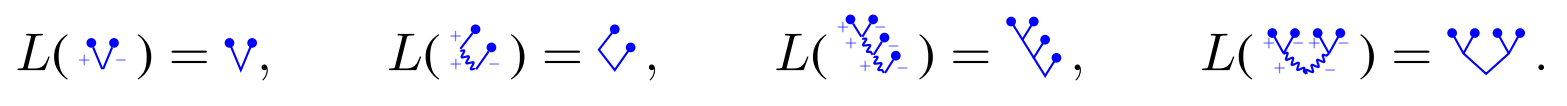

Then we write $L^{-1}$ for the map, which for every element $\tau \in \mathcal{W}$ returns the set $\{\bar\tau \in \widetilde{\mathcal{W}} : L(\bar\tau) = \tau\}$. For example,

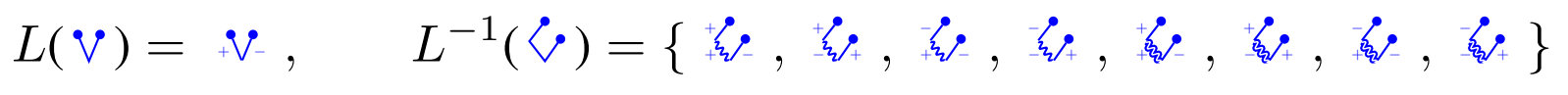

and

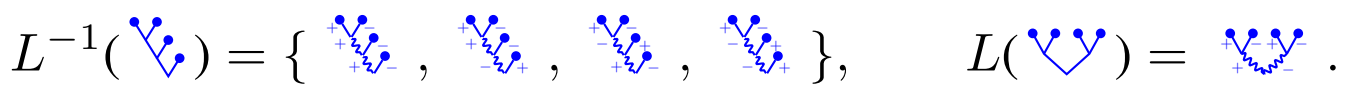

Let us denote the element $\mathcal{E}$ as a zigzag edge 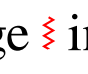 in the graphical representation of the elements of $\widetilde{\mathcal{W}}$. We draw respectively 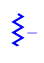 and 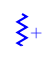 for $\partial_-\mathcal{E}$ and $\partial_+\mathcal{E}$. Then we have for example

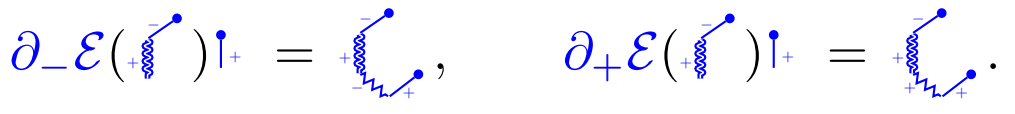

Using the set (3.12) we define

$$\widetilde{\mathcal{W}}_- := \Big\{$$ 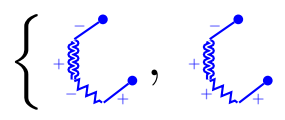 $$\Big\} \cup \bigcup_{\tau\in\mathcal{W}_-} L^{-1}(\tau)$$

and denote by $\widetilde{\mathcal{T}}_-$ a free unital algebra generated by the elements of this set. Then by analogy with the map (3.13) we define $\widetilde{\Delta}_- : \widetilde{\mathcal{T}} \to \widetilde{\mathcal{T}}_- \otimes \widetilde{\mathcal{T}}$ by

$$\widetilde{\Delta}_- \tau := \sum_{A \in \widetilde{\mathcal{N}}(\tau)} A \otimes \widetilde{\mathsf{Cut}}_A(\tau),$$

where the set $\widetilde{\mathcal{N}}(\tau)$ and the map $\widetilde{\mathsf{Cut}}_A$ is defined by analogy with the respective objects in (3.13) but using the set $\widetilde{\mathcal{W}}_-$. This maps satisfies

$$\widetilde{\Delta}_- \tau = \mathbf{1} \otimes \tau + \tau \otimes \mathbf{1}, \qquad \tau \in \Big\{ \text{[tree]}, \text{[tree]} \Big\},$$

and it acts on the other elements of $\widetilde{\mathcal{W}}_-$ by returning respective labeled expressions from Table 3.

The renormalisation group $\widetilde{\mathfrak{R}}$ for the regularity structure $\widetilde{\mathscr{T}}$ is defined as a character group of $\widetilde{\mathcal{T}}_-$, and for every $\widetilde{g} \in \widetilde{\mathfrak{R}}$ we define a linear map $\widetilde{M}_{\widetilde{g}} : \widetilde{\mathcal{T}} \to \widetilde{\mathcal{T}}$, renormalising each element $\tau \in \widetilde{\mathcal{W}}$ as

$$\widetilde{M}_{\widetilde{g}} \tau := (\widetilde{g} \otimes I) \widetilde{\Delta}_- \tau.$$

We denote the renormalisation constants by

$$\widetilde{g}(\tau) = \widetilde{C}(\tau). \tag{6.8}$$

Furthermore, given $Z^\varepsilon = (\Pi^\varepsilon, \Gamma^\varepsilon)$, we define its renormalisaiton $Z^{\varepsilon, \widetilde{g}} = (\Pi^{\varepsilon, \widetilde{g}}, \Gamma^{\varepsilon, \widetilde{g}})$, as in (3.16), to be

$$\Pi_z^{\varepsilon, \widetilde{g}} := \Pi_z^\varepsilon \widetilde{M}_{\widetilde{g}}, \qquad \Gamma_{z\bar{z}}^{\varepsilon, \widetilde{g}} := \widetilde{M}_{\widetilde{g}}^{-1} \Gamma_{z\bar{z}}^\varepsilon \widetilde{M}_{\widetilde{g}}. \tag{6.9}$$

If $\Pi^\varepsilon$ and $\Gamma^\varepsilon$ satisfy the algebraic identities from Definition 5.3, the same is clearly true for the maps $\Pi^{\varepsilon, \widetilde{g}}$ and $\Gamma^{\varepsilon, \widetilde{g}}$.

Considering the discrete model $Z^{\varepsilon, \mathfrak{m}} = (\Pi^{\varepsilon, \mathfrak{m}}, \Gamma^{\varepsilon, \mathfrak{m}})$ defined in Section 6, we need to take an $\varepsilon$-dependent element $\widetilde{g}_\varepsilon \in \widetilde{\mathfrak{R}}$, whose values are given by $\varepsilon$-dependent constants $\widetilde{C}_\varepsilon(\tau)$, as in (6.8). Then we denote by $\widetilde{Z}^{\varepsilon, \mathfrak{m}} = (\widetilde{\Pi}^{\varepsilon, \mathfrak{m}}, \widetilde{\Gamma}^{\varepsilon, \mathfrak{m}}) := Z^{\varepsilon, \widetilde{g}_\varepsilon}$ the respective renormalised model, defined in (6.9). For example, we have

$$(\widetilde{\Pi}_z^{\varepsilon, \mathfrak{m}} \text{[tree]})(\bar{z}) = (\Pi_z^{\varepsilon, \mathfrak{m}} \text{[tree]})(\bar{z}) - \widetilde{C}_\varepsilon(\text{[tree]}), \tag{6.10}$$

$$(\widetilde{\Pi}_z^{\varepsilon, \mathfrak{m}} \text{[tree]})(\bar{z}) = (\Pi_z^{\varepsilon, \mathfrak{m}} \text{[tree]})(\bar{z}) - \widetilde{C}_\varepsilon(\text{[tree]}) \tag{6.11}$$

and

$$\begin{aligned}(\widetilde{\Pi}_z^{\varepsilon, \mathfrak{m}} \text{[tree]})(\bar{z}) = (\Pi_z^{\varepsilon, \mathfrak{m}} \text{[tree]})(\bar{z}) - \widetilde{C}_\varepsilon(\text{[tree]})(\Pi_z^{\varepsilon, \mathfrak{m}} \text{[tree]})(\bar{z}) \\ - \widetilde{C}_\varepsilon(\text{[tree]})(\Pi_z^{\varepsilon, \mathfrak{m}} \text{[tree]})(\bar{z}) - \widetilde{C}_\varepsilon(\text{[tree]})(\Pi_z^{\varepsilon, \mathfrak{m}} \text{[tree]})(\bar{z}) - \widetilde{C}_\varepsilon(\text{[tree]}).\end{aligned} \tag{6.12}$$

A special role is played by the element [tree] which requires extra renormalisation by a function $\widetilde{C}_{\varepsilon, \mathfrak{m}}(\bar{z})$. More precisely, we define

$$(\widehat{\Pi}_z^{\varepsilon, \mathfrak{m}} \text{[tree]})(\bar{z}) = (\widetilde{\Pi}_z^{\varepsilon, \mathfrak{m}} \text{[tree]})(\bar{z}) - \widetilde{C}_{\varepsilon, \mathfrak{m}}(\bar{z}). \tag{6.13}$$

For the other elements $\tau \in \widetilde{\mathcal{W}} \setminus \{ \text{[tree]} \}$ we set $(\widehat{\Pi}_z^{\varepsilon, \mathfrak{m}} \tau)(\bar{z}) = (\widetilde{\Pi}_z^{\varepsilon, \mathfrak{m}} \tau)(\bar{z})$. Furthermore, we set $\widehat{\Gamma}_{z\bar{z}}^{\varepsilon, \mathfrak{m}} \tau = \widetilde{\Gamma}_{z\bar{z}}^{\varepsilon, \mathfrak{m}} \tau$ for all $\tau \in \widetilde{\mathcal{W}}$, and denote $\widehat{Z}^{\varepsilon, \mathfrak{m}} = (\widehat{\Pi}^{\varepsilon, \mathfrak{m}}, \widehat{\Gamma}^{\varepsilon, \mathfrak{m}})$.

**Lemma 6.3** *The just defined maps $\widehat{\Pi}^{\varepsilon, \mathfrak{m}}$ and $\widehat{\Gamma}^{\varepsilon, \mathfrak{m}}$ satisfy the algebraic properties of a discrete model from Definition 5.3.*

*Proof.* Since the maps $\Pi^{\varepsilon,\mathfrak{m}}$ and $\Gamma^{\varepsilon,\mathfrak{m}}$ satisfy the algebraic properties from Definition 5.3, the same if true for the maps $\widetilde{\Pi}^{\varepsilon,\mathfrak{m}}$ and $\widetilde{\Gamma}^{\varepsilon,\mathfrak{m}}$. Hence, we only need to check property 2 from Definition 5.3 for the element . This follows readily because the transformation (6.13) is independent of the base point $z$. □

The reconstruction map $\widehat{\mathcal{R}}^{\varepsilon,\mathfrak{m}}$ for $\widehat{Z}^{\varepsilon,\mathfrak{m}}$ is defined by (5.19). The next result follows readily from Lemma 6.2, because the kernels $\widetilde{K}^\varepsilon$ and $\widetilde{\widetilde{K}}^\varepsilon$ integrate to 0.

**Lemma 6.4** *For the elements of the form $\partial_\pm \widetilde{\mathcal{I}}(\tau)$ and $\partial_\pm \widetilde{\widetilde{\mathcal{I}}}(\tau)$ with $1 < |\widetilde{\mathcal{I}}(\tau)|, |\widetilde{\widetilde{\mathcal{I}}}(\tau)| < \frac{3}{2}$ one has*

$$(\widehat{\mathcal{R}}^{\varepsilon,\mathfrak{m}} \partial_\pm \widetilde{\mathcal{I}}(\tau))(z) = \pm\frac{\varepsilon}{2} \int_{D_\varepsilon} \nabla_\varepsilon^+ \widetilde{K}^\varepsilon(\bar{z} - \tilde{z}) \nabla_\varepsilon^- (\widehat{\mathcal{R}}^{\varepsilon,\mathfrak{m}} \tau)(\tilde{z}) d\tilde{z},$$

$$(\widehat{\mathcal{R}}^{\varepsilon,\mathfrak{m}} \partial_\pm \widetilde{\widetilde{\mathcal{I}}}(\tau))(z) = \pm\frac{\varepsilon}{2} \int_{D_\varepsilon} \nabla_\varepsilon^+ \widetilde{\widetilde{K}}^\varepsilon(\bar{z} - \tilde{z}) \nabla_\varepsilon^- (\widehat{\mathcal{R}}^{\varepsilon,\mathfrak{m}} \tau)(\tilde{z}) d\tilde{z},$$

*while for other elements $\tau \in \widetilde{\mathcal{W}}$, such that $|\tau| > 0$, one has $(\widehat{\mathcal{R}}^{\varepsilon,\mathfrak{m}} \tau)(z) = 0$.*

As we show below, the maps $\widehat{Z}^{\varepsilon,\mathfrak{m}} = (\widehat{\Pi}^{\varepsilon,\mathfrak{m}}, \widehat{\Gamma}^{\varepsilon,\mathfrak{m}})$ converge as $\varepsilon \to 0$ in respective topologies, if we make the following choice of the renormalisation constants and function. The random function in (6.13) is taken to be

$$\widetilde{C}_{\varepsilon,\mathfrak{m}}(t,x) = \int_{D_\varepsilon} \nabla_\varepsilon^- K^\varepsilon(t-s, x-y) \nabla_\varepsilon^+ K^\varepsilon(t-s, x-y) (\widetilde{\mathbf{C}}_{\varepsilon,\mathfrak{m}}(s,y) - \nu^\varepsilon) ds, \qquad (6.14)$$

which is defined by analogy with (2.18) but through the function $\widetilde{\mathbf{C}}_{\varepsilon,\mathfrak{m}}$ defined in (2.31) and (6.1). As in (2.18), we remove the constant term $\nu^\varepsilon$ from this function because it will be subtracted as the renormalisation constant $\widetilde{C}_\varepsilon(\ )$ in the definition of the renormalised model (6.10).

To define the renormalisation constants, we are going to use a graphical representation of convolutions of kernels. More precisely, the vertex "●" labelled with $z$ represents a fixed point $z \in D_\varepsilon$. In the following diagrams we take $z = 0$ and we do not write this label. The arrow "⟶" represents the discrete kernel $K^\varepsilon$, and this arrow labeled by "+" or "−" corresponds to $\nabla_\varepsilon^+ K^\varepsilon$ or $\nabla_\varepsilon^- K^\varepsilon$ respectively. The kernels $\widetilde{K}^\varepsilon$ and $\widetilde{\widetilde{K}}^\varepsilon$ are represented by the arrows "⇝" and "⇝⇝" respectively, and we label them in the same way to represent the discrete derivatives. By the node "•" we denote a variable integrated out in $D_\varepsilon$.

For each element $\tau$ in the basis of the regularity structure, let $[\tau]$ be the number of occurrences of $\Xi$ in the symbolic representation of $\tau$. Then we choose the renormalisation constants to be

$$\widetilde{C}_\varepsilon(\tau) = (\nu^\varepsilon)^{[\tau]/2} C_\varepsilon(\tau), \qquad (6.15)$$

where $\nu^\varepsilon$ is the constant term in the bracket process (2.13), which appears exactly $[\tau]/2$ times in the renormalisation constant because there exactly $[\tau]/2$ pairings of the noises, and

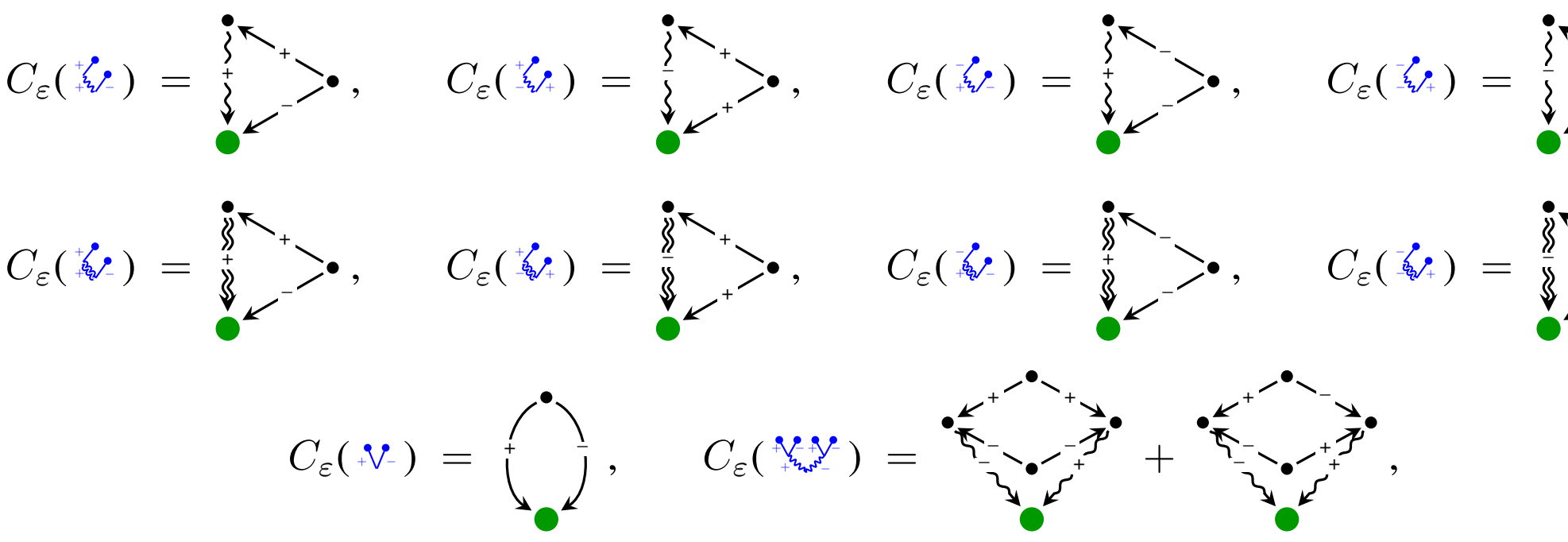

(6.16)

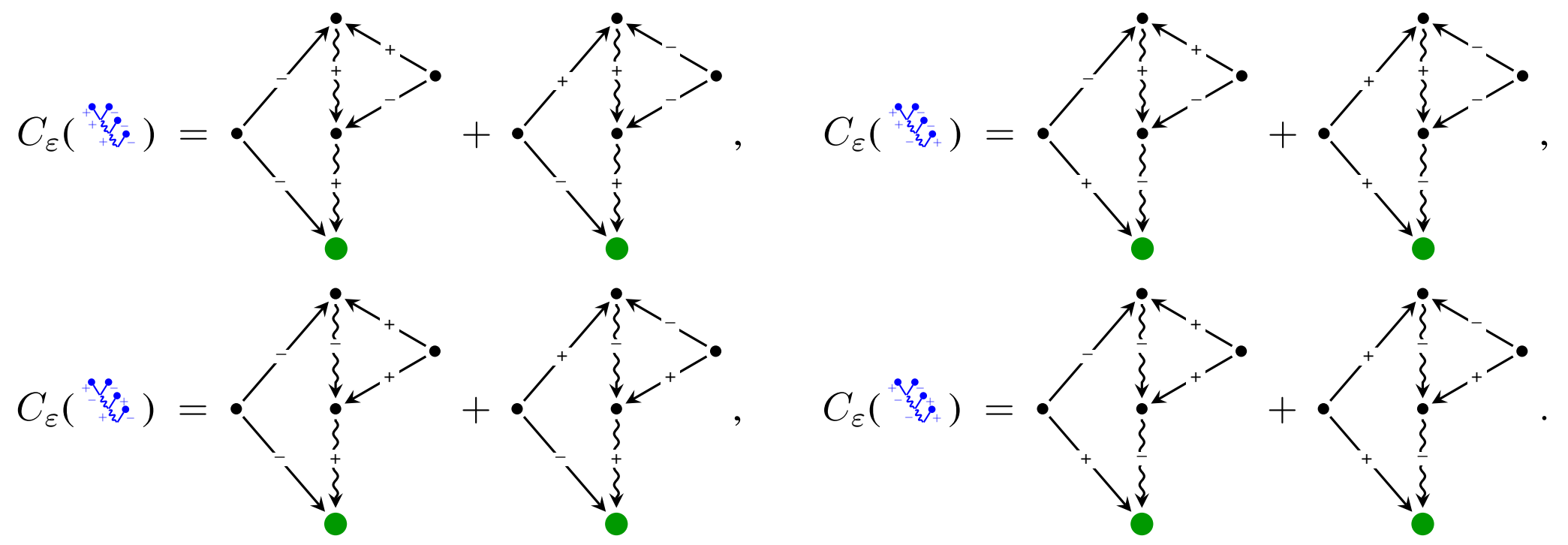

Finally, we define the renormalisation constants for the last two elements of $\widetilde{\mathcal{W}}$ to be

$$\begin{aligned} C_\varepsilon(\;) &= \varepsilon \iint_{D_\varepsilon^2} \nabla_\varepsilon^- \nabla_\varepsilon^+ \widetilde{\widetilde{K}}^\varepsilon(z_1 - z) \nabla_\varepsilon^- K^\varepsilon(z_2 - z_1) \nabla_\varepsilon^+ K^\varepsilon(z_2 - z) dz_1 dz_2, \\ C_\varepsilon(\;) &= \varepsilon \iint_{D_\varepsilon^2} \nabla_\varepsilon^+ \nabla_\varepsilon^+ \widetilde{\widetilde{K}}^\varepsilon(z_1 - z) \nabla_\varepsilon^- K^\varepsilon(z_2 - z_1) \nabla_\varepsilon^+ K^\varepsilon(z_2 - z) dz_1 dz_2, \end{aligned} \tag{6.17}$$

where the discrete derivatives are applied to the spatial variables. We note that these constant are independent of $z$, by a simple variable change in the integrals.

## 6.2 Properties of the renormalisation constants

We state relations between the renormalisation constants (6.16) in the following lemmas.

**Lemma 6.5** *The constant $C_\varepsilon(\;)$ is bounded uniformly in $\varepsilon \in (0,1)$.*

*Proof.* If we replace at least one of the kernels $K^\varepsilon$ in the definition of this renormalisation constant by $K^\varepsilon - G^\varepsilon$, then the resultant constant is uniformly bounded in $\varepsilon$. Hence, we can write

$$C_\varepsilon(\;) = \int_{D_\varepsilon} \nabla_\varepsilon^- G^\varepsilon(z) \nabla_\varepsilon^+ G^\varepsilon(z)\, dz + \mathcal{O}(1).$$

The integrand is the function $Q^\varepsilon(z)$, defined in (2.16), and by Lemma 4.4 the integral vanishes. □

Throughout the following proofs we use the functions

$$f_\varepsilon(k) := \varepsilon^{-2}(1 - \cos(\pi\varepsilon k)), \qquad m_\varepsilon^+(k) := \varepsilon^{-1}(e^{\pi i k \varepsilon} - 1), \qquad m_\varepsilon^-(k) := \varepsilon^{-1}(1 - e^{-\pi i k \varepsilon}) \tag{6.18}$$

and the Dirichlet kernel

$$D_N(x) = \sum_{k=-N}^{N} e^{ikx} = \frac{\sin((N+1/2)x)}{\sin(x/2)}, \tag{6.19}$$

where the last identity holds at $x = 0$ by the L'Hopital's rule.

**Lemma 6.6** *The renormalisation constant*

$$C_\varepsilon(\;) + C_\varepsilon(\;) + C_\varepsilon(\;) + C_\varepsilon(\;) + C_\varepsilon(\;) \tag{6.20}$$

*is bounded uniformly in $\varepsilon$.*

**Remark 6.7** In [Hai13, Lem. 6.5], such a cancellation property is observed for the first time and shown for the case where the regularisation is done by cutting off high Fourier modes of a Gaussian space-time white noise in the KPZ equation; similarly in [GP17, Sec. 9.6]. Since we have a very different regularisation (i.e. lattice discretisation of the heat kernel, and jump martingales instead of Gaussian white noise), we have to prove this cancellation property anew.

*Proof of Lemma 6.6.* From the definitions (6.16) and the results from Section 7.1 in [HM18] we conclude that each constant in (6.20) diverges at most logarithmically as $\varepsilon \to 0$. If we moreover replace at least one of the kernels $K^\varepsilon$ or $\widetilde{K}^\varepsilon$ by $K^\varepsilon - P^\varepsilon$ or $\widetilde{K}^\varepsilon - \widetilde{P}^\varepsilon$, where the periodic kernels are defined in (4.11) and (4.42), then the respective renormalisation constant becomes uniformly bounded in $\varepsilon$. Hence, it is sufficient to consider (6.20) where the renormalisation constants are defined as before but via the kernels $P^\varepsilon$ and $\widetilde{P}^\varepsilon$.

Let us first analyse the described modification of the constant $\widetilde{C}_\varepsilon($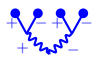$)$. For two signs $\mathfrak{l}_1, \mathfrak{l}_2 \in \{-,+\}$ we define the function

$$Q^\varepsilon_{\mathfrak{l}_1,\mathfrak{l}_2}(z_1, z_2) := \int_{\mathbb{R}_+ \times \mathbb{T}_\varepsilon} \nabla^{\mathfrak{l}_1}_\varepsilon P^\varepsilon(z_1 - z) \nabla^{\mathfrak{l}_2}_\varepsilon P^\varepsilon(z_2 - z) dz. \tag{6.21}$$

This function depends on the difference $z_1 - z_2$ and we denote $Q^\varepsilon_{\mathfrak{l}_1,\mathfrak{l}_2}(z_1 - z_2) := Q^\varepsilon_{\mathfrak{l}_1,\mathfrak{l}_2}(z_1, z_2)$. Then

$$\widetilde{C}_\varepsilon(\text{tree}) = \iint_{(\mathbb{R}_+ \times \mathbb{T}_\varepsilon)^2} \nabla^-_\varepsilon \widetilde{P}^\varepsilon(-z_1) Q^{\varepsilon,1}(z_1 - z_2) \nabla^+_\varepsilon \widetilde{P}^\varepsilon(-z_2) dz_1 dz_2 + \iint_{(\mathbb{R}_+ \times \mathbb{T}_\varepsilon)^2} \nabla^-_\varepsilon \widetilde{P}^\varepsilon(-z_1) Q^{\varepsilon,2}(z_1 - z_2) \nabla^+_\varepsilon \widetilde{P}^\varepsilon(-z_2) dz_1 dz_2 \tag{6.22}$$

with the kernels $Q^{\varepsilon,1}(z_1 - z_2) := Q^\varepsilon_{+,+}(z_1 - z_2) Q^\varepsilon_{-,-}(z_1 - z_2)$ and $Q^{\varepsilon,2}(z_1 - z_2) := Q^\varepsilon_{+,-}(z_1 - z_2) Q^\varepsilon_{-,+}(z_1 - z_2)$.

One can conclude from Lemma 4.8 that replacing each $\widetilde{P}^\varepsilon$ by $P^\varepsilon$ in (6.22) gives an error term bounded uniformly in $\varepsilon$. We denote by $C^{(1)}_\varepsilon$ and $C^{(2,i)}_\varepsilon$, for $i = 1, \ldots, 4$, the five renormalisation constants in (6.20) in the definitions of which (6.16) the kernels $K^\varepsilon$ and $\widetilde{K}^\varepsilon$ are replaced by $P^\varepsilon$.

Let us start with preliminary computations. We will compute convolutions (6.21), which appears in every definition of the renormalisation constants. Using the property that $\widehat{g}(k) = \widehat{f}(-k)$ when $g(x) = f(-x)$, we have that

$$\widehat{Q^\varepsilon_{\mathfrak{l}_1,\mathfrak{l}_2}}(t_1 - t_2, k) = \int_{-\infty}^{t_1 \wedge t_2} \widehat{\nabla^{\mathfrak{l}_1}_\varepsilon P^\varepsilon}(t_1 - s, k) \widehat{\nabla^{\mathfrak{l}_2}_\varepsilon P^\varepsilon}(t_2 - s, -k) ds.$$

We denote throughout the proof by $\widehat{g}(k)$ the discrete Fourier transform. We readily get from (4.12)

$$\widehat{\nabla^\pm_\varepsilon P^\varepsilon}(t, k) = m^\pm_\varepsilon(k) e^{-t f_\varepsilon(k)}, \qquad k \in \mathbb{T}_N, \tag{6.23}$$

where we use the functions (6.18). We note that $m^\pm_\varepsilon(0) = 0$ and hence

$$\begin{aligned} \widehat{Q^\varepsilon_{\mathfrak{l}_1,\mathfrak{l}_2}}(t_1 - t_2, k) &= \mathbf{1}_{k \neq 0} \int_{-\infty}^{t_1 \wedge t_2} m^{\mathfrak{l}_1}_\varepsilon(k) m^{\mathfrak{l}_2}_\varepsilon(-k) e^{-(t_1 + t_2 - 2s) f_\varepsilon(k)} ds \\ &= \mathbf{1}_{k \neq 0} \frac{m^{\mathfrak{l}_1}_\varepsilon(k) m^{\mathfrak{l}_2}_\varepsilon(-k)}{2 f_\varepsilon(k)} e^{-|t_1 - t_2| f_\varepsilon(k)}. \end{aligned} \tag{6.24}$$

We have $m^{\mathfrak{l}}_\varepsilon(k) m^{\mathfrak{l}}_\varepsilon(-k) = 2 f_\varepsilon(k)$ for $\mathfrak{l} \in \{-,+\}$, which yields

$$\widehat{Q^\varepsilon_{\mathfrak{l},\mathfrak{l}}}(t_1 - t_2, k) = \mathbf{1}_{k \neq 0} e^{-|t_1 - t_2| f_\varepsilon(k)}. \tag{6.25}$$

Let us start with the renormalisation constant $C^{(1)}_\varepsilon$. It is defined in (6.16) as a sum of two constants which we denote by $C^{(1)}_{\varepsilon,1}$ and $C^{(1)}_{\varepsilon,2}$ respectively. We then have

$$C^{(1)}_{\varepsilon,1} = \iint_{(\mathbb{R}_+\times\mathbb{T}_\varepsilon)^2} \nabla^-_\varepsilon \widetilde{P}^\varepsilon(-z_1) Q^{\varepsilon,1}(z_1-z_2)\nabla^+_\varepsilon \widetilde{P}^\varepsilon(-z_2) dz_1 dz_2,$$

where $Q^{\varepsilon,1}(z_1-z_2) := Q^\varepsilon_{+,+}(z_1-z_2)Q^\varepsilon_{-,-}(z_1-z_2)$. Applying Parseval's identity to the sum with respect to $x_1$ and using (6.23) yields

$$\begin{aligned}C^{(1)}_{\varepsilon,1} &= \frac{1}{2}\int_{-\infty}^0\int_{-\infty}^0 \sum_{k\in\mathbb{T}_N} \overline{\widehat{\nabla^-_\varepsilon P^\varepsilon}(-t_1,-k)}\widehat{Q^{\varepsilon,1}}(t_1-t_2,k)\widehat{\nabla^+_\varepsilon P^\varepsilon}(-t_2,-k)dt_1dt_2\\ &= \frac{1}{2}\int_{-\infty}^0\int_{-\infty}^0 \sum_{k\in\mathbb{T}_N} m^-_\varepsilon(k)m^+_\varepsilon(-k)\widehat{Q^{\varepsilon,1}}(t_1-t_2,k)e^{(t_1+t_2)f_\varepsilon(k)}dt_1dt_2,\end{aligned}$$

where we used $\overline{m^-_\varepsilon(-k)} = m^-_\varepsilon(k)$. Using (6.25) we can compute

$$\widehat{Q^{\varepsilon,1}}(t_1-t_2,k) = \frac{1}{2}\sum_{\substack{\ell\in\mathbb{T}^\star_N\\ \ell\neq k}} \widehat{Q^\varepsilon_{-,-}}(t_1-t_2,\ell)\widehat{Q^\varepsilon_{+,+}}(t_1-t_2,k-\ell) = \frac{1}{2}\sum_{\substack{\ell\in\mathbb{T}^\star_N\\ \ell\neq k}} e^{-|t_1-t_2|(f_\varepsilon(\ell)+f_\varepsilon(k-\ell))},$$

where $\mathbb{T}^\star_N := \mathbb{T}_N\setminus\{0\}$. Plugging this into the preceding identity, we write $C^{(1)}_{\varepsilon,1}$ as

$$\frac{1}{4}\int_{-\infty}^0\int_{-\infty}^0 \sum_{k\neq\ell\in\mathbb{T}^\star_N} m^-_\varepsilon(k)m^+_\varepsilon(-k)e^{-|t_1-t_2|(f_\varepsilon(\ell)+f_\varepsilon(k-\ell))+(t_1+t_2)f_\varepsilon(k)}dt_1dt_2,$$

where the sum over $k$ can be restricted to $\mathbb{T}^\star_N$ because at $k=0$ the functions $m^\pm_\varepsilon(k)$ vanish. The integrand is symmetric with respect to $t_1$ and $t_2$, and this expression equals

$$\begin{aligned}C^{(1)}_{\varepsilon,1} &= \frac{1}{2}\int_{-\infty}^0\int_{-\infty}^{t_2} \sum_{k\neq\ell\in\mathbb{T}^\star_N} m^-_\varepsilon(k)m^+_\varepsilon(-k)e^{-(t_2-t_1)(f_\varepsilon(\ell)+f_\varepsilon(k-\ell))+(t_1+t_2)f_\varepsilon(k)}dt_1dt_2\\ &= \frac{1}{2}\sum_{k\neq\ell\in\mathbb{T}^\star_N} \frac{m^-_\varepsilon(k)m^+_\varepsilon(-k)}{f_\varepsilon(k)+f_\varepsilon(k-\ell)+f_\varepsilon(\ell)}\int_{-\infty}^0 e^{2t_2f_\varepsilon(k)}dt_2\\ &= \frac{1}{4}\sum_{k\neq\ell\in\mathbb{T}^\star_N} \frac{m^-_\varepsilon(k)m^+_\varepsilon(-k)}{f_\varepsilon(k)(f_\varepsilon(k)+f_\varepsilon(k-\ell)+f_\varepsilon(\ell))}.\end{aligned} \tag{6.26}$$

Recalling that $\varepsilon = \frac{2}{2N+1}$, we have $\varepsilon|k|<1$ for $|k|\leq N$. Then from the definitions of the functions $f_\varepsilon$ and $m^\pm_\varepsilon$ we get

$$|f_\varepsilon(k) - \tfrac{\pi^2k^2}{2}| \leq c_0\varepsilon^\theta|k|^{\theta+2}, \qquad |m^\pm_\varepsilon(k) - \pi i k| \leq c_0\varepsilon^\theta|k|^{\theta+1}, \tag{6.27}$$

for all $|k|\leq N$, any $\theta\in(0,1]$ and some $c_0>0$ depending on $\theta$. Applying these bounds to (6.26), we get

$$C^{(1)}_{\varepsilon,1} = \frac{1}{2\pi^2}\sum_{k\neq\ell\in\mathbb{T}^\star_N}\frac{1}{k^2+\ell^2-k\ell} + R^{(1)}_{\varepsilon,1} \tag{6.28}$$

with a remainder $R^{(1)}_{\varepsilon,1}$ bounded uniformly in $\varepsilon$.

Now, we will compute the constant $C^{(1)}_{\varepsilon,2}$. We have

$$C^{(1)}_{\varepsilon,2} = \iint_{(\mathbb{R}_+\times\mathbb{T}_\varepsilon)^2} \nabla^-_\varepsilon P^\varepsilon(-z_1)Q^{\varepsilon,2}(z_1-z_2)\nabla^+_\varepsilon P^\varepsilon(-z_2)dz_1dz_2,$$

where $Q^{\varepsilon,2}(z_1-z_2) := Q^{\varepsilon}_{+,-}(z_1-z_2)Q^{\varepsilon}_{-,+}(z_1-z_2)$. Using (6.24) we can compute

$$\widehat{Q^{\varepsilon,2}}(t_1-t_2,k) = \frac{1}{2}\sum_{\ell\in\mathbb{T}_N} \widehat{Q^{\varepsilon}_{+,-}}(t_1-t_2,\ell)\widehat{Q^{\varepsilon}_{-,+}}(t_1-t_2,k-\ell)$$
$$= \frac{1}{8}\sum_{\substack{\ell\in\mathbb{T}_N^\star\\ \ell\neq k}} \frac{m_\varepsilon^+(\ell)m_\varepsilon^-(-\ell)m_\varepsilon^-(k-\ell)m_\varepsilon^+(\ell-k)}{f_\varepsilon(\ell)f_\varepsilon(k-\ell)} e^{-|t_1-t_2|(f_\varepsilon(\ell)+f_\varepsilon(k-\ell))}.$$

Repeating computations as in (6.28), we get

$$C^{(1)}_{\varepsilon,2} = \frac{1}{16}\sum_{k\neq\ell\in\mathbb{T}_N^\star} \frac{m_\varepsilon^-(k)m_\varepsilon^+(-k)m_\varepsilon^+(\ell)m_\varepsilon^-(-\ell)m_\varepsilon^-(k-\ell)m_\varepsilon^+(\ell-k)}{f_\varepsilon(k)f_\varepsilon(\ell)f_\varepsilon(k-\ell)(f_\varepsilon(k)+f_\varepsilon(k-\ell)+f_\varepsilon(\ell))}.$$

Using bounds (6.27), we can write it as

$$C^{(1)}_{\varepsilon,2} = \frac{1}{2\pi^2}\sum_{k\neq\ell\in\mathbb{T}_N^\star} \frac{1}{k^2+\ell^2-k\ell} + R^{(1)}_{\varepsilon,2}, \tag{6.29}$$

where the remainder $R^{(1)}_{\varepsilon,2}$ is bounded uniformly in $\varepsilon$.

We write each of the first four constants $C^{(2,i)}_\varepsilon$ in (6.20) as $C^{(2,i)}_{\varepsilon,1} + C^{(2,i)}_{\varepsilon,2}$ for $i=1,\dots,4$. We have

$$C^{(2,1)}_{\varepsilon,1} = \iint_{(\mathbb{R}_+\times\mathbb{T}_\varepsilon)^2} \nabla_\varepsilon^+ P^\varepsilon(-z_1)Q^{\varepsilon,3}(z_1-z_2)Q^{\varepsilon}_{-,-}(-z_2)dz_1dz_2,$$

where $Q^{\varepsilon,3}(z_1-z_2) := \nabla_\varepsilon^+ P^\varepsilon(z_1-z_2)Q^{\varepsilon}_{-,+}(z_1-z_2)$. Applying Parseval's identity to the sum with respect to $x_1$, and using (6.23) and (6.24) yields

$$C^{(2,1)}_{\varepsilon,1} = \frac{1}{2}\int_{-\infty}^0\int_{-\infty}^0 \sum_{k\in\mathbb{T}_N} \overline{\widehat{\nabla_\varepsilon^+ P^\varepsilon}(-t_1,-k)}\,\widehat{Q^{\varepsilon,3}}(t_1-t_2,k)\,\widehat{Q^{\varepsilon}_{-,-}}(-t_2,-k)\,dt_1dt_2$$
$$= \frac{1}{4}\int_{-\infty}^0\int_{-\infty}^0 \sum_{k\in\mathbb{T}_N^\star} \frac{m_\varepsilon^+(k)m_\varepsilon^-(-k)m_\varepsilon^-(k)}{f_\varepsilon(k)}\widehat{Q^{\varepsilon,3}}(t_1-t_2,k)\,e^{(t_1+t_2)f_\varepsilon(k)}\,dt_1dt_2,$$

where we used $\overline{m_\varepsilon^+(-k)} = m_\varepsilon^+(k)$. Furthermore, using (6.23) and (6.24) we get

$$\widehat{Q^{\varepsilon,3}}(t_1-t_2,k) = \frac{1}{2}\sum_{\ell\in\mathbb{T}_N} \widehat{\nabla_\varepsilon^+ P^\varepsilon}(t_1-t_2,\ell)\widehat{Q^{\varepsilon}_{-,+}}(t_1-t_2,k-\ell)$$
$$= \frac{1}{4}\sum_{\ell\in\mathbb{T}_N^\star,\ell\neq k} \frac{m_\varepsilon^+(\ell)m_\varepsilon^-(k-\ell)m_\varepsilon^+(\ell-k)}{f_\varepsilon(k-\ell)} e^{-(t_1-t_2)(f_\varepsilon(\ell)+f_\varepsilon(k-\ell))}$$

for $t_1-t_2\geq 0$, and this function vanishes for $t_1-t_2<0$. Then the constant $C^{(2,1)}_{\varepsilon,1}$ equals

$$\frac{1}{16}\int_{-\infty}^0\int_{-\infty}^{t_1} \sum_{k\neq\ell\in\mathbb{T}_N^\star} \frac{m_\varepsilon^+(k)m_\varepsilon^-(-k)m_\varepsilon^-(k)m_\varepsilon^+(\ell)m_\varepsilon^-(k-\ell)m_\varepsilon^+(\ell-k)}{f_\varepsilon(k)f_\varepsilon(k-\ell)}$$
$$\times\, e^{-(t_1-t_2)(f_\varepsilon(\ell)+f_\varepsilon(k-\ell))+(t_1+t_2)f_\varepsilon(k)}\,dt_2dt_1$$
$$= \frac{1}{16}\int_{-\infty}^0 \sum_{k\neq\ell\in\mathbb{T}_N^\star} \frac{m_\varepsilon^+(k)m_\varepsilon^-(-k)m_\varepsilon^-(k)m_\varepsilon^+(\ell)m_\varepsilon^-(k-\ell)m_\varepsilon^+(\ell-k)}{f_\varepsilon(k)f_\varepsilon(k-\ell)(f_\varepsilon(k)+f_\varepsilon(k-\ell)+f_\varepsilon(\ell))} e^{2t_1f_\varepsilon(k)}\,dt_1$$
$$= \frac{1}{32}\sum_{k\neq\ell\in\mathbb{T}_N^\star} \frac{m_\varepsilon^+(k)m_\varepsilon^-(-k)m_\varepsilon^-(k)m_\varepsilon^+(\ell)m_\varepsilon^-(k-\ell)m_\varepsilon^+(\ell-k)}{f_\varepsilon(k)^2f_\varepsilon(k-\ell)(f_\varepsilon(k)+f_\varepsilon(k-\ell)+f_\varepsilon(\ell))}.$$

When computing the other constants, we get the same expression but with different signs of the functions $m_\varepsilon$. Namely, we get

$$C^{(2,1)}_{\varepsilon,2} = \frac{1}{32} \sum_{k\neq\ell\in\mathbb{T}^\star_N} \frac{m^-_\varepsilon(-k)m^+_\varepsilon(k)^2 m^+_\varepsilon(\ell)m^-_\varepsilon(k-\ell)m^-_\varepsilon(\ell-k)}{f_\varepsilon(k)^2 f_\varepsilon(k-\ell)(f_\varepsilon(k)+f_\varepsilon(k-\ell)+f_\varepsilon(\ell))},$$

$$C^{(2,2)}_{\varepsilon,1} = \frac{1}{32} \sum_{k\neq\ell\in\mathbb{T}^\star_N} \frac{m^-_\varepsilon(k)^2 m^+_\varepsilon(-k)m^+_\varepsilon(\ell)m^-_\varepsilon(k-\ell)m^+_\varepsilon(\ell-k)}{f_\varepsilon(k)^2 f_\varepsilon(k-\ell)(f_\varepsilon(k)+f_\varepsilon(k-\ell)+f_\varepsilon(\ell))},$$

$$C^{(2,2)}_{\varepsilon,2} = \frac{1}{32} \sum_{k\neq\ell\in\mathbb{T}^\star_N} \frac{m^-_\varepsilon(k)m^+_\varepsilon(-k)m^+_\varepsilon(k)m^+_\varepsilon(\ell)m^-_\varepsilon(k-\ell)m^-_\varepsilon(\ell-k)}{f_\varepsilon(k)^2 f_\varepsilon(k-\ell)(f_\varepsilon(k)+f_\varepsilon(k-\ell)+f_\varepsilon(\ell))},$$

$$C^{(2,3)}_{\varepsilon,1} = \frac{1}{32} \sum_{k\neq\ell\in\mathbb{T}^\star_N} \frac{m^+_\varepsilon(k)m^-_\varepsilon(-k)m^-_\varepsilon(k)m^-_\varepsilon(\ell)m^+_\varepsilon(k-\ell)m^+_\varepsilon(\ell-k)}{f_\varepsilon(k)^2 f_\varepsilon(k-\ell)(f_\varepsilon(k)+f_\varepsilon(k-\ell)+f_\varepsilon(\ell))},$$

$$C^{(2,3)}_{\varepsilon,2} = \frac{1}{32} \sum_{k\neq\ell\in\mathbb{T}^\star_N} \frac{m^+_\varepsilon(k)^2 m^-_\varepsilon(-k)m^-_\varepsilon(\ell)m^+_\varepsilon(k-\ell)m^-_\varepsilon(\ell-k)}{f_\varepsilon(k)^2 f_\varepsilon(k-\ell)(f_\varepsilon(k)+f_\varepsilon(k-\ell)+f_\varepsilon(\ell))},$$

$$C^{(2,4)}_{\varepsilon,1} = \frac{1}{32} \sum_{k\neq\ell\in\mathbb{T}^\star_N} \frac{m^-_\varepsilon(k)^2 m^+_\varepsilon(-k)m^-_\varepsilon(\ell)m^+_\varepsilon(k-\ell)m^+_\varepsilon(\ell-k)}{f_\varepsilon(k)^2 f_\varepsilon(k-\ell)(f_\varepsilon(k)+f_\varepsilon(k-\ell)+f_\varepsilon(\ell))},$$

$$C^{(2,4)}_{\varepsilon,2} = \frac{1}{32} \sum_{k\neq\ell\in\mathbb{T}^\star_N} \frac{m^-_\varepsilon(k)m^+_\varepsilon(-k)m^+_\varepsilon(k)m^-_\varepsilon(\ell)m^+_\varepsilon(k-\ell)m^-_\varepsilon(\ell-k)}{f_\varepsilon(k)^2 f_\varepsilon(k-\ell)(f_\varepsilon(k)+f_\varepsilon(k-\ell)+f_\varepsilon(\ell))}.$$

Using the bounds (6.27), we can write these constants as

$$C^{(2,i)}_{\varepsilon,j} = -\frac{1}{4\pi^2} \sum_{k\neq\ell\in\mathbb{T}^\star_N} \frac{\ell}{k(k^2+\ell^2-k\ell)} + R^{(2,i)}_{\varepsilon,j}, \tag{6.30}$$

where the remainder $R^{(2,i)}_{\varepsilon,j}$ is bounded uniformly in $\varepsilon$.

From (6.28), (6.29) and (6.30) we get

$$C^{(1)}_{\varepsilon,1} + C^{(1)}_{\varepsilon,2} + \sum_{i=1}^{4}\sum_{j=1,2} C^{(2,i)}_{\varepsilon,j} = \frac{1}{\pi^2} \sum_{k\neq\ell\in\mathbb{T}^\star_N} \frac{1-2\ell/k}{k^2+\ell^2-k\ell} + R_\varepsilon, \tag{6.31}$$

where $R_\varepsilon$ is bounded uniformly in $\varepsilon$. Including the term $k=\ell$ into the sum, we get

$$\frac{1}{\pi^2}\sum_{k\in\mathbb{T}^\star_N}\frac{1}{k^2} - \frac{1}{\pi^2}\sum_{k,\ell\in\mathbb{T}^\star_N}\frac{1-2\ell/k}{k^2+\ell^2-k\ell} + R^\varepsilon,$$

We can write the first term in (6.31) as

$$\frac{1}{\pi^2}\sum_{k\in\mathbb{T}^\star_N}\frac{1}{k^2} = \frac{1}{\pi^2}\sum_{k\in\mathbb{Z}\setminus\{0\}}\frac{1}{k^2} - \frac{1}{\pi^2}\sum_{|k|>N}\frac{1}{k^2} = \frac{1}{3} + \mathcal{O}(\varepsilon),$$

because $\mathcal{O}(N^{-1}) = \mathcal{O}(\varepsilon)$. Furthermore, by renaming the variables $(k,\ell)\mapsto(\ell,k)$ we write the second term in (6.31) as

$$\frac{1}{2\pi^2}\sum_{k,\ell\in\mathbb{T}^\star_N}\frac{1-2\ell/k}{k^2+\ell^2-k\ell} + \frac{1}{2\pi^2}\sum_{k,\ell\in\mathbb{T}^\star_N}\frac{1-2k/\ell}{k^2+\ell^2-k\ell} = -\frac{1}{\pi^2}\sum_{k,\ell\in\mathbb{T}^\star_N}\frac{1}{k\ell} = -\frac{1}{\pi^2}\Big(\sum_{k\in\mathbb{T}^\star_N}\frac{1}{k}\Big)^2.$$

The sum in the parentheses vanishes and we conclude that (6.31) is bounded uniformly in $\varepsilon$. □

**Lemma 6.8** *One has*

$$C_\varepsilon([\text{tree}]) + C_\varepsilon([\text{tree}]) = C_\varepsilon([\text{tree}]) + C_\varepsilon([\text{tree}]) = C_\varepsilon([\text{tree}]) + C_\varepsilon([\text{tree}]) = C_\varepsilon([\text{tree}]) + C_\varepsilon([\text{tree}]) = 0 \quad (6.32)$$

*and moreover* $C_\varepsilon([\text{tree}]) + C_\varepsilon([\text{tree}])$ *and* $C_\varepsilon([\text{tree}]) + C_\varepsilon([\text{tree}])$ *are bounded uniformly in* $\varepsilon$.

*Proof.* Let us start with proving $C_\varepsilon([\text{tree}]) + C_\varepsilon([\text{tree}]) = 0$. We recall that both $K^\varepsilon(t,x)$ and $\widetilde{K}^\varepsilon(t,x)$ are even functions in the spatial variable $x$, which yields $\nabla_\varepsilon^+ K^\varepsilon(t,x) = -\nabla_\varepsilon^- K^\varepsilon(t,-x)$ and the same for $\widetilde{K}^\varepsilon$. Recalling the definition of the renormalisation constants (6.16), we get

$$\begin{aligned} C_\varepsilon([\text{tree}]) &= \iint_{D_\varepsilon^2} \nabla_\varepsilon^- \widetilde{K}^\varepsilon(-\bar z)\nabla_\varepsilon^- K^\varepsilon(\bar z - z)\nabla_\varepsilon^+ K^\varepsilon(-z)dz d\bar z \\ &= -\iint_{D_\varepsilon^2} \nabla_\varepsilon^+ \widetilde{K}^\varepsilon(-\bar t, \bar x)\nabla_\varepsilon^+ K^\varepsilon(\bar t - t, x - \bar x)\nabla_\varepsilon^- K^\varepsilon(-t,x)dz d\bar z \\ &= -\iint_{D_\varepsilon^2} \nabla_\varepsilon^+ \widetilde{K}^\varepsilon(-\bar z)\nabla_\varepsilon^+ K^\varepsilon(\bar z - z)\nabla_\varepsilon^- K^\varepsilon(-z)dz d\bar z, \end{aligned}$$

where in the last equality we made the change of variables $x \mapsto -x$ and $\bar x \mapsto -\bar x$. The resulting expression equals $-C_\varepsilon([\text{tree}])$, as claimed. The other identities in (6.32) can be proved in the same way.

Now, we will prove the second part of this lemma. As in the proof of Lemma 6.6, we can replace all instances of $K^\varepsilon$, $\widetilde{K}^\varepsilon$ and $\widetilde{\widetilde{K}}^\varepsilon$ with $P^\varepsilon$ in the definitions of the renormalisation constant $C_\varepsilon([\text{tree}])$, $C_\varepsilon([\text{tree}])$, $C_\varepsilon([\text{tree}])$ and $C_\varepsilon([\text{tree}])$, up to error terms bounded uniformly in $\varepsilon$. We denote such modified constants by $C_\varepsilon^{(i)}$, for $i = 1, \dots, 4$, respectively. Then we clearly have $C_\varepsilon^{(3)} = C_\varepsilon^{(1)}$ and $C_\varepsilon^{(4)} = C_\varepsilon^{(2)}$, and we should prove that $C_\varepsilon^{(1)} + C_\varepsilon^{(2)}$ is bounded uniformly in $\varepsilon$.

We will perform calculations similar to those in the proof of Lemma 6.6. It will be convenient to use the kernels (6.21) and we recall that they satisfy $Q^\varepsilon_{\mathfrak{l}_1,\mathfrak{l}_2}(z_1 - z_2) := Q^\varepsilon_{\mathfrak{l}_1,\mathfrak{l}_2}(z_1, z_2)$. Then Parceval's identity and the Fourier transforms (6.23) and (6.24) yield

$$\begin{aligned} C_\varepsilon^{(1)} &= \int_{\mathbb{R}_+\times\mathbb{T}_\varepsilon} \nabla_\varepsilon^- P^\varepsilon(-\bar z) Q^\varepsilon_{-,+}(\bar z) d\bar z = \frac{1}{2}\int_{-\infty}^0 \sum_{k\in\mathbb{T}_N} \overline{\widehat{\nabla_\varepsilon^- P^\varepsilon}(-\bar t, -k)} \widehat{Q^\varepsilon_{-,+}}(\bar t, k) d\bar t \quad (6.33) \\ &= \frac{1}{2}\int_{-\infty}^0 \sum_{k\in\mathbb{T}_N^\star} \frac{m_\varepsilon^-(k)^2 m_\varepsilon^+(-k)}{2 f_\varepsilon(k)} e^{2\bar t f_\varepsilon(k)} d\bar t = -\sum_{k\in\mathbb{T}_N^\star} \frac{m_\varepsilon^-(k)^2 m_\varepsilon^+(-k)}{8 f_\varepsilon(k)^2}, \end{aligned}$$

where we used $\overline{m_\varepsilon^-(-k)} = m_\varepsilon^-(k)$ in the second line. In the same way we compute

$$\begin{aligned} C_\varepsilon^{(2)} &= \int_{\mathbb{R}_+\times\mathbb{T}_\varepsilon} \nabla_\varepsilon^- P^\varepsilon(-\bar z) Q^\varepsilon_{+,-}(\bar z) d\bar z = \frac{1}{2}\int_{-\infty}^0 \sum_{k\in\mathbb{T}_N} \overline{\widehat{\nabla_\varepsilon^- P^\varepsilon}(-\bar t, -k)} \widehat{Q^\varepsilon_{+,-}}(\bar t, k) d\bar t \\ &= \frac{1}{2}\int_{-\infty}^0 \sum_{k\in\mathbb{T}_N^\star} \frac{m_\varepsilon^-(k) m_\varepsilon^+(k) m_\varepsilon^-(-k)}{2 f_\varepsilon(k)} e^{2\bar t f_\varepsilon(k)} d\bar t = -\sum_{k\in\mathbb{T}_N^\star} \frac{m_\varepsilon^-(k) m_\varepsilon^+(k) m_\varepsilon^-(-k)}{8 f_\varepsilon(k)^2}. \end{aligned}$$

Hence, we have

$$C_\varepsilon^{(1)} + C_\varepsilon^{(2)} = -\sum_{k\in\mathbb{T}_N^\star} \frac{m_\varepsilon^-(k)(m_\varepsilon^-(k) m_\varepsilon^+(-k) + m_\varepsilon^+(k) m_\varepsilon^-(-k))}{8 f_\varepsilon(k)^2}.$$

We note that $m_\varepsilon^-(k) m_\varepsilon^+(-k) = 2 f_\varepsilon(k) e^{-\pi i k \varepsilon}$ and $m_\varepsilon^+(k) m_\varepsilon^-(-k) = 2 f_\varepsilon(k) e^{\pi i k \varepsilon}$. Moreover, $m_\varepsilon^-(k) m_\varepsilon^-(-k) = 2 f_\varepsilon(k)$. Then the preceding expression equals

$$-\sum_{k\in\mathbb{T}_N^\star} \frac{m_\varepsilon^-(k)(e^{-\pi i k\varepsilon} + e^{\pi i k \varepsilon})}{4 f_\varepsilon(k)} = -\sum_{k\in\mathbb{T}_N^\star} \frac{e^{-\pi i k\varepsilon} + e^{\pi i k \varepsilon}}{2 m_\varepsilon^-(-k)} = -\frac{\varepsilon i}{2} \sum_{k\in\mathbb{T}_N^\star} e^{-\pi i k \varepsilon/2} \frac{\cos(\pi k \varepsilon)}{\sin(\pi k \varepsilon/2)},$$

where we used $e^{-\pi ik\varepsilon} + e^{\pi ik\varepsilon} = 2\cos(\pi k\varepsilon)$ and $m_\varepsilon^-(-k) = -2i\varepsilon^{-1}e^{\pi ik\varepsilon/2}\sin(\pi k\varepsilon/2)$. Writing $e^{-\pi ik\varepsilon/2} = \cos(\pi k\varepsilon/2) - i\sin(\pi k\varepsilon/2)$, we note that the real part of the sum over $k$ vanishes because the summands are odd with respect to $k$. Hence, the preceding expression equals $-\frac{\varepsilon}{2}\sum_{k\in\mathbb{T}_N^\star}\cos(\pi k\varepsilon)$. This expression can be written in terms of the Dirichlet kernel (6.19) as

$$-\frac{\varepsilon}{2}(D_N(\pi\varepsilon) - 1) = \frac{\varepsilon}{2},$$

where we have used $D_N(\pi\varepsilon) = 0$. □

**Lemma 6.9** *The renormalisation constants* (6.17) *are bounded uniformly in* $\varepsilon$.

*Proof.* As in the proof of the previous lemma, we can replace $K^\varepsilon$ and $\widetilde{K}^\varepsilon$ in (6.17) with $P^\varepsilon$, up to an error bounded uniformly in $\varepsilon$. We denote the two modified constants by $C_\varepsilon^{(1)}$ and $C_\varepsilon^{(2)}$. Then we get as in (6.33)

$$\begin{aligned}C_\varepsilon^{(1)} &= \varepsilon\int_{\mathbb{R}_+\times\mathbb{T}_\varepsilon}(\nabla_\varepsilon^+)^2P^\varepsilon(-\bar z)Q_{-,+}^\varepsilon(\bar z)d\bar z = \frac{\varepsilon}{2}\int_{-\infty}^0\sum_{k\in\mathbb{T}_N}\overline{\widehat{(\nabla_\varepsilon^+)^2P^\varepsilon}(-\bar t,-k)\widehat{Q_{-,+}^\varepsilon}(\bar t,k)}d\bar t\\ &= \frac{\varepsilon}{4}\int_{-\infty}^0\sum_{k\in\mathbb{T}_N^\star}\frac{m_\varepsilon^+(k)^2m_\varepsilon^-(k)m_\varepsilon^+(-k)}{f_\varepsilon(k)}e^{2\bar t f_\varepsilon(k)}d\bar t = -\frac{\varepsilon}{8}\sum_{k\in\mathbb{T}_N^\star}\frac{m_\varepsilon^+(k)^2m_\varepsilon^-(k)m_\varepsilon^+(-k)}{f_\varepsilon(k)^2}.\end{aligned}$$

Using $m_\varepsilon^-(k)m_\varepsilon^+(-k) = 2f_\varepsilon(k)e^{-\pi ik\varepsilon}$ and $m_\varepsilon^+(k)^2 = -2f_\varepsilon(k)e^{\pi ik\varepsilon}$, the expression inside the sum is exactly 1 and the preceding expression equals $\frac{\varepsilon}{2}|\mathbb{T}_N^\star| = \frac{2N}{2N+1} < 1$.

Performing similar computations, we get

$$C_\varepsilon^{(2)} = \varepsilon\int_{\mathbb{R}_+\times\mathbb{T}_\varepsilon}\nabla_\varepsilon^-\nabla_\varepsilon^+P^\varepsilon(-\bar z)Q_{-,+}^\varepsilon(\bar z)d\bar z = -\frac{\varepsilon}{8}\sum_{k\in\mathbb{T}_N^\star}\frac{m_\varepsilon^-(k)m_\varepsilon^+(k)m_\varepsilon^-(k)m_\varepsilon^+(-k)}{f_\varepsilon(k)^2}.$$

We use $m_\varepsilon^-(k)m_\varepsilon^+(-k) = 2f_\varepsilon(k)e^{-\pi ik\varepsilon}$ and $m_\varepsilon^-(k)m_\varepsilon^+(k) = -2f_\varepsilon(k)$ to write this expression as

$$C_\varepsilon^{(2)} = \frac{\varepsilon}{2}\sum_{k\in\mathbb{T}_N^\star}e^{-\pi ik\varepsilon} = \frac{\varepsilon}{2}(D_N(-\pi\varepsilon) - 1) = -\frac{\varepsilon}{2},$$

where we used the Dirichlet kernel (6.19) and $D_N(-\pi\varepsilon) = 0$. □

## 6.3 Properties of the martingales

The renormalised discrete models is defined in terms of stochastic integrals with respect to the martingales $\widetilde{M}^{\varepsilon,\mathfrak{m}}$, and we are going to use the results of [GMW24] to bound the model. For this, we need the following result:

**Lemma 6.10** *The martingales* $\widetilde{M}^{\varepsilon,\mathfrak{m}}$ *satisfy [GMW24, Assm. 1].*

*Proof.* From the definition of the height function we conclude that the function $\widetilde{\mathbf{C}}_{\varepsilon,\mathfrak{m}}$, defined in Section 2.2, is absolutely bounded by a constant independent of $\varepsilon \in (0,1]$. This yields property (1) in [GMW24, Assm. 1]. Property (2) in [GMW24, Assm. 1] follows from the definition of the martingales, and property (3) in [GMW24, Assm. 1] is satisfied with parameter $\mathbf{k} = \frac{1}{2}$ (the size of a jump of the martingale is proportional to $\varepsilon^{\frac{1}{2}}$). Finally, from equation (2.33) we get

$$\widetilde{M}^{\varepsilon,\mathfrak{m}}(t,x) = \widetilde{\mathfrak{J}}^{\varepsilon,\mathfrak{m}}(t,x) - \varepsilon^{-\frac{3}{2}}\int_0^t\widetilde{\mathfrak{C}}_{\varepsilon,\mathfrak{m}}(s,x)ds, \qquad (6.34)$$

where $t \mapsto \widetilde{\mathsf{J}}^{\varepsilon,\mathfrak{m}}(t,x)$ is a pure jump process, and $\widetilde{\mathsf{C}}_{\varepsilon,\mathfrak{m}}(t,x) = \widetilde{\mathsf{C}}_{\varepsilon}(t,x)$ for $t < \tau_{\varepsilon,\mathfrak{m}}$, where

$$\begin{aligned}\widetilde{\mathsf{C}}_{\varepsilon}(t,x) := \frac{1}{2}\varepsilon^{\frac{3}{2}}\big((1+\sqrt{\varepsilon})\Delta_\varepsilon \tilde{h}^\varepsilon(t,x) - \nabla_\varepsilon^- \tilde{h}^\varepsilon(t,x)\nabla_\varepsilon^+ \tilde{h}^\varepsilon(t,x) \\ - 2\varrho_\varepsilon \varepsilon^{-\frac{1}{2}}\nabla_\varepsilon \tilde{h}^\varepsilon(t,x) - \varrho_\varepsilon^2 \varepsilon^{-1}\big),\end{aligned} \tag{6.35}$$

and

$$\widetilde{\mathsf{C}}_{\varepsilon,\mathfrak{m}}(t,x) = \frac{1}{2}\varepsilon^{\frac{3}{2}}\nu^\varepsilon \Delta_\varepsilon \tilde{h}^{\varepsilon,\mathfrak{m}}_{\mathrm{sym}}(t,x) \tag{6.36}$$

for $t \geq \tau_{\varepsilon,\mathfrak{m}}$ where we use the constant (2.8). The function $\widetilde{\mathsf{C}}_{\varepsilon,\mathfrak{m}}(t,x)$ is absolutely bounded uniformly in $\varepsilon$. This yields property (4) in [GMW24, Assm. 1]. □

We also know from [GMW24, Lem. 2.3] that the martingales $\bar{M}^{\varepsilon,\mathfrak{m}}(t,x) = \varepsilon^{-\frac{1}{2}}([\widetilde{M}^{\varepsilon,\mathfrak{m}}(x)]_t - \langle \widetilde{M}^{\varepsilon,\mathfrak{m}}(x)\rangle_t)$ satisfy the same [GMW24, Assm. 1] with $\bar{\mathsf{C}}_{\varepsilon,\mathfrak{m}}(t,x)$ replacing $\widetilde{\mathsf{C}}_{\varepsilon,\mathfrak{m}}(t,x)$. Hence, we can apply the moment bounds proved in [GMW24] to multiple stochastic integrals driven by the martingales $\widetilde{M}^{\varepsilon,\mathfrak{m}}$ and $\bar{M}^{\varepsilon,\mathfrak{m}}$.

### 6.4 Moment bounds for the renormalised discrete model

The main result of this section is the following moment bounds for the renormalised discrete model $\widehat{Z}^{\varepsilon,\mathfrak{m}}$ defined in Section 6.1. Moreover, we introduce a new discrete model $\widehat{Z}^{\varepsilon,\mathfrak{m},\delta}$, defined as $\widehat{Z}^{\varepsilon,\mathfrak{m}}$ but via mollified martingales. This new model will be important in our proof of Theorem 1.1 via intermediate processes.

Let $\varrho : \mathbb{R}^2 \to \mathbb{R}$ be a symmetric smooth function, supported on the ball of radius 1 (with respect to the parabolic distance $\|\bullet\|_{\mathfrak{s}}$) and satisfying $\int_{\mathbb{R}^2}\varrho(z)dz = 1$. For any $\delta \in (0,1)$ we define its rescaling

$$\varrho_\delta(t,x) := \frac{1}{\delta^3}\varrho\Big(\frac{t}{\delta^2}, \frac{x}{\delta}\Big).$$

We need to modify this function in a way that its integral over $D_\varepsilon$ becomes 1. For this, we approximate the function by its local averages as

$$\varrho_{\delta,\varepsilon}(t,x) := \varepsilon^{-1}\int_{y\in\mathbb{R}:|y-x|_\infty\leq\varepsilon/2}\varrho_\delta(t,y)dy, \tag{6.37}$$

which satisfies $\int_{D_\varepsilon}\varrho_{\delta,\varepsilon}(z)dz = 1$. Then we regularise the martingales $\widetilde{M}^{\varepsilon,\mathfrak{m}}$ in the following way:

$$\zeta^{\varepsilon,\mathfrak{m},\delta}(t,x) := \varepsilon\sum_{y\in\Lambda_\varepsilon}\int_{\mathbb{R}}\varrho_{\delta,\varepsilon}(t-s,x-y)\,d\widetilde{M}^{\varepsilon,\mathfrak{m}}(s,y). \tag{6.38}$$

The process $\zeta^{\varepsilon,\mathfrak{m},\delta}(t,x)$ is defined on $(t,x)\in\mathbb{R}\times\Lambda_\varepsilon$ and not a martingale anymore, but a convolution with this process can be interpreted as a stochastic integral. For example, a space-time convolution with the kernel $K^\varepsilon$ may be written as

$$\frac{1}{\sqrt{2}}\varepsilon\sum_{y\in\Lambda_\varepsilon}\int_{\mathbb{R}}K^{\varepsilon,\delta}_{t-s}(x-y)\,d\widetilde{M}^{\varepsilon,\mathfrak{m}}(s,y),$$

where

$$K^{\varepsilon,\delta} := K^\varepsilon \star_\varepsilon \varrho_{\delta,\varepsilon}. \tag{6.39}$$

We define a renormalised model $\widehat{Z}^{\varepsilon,\mathfrak{m},\delta}$ in the same way as we defined $\widehat{Z}^{\varepsilon,\mathfrak{m}}$ in Section 6.1, but we use the random process $\zeta^{\varepsilon,\mathfrak{m},\delta}(t,x)dt$ everywhere in place of the martingale $d\widetilde{M}^{\varepsilon,\mathfrak{m}}(t,x)$. For example, we have

$$(\widehat{\Pi}^{\varepsilon,\mathfrak{m},\delta}_z \curlyvee)(\bar{z}) = (\widetilde{\Pi}^{\varepsilon,\mathfrak{m},\delta}_z \curlyvee)(\bar{z}) - \widetilde{C}_{\varepsilon,\delta}(\curlyvee) - \widetilde{C}_{\varepsilon,\mathfrak{m},\delta}(\bar{z}),$$

with the renormalisation function

$$\widetilde{C}_{\varepsilon,\mathfrak{m},\delta}(t,x) := \int_{D_\varepsilon} \nabla_\varepsilon^- K^{\varepsilon,\delta}(t-s,x-y)\nabla_\varepsilon^+ K^{\varepsilon,\delta}(t-s,x-y)(\widetilde{\mathbf{C}}_{\varepsilon,\mathfrak{m}}(s,y)-\nu^\varepsilon)ds \tag{6.40}$$

and the constant

$$\widetilde{C}_{\varepsilon,\delta}(\cdot) := \nu^\varepsilon \int_{D_\varepsilon} \nabla_\varepsilon^+ K^{\varepsilon,\delta}(z)\nabla_\varepsilon^- K^{\varepsilon,\delta}(z)\,dz, \tag{6.41}$$

where we use the constant (2.8). The renormalisation constants $\widetilde{C}_{\varepsilon,\mathfrak{m},\delta}(\tau)$ involved in this definition are defined by (6.15), (6.16) and (6.17), where each kernel $K^\varepsilon$ is replaced by $K^{\varepsilon,\delta}$.

The following proposition provides moment bounds on these renormalised discrete models.

**Proposition 6.11** *There exists $\varepsilon_0 > 0$ and $\theta > 0$ such that for any $p \geq 1$ and $T > 0$ there is $C > 0$ such that the renormalised discrete models satisfy*

$$\sup_{\varepsilon\in(0,\varepsilon_0)} \mathbb{E}\Big[(\|\widehat{Z}^{\varepsilon,\mathfrak{m}}\|_T^{(\varepsilon)})^p\Big] \leq C, \qquad \sup_{\varepsilon\in(0,\varepsilon_0)} \mathbb{E}\Big[(\|\widehat{Z}^{\varepsilon,\mathfrak{m}};\widehat{Z}^{\varepsilon,\mathfrak{m},\delta}\|_T^{(\varepsilon)})^p\Big] \leq C\delta^{\theta p}, \tag{6.42}$$

*for any $\delta \in (0,1)$. The constant $C$ depends polynomially on $\mathfrak{m}$.*

We start with proving moment bounds on the map $\widehat{\Pi}^{\varepsilon,\mathfrak{m}}\tau$, where $\tau$ are the basis elements of $\widetilde{\mathcal{W}}$ satisfying $|\tau| < 0$. We denote the set of such $\tau$ by $\widetilde{\mathcal{W}}_-$. A proof of moment bound for all basis elements of $\widetilde{\mathcal{T}}$ is then a standard procedure and we prefer not to go into details. One can find an analogous proof in a discrete setting in [GMW25, Prop. 7.1].

**Proposition 6.12** *There are constants $\bar\kappa > 0$, $\varepsilon_0 > 0$ and $\theta > 0$, such that for any $\tau \in \widetilde{\mathcal{W}}_-$, $p \geq 1$ and $T > 0$ there is $C > 0$ for which we have the bounds*

$$\Big(\mathbb{E}|\iota_\varepsilon(\widehat{\Pi}_z^{\varepsilon,\mathfrak{m}}\tau)(\varphi_z^\lambda)|^p\Big)^{\frac{1}{p}} \leq C(\lambda\vee\varepsilon)^{|\tau|+\bar\kappa}, \tag{6.43}$$

$$\Big(\mathbb{E}|\iota_\varepsilon(\widehat{\Pi}_z^{\varepsilon,\mathfrak{m}}\tau - \widehat{\Pi}_z^{\varepsilon,\mathfrak{m},\delta}\tau)(\varphi_z^\lambda)|^p\Big)^{\frac{1}{p}} \leq C\delta^\theta(\lambda\vee\varepsilon)^{|\tau|+\bar\kappa-\theta}, \tag{6.44}$$

*uniformly in $z \in D_\varepsilon$, $\lambda \in (0,1]$, $\varphi \in \mathcal{B}^2_{\mathfrak{s}}$, $\delta \in (0,1)$ and $\varepsilon \in (0,\varepsilon_0)$. The constant $C$ depends polynomially on $\mathfrak{m}$.*

We prove these bounds for any $p$ sufficiently large in the rest of this section. The bounds for any $p \geq 1$ follow then by Hölder's inequality. These bounds will be proved using the following method. For every symbol $\tau \in \widetilde{\mathcal{W}}_-$, we use the definition of the renormalised discrete model in Section 6.1 and the expansion [GMW24, Eq. 2.16] to write $\iota_\varepsilon(\widehat{\Pi}_z^{\varepsilon,\mathfrak{m}}\tau)(\varphi_z^\lambda)$ as a sum of terms of the form

$$\begin{aligned}\int_{D_\varepsilon} \varphi_z^\lambda(\bar z)\Big(\int_{D_\varepsilon^n} F_{\bar z}(z_1,\ldots,z_n)\,d\widetilde{\mathbf{M}}^n_{\varepsilon,\mathfrak{m}}(z_1,\ldots,z_n)\Big)d\bar z \\ = \int_{D_\varepsilon^n}\Big(\int_{D_\varepsilon} \varphi_z^\lambda(\bar z)F_{\bar z}(z_1,\ldots,z_n)\,d\bar z\Big)\,d\widetilde{\mathbf{M}}^n_{\varepsilon,\mathfrak{m}}(z_1,\ldots,z_n),\end{aligned} \tag{6.45}$$

where the measure $\widetilde{\mathbf{M}}^n_{\varepsilon,\mathfrak{m}}$ is defined in Section 2.1 in [GMW24] for the martingales $\widetilde{M}^{\varepsilon,\mathfrak{m}}$, and a function $F$ of $n$ space-time variables. Similarly, we write $\iota_\varepsilon(\widehat{\Pi}_z^{\varepsilon,\mathfrak{m},\delta}\tau)(\varphi_z^\lambda)$ as a sum of terms of the form

$$\int_{D_\varepsilon} \varphi_z^\lambda(\bar z)\Big(\int_{D_\varepsilon^n} F_{\bar z}(z_1,\ldots,z_n)\,d\widetilde{\mathbf{M}}^n_{\varepsilon,\mathfrak{m},\delta}(z_1,\ldots,z_n)\Big)d\bar z \tag{6.46}$$

$$= \int_{D_\varepsilon^n} \left( \int_{D_\varepsilon} \varphi_z^\lambda(\bar z)(F_{\bar z} \star_\varepsilon \varrho_{\delta,\varepsilon})(z_1, \dots, z_n)\, d\bar z \right) d\widetilde{\mathbf{M}}^n_{\varepsilon,\mathfrak{m}}(z_1, \dots, z_n),$$

where $\widetilde{\mathbf{M}}^n_{\varepsilon,\mathfrak{m},\delta}(z_1, \dots, z_n)$ stays for the product measure associated to the regularised martingales $\zeta^{\varepsilon,\mathfrak{m},\delta}$, defined in (6.38), and $\varrho_{\delta,\varepsilon}$ is the mollifier from (6.38). The functions $F$ will be defined in terms of the kernels $K^\varepsilon$ and $\widetilde{K}^\varepsilon$, so that the function $F_{\bar z} \star_\varepsilon \varrho_{\delta,\varepsilon}$ will be written in terms of $K^{\varepsilon,\delta}$ and $\widetilde{K}^\varepsilon$, where we use the kernel (6.39). To bound the terms (6.45) and their difference with those in (6.46), we are going to use a general moment bound for multiple stochastic integrals proved in Theorem C.1.

We are going to provide complete proofs of the bounds (6.43)-(6.44) only for the elements

$\tau \in \{$ 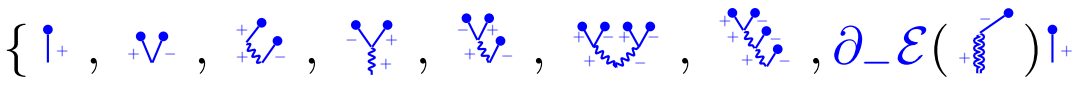 $\}$ (6.47)

while the bounds for the other elements in $\widetilde{\mathcal{W}}_-$ can be proved similarly, since the variation is only in $\pm$ signs which do not make much difference for these bounds. Each of the following section is devoted to one of these elements. To represent the map $\widehat{\Pi}^{\varepsilon,\mathfrak{m}}_z\tau$ as a stochastic integral we use the graphical notation introduced in Section 6.1. Moreover, in what follows, the vertex "●" labelled with $z$ represents the basis point $z \in D_\varepsilon$ and the arrow " ⟶ " represents a test function $\varphi_z^\lambda$. Each edge corresponding to a kernel is labeled by a triplet $(a, r, s)$, where the pair $(a, r)$ corresponds to the labels on graphs as described in [GMW24, Sec. 4] (more precisely, $a$ is the order of singularity of the kernel and $r$ is the order of its positive or negative renormalisation) and $s$ is either "+" or "−", corresponding to one of the discrete derivatives (2.1) applied to the kernel. In particular, Propositions 4.2 and 4.12 imply that the kernels $K^\varepsilon$ and $\widetilde{K}^\varepsilon$ have orders of singularities 1, and hence $\nabla_\varepsilon^\pm K^\varepsilon$ and $\nabla_\varepsilon^\pm \widetilde{K}^\varepsilon$ will be labeled by the triplets $(2, 0, \pm)$. Each variable $z_i$, integrated with respect to the measure $\widetilde{\mathbf{M}}^n_{\varepsilon,\mathfrak{m}}$ with $n \geq 2$ is denoted by a node "●"; the variable integrated with respect to the martingale $\widetilde{M}^{\varepsilon,\mathfrak{m}}$ we denote by "○".

The bound (6.44) on the difference of two maps can be proved in exactly the same way as we prove (6.43), but taking into account that the bound on the kernel $K^{\varepsilon,\delta} - K^\varepsilon$ is of order $\delta^\theta$, for $\theta > 0$ small. We prefer not to go into more details.

### 6.4.1 Bound on the first element in (6.47)

For the element $\tau =$ 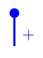 we have

$\iota_\varepsilon(\widehat{\Pi}^{\varepsilon,\mathfrak{m}}_z\tau)(\varphi_z^\lambda) =$ 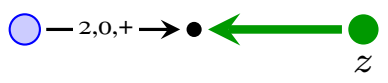 .

For this case, the constant [GMW24, Eq. 4.15] is $\nu_\gamma = -\frac{1}{2}$. The only choice in the first term on the right of (C.1) is $\mathbf{p}(1) = 2$ with $m = 1$, $\alpha_\gamma(\mathbf{p}) = \delta_\gamma(\mathbf{p}) = 0$; and the only choice in the second term on the right of (C.1) is $\mathbf{p}(1) = \infty$ with $\alpha_\gamma(\mathbf{p}) = \delta_\gamma(\mathbf{p}) = \frac{3}{2}$. Hence, for any $p \geq 2$ and $\theta > 0$ sufficiently small, by Theorem C.1 (taking $\theta_1 = \theta_2 = \theta$, $\mathfrak{e} = \varepsilon$),

$$\left( \mathbb{E} |\iota_\varepsilon(\widehat{\Pi}^{\varepsilon,\mathfrak{m}}_z\tau)(\varphi_z^\lambda)|^p \right)^{\frac{1}{p}} \lesssim (\lambda \vee \varepsilon)^{-\frac{1}{2}} + (\lambda \vee \varepsilon)^{-\frac{1}{2}-\theta} \lesssim (\lambda \vee \varepsilon)^{-\frac{1}{2}-\theta}.$$

This is the required bound (6.43) if $\theta$ is small enough because $|\tau| = -\frac{1}{2} - \kappa$.

We are going to use the constant $\nu_\gamma$, defined in [GMW24, Eq. 4.15], without reference.

### 6.4.2 Bound on the second element in (6.47)

Analysis of the element $\tau =$ 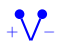 is more involved because of the renormalisation (6.13). For this, we need to introduce more notation for the diagrams representing this element. Let the vertex "▼" correspond to a variable $(t, x)$ integration against the predictable quadratic variation

$(t,x) \mapsto \langle \widetilde{M}^{\varepsilon,\mathfrak{m}}(x)\rangle_t$, and let the vertex "◇" correspond to a variable $(t,x)$ integrated agains the family of martingales $(t,x) \mapsto [\widetilde{M}^{\varepsilon,\mathfrak{m}}(x)]_t - \langle \widetilde{M}^{\varepsilon,\mathfrak{m}}(x)\rangle_t$. Then our definitions yield

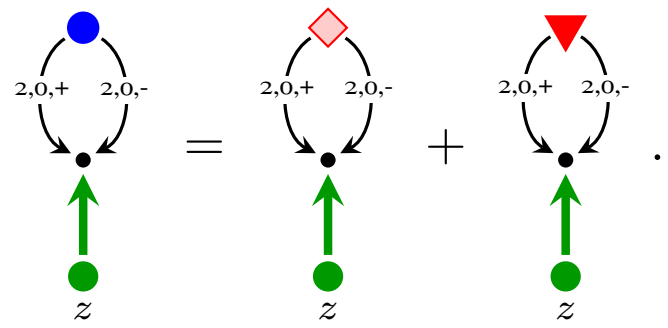

(6.48)

The renormalisation function (6.14) can be represented by the diagram

$$(\iota_\varepsilon C_{\varepsilon,\mathfrak{m}})(\varphi_z^\lambda) =$$

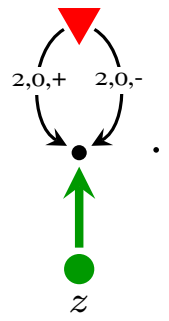

Recall that in (6.13) the model is renormalised twice: first time the model is renormalised by the constant $\widetilde{C}_\varepsilon($ 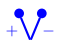 $)$ (this is a part of the definition of the map $\widetilde{\Pi}_z^{\varepsilon,\mathfrak{m}}$ 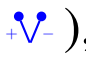), and second time it is renormalised by the random function $\widetilde{C}_{\varepsilon,\mathfrak{m}}$. The definitions (6.16) and (6.14) of these constant and function yield

$(\widehat{\Pi}_z^{\varepsilon,\mathfrak{m}}$ 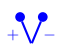 $)(\bar{z}) = (\Pi_z^{\varepsilon,\mathfrak{m}}$ 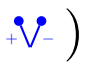 $)(\bar{z}) - \widetilde{C}_\varepsilon($ 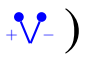 $) - \widetilde{C}_{\varepsilon,\mathfrak{m}}(\bar{z})$

where

$\widetilde{C}_\varepsilon($ 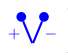 $) + \widetilde{C}_{\varepsilon,\mathfrak{m}}(t,x) = \int_{D_\varepsilon} \nabla_\varepsilon^- K^\varepsilon(t-s,x-y)\nabla_\varepsilon^+ K^\varepsilon(t-s,x-y)\widetilde{\mathbf{C}}_{\varepsilon,\mathfrak{m}}(s,y)\,ds,$

The latter is exactly the last term in (6.48), tested agains the function $\varphi_z^\lambda$. Hence, our definition (6.13) and [GMW24, Eq. 2.16] yield

$$\iota_\varepsilon(\widehat{\Pi}_z^{\varepsilon,\mathfrak{m}}\tau)(\varphi_z^\lambda) =$$

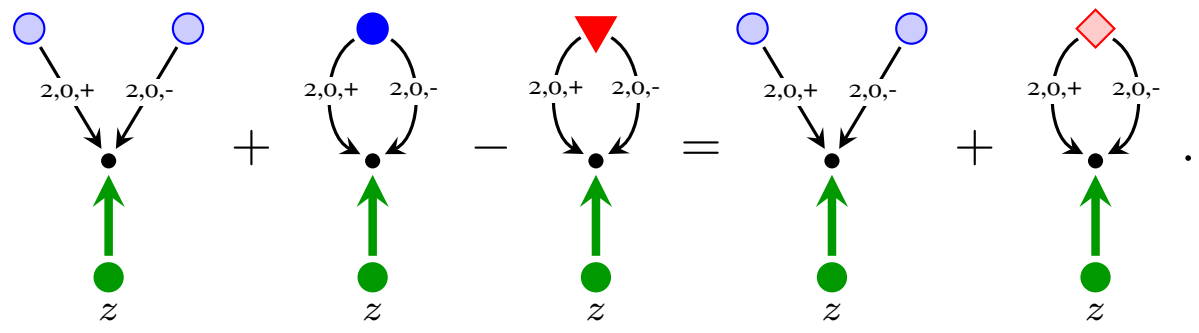

(6.49)

We denote the last two diagrams by $\iota_\varepsilon(\widehat{\Pi}_z^{\varepsilon,\mathfrak{m},1}\tau)(\varphi_z^\lambda)$ and $\iota_\varepsilon(\widehat{\Pi}_z^{\varepsilon,\mathfrak{m},2}\tau)(\varphi_z^\lambda)$ respectively.

For the first diagram in (6.49) we have $\nu_\gamma = -1$. With $m = 2$, the only choice in the first term on the right of (C.1) is $\mathbf{p}(1) = \mathbf{p}(2) = 2$ with $\alpha_\gamma(\mathbf{p}) = \delta_\gamma(\mathbf{p}) = 0$; and there are three choices in the second term on the right of (C.1), namely $\mathbf{p}(1) = \mathbf{p}(2) = \infty$ with $\alpha_\gamma(\mathbf{p}) = \delta_\gamma(\mathbf{p}) = 3$; $\mathbf{p}(1) = 2$, $\mathbf{p}(2) = \infty$ with $\alpha_\gamma(\mathbf{p}) = \delta_\gamma(\mathbf{p}) = \frac{3}{2}$; $\mathbf{p}(1) = \infty$, $\mathbf{p}(2) = 2$ with $\alpha_\gamma(\mathbf{p}) = \delta_\gamma(\mathbf{p}) = \frac{3}{2}$. Hence, for any $p \geq 2$ and $\theta > 0$ sufficiently small, Theorem C.1 yields

$$\Big(\mathbb{E}|\iota_\varepsilon(\widehat{\Pi}_z^{\varepsilon,\mathfrak{m},1}\tau)(\varphi_z^\lambda)|^p\Big)^{\frac{1}{p}} \lesssim (\lambda \vee \varepsilon)^{-1} + (\lambda \vee \varepsilon)^{-1-\theta} \lesssim (\lambda \vee \varepsilon)^{-1-\theta}. \tag{6.50}$$

The second diagram in (6.49) does not satisfy [GMW24, Assm. 3(1)], and we need to modify it by multiplying the singular kernels by powers of $\varepsilon$, as follows:

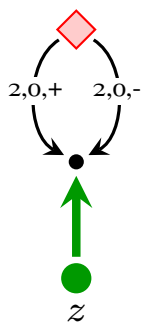

$$= \varepsilon^{2(a-2)}$$

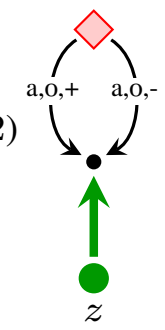

(6.51)

For any $a \in (\frac{3}{4}, \frac{3}{2})$, the new diagram satisfies [GMW24, Assm. 3] and we have $\nu_\gamma = \frac{3}{2} - 2a$. Since $m = 1$, the only choice in the first term on the right of (C.1) is $\mathbf{p}(1) = 2$ with $\alpha_\gamma(\mathbf{p}) = \frac{3}{2}$, $\delta_\gamma(\mathbf{p}) = 0$;

and the only choice in the second term on the right of (C.1) is $\mathbf{p}(1) = \infty$ with $\alpha_\gamma(\mathbf{p}) = 3$, $\delta_\gamma(\mathbf{p}) = \frac{3}{2}$. Hence, for any $p \geq 2$ and $\theta > 0$ sufficiently small, Theorem C.1 yields

$$\Big(\mathbb{E}|\iota_\varepsilon(\widehat{\Pi}_z^{\varepsilon,\mathtt{m},2}\tau)(\varphi_z^\lambda)|^p\Big)^{\frac{1}{p}} \lesssim \varepsilon^{2(a-2)}((\lambda\vee\varepsilon)^{\frac{3}{2}-2a}\varepsilon^{\frac{3}{2}} + (\lambda\vee\varepsilon)^{\frac{3}{2}-2a-\theta}\varepsilon^{\frac{3}{2}}) \lesssim \varepsilon^{2a-\frac{5}{2}}(\lambda\vee\varepsilon)^{\frac{3}{2}-2a-\theta}.$$

Choosing $a$ close enough to $\frac{3}{2}$ such that $2a - \frac{5}{2} > 0$, we get

$$\Big(\mathbb{E}|\iota_\varepsilon(\widehat{\Pi}_z^{\varepsilon,\mathtt{m},2}\tau)(\varphi_z^\lambda)|^p\Big)^{\frac{1}{p}} \lesssim \varepsilon^{\theta}(\lambda\vee\varepsilon)^{-1-2\theta} \tag{6.52}$$

for any $\theta > 0$ small enough. The two bounds (6.50) and (6.52) yield (6.43) if we take $\theta$ sufficiently small because $|\tau| = -1-2\kappa$. Moreover, (6.52) vanishes as $\varepsilon \to 0$ in the sense of (6.43).

### 6.4.3 Bound on the third element in (6.47)

For the element $\tau =$ 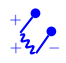, [GMW24, Eq. 2.16] yields

$$\iota_\varepsilon(\widehat{\Pi}_z^{\varepsilon,\mathtt{m}}\tau)(\varphi_z^\lambda) =$$

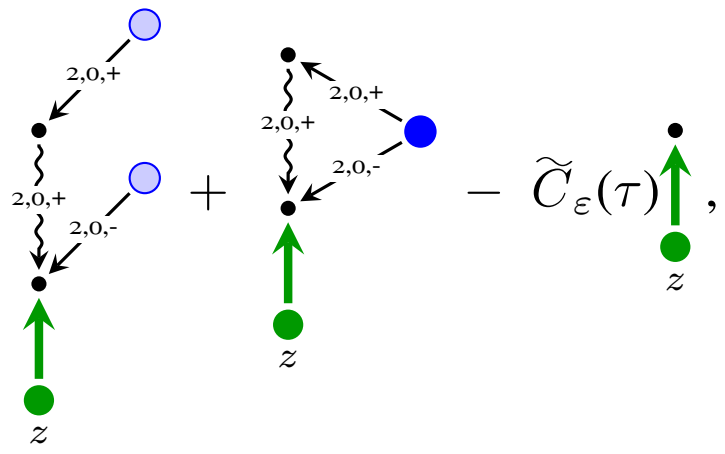

where $\widetilde{C}_\varepsilon(\tau)$ is the respective renormalisation constant defined in (6.15)-(6.16). As in (6.48), we can write furthermore

$$\iota_\varepsilon(\widehat{\Pi}_z^{\varepsilon,\mathtt{m}}\tau)(\varphi_z^\lambda) =$$

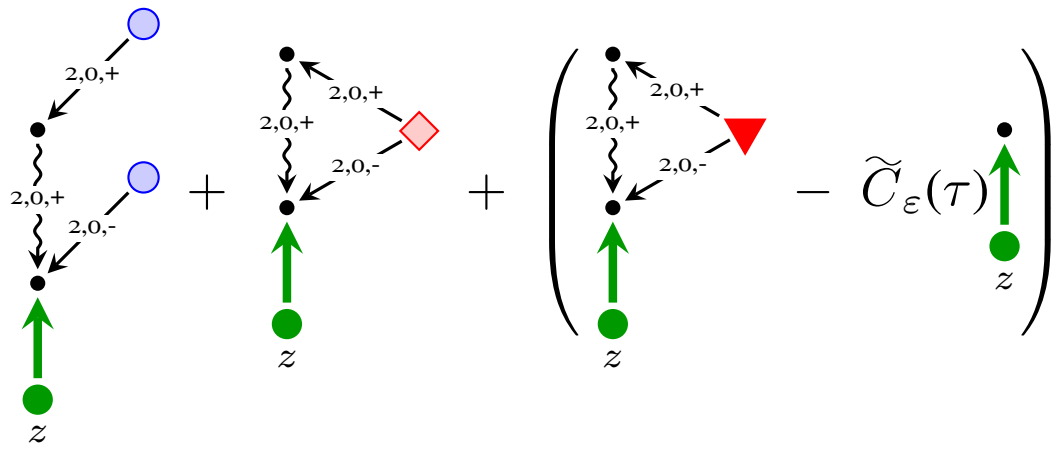

(6.53)

We denote the first two diagrams of (6.53) by $\iota_\varepsilon(\widehat{\Pi}_z^{\varepsilon,\mathtt{m},1}\tau)(\varphi_z^\lambda)$ and $\iota_\varepsilon(\widehat{\Pi}_z^{\varepsilon,\mathtt{m},2}\tau)(\varphi_z^\lambda)$ respectively, and the difference of the last two diagrams in parenthesis by $\iota_\varepsilon(\widehat{\Pi}_z^{\varepsilon,\mathtt{m},3}\tau)(\varphi_z^\lambda)$.

For the first diagram in (6.53) we have $\nu_\gamma = 0$. With $m = 2$, the only choice in the first term on the right of (C.1) is $\mathbf{p}(1) = \mathbf{p}(2) = 2$ with $\alpha_\gamma(\mathbf{p}) = \delta_\gamma(\mathbf{p}) = 0$; and there are three choices in the second term on the right of (C.1), namely $\mathbf{p}(1) = \mathbf{p}(2) = \infty$ with $\alpha_\gamma(\mathbf{p}) = \delta_\gamma(\mathbf{p}) = 3$; $\mathbf{p}(1) = 2$, $\mathbf{p}(2) = \infty$ with $\alpha_\gamma(\mathbf{p}) = \delta_\gamma(\mathbf{p}) = \frac{3}{2}$; $\mathbf{p}(1) = \infty$, $\mathbf{p}(2) = 2$ with $\alpha_\gamma(\mathbf{p}) = \delta_\gamma(\mathbf{p}) = \frac{3}{2}$. Hence, for any $p \geq 2$ and $\theta > 0$ sufficiently small,

$$\Big(\mathbb{E}|\iota_\varepsilon(\widehat{\Pi}_z^{\varepsilon,\mathtt{m},1}\tau)(\varphi_z^\lambda)|^p\Big)^{\frac{1}{p}} \lesssim (\lambda\vee\varepsilon)^{-\theta},$$

which is a bound of the desired order (6.43) because $|\tau| = -2\kappa$.

For the second diagram in (6.53) we have $\nu_\gamma = -\frac{3}{2}$. Since $m = 1$, the only choice in the first term on the right of (C.1) is $\mathbf{p}(1) = 2$ with $\alpha_\gamma(\mathbf{p}) = \frac{3}{2}$, $\delta_\gamma(\mathbf{p}) = 0$; and the only choice in the second term on the right of (C.1) is $\mathbf{p}(1) = \infty$ with $\alpha_\gamma(\mathbf{p}) = 3$, $\delta_\gamma(\mathbf{p}) = \frac{3}{2}$. Hence, for any $p \geq 2$ and $\theta > 0$ sufficiently small

$$\Big(\mathbb{E}|\iota_\varepsilon(\widehat{\Pi}_z^{\varepsilon,\mathtt{m},2}\tau)(\varphi_z^\lambda)|^p\Big)^{\frac{1}{p}} \lesssim (\lambda\vee\varepsilon)^{-\frac{3}{2}}\varepsilon^{\frac{3}{2}} + (\lambda\vee\varepsilon)^{-\frac{3}{2}-\theta}\varepsilon^{\frac{3}{2}} \lesssim \varepsilon^{\theta}(\lambda\vee\varepsilon)^{-2\theta}.$$

This is a bound of the desired order (6.43), and this expression moreover vanishes as $\varepsilon \to 0$.

To analyse the expression in parenthesis in (6.53) we introduce a new type of vertex "$\triangledown$", which corresponds to integration against

$$(t,x) \quad \mapsto \quad \varepsilon\langle \widetilde{M}^{\varepsilon,\mathfrak{m}}(x)\rangle_t - \nu^\varepsilon t \tag{6.54}$$

with the constant $\nu^\varepsilon$ defined in (2.8). Then we have

$$\iota_\varepsilon(\widehat{\Pi}_z^{\varepsilon,\mathfrak{m},3}\tau)(\varphi_z^\lambda) \;=$$ 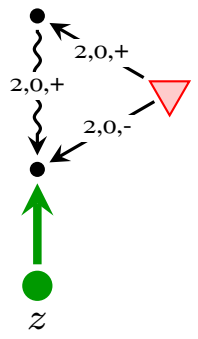 .

Recalling the definitions of the renormalisation constants (6.15)-(6.16), we can write

$$\begin{aligned}\iota_\varepsilon(\widehat{\Pi}_z^{\varepsilon,\mathfrak{m},3}\tau)(\varphi_z^\lambda) &= \int_{(D_\varepsilon)^3} dz_1 dz_2 d\bar z\; \varphi_z^\lambda(\bar z)\nabla_\varepsilon^+ K^\varepsilon(\bar z - z_1)\nabla_\varepsilon^+ K^\varepsilon(z_1 - z_2)\nabla_\varepsilon^- K^\varepsilon(\bar z - z_2) \\ &\qquad\qquad\qquad\qquad\qquad\qquad \times \big(\widetilde{\mathbf{C}}_{\varepsilon,\mathfrak{m}}(z_2) - \nu^\varepsilon\big) \\ &= \int_{(D_\varepsilon)^2} d\bar z dz_2\; \mathfrak{K}(\bar z)\varphi^\lambda_{z-\bar z}(z_2)\big(\widetilde{\mathbf{C}}_{\varepsilon,\mathfrak{m}}(z_2) - \nu^\varepsilon\big),\end{aligned}$$

with the new kernel

$$\mathfrak{K}(\bar z - z_2) := \nabla_\varepsilon^- K^\varepsilon(\bar z - z_2)\int_{D_\varepsilon} \nabla_\varepsilon^+ K^\varepsilon(\bar z - z_1)\nabla_\varepsilon^+ K^\varepsilon(z_1 - z_2) dz_1,$$

which satisfies $\int_{D_\varepsilon} |\mathfrak{K}(\bar z)| d\bar z \lesssim |\log \varepsilon|$. Lemma B.3 yields for any $\theta > 1 - 2\alpha + \underline{\kappa}$

$$\Big(\mathbb{E}|\iota_\varepsilon(\widehat{\Pi}_z^{\varepsilon,\mathfrak{m},3}\tau)(\varphi_z^\lambda)|^p\Big)^{\frac{1}{p}} \lesssim \varepsilon^{\theta/2}\mathfrak{m}^2(\lambda\vee\varepsilon)^{-\theta},$$

which is a bound of the required order (6.43). This quantity vanishes as $\varepsilon \to 0$.

We are going to bound integrals agains (6.54) in a similar manner.

### 6.4.4 Bound on the fourth element in (6.47)

Let $\tau =$ 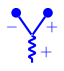 . [GMW24, Eq. 2.16] yields

$$\iota_\varepsilon(\widehat{\Pi}_z^{\varepsilon,\mathfrak{m}}\tau)(\varphi_z^\lambda) \;=$$ 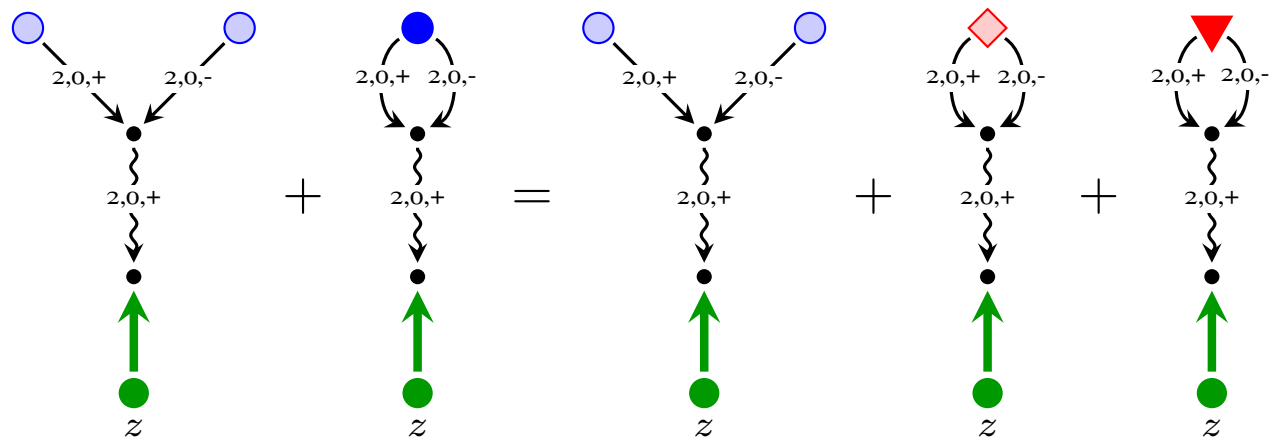 (6.55)

We denote the three diagrams of (6.55) by $\iota_\varepsilon(\widehat{\Pi}_z^{\varepsilon,\mathfrak{m},i}\tau)(\varphi_z^\lambda)$, $i = 1,2,3$, respectively.

For the first diagram in (6.55) we have $\nu_\gamma = 0$. With $m = 2$, the only choice in the first term on the right of (C.1) is $\mathbf{p}(1) = \mathbf{p}(2) = 2$ with $\alpha_\gamma(\mathbf{p}) = \delta_\gamma(\mathbf{p}) = 0$; and there are three choices in the second term on the right of (C.1), namely $\mathbf{p}(1) = \mathbf{p}(2) = \infty$ with $\alpha_\gamma(\mathbf{p}) = \delta_\gamma(\mathbf{p}) = 3$; $\mathbf{p}(1) = 2$, $\mathbf{p}(2) = \infty$ with $\alpha_\gamma(\mathbf{p}) = \delta_\gamma(\mathbf{p}) = \frac{3}{2}$; $\mathbf{p}(1) = \infty$, $\mathbf{p}(2) = 2$ with $\alpha_\gamma(\mathbf{p}) = \delta_\gamma(\mathbf{p}) = \frac{3}{2}$. Hence, for any $p \geq 2$ and $\theta > 0$ sufficiently small,

$$\Big(\mathbb{E}|\iota_\varepsilon(\widehat{\Pi}_z^{\varepsilon,\mathfrak{m},1}\tau)(\varphi_z^\lambda)|^p\Big)^{\frac{1}{p}} \lesssim (\lambda\vee\varepsilon)^{-\theta}.$$

This is a bound of the required order (6.43) with $|\tau| = -2\kappa$.

The second diagram in (6.55) does not satisfy [GMW24, Assm. 3(1)&(4)], and hence we need to modify it by multiplying the singular kernels by powers of $\varepsilon$, in the same way as we did in (6.51). We write

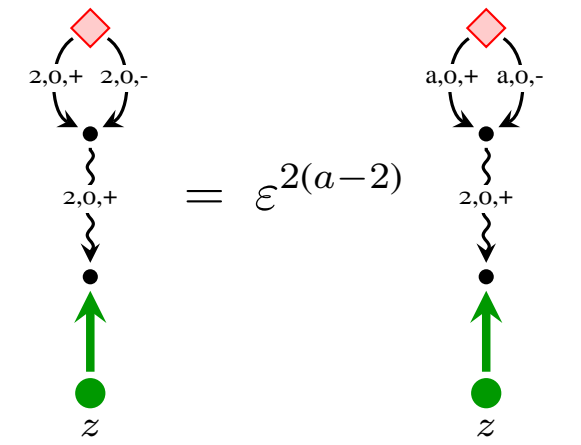

For any $a \in (\frac{5}{4}, \frac{3}{2})$, the new diagram satisfies [GMW24, Assm. 3], and the constant of [GMW24, Eq. 4.15] $\nu_\gamma = -\frac{1}{2} - 2a$. With $m = 1$, there is only one choice in the first term on the right of (C.1) $\mathbf{p}(1) = 2$ with $\alpha_\gamma(\mathbf{p})\frac{3}{2}$, $\delta_\gamma(\mathbf{p}) = 0$; and there is only one choice in the second term on the right of (C.1), namely $\mathbf{p}(1) = \infty$ with $\alpha_\gamma(\mathbf{p}) = 3$, $\delta_\gamma(\mathbf{p}) = \frac{3}{2}$. Hence, for any $p \geq 2$ and $\theta > 0$ sufficiently small,

$$\Big(\mathbb{E}|\iota_\varepsilon(\widehat{\Pi}_z^{\varepsilon,\mathfrak{m},2}\tau)(\varphi_z^\lambda)|^p\Big)^{\frac{1}{p}} \lesssim \varepsilon^{2(a-2)}((\lambda \vee \varepsilon)^{-\frac{1}{2}-2a}\varepsilon^{\frac{3}{2}} + (\lambda \vee \varepsilon)^{-\frac{1}{2}-2a-\theta}\varepsilon^{\frac{3}{2}}) \lesssim \varepsilon^\theta(\lambda \vee \varepsilon)^{-2\theta}.$$

This is a bound of the required order (6.43), and the expression vanishes as $\varepsilon \to 0$.

For the last diagram in (6.55), the analysis is similar to that of the last two diagrams in parenthesis in (6.53). Recalling from Propositions 4.2 and 4.12 that the kernels "kill" polynomials, the constant term of the bracket process (2.13) vanishes after convolving with the kernels. Hence, we get a renormalised triangular vertex as in (6.53). And we gat a bound on the last diagram in (6.55) of the required order.

### 6.4.5 Bound on the fifth element in (6.47)

Let $\tau =$ 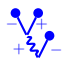 . Recalling the definition of the renormalised model, we get from [GMW24, Eq. 2.16]

$$\iota_\varepsilon(\widehat{\Pi}_z^{\varepsilon,\mathfrak{m}}\tau)(\varphi_z^\lambda) = \tag{6.56}$$

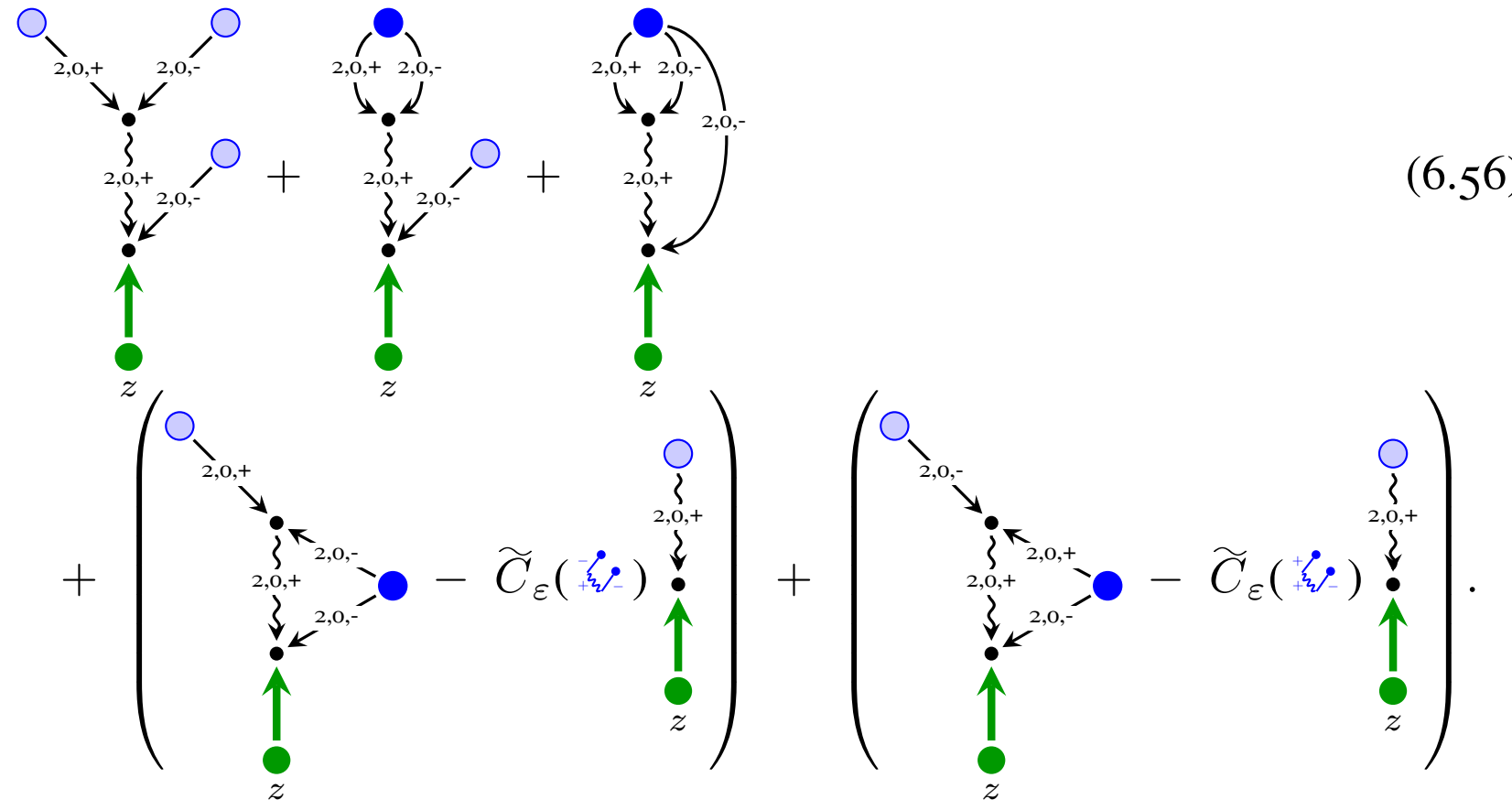

We denote the first three diagrams and the two expressions in the parentheses by $\iota_\varepsilon(\widehat{\Pi}_z^{\varepsilon,\mathfrak{m},i}\tau)(\varphi_z^\lambda)$, $i = 1, \ldots, 5$, respectively.

For the first diagram in (6.56) we have $\nu_\gamma = -\frac{1}{2}$. With $m = 3$, the only choice in the first term on the right of (C.1) is $\mathbf{p}(1) = \mathbf{p}(2) = \mathbf{p}(3) = 2$ with $\alpha_\gamma(\mathbf{p}) = \delta_\gamma(\mathbf{p}) = 0$; and there are seven choices in the second term on the right of (C.1), namely $\mathbf{p}(1) = \mathbf{p}(2) = \mathbf{p}(3) = \infty$ with $\alpha_\gamma(\mathbf{p}) = \delta_\gamma(\mathbf{p}) = \frac{9}{2}$; $\mathbf{p}(1) = \mathbf{p}(2) = \infty$, $\mathbf{p}(3) = 2$ with $\alpha_\gamma(\mathbf{p}) = \delta_\gamma(\mathbf{p}) = 3$ (and similarly for $\mathbf{p}(1) = \mathbf{p}(3) = \infty$, $\mathbf{p}(2) = 2$ and $\mathbf{p}(2) = \mathbf{p}(3) = \infty$, $\mathbf{p}(1) = 2$); $\mathbf{p}(1) = \infty$, $\mathbf{p}(2) = \mathbf{p}(3) = 2$ with $\alpha_\gamma(\mathbf{p}) = \delta_\gamma(\mathbf{p}) = \frac{3}{2}$ (and similarly for $\mathbf{p}(1) = \mathbf{p}(3) = 2$, $\mathbf{p}(2) = \infty$ and $\mathbf{p}(2) = \mathbf{p}(3) = 2$, $\mathbf{p}(1) = \infty$). Hence, for any $p \geq 2$ and $\theta > 0$ sufficiently small,

$$\Big(\mathbb{E}|\iota_\varepsilon(\widehat{\Pi}_z^{\varepsilon,\mathfrak{m},1}\tau)(\varphi_z^\lambda)|^p\Big)^{\frac{1}{p}} \lesssim (\lambda \vee \varepsilon)^{-\frac{1}{2}} + (\lambda \vee \varepsilon)^{-\frac{1}{2}-\theta} \lesssim (\lambda \vee \varepsilon)^{-\frac{1}{2}-\theta}.$$

The second diagram in (6.56) can be written as

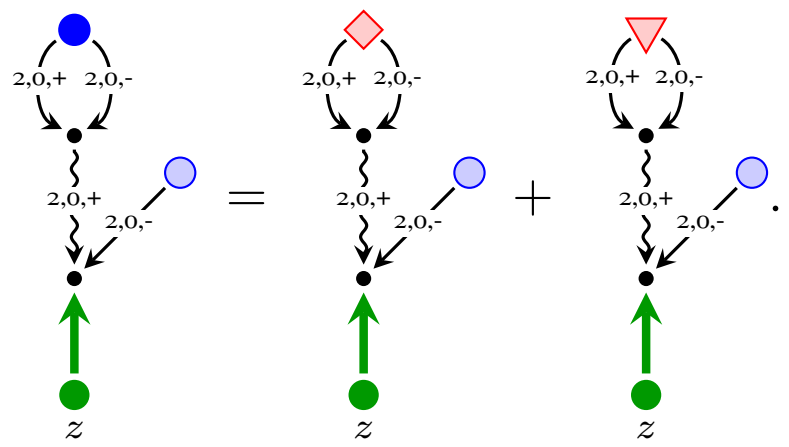

(6.57)

We denote these two diagrams by $\iota_\varepsilon(\widehat{\Pi}_z^{\varepsilon,\mathfrak{m},2,1}\tau)(\varphi_z^\lambda)$ and $\iota_\varepsilon(\widehat{\Pi}_z^{\varepsilon,\mathfrak{m},2,2}\tau)(\varphi_z^\lambda)$ respectively.

The first diagram in (6.57) does not satisfy [GMW24, Assm. 3(1)&(4)], and hence we need to modify it by multiplying the singular kernels by powers of $\varepsilon$, as we did in (6.51). We write

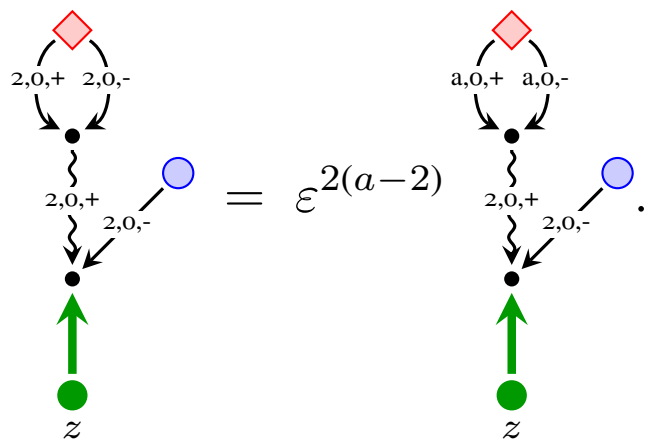

For $a \in (\frac{5}{4}, \frac{3}{2})$, the new diagram satisfies [GMW24, Assm. 3], and the constant of [GMW24, Eq. 4.15] $\nu_\gamma = 2 - 2a$. With $m = 2$, the only choice in the first term on the right of (C.1) is $\mathbf{p}(1) = \mathbf{p}(2) = 2$ with $\alpha_\gamma(\mathbf{p}) = \frac{3}{2}$, $\delta_\gamma(\mathbf{p}) = 0$; and there are three choices in the second term on the right of (C.1), namely $\mathbf{p}(1) = \mathbf{p}(2) = \infty$ with $\alpha_\gamma(\mathbf{p}) = \frac{9}{2}$, $\delta_\gamma(\mathbf{p}) = 3$; $\mathbf{p}(1) = \infty$, $\mathbf{p}(2) = 2$ with $\alpha_\gamma(\mathbf{p}) = 3$, $\delta_\gamma(\mathbf{p}) = \frac{3}{2}$; $\mathbf{p}(1) = 2$, $\mathbf{p}(2) = \infty$ with $\alpha_\gamma(\mathbf{p}) = 3$, $\delta_\gamma(\mathbf{p}) = \frac{3}{2}$. Hence, for any $p \geq 2$ and $\theta > 0$ sufficiently small,

$$\Big(\mathbb{E}|\iota_\varepsilon(\widehat{\Pi}_z^{\varepsilon,\mathfrak{m},2,1}\tau)(\varphi_z^\lambda)|^p\Big)^{\frac{1}{p}} \lesssim \varepsilon^{2(a-2)}((\lambda \vee \varepsilon)^{2-2a}\varepsilon^{\frac{3}{2}} + (\lambda \vee \varepsilon)^{2-2a-\theta}\varepsilon^{\frac{3}{2}}) \lesssim \varepsilon^\theta(\lambda \vee \varepsilon)^{-\frac{1}{2}-2\theta},$$

since $2a - \frac{5}{2} > 0$. This quantity vanishes as $\varepsilon \to 0$. We analyse the second diagram in (6.57) similarly to how we analysed the last diagram in (6.55).

Now, we will analyse the third diagram in (6.56). We recall that the contraction of three vertices corresponds to integration against the process $(t, x) \mapsto \varepsilon^2(\Delta_t \widetilde{M}^{\varepsilon,\mathfrak{m}}(x))^3$. Recalling that the jump size of the martingale is $\varepsilon^{1/2}$, this process equals $(t, x) \mapsto \varepsilon^3 \Delta_t \widetilde{M}^{\varepsilon,\mathfrak{m}}(x)$. Finally, using equation (6.34), we write the latter as $(t, x) \mapsto \varepsilon^3 d\widetilde{M}^{\varepsilon,\mathfrak{m}}(t, x) + \varepsilon^{\frac{3}{2}} \widetilde{\mathfrak{C}}_{\varepsilon,\mathfrak{m}}(t, x)$. Hence, we can write

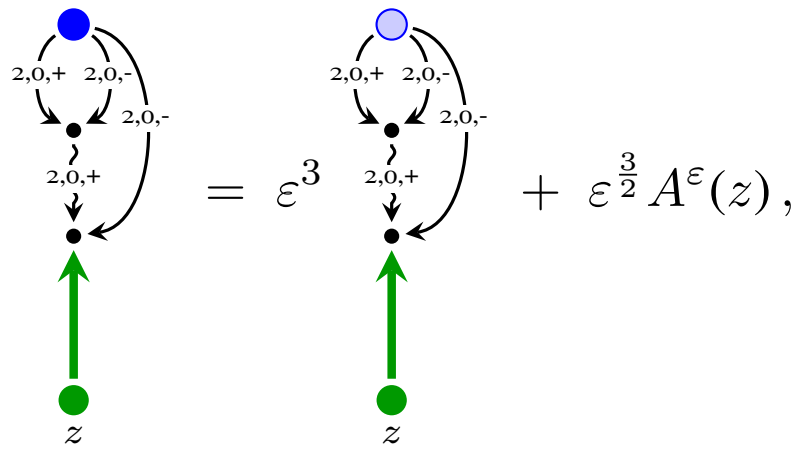

(6.58)

where $\widetilde{\mathfrak{C}}_{\varepsilon,\mathfrak{m}}$ is defined in Section 2.2, and

$$A^\varepsilon(z) = \int_{(D_\varepsilon)^3} \varphi_z^\lambda(\bar{z})\nabla_\varepsilon^+ K^\varepsilon(\bar{z} - z_1)\nabla_\varepsilon^+ K^\varepsilon(z_1 - \tilde{z})\nabla_\varepsilon^- K^\varepsilon(z_1 - \tilde{z})\nabla_\varepsilon^- K^\varepsilon(\bar{z} - \tilde{z})\widetilde{\mathfrak{C}}_{\varepsilon,\mathfrak{m}}(\tilde{z})\, dz_1 d\tilde{z} d\bar{z}.$$

The first diagram in (6.58), which we denote by $\iota_\varepsilon(\widehat{\Pi}_z^{\varepsilon,\mathfrak{m},5,1}\tau)(\varphi_z^\lambda)$, does not satisfy [GMW24, Assm. 3(1)], and hence we need to modify it by multiplying the singular kernels by powers of $\varepsilon$, as follows:

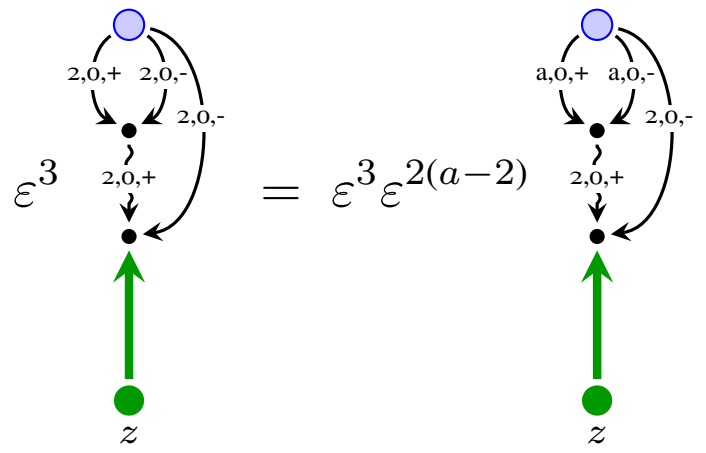

For $a \in (\frac{1}{2}, 1)$, the new diagram satisfies [GMW24, Assm. 3], and we have $\nu_\gamma = \frac{1}{2} - 2a$. Since $m = 1$, the only choice in the first term on the right of (C.1) is $\mathbf{p}(1) = 2$ with $\alpha_\gamma(\mathbf{p}) = 3$, $\delta_\gamma(\mathbf{p}) = 0$; and the only choice in the second term on the right of (C.1) is $\mathbf{p}(1) = \infty$ with $\alpha_\gamma(\mathbf{p}) = \frac{9}{2}$, $\delta_\gamma(\mathbf{p}) = \frac{3}{2}$. Hence, for any $p \geq 2$ and $\theta > 0$ sufficiently small

$$\left(\mathbb{E}|\iota_\varepsilon(\widehat{\Pi}_z^{\varepsilon,\mathfrak{m},5,1}\tau)(\varphi_z^\lambda)|^p\right)^{\frac{1}{p}} \lesssim \varepsilon^{2a-1}((\lambda \vee \varepsilon)^{\frac{1}{2}-2a}\varepsilon^3 + (\lambda \vee \varepsilon)^{\frac{1}{2}-2a-\theta}\varepsilon^3) \lesssim \varepsilon^{2a+2}(\lambda \vee \varepsilon)^{\frac{1}{2}-2a-\theta}.$$

Since $2a - 1 > 0$, we get a bound of order $\varepsilon^\theta(\lambda \vee \varepsilon)^{-\frac{1}{2}-2\theta}$ as required. This diagram vanishes in the limit $\varepsilon \to 0$.

Defining the new kernel

$$\mathfrak{K}(\bar{z} - \tilde{z}) := \nabla_\varepsilon^- K^\varepsilon(\bar{z} - \tilde{z}) \int_{D_\varepsilon} \nabla_\varepsilon^+ K^\varepsilon(\bar{z} - z_1)\nabla_\varepsilon^+ K^\varepsilon(z_1 - \tilde{z})\nabla_\varepsilon^- K^\varepsilon(z_1 - \tilde{z})dz_1\,,$$

which satisfies $|\mathfrak{K}(\bar{z})| \lesssim \|\bar{z}\|_\mathfrak{s}^{-5}$, we can write $\varepsilon^{\frac{3}{2}}A^\varepsilon(z)$ as

$$\varepsilon^{\frac{3}{2}}A^\varepsilon(z) = \varepsilon^{\frac{3}{2}} \int_{(D_\varepsilon)^2} d\tilde{z}d\bar{z}\; \mathfrak{K}(\bar{z})\varphi_{z-\bar{z}}^\lambda(\tilde{z})\widetilde{\mathsf{C}}_{\varepsilon,\mathfrak{m}}(\tilde{z}).$$

Lemma B.4 yields

$$\left(\mathbb{E}|\varepsilon^{\frac{3}{2}}A^\varepsilon(z)|^p\right)^{\frac{1}{p}} \lesssim \mathfrak{m}\varepsilon^\theta(\lambda \vee \varepsilon)^{-\frac{1}{2}-\theta},$$

for any $\theta > 0$ small enough. This is a bound of the required order and this expression vanishes as $\varepsilon \to 0$.

For the expressions in the last two parentheses in (6.56) we get bounds of the same order in the same way as we did above.

### 6.4.6 Bound on the sixth element in (6.47)

Now we consider the element $\tau =$ 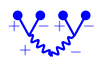 . We have

$$\iota_\varepsilon(\widehat{\Pi}_z^{\varepsilon,\mathfrak{m}}\tau)(\varphi_z^\lambda) =$$

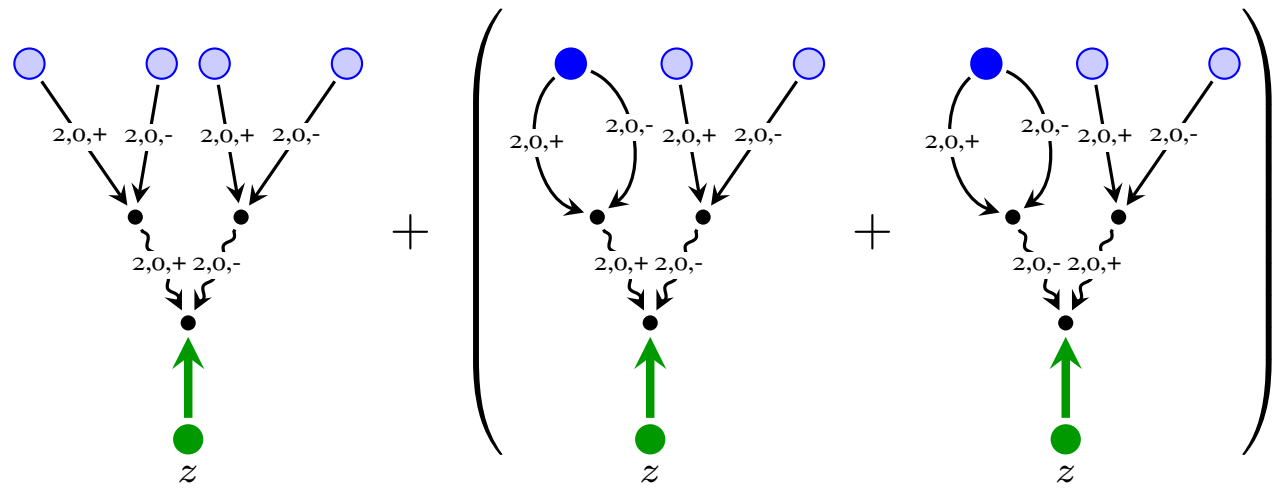

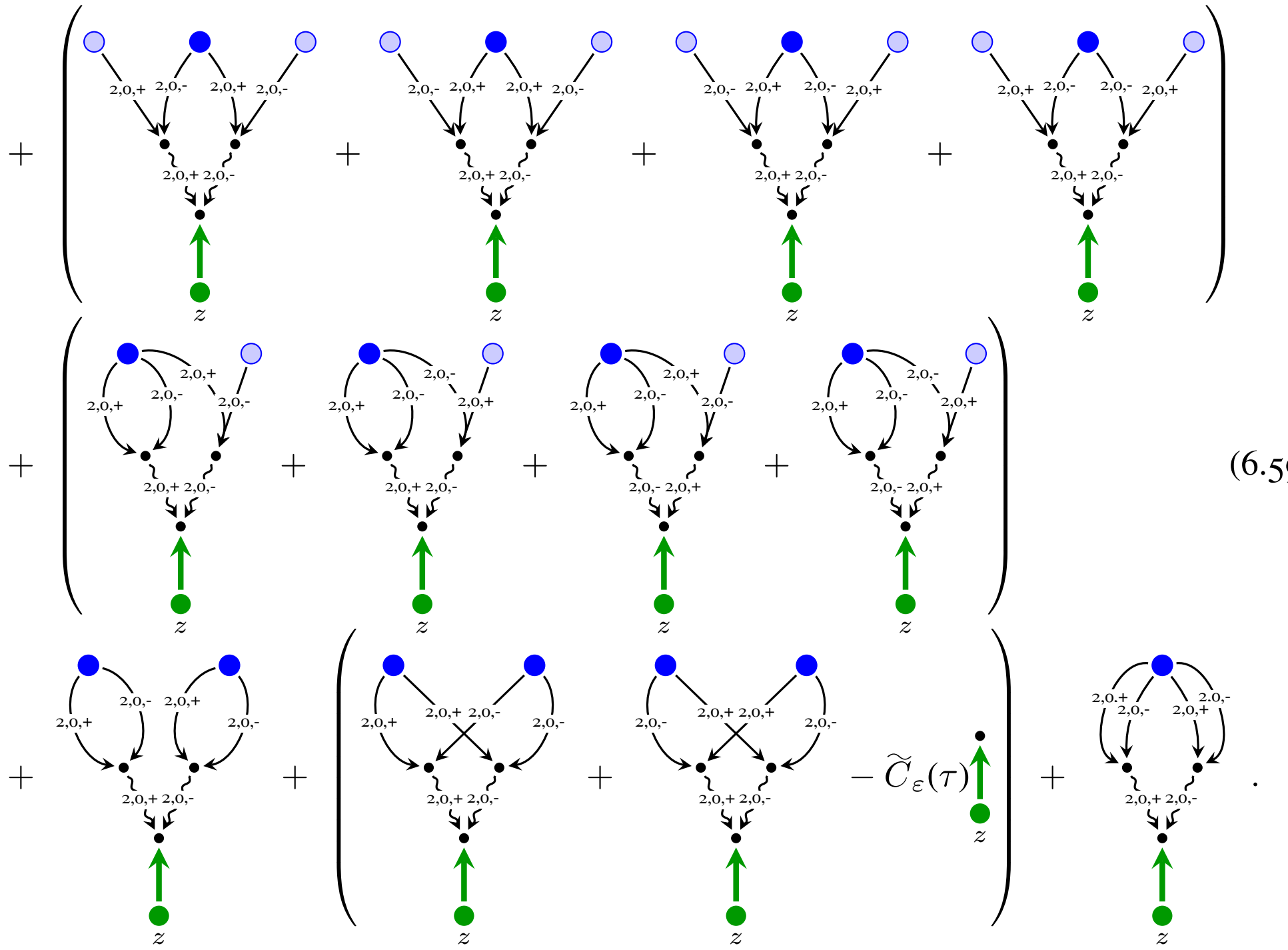

(6.59)

We denote the sixteen diagrams of (6.59) by $\iota_\varepsilon(\widehat{\Pi}_z^{\varepsilon,\mathrm{m},i}\tau)(\varphi_z^\lambda)$ for $i = 1, \ldots, 16$ respectively.

For the first diagram in (6.59) we have $\nu_\gamma = 0$. With $m = 4$, the only choice in the first term on the right of (C.1) is $\mathbf{p}(1) = \mathbf{p}(2) = \mathbf{p}(3) = \mathbf{p}(4) = 2$ with $\alpha_\gamma(\mathbf{p}) = \delta_\gamma(\mathbf{p}) = 0$; and there are fifteen choices in the second term on the right of (C.1), namely $\mathbf{p}(1) = \mathbf{p}(2) = \mathbf{p}(3) = \mathbf{p}(4) = \infty$ with $\alpha_\gamma(\mathbf{p}) = \delta_\gamma(\mathbf{p}) = 6$; $\mathbf{p}(1) = \mathbf{p}(2) = \mathbf{p}(3) = \infty$, $\mathbf{p}(4) = 2$ with $\alpha_\gamma(\mathbf{p}) = \delta_\gamma(\mathbf{p}) = \frac{9}{2}$ (and similarly for three other permutations); $\mathbf{p}(1) = \mathbf{p}(2) = \infty$, $\mathbf{p}(3) = \mathbf{p}(4) = 2$ with $\alpha_\gamma(\mathbf{p}) = \delta_\gamma(\mathbf{p}) = 3$ (and similarly for five other permutations); $\mathbf{p}(1) = \mathbf{p}(2) = \mathbf{p}(3) = 2$, $\mathbf{p}(4) = \infty$ with $\alpha_\gamma(\mathbf{p}) = \delta_\gamma(\mathbf{p}) = \frac{3}{2}$ (and similarly for three other permutations). Hence, for any $p \geq 2$ and $\theta > 0$ sufficiently small

$$\Big(\mathbb{E}|\iota_\varepsilon(\widehat{\Pi}_z^{\varepsilon,\mathrm{m},1}\tau)(\varphi_z^\lambda)|^p\Big)^{\frac{1}{p}} \lesssim (\lambda \vee \varepsilon)^{-\theta},$$

as required.

The moments of the other diagrams are bounded by a multiple of $\varepsilon^\theta(\lambda \vee \varepsilon)^{-2\theta}$, for $\theta > 0$ sufficiently small, using the same approach as above. Moreover, they vanish as $\varepsilon \to 0$.

### 6.4.7 Bound on the seventh element in (6.47)

For the element $\tau =$ 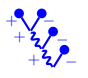 we have

$$\iota_\varepsilon(\widehat{\Pi}_z^{\varepsilon,\mathrm{m}}\tau)(\varphi_z^\lambda) =$$

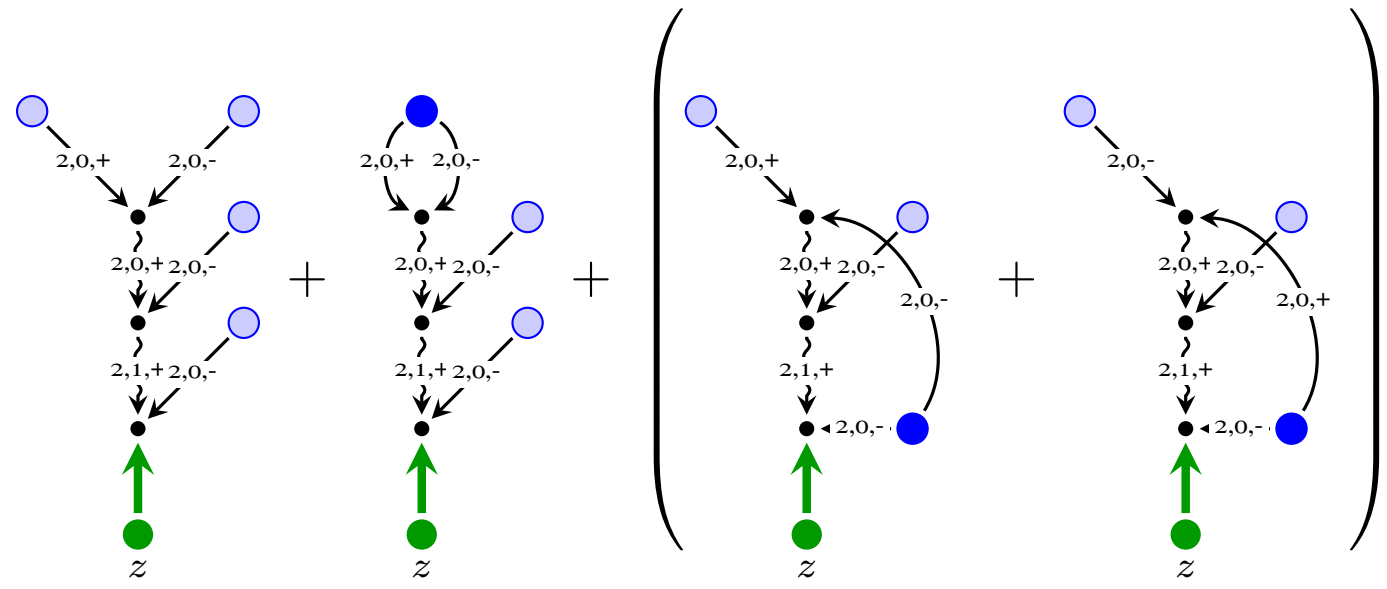

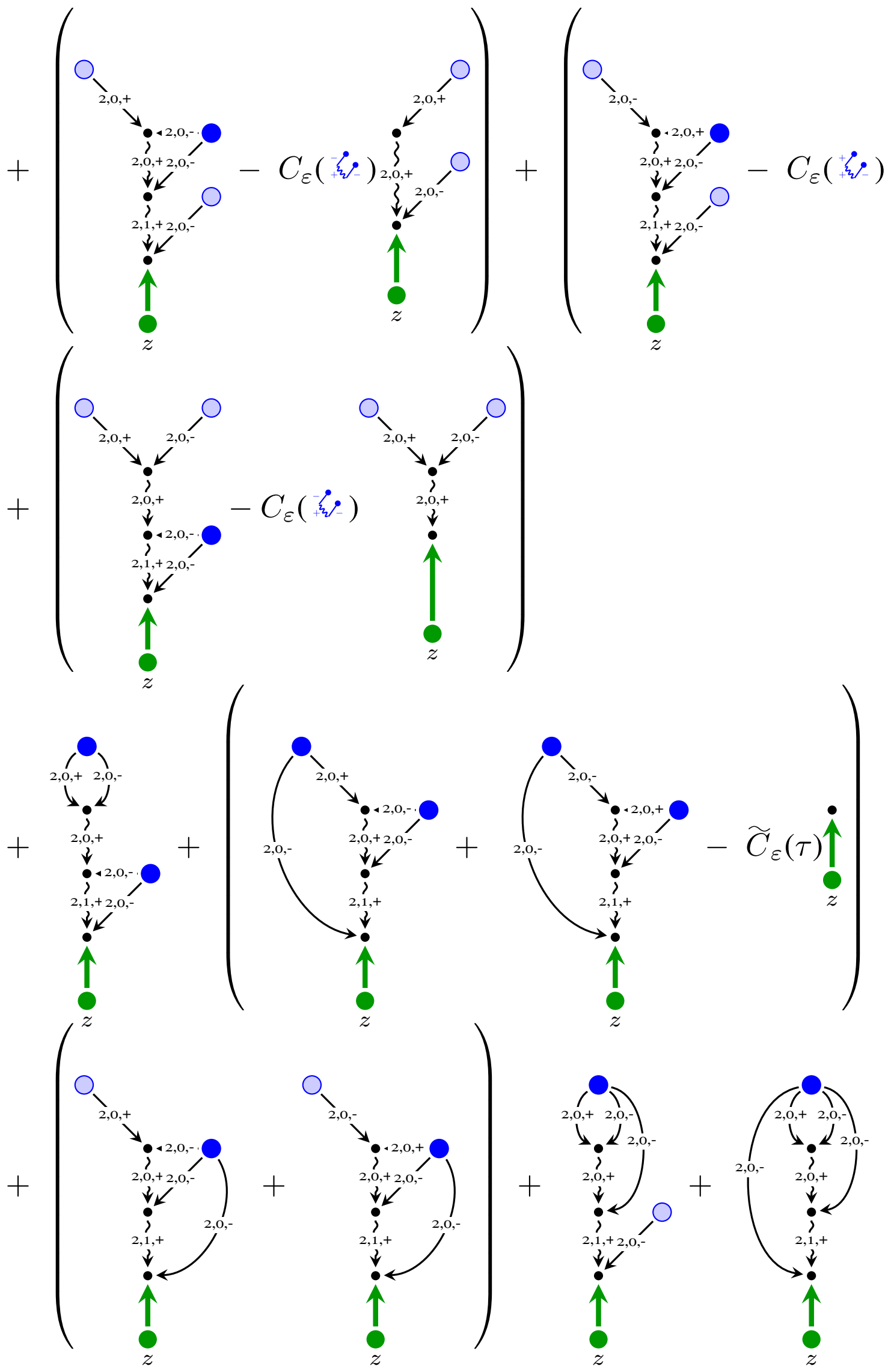

The moment bounds on these diagrams can be proved by analogy with our previous bounds. So, we obtain the required bounds (6.43) and only the first diagram does not vanish in the limit $\varepsilon \to 0$.

### 6.4.8 Bound on the last element in (6.47)

Now, we are going to prove moment bounds for the last element in (6.47), which we denote by $\tau = \partial_\mathcal{E}($ 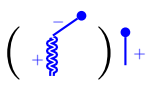 . The definitions of the discrete model and the renormalisation constant (6.17) yield

$$\iota_\varepsilon(\widehat{\Pi}^{\varepsilon,\mathfrak{m}}_z\tau)(\varphi^\lambda_z) =$$

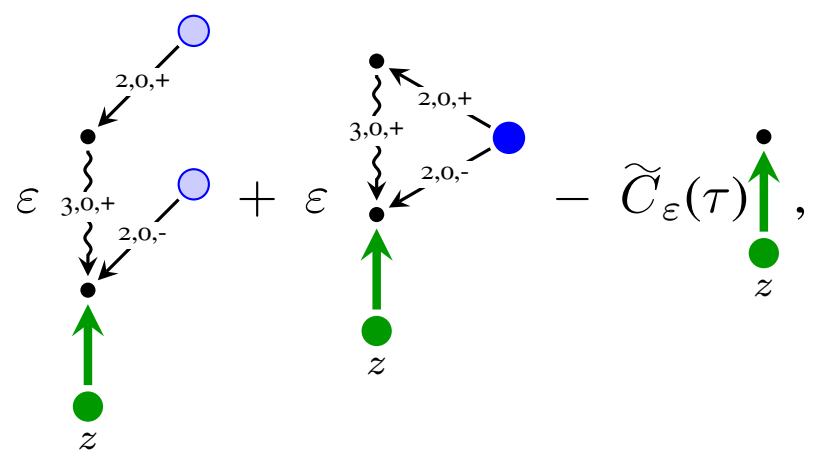

where we used the discrete Leibniz rule in (6.2). We can use the multiplier $\varepsilon$ to "improve" the singularity of the kernels labeled with 3. Then this element is bounded in the same way as the element in Section 6.4.3. Namely, we get

$$\Big(\mathbb{E}|\iota_\varepsilon(\widehat{\Pi}_z^{\varepsilon,\mathrm{m}}\tau)(\varphi_z^\lambda)|^p\Big)^{\frac{1}{p}} \lesssim \varepsilon^\theta(\lambda \vee \varepsilon)^{-2\theta}$$

for any $\theta > 0$ small enough, which vanishes as $\varepsilon \to 0$.

# 7 Discrete solution maps

We are going to define an equation on a suitable space $\mathcal{D}_\varepsilon^{\gamma,\eta}(\widehat{Z}^{\varepsilon,\mathrm{m}})$ with respect to the renormalised model $\widehat{Z}^{\varepsilon,\mathrm{m}}$, defined in Section 6.1, so that after reconstruction we obtain the discrete equation (2.33).

## 7.1 Operators on modelled distributions

In order to describe equation (2.33) in the framework of regularity structures, we need to define the derivatives and multiplication by $\varepsilon$ as operators on modelled distributions. More precisely, let us take a discrete modelled distribution $H^\varepsilon \in \mathcal{D}_\varepsilon^{\gamma,\eta}$, with respect to a discrete model $Z^\varepsilon$. Let $H^\varepsilon$ be in the domain of the map $\partial_+$, defined on a subset of the model space in Section 5. Then the definitions (5.2e) and (5.18) yield $\partial_+ H^\varepsilon \in \mathcal{D}_\varepsilon^{\gamma-1,\eta-1}$. However, Definition 5.3 of a model and the definition of the reconstruction map (5.19) do not imply that $(\mathcal{R}^\varepsilon\partial_+ H^\varepsilon)(z)$ is equal to $\nabla_\varepsilon^+(\mathcal{R}^\varepsilon H^\varepsilon)(z)$. Indeed, one can see it because $\Pi_z^\varepsilon\partial_+\tau$ does not have to be equal to $\nabla_\varepsilon^+\Pi_z^\varepsilon\tau$ for $\tau \notin \{\mathbf{1}, X_1\}$. We define in (7.1) and (7.2) the differential maps on modelled distributions which yield the respective discrete derivatives of $\mathcal{R}^\varepsilon H^\varepsilon$ after reconstruction (see Lemma 7.1). Using these maps in the abstract equation (7.14)-(7.15) below will guarantee that we get the required equation (2.33) after reconstruction.

We also have a problem with multiplication by $\varepsilon$. Namely, for a suitable modelled distribution $U^\varepsilon$, $\mathcal{R}^\varepsilon\mathcal{E}(U^\varepsilon)$ is not necessarily equal to $\varepsilon\mathcal{R}^\varepsilon U^\varepsilon$, because $U^\varepsilon$ usually has a polynomial part and $\mathcal{E}$ has not been defined on polynomials. We define in (7.6) a map on modelled distributions which yields multiplication by $\varepsilon$ after reconstruction (see Lemma 7.2).

We work throughout this section with a general discrete model $Z^\varepsilon = (\Pi^\varepsilon, \Gamma^\varepsilon)$.

For a modelled distribution $H^\varepsilon : D_\varepsilon \to \mathrm{span}\{\widetilde{\mathcal{U}}\}$ we set

$$(\mathscr{D}_\varepsilon^\pm H^\varepsilon)(z) := \partial_\pm H^\varepsilon(z) \pm \frac{1}{\varepsilon}(\mathcal{R}^\varepsilon H^\varepsilon - \Pi_z^\varepsilon H^\varepsilon(z))(z \pm \varepsilon)\mathbf{1}, \tag{7.1}$$

where the first term in a natural differentiation of the modelled distribution and the last term is a correction term which guarantees the desired identities (7.3) after reconstruction. We will also need an analogue of the symmetric discrete derivative (2.5), defined as

$$\mathscr{D}_\varepsilon H^\varepsilon := \frac{1}{2}(\mathscr{D}_\varepsilon^- H^\varepsilon + \mathscr{D}_\varepsilon^+ H^\varepsilon). \tag{7.2}$$

**Lemma 7.1** *The reconstructions of* (7.1) *and* (7.2) *satisfy*

$$\mathcal{R}^\varepsilon\mathscr{D}_\varepsilon^\pm H^\varepsilon = \nabla_\varepsilon^\pm(\mathcal{R}^\varepsilon H^\varepsilon), \qquad \mathcal{R}^\varepsilon\mathscr{D}_\varepsilon H^\varepsilon = \nabla_\varepsilon(\mathcal{R}^\varepsilon H^\varepsilon). \tag{7.3}$$

*Moreover, for any $\gamma > 1$ and $\eta \leq 1$ one has*

$$\|\mathscr{D}_\varepsilon^\pm H^\varepsilon\|_{\gamma-1,\eta-1;T}^{(\varepsilon)} \lesssim (1 + \|\Pi^\varepsilon\|_T^{(\varepsilon)})\|H^\varepsilon\|_{\gamma,\eta;T}^{(\varepsilon)}, \quad \|\mathscr{D}_\varepsilon H^\varepsilon\|_{\gamma-1,\eta-1;T}^{(\varepsilon)} \lesssim (1 + \|\Pi^\varepsilon\|_T^{(\varepsilon)})\|H^\varepsilon\|_{\gamma,\eta;T}^{(\varepsilon)},$$

*where the proportionality constants are independent of $\varepsilon$ and $T$.*

*Proof.* Throughout the proof we will use $z, \bar{z} \in [-T, T] \times [-1, 1]$, so that we use the norms from Remarks 5.6 and 5.7. Using the definition of the reconstruction map (5.19) we get

$$(\mathcal{R}^\varepsilon\mathscr{D}_\varepsilon^\pm H^\varepsilon)(z) = (\Pi_z^\varepsilon\partial_\pm H^\varepsilon(z))(z) \mp \frac{1}{\varepsilon}(\Pi_z^\varepsilon H^\varepsilon(z) - \Pi_{z\pm\varepsilon}^\varepsilon H^\varepsilon(z \pm \varepsilon))(z \pm \varepsilon).$$

Definition 5.3 yields $(\Pi^\varepsilon_z \partial_\pm H^\varepsilon(z))(z) = \nabla^\pm_\varepsilon (\Pi^\varepsilon_z H^\varepsilon(z))(\bar z)|_{\bar z = z}$, where the discrete derivatives are applied to the spatial component of $\bar z$. Writing explicitly, we get

$$(\Pi^\varepsilon_z \partial_\pm H^\varepsilon(z))(z) = \pm\frac{1}{\varepsilon}((\Pi^\varepsilon_z H^\varepsilon(z))(z \pm \varepsilon) - (\Pi^\varepsilon_z H^\varepsilon(z))(z)).$$

Plugging it into the preceding expression, we get

$$(\mathcal{R}^\varepsilon \mathscr{D}^\pm_\varepsilon H^\varepsilon)(z) = \pm\frac{1}{\varepsilon}((\Pi^\varepsilon_{z\pm\varepsilon} H^\varepsilon(z \pm \varepsilon))(z \pm \varepsilon) - (\Pi^\varepsilon_z H^\varepsilon(z))(z)) = \nabla^\pm_\varepsilon (\mathcal{R}^\varepsilon H^\varepsilon)(z).$$

The last identity in (7.3) follows from the first two.

A bound on $\mathscr{D}_\varepsilon H^\varepsilon$ follows from the bounds on $\mathscr{D}^\pm_\varepsilon H^\varepsilon$. Let us now prove the bound on $\mathscr{D}^-_\varepsilon H^\varepsilon$, and the bound on $\mathscr{D}^+_\varepsilon H^\varepsilon$ can be proved in the same way. We write $\mathscr{D}^-_\varepsilon H^\varepsilon = \partial_- H^\varepsilon + \bar H^\varepsilon$ where $\bar H^\varepsilon$ is the last term in the definition (7.1). Then

$$\| \mathscr{D}^-_\varepsilon H^\varepsilon \|^{(\varepsilon)}_{\gamma-1,\eta-1;T} \leq \| \partial_- H^\varepsilon \|^{(\varepsilon)}_{\gamma-1,\eta-1;T} + \| \bar H^\varepsilon \|^{(\varepsilon)}_{\gamma-1,\eta-1;T},$$

and we readily get $\| \partial_- H^\varepsilon \|^{(\varepsilon)}_{\gamma-1,\eta-1;T} \leq \| H^\varepsilon \|^{(\varepsilon)}_{\gamma,\eta;T}$ from the definitions (5.2e) and (5.18). We will show that $\| \bar H^\varepsilon \|^{(\varepsilon)}_{\gamma-1,\eta-1;T} \lesssim \|\Pi^\varepsilon\|^{(\varepsilon)}_T \| H^\varepsilon \|^{(\varepsilon)}_{\gamma,\eta;T}$. The definitions of the discrete model and the space $\mathcal{D}^{\gamma,\eta}_{\varepsilon,T}$ yield

$$\begin{aligned}
\|\bar H^\varepsilon(z)\|_0 &= \frac{1}{\varepsilon}|\Pi^\varepsilon_z (H^\varepsilon(z) - \Gamma^\varepsilon_{z,z-\varepsilon} H^\varepsilon(z-\varepsilon))(z-\varepsilon)| \\
&\leq \sum_{\beta \leq 0} \frac{1}{\varepsilon}|\Pi^\varepsilon_z \mathcal{Q}_\beta (H^\varepsilon(z) - \Gamma^\varepsilon_{z,z-\varepsilon} H^\varepsilon(z-\varepsilon))(z-\varepsilon)| \\
&\lesssim \sum_{\beta \leq 0} \varepsilon^{\beta-1} \|\Pi^\varepsilon\|^{(\varepsilon)}_T \|H^\varepsilon(z) - \Gamma^\varepsilon_{z,z-\varepsilon} H^\varepsilon(z-\varepsilon)\|_\beta \\
&\lesssim \varepsilon^{\gamma-1}(\|z\|_\mathfrak{s} + \varepsilon)^{\eta-\gamma} \|\Pi^\varepsilon\|^{(\varepsilon)}_T \| H^\varepsilon \|^{(\varepsilon)}_{\gamma,\eta;T}. \qquad (7.4)
\end{aligned}$$

Using the assumption $\gamma > 1$, we can bound it by a constant multiple of $(\|z\|_\mathfrak{s} + \varepsilon)^{\eta-1} \|\Pi^\varepsilon\|^{(\varepsilon)}_T \| H^\varepsilon \|^{(\varepsilon)}_{\gamma,\eta;T}$. Since $\eta \leq 1$, this expression equals $(\|z\|_\mathfrak{s} + \varepsilon)^{(\eta-1)\wedge 0} \|\Pi^\varepsilon\|^{(\varepsilon)}_T \| H^\varepsilon \|^{(\varepsilon)}_{\gamma,\eta;T}$, which yields the required bound on the first terms in (5.16) and (5.18).

For $\bar z \neq z$ we get

$$\|\bar H^\varepsilon(z) - \Gamma^\varepsilon_{z,\bar z} \bar H^\varepsilon(\bar z)\|_0 = \|\bar H^\varepsilon(z) - \bar H^\varepsilon(\bar z)\|_0 \leq \|\bar H^\varepsilon(z)\|_0 + \|\bar H^\varepsilon(\bar z)\|_0.$$

Using (7.4), we bound it by a constant times $\varepsilon^{\gamma-1}((\|z\|_\mathfrak{s} + \varepsilon)^{\eta-\gamma} + (\|\bar z\|_\mathfrak{s} + \varepsilon)^{\eta-\gamma}) \|\Pi^\varepsilon\|^{(\varepsilon)}_T \| H^\varepsilon \|^{(\varepsilon)}_{\gamma,\eta;T}$. Since $\eta - \gamma < 0$, we have

$$(\|z\|_\mathfrak{s} + \varepsilon)^{\eta-\gamma} + (\|\bar z\|_\mathfrak{s} + \varepsilon)^{\eta-\gamma} \leq 2((\|z\|_\mathfrak{s} + \varepsilon) \wedge (\|\bar z\|_\mathfrak{s} + \varepsilon))^{\eta-\gamma} = 2(\|z, \bar z\|_\mathfrak{s} + \varepsilon)^{\eta-\gamma}.$$

Furthermore, since $\gamma > 1$ we bound $\varepsilon^{\gamma-1} \leq (\|z - \bar z\|_\mathfrak{s} + \varepsilon)^{\gamma-1}$ and hence

$$\|\bar H^\varepsilon(z) - \Gamma^\varepsilon_{z,\bar z} \bar H^\varepsilon(\bar z)\|_0 \lesssim (\|z - \bar z\|_\mathfrak{s} + \varepsilon)^{\gamma-1} (\|z, \bar z\|_\mathfrak{s} + \varepsilon)^{\eta-\gamma} \|\Pi^\varepsilon\|^{(\varepsilon)}_T \| H^\varepsilon \|^{(\varepsilon)}_{\gamma,\eta;T}.$$

This gives the required bound on the second terms in (5.16) and (5.18). We have just proved $\| \bar H^\varepsilon \|^{(\varepsilon)}_{\gamma-1,\eta-1;T} \lesssim \|\Pi^\varepsilon\|^{(\varepsilon)}_T \| H^\varepsilon \|^{(\varepsilon)}_{\gamma,\eta;T}$, which yields $\| \mathscr{D}^-_\varepsilon H^\varepsilon \|^{(\varepsilon)}_{\gamma-1,\eta-1;T} \lesssim \|\Pi^\varepsilon\|^{(\varepsilon)}_T \| H^\varepsilon \|^{(\varepsilon)}_{\gamma,\eta;T}$. □

We define a map $\mathscr{E}_\varepsilon$ which corresponds to multiplication of a modelled distribution by $\varepsilon$. This operator will be used to describe $\varepsilon(\nabla^+_\varepsilon \widetilde{G}^\varepsilon \star^+_\varepsilon \nabla^-_\varepsilon \tilde h^\varepsilon)$ in the error term (2.24b) on the abstract level. As we explain in (5.5), the only element in $\widetilde{\mathcal{W}}$ of the form $\mathcal{E}(\tau)$ is $\mathcal{E}(\partial_+ \widetilde{\mathcal{I}}(\partial_- \mathcal{I}(\Xi)))$. Hence, the domain of $\mathscr{E}_\varepsilon$ should be the modelled distributions $U^\varepsilon : D_\varepsilon \to \text{span}\{\mathbf{1}, \mathbf{\Phi}\}$ with $\mathbf{\Phi} := \partial_+ \widetilde{\widetilde{\mathcal{I}}}(\partial_- \mathcal{I}(\Xi))$. Such modelled distribution can be written as

$$U^\varepsilon(z) = u^\varepsilon_{\mathbf{1}}(z)\mathbf{1} + u^\varepsilon_{\mathbf{\Phi}}(z)\mathbf{\Phi}, \tag{7.5}$$

for real-valued coefficients $u^\varepsilon_{\mathbf{1}}(z)$ and $u^\varepsilon_{\mathbf{\Phi}}(z)$. We set

$$(\mathscr{E}_\varepsilon U^\varepsilon)(z) := \varepsilon(\mathcal{R}^\varepsilon U^\varepsilon)(z)\mathbf{1} + \varepsilon\nabla^-_\varepsilon(\mathcal{R}^\varepsilon U^\varepsilon)(z)X_1 + u^\varepsilon_{\mathbf{\Phi}}(z)\mathcal{E}(\mathbf{\Phi}). \tag{7.6}$$

Now we can prove some properties of the map (7.6).

**Lemma 7.2** *Let the discrete model $Z^\varepsilon$ be such that $\mathcal{R}^\varepsilon\mathcal{E}(\mathbf{\Phi}) \equiv 0$. Then, in the setting of* (7.6),

$$\mathcal{R}^\varepsilon(\mathscr{E}_\varepsilon U^\varepsilon)(z) = \varepsilon(\mathcal{R}^\varepsilon U^\varepsilon)(z). \tag{7.7}$$

*Let furthermore the model satisfy*

$$\Gamma^\varepsilon_{z\bar z}\mathcal{E}(\mathbf{\Phi}) = \mathcal{E}(\mathbf{\Phi}) + g^\varepsilon_{\bar z}(z)\mathbf{1} + \nabla^-_\varepsilon g^\varepsilon_{\bar z}(z)X_1, \tag{7.8}$$

*for any $z, \bar z \in D_\varepsilon$, where $g^\varepsilon_{\bar z}(z) = \varepsilon((\Pi^\varepsilon_{\bar z}\mathbf{\Phi})(z) - (\mathcal{R}^\varepsilon\mathbf{\Phi})(\bar z) + \nabla^-_\varepsilon(\mathcal{R}^\varepsilon\mathbf{\Phi})(\bar z)\,(\bar x - x))$. Then*

$$\|\!|\mathscr{E}_\varepsilon U^\varepsilon|\!\|^{(\varepsilon)}_{\gamma+1,\eta+1;T} \lesssim (1+\sqrt{T})(1+\|\Pi^\varepsilon\|^{(\varepsilon)}_T)\|\!|U^\varepsilon|\!\|^{(\varepsilon)}_{\gamma,\eta;T}, \tag{7.9}$$

*for any $\gamma > |\mathbf{\Phi}|$ and $\eta \in [0, |\mathbf{\Phi}|]$, where the proportionality constant is independent of $\varepsilon$ and $T$.*

*Proof.* By our assumptions (see Definition 5.3), the reconstruction map $\mathcal{R}^\varepsilon$ maps the last two terms in (7.6) to 0, and the desired identity (7.7) follows.

Now, we will prove the bound (7.9). Throughout the proof we will use $z, \bar z \in [-T, T] \times [-1, 1]$, so that we use the norms from Remarks 5.6 and 5.7. From the definition (7.5) and the first terms in (5.16) and (5.18) we get

$$|u^\varepsilon_{\mathbf{1}}(z)| = |U^\varepsilon(z)|_0 \le \|\!|U^\varepsilon|\!\|^{(\varepsilon)}_{\gamma,\eta;T}$$

and

$$|u^\varepsilon_{\mathbf{\Phi}}(z)| = |U^\varepsilon(z)|_{|\mathbf{\Phi}|} \le (\|z\|_{\mathfrak{s}} + \varepsilon)^{\eta-|\mathbf{\Phi}|}\|\!|U^\varepsilon|\!\|^{(\varepsilon)}_{\gamma,\eta;T}, \tag{7.10}$$

where we used our assumption $\eta \in [0, |\mathbf{\Phi}|]$. Moreover, from the definition of the reconstruction map (5.19) and the second bound in (5.15a) on the model (see also (5.12)) we have

$$|(\mathcal{R}^\varepsilon\mathbf{\Phi})(z)| \le \varepsilon^{|\mathbf{\Phi}|}\|\Pi^\varepsilon\|^{(\varepsilon)}_T. \tag{7.11}$$

Then (7.5) yields

$$\begin{aligned}|(\mathcal{R}^\varepsilon U^\varepsilon)(z)| &\le |u^\varepsilon_{\mathbf{1}}(z)| + |u^\varepsilon_{\mathbf{\Phi}}(z)||(\mathcal{R}^\varepsilon\mathbf{\Phi})(z)| \\ &\le \Big(1 + (\|z\|_{\mathfrak{s}} + \varepsilon)^{\eta-|\mathbf{\Phi}|}\varepsilon^{|\mathbf{\Phi}|}\|\Pi^\varepsilon\|^{(\varepsilon)}_T\Big)\|\!|U^\varepsilon|\!\|^{(\varepsilon)}_{\gamma,\eta;T}.\end{aligned}$$

Since $z \in [-T, T] \times [-1, 1]$ and $\eta \in [0, 1]$, we get $(\|z\|_{\mathfrak{s}} + \varepsilon)^\eta \lesssim 1 + \sqrt{T}$. Then the preceding expression is bounded by a constant times $(1+\sqrt{T})(1+\|\Pi^\varepsilon\|^{(\varepsilon)}_T)\|\!|U^\varepsilon|\!\|^{(\varepsilon)}_{\gamma,\eta;T}$. Furthermore, using this bound we get

$$|\varepsilon\nabla^-_\varepsilon(\mathcal{R}^\varepsilon U^\varepsilon)(z)| \lesssim \max_{\bar z \in \{z, z-\varepsilon\}} |(\mathcal{R}^\varepsilon U^\varepsilon)(\bar z)| \lesssim (1+\sqrt{T})(1+\|\Pi^\varepsilon\|^{(\varepsilon)}_T)\|\!|U^\varepsilon|\!\|^{(\varepsilon)}_{\gamma,\eta;T},$$

where we use the shorthand $z - \varepsilon$ for $z - (0, \varepsilon)$.

We are going to use the just derived bound to estimate the coefficients in (7.6). We have

$$|\mathscr{E}_\varepsilon U^\varepsilon(z)|_0 = |\varepsilon(\mathcal{R}^\varepsilon U^\varepsilon)(z)| \lesssim (1+\sqrt{T})(1+\|\Pi^\varepsilon\|^{(\varepsilon)}_T)\|\!|U^\varepsilon|\!\|^{(\varepsilon)}_{\gamma,\eta;T},$$

$$|\mathscr{E}_\varepsilon U^\varepsilon(z)|_{X_1} = |\varepsilon\nabla_\varepsilon^-(\mathcal{R}^\varepsilon U^\varepsilon)(z)| \lesssim (1+\sqrt{T})(1+\|\Pi^\varepsilon\|_T^{(\varepsilon)})\|\!|U^\varepsilon|\!\|_{\gamma,\eta;T}^{(\varepsilon)}$$

and

$$|\mathscr{E}_\varepsilon U^\varepsilon(z)|_{\mathcal{E}(\mathbf{\Phi})} = |u_{\mathbf{\Phi}}^\varepsilon(z)| \leq (\|z\|_\mathfrak{s}+\varepsilon)^{\eta-|\mathbf{\Phi}|}\|\!|U^\varepsilon|\!\|_{\gamma,\eta;T}^{(\varepsilon)} = (\|z\|_\mathfrak{s}+\varepsilon)^{(\eta+1)-|\mathcal{E}(\mathbf{\Phi})|}\|\!|U^\varepsilon|\!\|_{\gamma,\eta;T}^{(\varepsilon)}.$$

We have just prove the part of the bound (7.9) on the first terms in (5.16) and (5.18).

For $z \neq \bar{z}$, we are now going to bound $\mathscr{E}_\varepsilon U^\varepsilon(z) - \Gamma_{z\bar{z}}^\varepsilon \mathscr{E}_\varepsilon U^\varepsilon(\bar{z})$. From assumption (7.8) and algebraic properties of the model we get

$$\Gamma_{z\bar{z}}^\varepsilon \mathscr{E}_\varepsilon U^\varepsilon(\bar{z}) = u_{\mathbf{\Phi}}^\varepsilon(\bar{z})\mathcal{E}(\mathbf{\Phi}) + \varepsilon(\mathcal{R}^\varepsilon \Gamma_{z\bar{z}}^\varepsilon U^\varepsilon(\bar{z}))(z)\mathbf{1} + \varepsilon\nabla_\varepsilon^-(\mathcal{R}^\varepsilon \Gamma_{z\bar{z}}^\varepsilon U^\varepsilon(\bar{z}))(z)X_1$$

and hence

$$\begin{aligned}\mathscr{E}_\varepsilon U^\varepsilon(z) - \Gamma_{z\bar{z}}^\varepsilon \mathscr{E}_\varepsilon U^\varepsilon(\bar{z}) = \varepsilon\mathcal{R}^\varepsilon(U^\varepsilon - \Gamma_{z\bar{z}}^\varepsilon U^\varepsilon(\bar{z}))(z)\mathbf{1} + \varepsilon\nabla_\varepsilon^-\mathcal{R}^\varepsilon(U^\varepsilon - \Gamma_{z\bar{z}}^\varepsilon U^\varepsilon(\bar{z}))(z)X_1 \\ + (u_{\mathbf{\Phi}}^\varepsilon(z) - u_{\mathbf{\Phi}}^\varepsilon(\bar{z}))\mathcal{E}(\mathbf{\Phi}). \qquad (7.12)\end{aligned}$$

We are going to bound the three terms on the right-hand side. For the first term, we use the projection (3.8) to get

$$\begin{aligned}|\mathscr{E}_\varepsilon U^\varepsilon(z) - \Gamma_{z\bar{z}}^\varepsilon \mathscr{E}_\varepsilon U^\varepsilon(\bar{z})|_0 &= \varepsilon|\mathcal{R}^\varepsilon(U^\varepsilon - \Gamma_{z\bar{z}}^\varepsilon U^\varepsilon(\bar{z}))(z)| \\ &\leq \varepsilon|\mathrm{Proj}_{\mathbf{1}}(U^\varepsilon(z) - \Gamma_{z\bar{z}}^\varepsilon U^\varepsilon(\bar{z}))| + \varepsilon|(\mathcal{R}^\varepsilon\mathbf{\Phi})(z)||\mathrm{Proj}_{\mathbf{\Phi}}(U^\varepsilon(z) - \Gamma_{z\bar{z}}^\varepsilon U^\varepsilon(\bar{z}))|.\end{aligned}$$

Using (7.11) we bound it by

$$\begin{aligned}&\varepsilon(\|z-\bar{z}\|_\mathfrak{s}+\varepsilon)^\gamma(\|z,\bar{z}\|_\mathfrak{s}+\varepsilon)^{\eta-\gamma}\|\!|U^\varepsilon|\!\|_{\gamma,\eta;T}^{(\varepsilon)} \\ &\quad + \varepsilon^{1+|\mathbf{\Phi}|}(\|z-\bar{z}\|_\mathfrak{s}+\varepsilon)^{\gamma-|\mathbf{\Phi}|}(\|z,\bar{z}\|_\mathfrak{s}+\varepsilon)^{\eta-\gamma}\|\Pi^\varepsilon\|_T^{(\varepsilon)}\|\!|U^\varepsilon|\!\|_{\gamma,\eta;T}^{(\varepsilon)} \\ &\lesssim (\|z-\bar{z}\|_\mathfrak{s}+\varepsilon)^{\gamma+1}(\|z,\bar{z}\|_\mathfrak{s}+\varepsilon)^{\eta-\gamma}(1+\|\Pi^\varepsilon\|_T^{(\varepsilon)})\|\!|U^\varepsilon|\!\|_{\gamma,\eta;T}^{(\varepsilon)}.\end{aligned}$$

We write the second term in (7.12) as

$$|\mathscr{E}_\varepsilon U^\varepsilon(z) - \Gamma_{z\bar{z}}^\varepsilon \mathscr{E}_\varepsilon U^\varepsilon(\bar{z})|_1 = \varepsilon|\nabla_\varepsilon^-\mathcal{R}^\varepsilon(U^\varepsilon - \Gamma_{z\bar{z}}^\varepsilon U^\varepsilon(\bar{z}))(z)|$$

and bound it in the same way by the same quantity. The last term in (7.12) can be written as

$$\begin{aligned}|\mathscr{E}_\varepsilon U^\varepsilon(z) - \Gamma_{z\bar{z}}^\varepsilon \mathscr{E}_\varepsilon U^\varepsilon(\bar{z})|_{|\mathcal{E}(\mathbf{\Phi})|} &= |u_{\mathbf{\Phi}}^\varepsilon(z) - u_{\mathbf{\Phi}}^\varepsilon(\bar{z})| = |U^\varepsilon(z) - \Gamma_{z\bar{z}}^\varepsilon U^\varepsilon(\bar{z})|_{|\mathbf{\Phi}|} \\ &\leq (\|z-\bar{z}\|_\mathfrak{s}+\varepsilon)^{\gamma-|\mathbf{\Phi}|}(\|z,\bar{z}\|_\mathfrak{s}+\varepsilon)^{\eta-\gamma}\|\!|U^\varepsilon|\!\|_{\gamma,\eta;T}^{(\varepsilon)} \\ &= (\|z-\bar{z}\|_\mathfrak{s}+\varepsilon)^{\gamma+1-|\mathcal{E}(\mathbf{\Phi})|}(\|z,\bar{z}\|_\mathfrak{s}+\varepsilon)^{\eta-\gamma}\|\!|U^\varepsilon|\!\|_{\gamma,\eta;T}^{(\varepsilon)},\end{aligned}$$

where we used (5.2e). This finishes the proof of (7.9). □

### 7.2 Abstract discrete equation

Let us first consider the non-renormalised lift $Z^{\varepsilon,\mathfrak{m}}$ defined in (6.3). One can readily see from Lemma 6.2 and the definition of $\Gamma^{\varepsilon,\mathfrak{m}}$ that the model satisfies the assumptions of Lemma 7.2, which we are going to use below. Thought this section we are using $\gamma = \frac{3}{2} + 2\kappa$ with $\kappa$ the same as in (3.3b).

We define the integration operators on the space of modelled distributions via the kernel $G^\varepsilon$ as in [EH19, Sec. 4]. More precisely, we write $G^\varepsilon = K^\varepsilon + R^\varepsilon$ as in Propositions 4.2. Then we use the singular part $K^\varepsilon$ to define the map $\mathscr{K}_{\gamma_1}^\varepsilon$ as in [EH19, Eq. 4.6] for the value $\beta = 2$ and $\gamma_1 > 0$. We note that we do not need to consider the map $\mathcal{A}^\varepsilon$ from [EH19, Eq. 4.26], since it vanishes in our case (see [EH19, Rem. 4.19]). We lift the smooth part $R^\varepsilon$ to a modelled distribution $\mathscr{R}_{\gamma_2}^\varepsilon$, for $\gamma_2 \in [\gamma_1, \gamma_1 + 2]$, by a Taylor's expansion as in [HM18, Eq. 5.17]. Then we define by analogy with (3.18) the map

$$\mathscr{P}_{\gamma_2}^\varepsilon := \mathscr{K}_{\gamma_1}^\varepsilon + \mathscr{R}_{\gamma_2}^\varepsilon \mathcal{R}^{\varepsilon,\mathfrak{m}} \qquad (7.13)$$

from $\mathcal{D}_\varepsilon^{\gamma_1,\eta_1}(Z^{\varepsilon,\mathfrak{m}})$ to $\mathcal{D}_\varepsilon^{\gamma_2,\eta_2}(Z^{\varepsilon,\mathfrak{m}})$ with suitable $\eta_1$ and $\eta_2$, where $\mathcal{R}^{\varepsilon,\mathfrak{m}}$ is the reconstruction map associated to the model $Z^{\varepsilon,\mathfrak{m}}$. Since the domain of this map will be always clear from the context, we do not indicate $\gamma_1$ in the notation. One of the key properties of this integration map is

$$\mathcal{R}^{\varepsilon,\mathfrak{m}} \mathscr{P}^\varepsilon_{\gamma_2} U^\varepsilon = G^\varepsilon \star^+_\varepsilon \mathcal{R}^{\varepsilon,\mathfrak{m}} U^\varepsilon,$$

for any modelled distirbution $U^\varepsilon$ of regularity at most $\gamma_1$. Similarly, we use the expansions $\widetilde{G}^\varepsilon = \widetilde{K}^\varepsilon + \widetilde{R}^\varepsilon$ and $\widetilde{\widetilde{G}}^\varepsilon = \widetilde{\widetilde{K}}^\varepsilon + \widetilde{\widetilde{R}}^\varepsilon$ from Propositions 4.12 and 4.15 and define the respective maps $\widetilde{\mathscr{P}}^\varepsilon_{\gamma_2}$ and $\widetilde{\widetilde{\mathscr{P}}}^\varepsilon_{\gamma_2}$. Finally, to define the error terms (2.26)-(2.27), we introduce the integration map $\bar{\mathscr{P}}^\varepsilon_{\gamma_2}$, defined as in (7.13) but on the space-time domain $\bar{\Lambda}_\varepsilon \times \Lambda_\varepsilon$ via the convolution (2.25).

Then, for $\gamma = \frac{3}{2} + 2\kappa$ and $\eta \in \mathbb{R}$ we define the abstract equation

$$H^{\varepsilon,\mathfrak{m}} = \mathcal{Q}_{<\gamma} \Big( G^\varepsilon \tilde{h}^\varepsilon_0 - \frac{1}{2} \widetilde{\mathscr{P}}^\varepsilon_\gamma \mathbf{1}_+ F_\varepsilon(H^{\varepsilon,\mathfrak{m}}) + \varrho_\varepsilon \widetilde{\widetilde{\mathscr{P}}}^\varepsilon_\gamma \mathbf{1}_+ (\mathscr{D}_\varepsilon H^{\varepsilon,\mathfrak{m}}) + \widetilde{E}^{\varepsilon,\mathfrak{m}}(H^{\varepsilon,\mathfrak{m}}) + \mathscr{P}^\varepsilon_\gamma \mathbf{1}_+ \Xi \Big), \tag{7.14}$$

for a modelled distribution $H^{\varepsilon,\mathfrak{m}} \in \mathcal{D}_\varepsilon^{\gamma,\eta}(Z^{\varepsilon,\mathfrak{m}})$, where $G^\varepsilon h^\varepsilon_0$ is the polynomial lift of the discrete heat semigroup applied to $h^\varepsilon_0$ and where

$$F_\varepsilon(H^{\varepsilon,\mathfrak{m}})(z) = \mathcal{Q}_{\leq 0}(\mathscr{D}^-_\varepsilon H^{\varepsilon,\mathfrak{m}})(z)(\mathscr{D}^+_\varepsilon H^{\varepsilon,\mathfrak{m}})(z) - \widetilde{C}_{\varepsilon,\mathfrak{m}}(z)\mathbf{1} \tag{7.15}$$

with the random function $\widetilde{C}_{\varepsilon,\mathfrak{m}}$ defined in (2.18). We need to subtract this function to get a renormalised product (2.17) after reconstruction. The projections $\mathbf{1}_+$, $\mathcal{Q}_{<\gamma}$ and $\mathcal{Q}_{\leq 0}$ are the same as in (3.17), and the error term

$$\widetilde{E}^{\varepsilon,\mathfrak{m}}(H^{\varepsilon,\mathfrak{m}}) = \widetilde{E}^{\varepsilon,\mathfrak{m}}_1(H^{\varepsilon,\mathfrak{m}}) + \widetilde{E}^{\varepsilon,\mathfrak{m}}_2(H^{\varepsilon,\mathfrak{m}}) \tag{7.16}$$

with

$$\widetilde{E}^\varepsilon_2(H^{\varepsilon,\mathfrak{m}}) = \varrho_\varepsilon(1 + \varepsilon^{-\frac{1}{2}})(\widetilde{\widetilde{\mathscr{P}}}^\varepsilon_\gamma \mathbf{1}_+ \mathscr{D}_\varepsilon H^{\varepsilon,\mathfrak{m}}) - \frac{1}{2}\mathscr{E}_\varepsilon(\mathscr{D}^+_\varepsilon \widetilde{\widetilde{\mathscr{P}}}^\varepsilon_\gamma \mathbf{1}_+ (\mathscr{D}^-_\varepsilon H^{\varepsilon,\mathfrak{m}})) \tag{7.17}$$

and

$$\widetilde{E}^\varepsilon_1(H^{\varepsilon,\mathfrak{m}}) = \varepsilon^{-\frac{1}{2}} \varrho_\varepsilon (\mathscr{D}^-_\varepsilon \bar{\mathscr{P}}^\varepsilon_{\gamma+2} - \mathscr{D}_\varepsilon \mathscr{P}^\varepsilon_{\gamma+2} + \mathscr{D}_\varepsilon \widetilde{\widetilde{\mathscr{P}}}^\varepsilon_{\gamma+2}) \mathbf{1}_+ H^{\varepsilon,\mathfrak{m}} \tag{7.18}$$

if $\varrho_\varepsilon > 0$, and

$$\widetilde{E}^\varepsilon_1(H^{\varepsilon,\mathfrak{m}}) = \varepsilon^{-\frac{1}{2}} \varrho_\varepsilon (\mathscr{D}^+_\varepsilon \bar{\mathscr{P}}^\varepsilon_{\gamma+2} - \mathscr{D}_\varepsilon \mathscr{P}^\varepsilon_{\gamma+2} + \mathscr{D}_\varepsilon \widetilde{\widetilde{\mathscr{P}}}^\varepsilon_{\gamma+2}) \mathbf{1}_+ H^{\varepsilon,\mathfrak{m}} \tag{7.19}$$

if $\varrho_\varepsilon < 0$. The abstract equation (7.14) is obtained from (2.33) by replacing all operators by their analogues in the framework of regularity structure. Formula for the error term $\widetilde{E}^\varepsilon_1$ follows from Lemma 2.2.

Let us explain the different projections and regularities of the integration maps, appearing in different terms. The projection $\mathcal{Q}_{\leq 0}$ in the non-linearity (7.15) is necessary to have the minimal elements in the regularity structure. This projection makes sure that the non-linearity is a modelled distribution of a very small regularity $\kappa > 0$, which is sufficient to reconstruct it. The most non-trivial term is (7.18)-(7.19), where the integration map should significantly increase the regularity of $H^{\varepsilon,\mathfrak{m}}$. For example, $\mathscr{P}^\varepsilon_{\gamma+2}\mathbf{1}_+ H^{\varepsilon,\mathfrak{m}}$ is a modelled distribution of regularity $\gamma + 2$. Hence, $\mathscr{D}_\varepsilon \mathscr{P}^\varepsilon_{\gamma+2}\mathbf{1}_+ H^{\varepsilon,\mathfrak{m}}$ has regularity $\gamma + 1$, and the difference between this regularity and the required regularity $\gamma$ will be used to compensate the factor $\varepsilon^{-\frac{1}{2}}$. Similarly, such factor will be compensated in the first term in (7.17). This will be done in the proof of Proposition 7.7 below, where we solve this equation.

One can see that the error term (7.17) is well-defined in the sense of (7.6) because the modelled distributions inside $\mathscr{E}_\varepsilon$ are of the form (7.5).

We can show that the reconstruction of the solution $H^{\varepsilon,\mathfrak{m}}$ gives a solution of a discrete equation.

**Lemma 7.3** *If $H^{\varepsilon,\mathfrak{m}} \in \mathcal{D}_\varepsilon^{\gamma,\eta}(Z^{\varepsilon,\mathfrak{m}})$ is a local solution of (7.14) on a time interval $[0,T]$, then the function $\tilde h^{\varepsilon,\mathfrak{m}} = \mathcal{R}^{\varepsilon,\mathfrak{m}} H^{\varepsilon,\mathfrak{m}}$ solves*

$$\tilde h^{\varepsilon,\mathfrak{m}}(t,x) = (G_t^\varepsilon *_\varepsilon \tilde h_0^\varepsilon)(x) - \frac{1}{2}\widetilde{G}^\varepsilon \star_\varepsilon^+ (\nabla_\varepsilon^- \tilde h^{\varepsilon,\mathfrak{m}} \nabla_\varepsilon^+ \tilde h^{\varepsilon,\mathfrak{m}})(t,x) + \varrho_\varepsilon(\widetilde{\widetilde{G}}^\varepsilon \star_\varepsilon^+ \nabla_\varepsilon \tilde h^{\varepsilon,\mathfrak{m}})(t,x) + \widetilde{E}^{\varepsilon,\mathfrak{m}}(\tilde h^{\varepsilon,\mathfrak{m}})(t,x) + (G^\varepsilon \star_\varepsilon^+ d\widetilde{M}^{\varepsilon,\mathfrak{m}})(t,x), \tag{7.20}$$

*on this time interval, which is different from (2.33) only by the renormalised product $\nabla_\varepsilon^- \tilde h^{\varepsilon,\mathfrak{m}} \nabla_\varepsilon^+ \tilde h^{\varepsilon,\mathfrak{m}}$.*

*Proof.* Equation (7.14) is in the form of a fixed-point problem $H^{\varepsilon,\mathfrak{m}} = \mathscr{M}^{\varepsilon,\mathfrak{m}}(H^{\varepsilon,\mathfrak{m}})$, where $\mathscr{M}^{\varepsilon,\mathfrak{m}}(H^{\varepsilon,\mathfrak{m}})$ is the right-hand side, which can be written, recalling the definitions of the involved maps, as

$$\mathscr{M}^{\varepsilon,\mathfrak{m}}(H^{\varepsilon,\mathfrak{m}}) = \mathcal{I}(\Xi) - \frac{1}{2}\mathcal{Q}_{<\gamma}\widetilde{\mathcal{I}}(\mathbf{1}_+\mathcal{Q}_{\le 0}(\partial_- H^{\varepsilon,\mathfrak{m}})(\partial_+ H^{\varepsilon,\mathfrak{m}})) + \varrho_\varepsilon \widetilde{\widetilde{\mathcal{I}}}(\mathbf{1}_+(\partial H^{\varepsilon,\mathfrak{m}})) - \frac{1}{2}\mathcal{E}(\partial_+\widetilde{\widetilde{\mathcal{I}}}(\mathbf{1}_+\partial_- H^{\varepsilon,\mathfrak{m}})) + \mathscr{N}^{\varepsilon,\mathfrak{m}}(H^{\varepsilon,\mathfrak{m}}),$$

where $\mathscr{N}^{\varepsilon,\mathfrak{m}}(H^{\varepsilon,\mathfrak{m}})$ takes values in a linear span of $\mathbf{1}$ and $X_1$ because we are solving the equation in $\mathcal{D}_\varepsilon^{\gamma,\eta}$ with $\gamma$ close to $\frac{3}{2}$. The last assumption implies that only the monomials of degrees smaller then $\gamma$ appear in the polynomial part of $\mathscr{M}^{\varepsilon,\mathfrak{m}}(H^{\varepsilon,\mathfrak{m}})$. The term containing $\mathcal{E}$ comes from the error term (7.17), while the other error term makes contributions only to the polynomial part $\mathscr{N}^{\varepsilon,\mathfrak{m}}(H^{\varepsilon,\mathfrak{m}})$.

We can perform Picard iterations to get a precise expression for $H^{\varepsilon,\mathfrak{m}}$, as a solution of this fixed-point problem. Namely, for two functions $h^{(0)}$ and $\bar h^{(0)}$ on the domain $D_\varepsilon$, the first approximation of $H^{\varepsilon,\mathfrak{m}}$ is

$$H^{(0)}(z) = h^{(0)}(z)\mathbf{1} + \bar h^{(0)}(z)X_1.$$

Computing $H^{(1)} = \mathscr{M}^{\varepsilon,\mathfrak{m}}(H^{(0)})$, we get the next approximation of the solution

$$H^{(1)}(z) = h^{(1)}(z)\mathbf{1} + \treeI + \bar h^{(1)}(z)X_1,$$

for some new functions $h^{(1)}$ and $\bar h^{(1)}$. We keep computing $H^{(i)} = \mathscr{M}^{\varepsilon,\mathfrak{m}}(H^{(i-1)})$ for $i \ge 2$ to get a sequence of approximate solutions:

$$\begin{aligned} H^{(2)}(z) &= h^{(2)}(z)\mathbf{1} + \treeI - \frac{1}{2}\treeY + \bar h^{(2)}(z)X_1 - \frac{1}{2}\bar h^{(1)}(z)(\treeL + \treeLb) \\ &\qquad + \frac{1}{2}\varrho_\varepsilon(\treeL + \treeLb) - \frac{1}{2}\mathcal{E}(\treeLc), \\ H^{(3)}(z) &= h^{(3)}(z)\mathbf{1} + \treeI - \frac{1}{2}\treeY + \bar h^{(3)}(z)X_1 + \frac{1}{4}(\treeYY + \treeYYb) - \frac{1}{2}\bar h^{(2)}(z)(\treeL + \treeLb) \\ &\qquad + \frac{1}{2}\varrho_\varepsilon(\treeL + \treeLb) - \frac{1}{2}\mathcal{E}(\treeLc), \end{aligned}$$

where we used (5.7g). The approximations $H^{(i)}$ for $i \ge 4$ have the same form, and written in terms of different functions $h^{(i)}$, $\bar h^{(i)}$ and $\bar h^{(i-1)}$.

Since we assumed that $H^{\varepsilon,\mathfrak{m}} \in \mathcal{D}_\varepsilon^{\gamma,\eta}(Z^{\varepsilon,\mathfrak{m}})$ is a local solution, it is the limit of $H^{(i)}$ as $i \to \infty$. Hence, the involved functions $h^{(i)}$ and $\bar h^{(i)}$ have limits as $i \to \infty$, and we get

$$\begin{aligned} H^{\varepsilon,\mathfrak{m}}(z) &= \tilde h^{\varepsilon,\mathfrak{m}}(z)\mathbf{1} + \treeI - \frac{1}{2}\treeY + \bar h^{\varepsilon,\mathfrak{m}}(z)X_1 + \frac{1}{4}(\treeYY + \treeYYb) - \frac{1}{2}\bar h^{\varepsilon,\mathfrak{m}}(z)(\treeL + \treeLb) \\ &\qquad + \frac{1}{2}\varrho_\varepsilon(\treeL + \treeLb) - \frac{1}{2}\mathcal{E}(\treeLc), \end{aligned} \tag{7.21}$$

where $h^{\varepsilon,\mathfrak{m}}$ and $\bar h^{\varepsilon,\mathfrak{m}}$ are two functions on $D_\varepsilon$. By Lemma 6.2, the reconstruction map $\mathcal{R}^{\varepsilon,\mathfrak{m}}$ maps all elements on the right-hand side, except $\mathbf{1}$, to zero. Hence, we get $\tilde h^{\varepsilon,\mathfrak{m}} = \mathcal{R}^{\varepsilon,\mathfrak{m}} H^{\varepsilon,\mathfrak{m}}$ from this expansion.

Applying the reconstruction map to both sides of (7.14) and using definitions of the involved operators, we get

$$\tilde{h}^{\varepsilon,\mathfrak{m}} = G^{\varepsilon}\tilde{h}_0^{\varepsilon} - \frac{1}{2}\widetilde{G}^{\varepsilon} \star_{\varepsilon}^{+} \mathcal{R}^{\varepsilon,\mathfrak{m}} F_{\varepsilon}(H^{\varepsilon,\mathfrak{m}}) + \varrho_{\varepsilon}\widetilde{\widetilde{G}}^{\varepsilon} \star_{\varepsilon}^{+} (\mathcal{R}^{\varepsilon,\mathfrak{m}}\mathscr{D}_{\varepsilon}H^{\varepsilon,\mathfrak{m}}) \\ + \mathcal{R}^{\varepsilon,\mathfrak{m}}\widetilde{E}^{\varepsilon,\mathfrak{m}}(H^{\varepsilon,\mathfrak{m}}) + G^{\varepsilon} \star_{\varepsilon}^{+} d\widetilde{M}^{\varepsilon,\mathfrak{m}}. \tag{7.22}$$

By Lemma 6.2 the reconstruction map satisfies $\mathcal{R}^{\varepsilon,\mathfrak{m}}\tau = 0$ for elements $|\tau| > 0$ which are not of the form (6.5)-(6.6). Then the projection $\mathcal{Q}_{\leq 0}$ can be omitted in the second term on the right-hand side of (7.22). Using then multiplicativity of the reconstruction map and Lemma 7.1, we get

$$\mathcal{R}^{\varepsilon,\mathfrak{m}} F_{\varepsilon}(H^{\varepsilon,\mathfrak{m}}) = \mathcal{R}^{\varepsilon,\mathfrak{m}}(\mathscr{D}_{\varepsilon}^{-} H^{\varepsilon,\mathfrak{m}})(\mathscr{D}_{\varepsilon}^{+} H^{\varepsilon,\mathfrak{m}})(z) \\ = \mathcal{R}^{\varepsilon,\mathfrak{m}}(\mathscr{D}_{\varepsilon}^{-} H^{\varepsilon,\mathfrak{m}})(z)\mathcal{R}^{\varepsilon,\mathfrak{m}}(\mathscr{D}_{\varepsilon}^{+} H^{\varepsilon,\mathfrak{m}})(z) = \nabla_{\varepsilon}^{-}\tilde{h}^{\varepsilon,\mathfrak{m}}(z)\nabla_{\varepsilon}^{+}\tilde{h}^{\varepsilon,\mathfrak{m}}(z). \tag{7.23}$$

The third term on the right-hand side of (7.22) equals $\nabla_{\varepsilon}\tilde{h}^{\varepsilon,\mathfrak{m}}(z)$, as follows from the properties of the involved operators.

It is left to compute the reconstruction of the error term in (7.22). From (7.16) we get

$$\mathcal{R}^{\varepsilon,\mathfrak{m}}\widetilde{E}^{\varepsilon,\mathfrak{m}}(H^{\varepsilon,\mathfrak{m}}) = \mathcal{R}^{\varepsilon,\mathfrak{m}}\widetilde{E}_1^{\varepsilon,\mathfrak{m}}(H^{\varepsilon,\mathfrak{m}}) + \mathcal{R}^{\varepsilon,\mathfrak{m}}\widetilde{E}_2^{\varepsilon,\mathfrak{m}}(H^{\varepsilon,\mathfrak{m}}).$$

The properties of the involved operators yield $\mathcal{R}^{\varepsilon,\mathfrak{m}}\widetilde{E}_1^{\varepsilon,\mathfrak{m}}(H^{\varepsilon,\mathfrak{m}}) = \widetilde{E}_1^{\varepsilon,\mathfrak{m}}(\tilde{h}^{\varepsilon,\mathfrak{m}})$, where the latter is given in Lemma 2.2. Furthermore, the definition (7.17) and Lemmas 7.1 and 7.2 yield that $\mathcal{R}^{\varepsilon,\mathfrak{m}}\widetilde{E}_2^{\varepsilon,\mathfrak{m}}(H^{\varepsilon,\mathfrak{m}})$ equals $\widetilde{E}_2^{\varepsilon,\mathfrak{m}}(\tilde{h}^{\varepsilon,\mathfrak{m}})$, where the latter is defined in (2.24b).

We have just showed that equation (7.22) coincides with (7.20) which implies that the solution is $\tilde{h}^{\varepsilon,\mathfrak{m}} = \mathcal{R}^{\varepsilon,\mathfrak{m}}H^{\varepsilon,\mathfrak{m}}$. □

We do not solve equation (7.14) defined via the discrete model $Z^{\varepsilon,\mathfrak{m}}$ because the latter is not bounded uniformly in $\varepsilon$. It means that we do not have a solution of (7.14) which is bounded uniformly in $\varepsilon$ in a suitable space of modelled distributions. In order to have such a uniformity, we need to define the equation with respect to a renormalised model $\widehat{Z}^{\varepsilon,\mathfrak{m}}$ from Section 6.1.

## 7.3 Renormalised equation

Our next goal is to define an abstract equation with respect to the renormalised model $\widehat{Z}^{\varepsilon,\mathfrak{m}}$ from Section 6.1. It is not difficult to see that if we consider equation (7.14) with respect to the model $\widehat{Z}^{\varepsilon,\mathfrak{m}}$, then we do not get a solution of (2.33) after reconstruction, because the renormalisation terms will appear on the right-hand side of the reconstructed equation. In order to use the model $\widehat{Z}^{\varepsilon,\mathfrak{m}}$, we need to make respective modifications of equation (7.14).

One can see that the renormalisation does not affect the assumptions of Lemma 7.2 and it can be used also for the model $\widehat{Z}^{\varepsilon,\mathfrak{m}}$. Then we define a modified version of (7.14) with respect to this model:

$$H^{\varepsilon,\mathfrak{m}} = \mathcal{Q}_{<\gamma}\Big(G^{\varepsilon}\tilde{h}_0^{\varepsilon} - \frac{1}{2}\widetilde{\mathscr{P}}_{\gamma}^{\varepsilon}\mathbf{1}_{+}\mathcal{Q}_{\leq 0}(\mathscr{D}_{\varepsilon}^{-} H^{\varepsilon,\mathfrak{m}})(\mathscr{D}_{\varepsilon}^{+} H^{\varepsilon,\mathfrak{m}}) + \varrho_{\varepsilon}\widetilde{\widetilde{\mathscr{P}}}_{\gamma}^{\varepsilon}\mathbf{1}_{+}(\mathscr{D}_{\varepsilon}H^{\varepsilon,\mathfrak{m}}) \\ + \widetilde{E}^{\varepsilon,\mathfrak{m}}(H^{\varepsilon,\mathfrak{m}}) + \mathscr{P}_{\gamma}^{\varepsilon}\mathbf{1}_{+}\Xi + \widetilde{\mathscr{P}}_{\gamma}^{\varepsilon}\mathbf{1}_{+}\widetilde{E}_0^{\varepsilon}\Big), \tag{7.24}$$

where the involved operators are defined as in (7.14) but via the model $\widehat{Z}^{\varepsilon,\mathfrak{m}}$. The error term $\widetilde{E}^{\varepsilon,\mathfrak{m}}$ is defined as in (7.16) but via the renormalised model. The new error term in this equation is

$$\widetilde{E}_0^{\varepsilon} = \frac{1}{4}(\widetilde{C}_{\varepsilon}(\text{[tree]}) + \widetilde{C}_{\varepsilon}(\text{[tree]}))(\text{[tree]}_{+} - \text{[tree]}_{-}), \tag{7.25}$$

which is defined via the renormalisation constants introduced in Section 6.1.

We show in the proof of Lemma 7.4 below, that the application of the renormalised reconstruction map to equation (7.24) causes appearance of some error terms proportional to the renormalisation constants. The goal of $\widetilde{E}_0^{\varepsilon}$ is to cancel those error terms.

We can show that a reconstruction of a solution for (7.24) gives a solution to (2.33).

**Lemma 7.4** *If $H^{\varepsilon,\mathfrak{m}}$ is a local solution of* (7.24) *on* $[0,T]$, *then the function*

$$\tilde h^{\varepsilon,\mathfrak{m}}(t,x) = \widehat{\mathcal{R}}^{\varepsilon,\mathfrak{m}} H^{\varepsilon,\mathfrak{m}}(t,x) + \tilde c_0^\varepsilon t$$

*solves* (2.33) *on the same time interval, where*

$$\tilde c_0^\varepsilon = -\frac{1}{2}\widetilde C_\varepsilon(\ ) - \frac{1}{8}\Big(\widetilde C_\varepsilon(\ ) + \widetilde C_\varepsilon(\ ) + \widetilde C_\varepsilon(\ ) + \widetilde C_\varepsilon(\ ) + \widetilde C_\varepsilon(\ )\Big) + \frac{1}{4}\Big(\widetilde C_\varepsilon(\ ) + \widetilde C_\varepsilon(\ )\Big), \tag{7.26}$$

*where we use the renormalisation constants defined in* (6.15)-(6.17). *Moreover,* $\tilde c_0^\varepsilon$ *is bounded uniformly in* $\varepsilon$.

*Proof.* The fact that $\tilde c_0^\varepsilon$ is bounded uniformly in $\varepsilon$ follows from Lemma 6.5, 6.20 and 6.9.

As in the proof of Lemma 7.3, we get the expansion (7.21) for the solution $H^{\varepsilon,\mathfrak{m}}$ with some functions $\tilde h^{\varepsilon,\mathfrak{m}}$ and $\bar h^{\varepsilon,\mathfrak{m}}$. Applying the reconstruction map $\widehat{\mathcal{R}}^{\varepsilon,\mathfrak{m}}$ to the equations (7.24), we get

$$\tilde h^{\varepsilon,\mathfrak{m}} = G^\varepsilon \tilde h_0^\varepsilon - \frac{1}{2}\widetilde G^\varepsilon \star_\varepsilon^+ \widehat{\mathcal{R}}^{\varepsilon,\mathfrak{m}} \mathcal{Q}_{\le 0}(\mathscr{D}_\varepsilon^- H^{\varepsilon,\mathfrak{m}})(\mathscr{D}_\varepsilon^+ H^{\varepsilon,\mathfrak{m}}) + \varrho_\varepsilon \widetilde{\widetilde G}^\varepsilon \star_\varepsilon^+ (\widehat{\mathcal{R}}^{\varepsilon,\mathfrak{m}} \mathscr{D}_\varepsilon H^{\varepsilon,\mathfrak{m}}) + \widehat{\mathcal{R}}^{\varepsilon,\mathfrak{m}} \widetilde E^{\varepsilon,\mathfrak{m}}(H^{\varepsilon,\mathfrak{m}}) + G^\varepsilon \star_\varepsilon^+ d\widetilde M^{\varepsilon,\mathfrak{m}} + \frac{1}{2}\widetilde G^\varepsilon \star_\varepsilon^+ (\widehat{\mathcal{R}}^{\varepsilon,\mathfrak{m}} \widetilde E_0^\varepsilon). \tag{7.27}$$

The reconstructions of all the terms, except the non-linearity and the last error term, are exactly the same as in (7.22) and they give the terms of equation (2.33) except the renormalised product.

We get from the expansion (7.21) of $H^{\varepsilon,\mathfrak{m}}$ and the definitions (7.1)

$$(\mathscr{D}_\varepsilon^- H^{\varepsilon,\mathfrak{m}})(z) = \ - \frac{1}{2}\ + \bar h_-^{\varepsilon,\mathfrak{m}}(z)\mathbf{1} + \frac{1}{4}(\ + \ ) - \frac{1}{2}\bar h^{\varepsilon,\mathfrak{m}}(z)(\ + \ ) + \frac{1}{2}\varrho_\varepsilon(\ + \ ) - \frac{1}{2}$$

and

$$(\mathscr{D}_\varepsilon^+ H^{\varepsilon,\mathfrak{m}})(z) = \ - \frac{1}{2}\ + \bar h_+^{\varepsilon,\mathfrak{m}}(z)\mathbf{1} + \frac{1}{4}(\ + \ ) - \frac{1}{2}\bar h^{\varepsilon,\mathfrak{m}}(z)(\ + \ ) + \frac{1}{2}\varrho_\varepsilon(\ + \ ) - \frac{1}{2}\ ,$$

where $\bar h_\pm^{\varepsilon,\mathfrak{m}}(z) = \bar h^{\varepsilon,\mathfrak{m}}(z) \mp \frac{1}{\varepsilon}(\Pi_z^\varepsilon H^\varepsilon(z) - \Pi_{z\pm\varepsilon}^\varepsilon H^\varepsilon(z\pm\varepsilon))(z\pm\varepsilon)$. Then

$$\begin{aligned}\mathcal{Q}_{\le 0}(\mathscr{D}_\varepsilon^- H^{\varepsilon,\mathfrak{m}})(z)(\mathscr{D}_\varepsilon^+ H^{\varepsilon,\mathfrak{m}})(z) = &\ - \frac{1}{2}(\ + \ ) + \Big(\bar h_-^{\varepsilon,\mathfrak{m}}(z)\ + \bar h_+^{\varepsilon,\mathfrak{m}}(z)\ \Big)\\ &+ \frac{1}{4}\Big(\ + \ + \ + \ \Big) + \frac{1}{4}\ - \frac{1}{2}\bar h^{\varepsilon,\mathfrak{m}}(z)\Big(\ + \ + \ + \ \Big)\\ &+ \frac{1}{2}\varrho_\varepsilon\Big(\ + \ + \ + \ \Big) - \frac{1}{2}\Big(\ + \ \Big)\\ &- \frac{1}{2}\Big(\bar h_-^{\varepsilon,\mathfrak{m}}(z)\ + \bar h_+^{\varepsilon,\mathfrak{m}}(z)\ \Big) + \bar h_-^{\varepsilon,\mathfrak{m}}(z)\bar h_+^{\varepsilon,\mathfrak{m}}(z)\mathbf{1}.\end{aligned}$$

We note that the reconstruction map $\widehat{\mathcal{R}}^{\varepsilon,\mathfrak{m}}$ is not multiplicative because of the renormalisation defined at the beginning of Section 6.1. More precisely, we have due to (6.13)

$$(\widehat{\mathcal{R}}^{\varepsilon,\mathfrak{m}}\ )(z) = (\widehat{\mathcal{R}}^{\varepsilon,\mathfrak{m}}\ )(z)(\widehat{\mathcal{R}}^{\varepsilon,\mathfrak{m}}\ )(z) - \widetilde C_\varepsilon(\ ) - \widetilde C_{\varepsilon,\mathfrak{m}}(z).$$

Any other element $\tau \in \widetilde{\mathcal{W}} \setminus \{\vee\}$ is of the form $\tau = \tau_1\tau_2$ and we have

$$(\widehat{\mathcal{R}}^{\varepsilon,\mathfrak{m}}\tau)(z) = (\widehat{\mathcal{R}}^{\varepsilon,\mathfrak{m}}\tau_1)(z)(\widehat{\mathcal{R}}^{\varepsilon,\mathfrak{m}}\tau_2)(z) - \widetilde{C}_\varepsilon(\tau),$$

where the latter is the respective renormalisation constant from (6.8). For the elements $\tau = \tau_1\tau_2 \notin \widetilde{\mathcal{W}}$ we have $(\widehat{\mathcal{R}}^{\varepsilon,\mathfrak{m}}\tau)(z) = (\widehat{\mathcal{R}}^{\varepsilon,\mathfrak{m}}\tau_1\tau_2)(z)$. These identities yield

$$\begin{aligned}\widehat{\mathcal{R}}^{\varepsilon,\mathfrak{m}}(\mathscr{D}^-_\varepsilon H^{\varepsilon,\mathfrak{m}})(\mathscr{D}^+_\varepsilon H^{\varepsilon,\mathfrak{m}})(z) &= \widehat{\mathcal{R}}^{\varepsilon,\mathfrak{m}}(\mathscr{D}^-_\varepsilon H^{\varepsilon,\mathfrak{m}})(z)\widehat{\mathcal{R}}^{\varepsilon,\mathfrak{m}}(\mathscr{D}^+_\varepsilon H^{\varepsilon,\mathfrak{m}})(z) - \widetilde{C}_\varepsilon(\vee) - \widetilde{C}_{\varepsilon,\mathfrak{m}}(z)\\ &\quad - \frac14\Big(\widetilde{C}_\varepsilon(\ast) + \widetilde{C}_\varepsilon(\ast) + \widetilde{C}_\varepsilon(\ast) + \widetilde{C}_\varepsilon(\ast) + \widetilde{C}_\varepsilon(\ast)\Big)\\ &\quad + \frac12\Big((\widetilde{C}_\varepsilon(\ast) + \widetilde{C}_\varepsilon(\ast))(\widehat{\mathcal{R}}^{\varepsilon,\mathfrak{m}}\uparrow_+)(z) + (\widetilde{C}_\varepsilon(\ast) + \widetilde{C}_\varepsilon(\ast))(\widehat{\mathcal{R}}^{\varepsilon,\mathfrak{m}}\uparrow_-)(z)\Big)\\ &\quad + \frac12\bar{h}^{\varepsilon,\mathfrak{m}}(z)\Big(\widetilde{C}_\varepsilon(\ast) + \widetilde{C}_\varepsilon(\ast) + \widetilde{C}_\varepsilon(\ast) + \widetilde{C}_\varepsilon(\ast)\Big)\\ &\quad - \frac12\varrho_\varepsilon\Big(\widetilde{C}_\varepsilon(\ast) + \widetilde{C}_\varepsilon(\ast) + \widetilde{C}_\varepsilon(\ast) + \widetilde{C}_\varepsilon(\ast)\Big)\\ &\quad + \frac12\Big(\widetilde{C}_\varepsilon(\ast) + \widetilde{C}_\varepsilon(\ast)\Big),\end{aligned}\tag{7.28}$$

where we omitted the projection $\mathcal{Q}_{\leq 0}$ by the same reason as in (7.23). Lemma 7.1 yields $\widehat{\mathcal{R}}^{\varepsilon,\mathfrak{m}}(\mathscr{D}^-_\varepsilon H^{\varepsilon,\mathfrak{m}})(z)\widehat{\mathcal{R}}^{\varepsilon,\mathfrak{m}}(\mathscr{D}^+_\varepsilon H^{\varepsilon,\mathfrak{m}})(z) = \nabla^-_\varepsilon h^{\varepsilon,\mathfrak{m}}(z)\nabla^+_\varepsilon h^{\varepsilon,\mathfrak{m}}(z)$. We recall moreover our definition $\nabla^-_\varepsilon h^{\varepsilon,\mathfrak{m}} \diamond \nabla^+_\varepsilon h^{\varepsilon,\mathfrak{m}} = \nabla^-_\varepsilon h^{\varepsilon,\mathfrak{m}}\nabla^+_\varepsilon h^{\varepsilon,\mathfrak{m}} - \widetilde{C}_{\varepsilon,\mathfrak{m}}$ of the renormalised product in (2.17). Furthermore, Lemma 6.8 allows to write the terms in (7.28) involving the renormalisation constants as

$$\begin{aligned}&-\widetilde{C}_\varepsilon(\vee) - \frac14\Big(\widetilde{C}_\varepsilon(\ast) + \widetilde{C}_\varepsilon(\ast) + \widetilde{C}_\varepsilon(\ast) + \widetilde{C}_\varepsilon(\ast) + \widetilde{C}_\varepsilon(\ast)\Big)\\ &\qquad + \frac12(\widetilde{C}_\varepsilon(\ast) + \widetilde{C}_\varepsilon(\ast))\widehat{\mathcal{R}}^{\varepsilon,\mathfrak{m}}(\uparrow_+ - \uparrow_-)(z) + \frac12\Big(\widetilde{C}_\varepsilon(\ast) + \widetilde{C}_\varepsilon(\ast)\Big)\end{aligned}$$

Plugging all these into (7.27), the terms $\frac12(\widetilde{C}_\varepsilon(\ast) + \widetilde{C}_\varepsilon(\ast))\widehat{\mathcal{R}}^{\varepsilon,\mathfrak{m}}(\uparrow_+ - \uparrow_-)(z)$ cancel out with the last error term in (7.27). The other terms in the preceding expression cancel out by the constant $\tilde{c}^\varepsilon_0$ as in the statement of this lemma and we get equation (2.33) as desired. □

Let furthermore $\widehat{Z}^{\varepsilon,\mathfrak{m},\delta}$ be the discrete renormalised model, defined in Section 6.4, and let us consider the equation

$$H^{\varepsilon,\mathfrak{m},\delta} = \mathcal{Q}_{<\gamma}\Big(G^\varepsilon \tilde{h}^{\varepsilon,\delta}_0 - \frac12\widetilde{\mathscr{P}}^{\varepsilon,\delta}_\gamma \mathbf{1}_+\mathcal{Q}_{\leq 0}(\mathscr{D}^-_{\varepsilon,\delta}H^{\varepsilon,\mathfrak{m},\delta})(\mathscr{D}^+_{\varepsilon,\delta}H^{\varepsilon,\mathfrak{m},\delta}) + \mathscr{P}^{\varepsilon,\delta}_\gamma\mathbf{1}_+\Xi + \widetilde{\mathscr{P}}^{\varepsilon,\delta}_\gamma\mathbf{1}_+\widetilde{E}^{\varepsilon,\delta}_0\Big),\tag{7.29}$$

with respect to this model, where the error term is

$$\widetilde{E}^{\varepsilon,\delta}_0 = \frac14(\widetilde{C}_{\varepsilon,\delta}(\ast) + \widetilde{C}_{\varepsilon,\delta}(\ast))(\uparrow_+ - \uparrow_-)\tag{7.30}$$

and the involved renormalisation constants are defined in Section 6.4. The integration maps $\mathscr{P}^{\varepsilon,\delta}_\gamma$ and $\widetilde{\mathscr{P}}^{\varepsilon,\delta}_\gamma$ are defined as above but with respect to this model, what makes them depend on $\delta$. The maps $\mathscr{D}^\pm_{\varepsilon,\delta}$ are defined by (7.1) with respect to this model. The following result is proved similarly to Lemma 7.4.

**Lemma 7.5** *Let $H^{\varepsilon,\mathfrak{m},\delta}$ be a local solution of (7.29) on $[0,T]$, and let*

$$\tilde{c}^{\varepsilon,\delta}_0 = -\frac18\Big(\widetilde{C}_{\varepsilon,\delta}(\ast) + \widetilde{C}_{\varepsilon,\delta}(\ast) + \widetilde{C}_{\varepsilon,\delta}(\ast) + \widetilde{C}_{\varepsilon,\delta}(\ast) + \widetilde{C}_{\varepsilon,\delta}(\ast)\Big).\tag{7.31}$$

*Then the function* $\tilde h^{\varepsilon,\mathfrak{m},\delta}(t,x) = \widehat{\mathcal{R}}^{\varepsilon,\mathfrak{m},\delta} H^{\varepsilon,\mathfrak{m},\delta}(t,x) + \tilde c_0^{\varepsilon,\delta} t$ *solves*

$$\begin{aligned}\tilde h^{\varepsilon,\mathfrak{m},\delta}(t,x) = (G_t^\varepsilon *_\varepsilon \tilde h_0^{\varepsilon,\delta})(x) - \frac{1}{2}\widetilde G^\varepsilon \star_\varepsilon^+ (\nabla_\varepsilon^- \tilde h^{\varepsilon,\mathfrak{m},\delta} \diamond \nabla_\varepsilon^+ \tilde h^{\varepsilon,\mathfrak{m},\delta} - \widetilde C_{\varepsilon,\delta}(\vee))(t,x) \\ + (G^\varepsilon \star_\varepsilon^+ d\widetilde M^{\varepsilon,\mathfrak{m},\delta})(t,x)\end{aligned} \tag{7.32}$$

*on the same time interval. The renormalised product is defined here as*

$$(\nabla_\varepsilon^- \tilde h^{\varepsilon,\mathfrak{m},\delta} \diamond \nabla_\varepsilon^+ \tilde h^{\varepsilon,\mathfrak{m},\delta})(z) := (\nabla_\varepsilon^- \tilde h^{\varepsilon,\mathfrak{m},\delta} \nabla_\varepsilon^+ \tilde h^{\varepsilon,\mathfrak{m},\delta})(z) - \widetilde C_{\varepsilon,\mathfrak{m},\delta}(z), \tag{7.33}$$

*via the function* (6.40). *Moreover,* $\tilde c_0^{\varepsilon,\delta}$ *has a limit as* $\delta \to 0$ *uniformly in* $\varepsilon$.

Note that we do not include that error term $\widetilde E^{\varepsilon,\mathfrak{m}}(H^{\varepsilon,\mathfrak{m}})$ into equation (7.29), what makes the constant (7.31) slightly easier compared to (7.26). Moreover, we do not include that constant $\widetilde C_{\varepsilon,\delta}(\vee)$ into (7.31) because it is not bounded uniformly in $\delta$.

Before proceeding, we need to prove bounds on the integration maps involved in the dentition of the error term $\widetilde E_1^{\varepsilon,\mathfrak{m}}$ from (7.16).

**Lemma 7.6** *For any* $\alpha \in [0,1]$ *there is a constant* $C = C(\alpha) > 0$ *such that for any modelled distribution* $U^\varepsilon$ *from the domains of* $\mathscr{P}_{\gamma+2}^\varepsilon$ *and* $\bar{\mathscr{P}}_{\gamma+2}^\varepsilon$ *and for any* $T > 0$ *one has*

$$\|(\bar{\mathscr{P}}_{\gamma+2}^\varepsilon - \mathscr{P}_{\gamma+2}^\varepsilon) U^\varepsilon\|^{(\varepsilon)}_{\gamma+2-\alpha,\eta+2-\alpha;T} \le C\varepsilon^\alpha \|U^\varepsilon\|^{(\varepsilon)}_{\gamma,\eta;T}. \tag{7.34}$$

*Similarly, for any* $U^\varepsilon$ *from the domain of* $\widetilde{\mathscr{P}}_{\gamma+2}^\varepsilon$ *one has*

$$\|\widetilde{\mathscr{P}}_{\gamma+2}^\varepsilon U^\varepsilon\|^{(\varepsilon)}_{\gamma+2-\alpha,\eta+2-\alpha;T} \le C\varepsilon^\alpha \|U^\varepsilon\|^{(\varepsilon)}_{\gamma,\eta;T}. \tag{7.35}$$

*Proof.* The proofs of these bounds are analogous to the respective proof from [EH19, Thm. 4.9], combined with

$$|(f \bar\star_\varepsilon g)(z) - (f \star_\varepsilon^+ g)(z)| \lesssim \varepsilon^\alpha (\|f\|_{L^\infty}\|g\|_{\mathcal{C}^\alpha} + \|f\|_{\mathcal{C}^\alpha}\|g\|_{L^\infty})$$

in the case of (7.34) (see the definition of the convolutions in (1.18) and (2.25)). The bound (7.35) is proved similarly, using the estimates (4.41) on the discrete kernel. □

We can now study the solution of equation (7.24). In the statement of the following proposition we use the discrete norms defined in (1.21).

**Proposition 7.7** *Let* $\widehat Z^{\varepsilon,\mathfrak{m}}$ *be the renormalised discrete model defined in Section 6.1 and let the initial state* $\tilde h_0^\varepsilon$ *satisfy the assumptions of Theorem 1.1. Then there exists* $\varepsilon_\star \in (0,1)$ *such that the following holds: for every* $\varepsilon \in (0,\varepsilon_\star)$ *and almost every realisation of* $\widehat Z^{\varepsilon,\mathfrak{m}}$ *there exists a time* $T_{\varepsilon,\mathfrak{m}} \in (0,\infty]$ *such that equation* (7.24) *has a unique local solution* $H^{\varepsilon,\mathfrak{m}} \in \mathcal{D}_\varepsilon^{\gamma,\eta}(\widehat Z^{\varepsilon,\mathfrak{m}})$ *on the time interval* $[0, T_{\varepsilon,\mathfrak{m}})$*, with* $\gamma = \frac{3}{2} + 2\kappa$ *and* $\eta = \alpha$*, where* $\alpha$ *is as in the statement of Theorem 1.1 and* $\kappa$ *is from* (5.2b).

*Let moreover* $\tilde h^{\varepsilon,\mathfrak{m}} = \widehat{\mathcal{R}}^{\varepsilon,\mathfrak{m}} H^{\varepsilon,\mathfrak{m}} + \tilde c_0^\varepsilon t$*, where* $\widehat{\mathcal{R}}^{\varepsilon,\mathfrak{m}}$ *is the reconstruction map associated to the model and* $\tilde c_0^\varepsilon$ *is defined in* (7.26). *Then for every* $L > 0$ *there is* $T_{\varepsilon,\mathfrak{m}}^L \in (0, T_{\varepsilon,\mathfrak{m}})$*, such that* $\lim_{L\to\infty} T_{\varepsilon,\mathfrak{m}}^L = T_{\varepsilon,\mathfrak{m}}$ *almost surely, and*

$$\sup_{t\in[0,T\wedge T_{\varepsilon,\mathfrak{m}}^L]} \|\tilde h^{\varepsilon,\mathfrak{m}}(t) - \tilde c_0^\varepsilon t\|^{(\varepsilon)}_{\mathcal{C}^\alpha} \le C, \tag{7.36}$$

*for any* $T > 0$*, provided* $\|\tilde h_0^{\varepsilon,\mathfrak{m}}\|^{(\varepsilon)}_{\mathcal{C}^\alpha} \le L$ *and* $\|\widehat Z^{\varepsilon,\mathfrak{m}}\|^{(\varepsilon)}_{T+1} \le L$*. The constant* $C$ *depends polynomially on* $L$ *and* $\mathfrak{m}$ *and is independent of* $\varepsilon$ *and* $T$.

*Let furthermore $\widehat{Z}^{\varepsilon,\mathfrak{m},\delta}$ be the discrete renormalised model, defined in Section 6.4. Let $H^{\varepsilon,\mathfrak{m},\delta} \in \mathcal{D}_\varepsilon^{\gamma,\eta}(\widehat{Z}^{\varepsilon,\mathfrak{m},\delta})$ be a local solution of equation (7.29) on an interval $[0, T_{\varepsilon,\mathfrak{m},\delta})$. Let $\tilde{h}^{\varepsilon,\mathfrak{m},\delta} = \widehat{\mathcal{R}}^{\varepsilon,\mathfrak{m},\delta} H^{\varepsilon,\mathfrak{m},\delta} + \tilde{c}_0^{\varepsilon,\delta} t$, where $\widehat{\mathcal{R}}^{\varepsilon,\mathfrak{m},\delta}$ is the respective reconstruction map and the constant $\tilde{c}_0^{\varepsilon,\delta}$ is defined in (7.31). Then there exist $\delta_0 > 0$, $\theta > 0$ and $T^L_{\varepsilon,\mathfrak{m},\delta} \in (0, T_{\varepsilon,\mathfrak{m},\delta})$, such that $\lim_{L\to\infty} T^L_{\varepsilon,\mathfrak{m},\delta} = T_{\varepsilon,\mathfrak{m},\delta}$ almost surely and*

$$\sup_{t\in[0,T\wedge T^L_{\varepsilon,\mathfrak{m}}\wedge T^L_{\varepsilon,\mathfrak{m},\delta}]} \|(\tilde{h}^{\varepsilon,\mathfrak{m}} - \tilde{h}^{\varepsilon,\mathfrak{m},\delta})(t) - (\tilde{c}_0^{\varepsilon} - \tilde{c}_0^{\varepsilon,\delta})t\|^{(\varepsilon)}_{\mathcal{C}^\alpha} \le C\delta^\theta, \tag{7.37}$$

*uniformly over $\delta \in (0, \delta_0)$, provided $\|\tilde{h}_0^{\varepsilon,\mathfrak{m}} - \tilde{h}_0^{\varepsilon,\mathfrak{m},\delta}\|^{(\varepsilon)}_{\mathcal{C}^\alpha} \le \delta^\theta$ and $\|\widehat{Z}^{\varepsilon,\mathfrak{m}}; \widehat{Z}^{\varepsilon,\mathfrak{m},\delta}\|^{(\varepsilon)}_{T+1} \le \delta^\theta$.*

*Proof.* To prove existence of a local solution, we use a purely deterministic argument. For this, we take $T > 0$ and a realisation of the discrete model $\widehat{Z}^{\varepsilon,\mathfrak{m}}$ such that $\|\widehat{Z}^{\varepsilon,\mathfrak{m}}\|^{(\varepsilon)}_{T+1}$ is finite. Proposition 6.11 suggests that this is almost surely the case. Let $\mathscr{M}^{\varepsilon,\mathfrak{m}}_T(H^{\varepsilon,\mathfrak{m}})$ be the right-hand side of equation (7.24) restricted to the time interval $[0, T]$. We are going to solve the fixed point problem

$$H^{\varepsilon,\mathfrak{m}} = \mathscr{M}^{\varepsilon,\mathfrak{m}}_T(H^{\varepsilon,\mathfrak{m}}) \tag{7.38}$$

on the space $\mathcal{D}^{\gamma,\eta}_{\varepsilon,T} = \mathcal{D}^{\gamma,\eta}_{\varepsilon,T}(\widehat{Z}^{\varepsilon,\mathfrak{m}})$. For this, we will prove that $\mathscr{M}^{\varepsilon,\mathfrak{m}}_T$ is a contraction map on $\mathcal{D}^{\gamma,\eta}_{\varepsilon,T}$, uniformly in $\varepsilon$ and for $T > 0$ small enough. More precisely, let us take $H^{\varepsilon,\mathfrak{m}}, \bar{H}^{\varepsilon,\mathfrak{m}} \in \mathcal{D}^{\gamma,\eta}_{\varepsilon,T}$. Then we will prove that there exists $\upsilon > 0$ and a constant $C_\varepsilon > 0$, which depends on $\varepsilon$ and $T$ only through $\|\widehat{Z}^{\varepsilon,\mathfrak{m}}\|^{(\varepsilon)}_{T+1}$, such that

$$\|\mathscr{M}^{\varepsilon,\mathfrak{m}}_T(H^{\varepsilon,\mathfrak{m}})\|^{(\varepsilon)}_{\gamma,\eta;T} \le C_\varepsilon\Big(\|\tilde{h}_0^{\varepsilon,\mathfrak{m}}\|^{(\varepsilon)}_{\mathcal{C}^\alpha} + T^\upsilon(1 + \|H^{\varepsilon,\mathfrak{m}}\|^{(\varepsilon)}_{\gamma,\eta;T})^2\Big), \tag{7.39a}$$

$$\begin{aligned}&\|\mathscr{M}^{\varepsilon,\mathfrak{m}}_T(H^{\varepsilon,\mathfrak{m}}); \mathscr{M}^{\varepsilon,\mathfrak{m}}_T(\bar{H}^{\varepsilon,\mathfrak{m}})\|^{(\varepsilon)}_{\gamma,\eta;T} \\ &\qquad \le C_\varepsilon T^\upsilon(1 + \|H^{\varepsilon,\mathfrak{m}}\|^{(\varepsilon)}_{\gamma,\eta;T} + \|\bar{H}^{\varepsilon,\mathfrak{m}}\|^{(\varepsilon)}_{\gamma,\eta;T})\|H^{\varepsilon,\mathfrak{m}}; \bar{H}^{\varepsilon,\mathfrak{m}}\|^{(\varepsilon)}_{\gamma,\eta;T}.\end{aligned} \tag{7.39b}$$

Then for $T > 0$ small enough, $\mathscr{M}^{\varepsilon,\mathfrak{m}}_T$ is a contraction map on $\mathcal{D}^{\gamma,\eta}_{\varepsilon,T}$. The proportionality constants in the bounds below depend on $\|\widehat{Z}^{\varepsilon,\mathfrak{m}}\|^{(\varepsilon)}_{T+1}$ and we prefer not to repeat it every time.

We are going to prove only the bound (7.39a) and (7.39b) is proved by analogous computations. For this, we fix $H^{\varepsilon,\mathfrak{m}} \in \mathcal{D}^{\gamma,\eta}_{\varepsilon,T}$ and get

$$\begin{aligned}\|\mathscr{M}^{\varepsilon,\mathfrak{m}}_T(H^{\varepsilon,\mathfrak{m}})\|^{(\varepsilon)}_{\gamma,\eta;T} &\le \|G^\varepsilon\tilde{h}_0^\varepsilon\|^{(\varepsilon)}_{\gamma,\eta;T} + \frac{1}{2}\|\widetilde{\mathscr{P}}^\varepsilon_\gamma \mathbf{1}_+ \mathcal{Q}_{\le 0}(\mathscr{D}^-_\varepsilon H^{\varepsilon,\mathfrak{m}})(\mathscr{D}^+_\varepsilon H^{\varepsilon,\mathfrak{m}})\|^{(\varepsilon)}_{\gamma,\eta;T} \\ &\quad + \varrho_\varepsilon\|\widetilde{\mathscr{P}}^\varepsilon_\gamma \mathbf{1}_+(\mathscr{D}_\varepsilon H^{\varepsilon,\mathfrak{m}})\|^{(\varepsilon)}_{\gamma,\eta;T} + \|\widetilde{E}^{\varepsilon,\mathfrak{m}}(H^{\varepsilon,\mathfrak{m}})\|^{(\varepsilon)}_{\gamma,\eta;T} \\ &\quad + \|\mathscr{P}^\varepsilon_\gamma \mathbf{1}_+ \Xi\|^{(\varepsilon)}_{\gamma,\eta;T} + \frac{1}{2}\|\widetilde{\mathscr{P}}^\varepsilon_\gamma \mathbf{1}_+ \widetilde{E}^\varepsilon_0\|^{(\varepsilon)}_{\gamma,\eta;T}.\end{aligned}$$

Similarly to [Hai14, Lem. 7.5] we get $\|G^\varepsilon\tilde{h}_0^\varepsilon\|^{(\varepsilon)}_{\gamma,\eta;T} \lesssim \|\tilde{h}_0^{\varepsilon,\mathfrak{m}}\|^{(\varepsilon)}_{\mathcal{C}^\alpha}$, where we recall that $\eta = \alpha$. The elements in $\widetilde{\mathcal{U}}$ with lowest homogeneities are $\mathbf{1}$ and $\mathcal{I}(\Xi)$, and hence $\mathscr{D}^-_\varepsilon H^{\varepsilon,\mathfrak{m}} \in \mathcal{D}^{\gamma-1,\eta-1}_{\varepsilon,T}$ and $\mathscr{D}^+_\varepsilon H^{\varepsilon,\mathfrak{m}} \in \mathcal{D}^{\gamma-1,\eta-1}_{\varepsilon,T}$ are in the spans of elements with minimal homogeneities $|\partial_-\mathcal{I}(\Xi)| = -\frac{1}{2} - \kappa$ and $|\partial_+\mathcal{I}(\Xi)| = -\frac{1}{2} - \kappa$ respectively. Then, by the results of [EH19, Sec. 5.1] and Lemma 7.1, $(\mathscr{D}^-_\varepsilon H^{\varepsilon,\mathfrak{m}})(\mathscr{D}^+_\varepsilon H^{\varepsilon,\mathfrak{m}}) \in \mathcal{D}^{\bar{\gamma},\bar{\eta}}_{\varepsilon,T}$ with $\bar{\gamma} = \gamma - \frac{3}{2} - \kappa$ and $\bar{\eta} = (\eta - \frac{3}{2} - \kappa) \wedge 2(\eta - 1) = 2(\eta - 1)$, where the latter holds if we take $\kappa \le \frac{1}{2} - \eta$. Moreover, $\|(\mathscr{D}^-_\varepsilon H^{\varepsilon,\mathfrak{m}})(\mathscr{D}^+_\varepsilon H^{\varepsilon,\mathfrak{m}})\|^{(\varepsilon)}_{\bar{\gamma},\bar{\eta};T} \lesssim (\|H^{\varepsilon,\mathfrak{m}}\|^{(\varepsilon)}_{\gamma,\eta;T})^2$ and the modelled distribution $(\mathscr{D}^-_\varepsilon H^{\varepsilon,\mathfrak{m}})(\mathscr{D}^+_\varepsilon H^{\varepsilon,\mathfrak{m}})$ takes values in a sector of minimal homogeneity $-1 - 2\kappa$. Then from [EH19, Thm. 4.26] (see also [Hai14, Thm. 7.1] for the analysis of the projection and appearance of a power of $T$) we get

$$\|\widetilde{\mathscr{P}}^\varepsilon_\gamma \mathbf{1}_+ \mathcal{Q}_{\le 0}(\mathscr{D}^-_\varepsilon H^{\varepsilon,\mathfrak{m}})(\mathscr{D}^+_\varepsilon H^{\varepsilon,\mathfrak{m}})\|^{(\varepsilon)}_{\gamma,\eta;T} \lesssim T^\upsilon(\|H^{\varepsilon,\mathfrak{m}}\|^{(\varepsilon)}_{\gamma,\eta;T})^2,$$

for $0 < \upsilon \leq \frac{\eta}{2}$. Similarly, we use Lemma 7.6 to get

$$\|\widetilde{\mathscr{P}}^\varepsilon_\gamma \mathbf{1}_+ \mathscr{D}_\varepsilon H^{\varepsilon,\mathrm{m}}\|^{(\varepsilon)}_{\gamma,\eta;T} \lesssim \varepsilon^\beta \|H^{\varepsilon,\mathrm{m}}\|^{(\varepsilon)}_{\gamma,\eta;T} \tag{7.40}$$

for some $\beta > \frac{1}{2}$ and

$$\|\mathscr{P}^\varepsilon_\gamma \mathbf{1}_+ \Xi\|^{(\varepsilon)}_{\gamma,\eta;T} \lesssim T^\upsilon, \qquad \|\widetilde{\mathscr{P}}^\varepsilon_\gamma \mathbf{1}_+ \widetilde{E}^\varepsilon_0\|^{(\varepsilon)}_{\gamma,\eta;T} \lesssim T^\upsilon.$$

In the last two bounds we use the fact that $\Xi$ and $\widetilde{E}^\varepsilon_0$ can be interpreted as elements of $\mathcal{D}^{\bar{\gamma},\bar{\eta}}_{\varepsilon,T}$.

It is left to bound $\|\widetilde{E}^{\varepsilon,\mathrm{m}}(H^{\varepsilon,\mathrm{m}})\|^{(\varepsilon)}_{\gamma,\eta;T}$. From (7.16) we get

$$\|\widetilde{E}^{\varepsilon,\mathrm{m}}(H^{\varepsilon,\mathrm{m}})\|^{(\varepsilon)}_{\gamma,\eta;T} \leq \|\widetilde{E}^{\varepsilon,\mathrm{m}}_1(H^{\varepsilon,\mathrm{m}})\|^{(\varepsilon)}_{\gamma,\eta;T} + \|\widetilde{E}^{\varepsilon,\mathrm{m}}_2(H^{\varepsilon,\mathrm{m}})\|^{(\varepsilon)}_{\gamma,\eta;T},$$

where (7.17) yields

$$\|\widetilde{E}^{\varepsilon,\mathrm{m}}_2(H^{\varepsilon,\mathrm{m}})\|^{(\varepsilon)}_{\gamma,\eta;T} \leq \varrho_\varepsilon(1+\varepsilon^{-\frac{1}{2}})\|\widetilde{\mathscr{P}}^\varepsilon_\gamma \mathbf{1}_+ \mathscr{D}_\varepsilon H^{\varepsilon,\mathrm{m}}\|^{(\varepsilon)}_{\gamma,\eta;T} + \frac{1}{2}\|\mathscr{E}_\varepsilon(\mathscr{D}^+_\varepsilon \widetilde{\mathscr{P}}^\varepsilon_\gamma \mathbf{1}_+ (\mathscr{D}^-_\varepsilon H^{\varepsilon,\mathrm{m}}))\|^{(\varepsilon)}_{\gamma,\eta;T}.$$

(7.40) yields a bound of order $\varepsilon^{\beta-\frac{1}{2}}\|H^{\varepsilon,\mathrm{m}}\|^{(\varepsilon)}_{\gamma,\eta;T}$ on the first term and Lemmas 7.1 and 7.2 allow to bound the second term by a constant times $\|\widetilde{\mathscr{P}}^\varepsilon_\gamma \mathbf{1}_+(\mathscr{D}^-_\varepsilon H^{\varepsilon,\mathrm{m}})\|^{(\varepsilon)}_{\gamma,\eta;T}$. We bound it, as in (7.40), by a constant multiple of $\varepsilon^\beta\|H^{\varepsilon,\mathrm{m}}\|^{(\varepsilon)}_{\gamma,\eta;T}$.

Finally, we need to bound (7.18) and (7.19). Let us start with (7.18), where $\varrho_\varepsilon > 0$. We get

$$\begin{aligned}\|\widetilde{E}^{\varepsilon,\mathrm{m}}_1(H^{\varepsilon,\mathrm{m}})\|^{(\varepsilon)}_{\gamma,\eta;T} \lesssim \varepsilon^{-\frac{1}{2}}\|(\mathscr{D}^-_\varepsilon - \mathscr{D}_\varepsilon)\mathscr{P}^\varepsilon_{\gamma+2}\mathbf{1}_+ H^{\varepsilon,\mathrm{m}}\|^{(\varepsilon)}_{\gamma,\eta;T} \\ + \varepsilon^{-\frac{1}{2}}\|\mathscr{D}^-_\varepsilon(\bar{\mathscr{P}}^\varepsilon_{\gamma+2} - \mathscr{P}^\varepsilon_{\gamma+2})\mathbf{1}_+ H^{\varepsilon,\mathrm{m}}\|^{(\varepsilon)}_{\gamma,\eta;T} \\ + \varepsilon^{-\frac{1}{2}}\|\mathscr{D}_\varepsilon \widetilde{\mathscr{P}}^\varepsilon_{\gamma+2}\mathbf{1}_+ H^{\varepsilon,\mathrm{m}}\|^{(\varepsilon)}_{\gamma,\eta;T}.\end{aligned} \tag{7.41}$$

Lemma 7.6 allow to bound the last two terms by a constant times $\varepsilon^\upsilon\|H^{\varepsilon,\mathrm{m}}\|^{(\varepsilon)}_{\gamma,\eta;T}$ for some $\upsilon > 0$. Hence, the last two terms can be made arbitrarily small by taking $\varepsilon$ small. Now, we will bound the first term in (7.41). The model distribution $H^{\varepsilon,\mathrm{m}}$ is function-like and hence the modelled distribution $U^\varepsilon(z) = \mathscr{P}^\varepsilon_{\gamma+2}\mathbf{1}_+ H^{\varepsilon,\mathrm{m}}(z)$ is in the span of $X^\ell$ for $\ell \in \mathbf{N}^2$ such that $|\ell|_\mathfrak{s} < \gamma + 2$. We write it as

$$U^\varepsilon(z) = \sum_{\ell \in \mathbf{N}^2 : |\ell|_\mathfrak{s} \leq \gamma+2} \frac{1}{\ell!} U^\varepsilon_\ell(z) X^\ell.$$

Moreover, as above we get

$$\|U^\varepsilon\|^{(\varepsilon)}_{\gamma+2,\eta+2;T} \lesssim \|H^{\varepsilon,\mathrm{m}}\|^{(\varepsilon)}_{\gamma,\eta;T}. \tag{7.42}$$

Then the definitions (7.1), the action of the maps $\partial_\pm$ on polynomials and the action of the model on polynomials (see Definition 5.3) yield

$$\begin{aligned}(\mathscr{D}^-_\varepsilon U^\varepsilon)(z) &= U^\varepsilon_{(0,1)}(z)\mathbf{1} + U^\varepsilon_{(0,2)}(z)X_1 \\ &\quad - \frac{1}{\varepsilon}\Big(U^\varepsilon_{(0,0)}(z-\varepsilon) - U^\varepsilon_{(0,0)}(z) + \varepsilon U^\varepsilon_{(0,1)}(z) - \frac{\varepsilon^2}{2}U^\varepsilon_{(0,2)}(z)\Big)\mathbf{1} \\ &= \frac{1}{\varepsilon}(U^\varepsilon_{(0,0)}(z) - U^\varepsilon_{(0,0)}(z-\varepsilon))\mathbf{1} + U^\varepsilon_{(0,2)}(z)(\frac{\varepsilon}{2}\mathbf{1} + X_1)\end{aligned}$$

and

$$\begin{aligned}(\mathscr{D}^+_\varepsilon U^\varepsilon)(z) &= U^\varepsilon_{(0,1)}(z)\mathbf{1} + U^\varepsilon_{(0,2)}(z)X_1 \\ &\quad + \frac{1}{\varepsilon}\Big(U^\varepsilon_{(0,0)}(z+\varepsilon) - U^\varepsilon_{(0,0)}(z) - \varepsilon U^\varepsilon_{(0,1)}(z) - \frac{\varepsilon^2}{2}U^\varepsilon_{(0,2)}(z)\Big)\mathbf{1}\end{aligned}$$

$$= \frac{1}{\varepsilon}(U^\varepsilon_{(0,0)}(z+\varepsilon) - U^\varepsilon_{(0,0)}(z))\mathbf{1} + U^\varepsilon_{(0,2)}(z)(-\frac{\varepsilon}{2}\mathbf{1} + X_1).$$

Hence, we have

$$(\mathscr{D}^-_\varepsilon - \mathscr{D}^+_\varepsilon)U^\varepsilon(z) = \frac{1}{\varepsilon}\mathrm{Proj}_{\mathbf{1}}\Big(U^\varepsilon_{(0,0)}(z) - \widehat{\Gamma}^{\varepsilon,\mathfrak{m}}_{z,z-\varepsilon}U^\varepsilon(z-\varepsilon)\Big)\mathbf{1}$$
$$+ \frac{1}{\varepsilon}\mathrm{Proj}_{\mathbf{1}}\Big(U^\varepsilon_{(0,0)}(z) - \widehat{\Gamma}^{\varepsilon,\mathfrak{m}}_{z,z+\varepsilon}U^\varepsilon(z+\varepsilon)\Big)\mathbf{1}.$$

Using (7.42) we can bound this expression as

$$\varepsilon^{-\frac{1}{2}}|\!|\!|(\mathscr{D}^-_\varepsilon - \mathscr{D}^+_\varepsilon)U^\varepsilon|\!|\!|^{(\varepsilon)}_{\gamma,\eta;T} \lesssim \varepsilon^{\gamma+2-\frac{3}{2}}|\!|\!|H^{\varepsilon,\mathfrak{m}}|\!|\!|^{(\varepsilon)}_{\gamma,\eta;T},$$

which is a bound on the first term in (7.41). We note that $\gamma + 2 - \frac{3}{2} > 0$ and this expression can be made arbitrarily small by taking $\varepsilon$ small.

For any fixed $T > 0$ and any $\varepsilon$ sufficiently small the preceding computations give the bound (7.39a), as desired.

Having a local solution of (7.38), the other statements of Proposition 7.7 are the standard results which can be found in [Hai14, Cor. 7.12]. The bound (7.37) is derived from the bound on the two discrete models (6.42). □

## 8 Proof of the main convergence result

We have everything we need to prove our main convergence result, Theorem 1.1, which is the content of this section.

Let $\tilde{h}^\varepsilon$ be the process on $\mathbb{R}_+ \times \Lambda_\varepsilon$ defined in (1.8) which satisfies the (renormalised) semi-discrete equation (2.22), and $h$ be the solution of the KPZ equation (1.11) on $\mathbb{R}_+ \times \mathbb{R}$, provided by Theorem 3.1. As follows from Theorem 3.1, there exists a constant $c_0$ such that $h_{\mathrm{CH}}(t,x) = h(t,x) - c_0 t$ coincides with the Cole-Hopf solution. We fix this constant $c_0$ and use it in what follows.

Let $\mathcal{Y}_\varepsilon(\tilde{h}^\varepsilon)$ be the piece-wise linear extension of $\tilde{h}^\varepsilon$ to $\mathbb{R}_+ \times \mathbb{R}$, i.e. an extension off the grid in the spatial variable. Our goal is to prove that there is a constant $c_0^\varepsilon$, bounded uniformly in $\varepsilon$, such that

$$\lim_{\varepsilon\to 0}\mathbb{E}[F(\mathcal{Y}_\varepsilon(\tilde{h}^\varepsilon) - c_0^\varepsilon\,\bullet)] = \mathbb{E}[F(h - c_0\,\bullet)], \tag{8.1}$$

for any fixed $T > 0$ and any bounded, uniformly continuous function $F : \mathcal{D}([0,T], \mathcal{C}(\mathbb{R})) \to \mathbb{R}$. The function in the argument of $F$ on the left-hand side is $\mathcal{Y}_\varepsilon(\tilde{h}^\varepsilon)(t,x) - c_0^\varepsilon t$. This is exactly the convergence stated in Theorem 1.1. We note that the processes $\tilde{h}^\varepsilon$ and $h$ are not required to be coupled, and the expectations in (8.1) may be on different probability spaces. We fix the value $T > 0$ throughout this section and prove the convergence of the processes on the time interval $[0,T]$.

Following the standard approach (see e.g. [GMW25]), we prove the convergence (8.1) using some intermediate processes. More precisely, we use the process $\tilde{h}^{\varepsilon,\mathfrak{m}}$, defined in Section 2.2 via a constant $\mathfrak{m} \geq 1$. The limit (8.1) follows if we show that for some $\varepsilon_0 > 0$ we have

$$\lim_{\varepsilon\to 0}\mathbb{E}[F(\mathcal{Y}_\varepsilon(\tilde{h}^{\varepsilon,\mathfrak{m}}) - c_0^\varepsilon\,\bullet)] = \mathbb{E}[F(h - c_0\,\bullet)], \tag{8.2a}$$

$$\lim_{\mathfrak{m}\to\infty}\sup_{\varepsilon\in(0,\varepsilon_0)}\mathbb{E}|F(\mathcal{Y}_\varepsilon(\tilde{h}^\varepsilon) - c_0^\varepsilon\,\bullet) - F(\mathcal{Y}_\varepsilon(\tilde{h}^{\varepsilon,\mathfrak{m}}) - c_0^\varepsilon\,\bullet)| = 0, \tag{8.2b}$$

where (8.2a) holds for any fixed $\mathfrak{m} \geq 1$. Note that the two processes in (8.2b) are defined on the same probability space. Then the limit (8.2b) follows from the divergence of the stopping time (2.30) as $\mathfrak{m} \to \infty$. This can be done a similar way how it was done in [GMW25] and we prefer not to repeat the argument here. To prove the limit (8.2a) we need to introduce more intermediate processes. We listed all the auxiliary processes and constants in Table 4 with the respective references.

| Function | Equation | Def. of constant | Limits |
|---|---|---|---|
| $\tilde{h}^\varepsilon(t,\bullet) - c_0^\varepsilon t$ | (1.8) | (8.10) | (8.1), (8.2b) |
| $\tilde{h}^{\varepsilon,\mathtt{m}}(t,\bullet) - c_0^\varepsilon t$ | (2.33) | (8.10) | (8.2a), (8.2b), (8.5c) |
| $\tilde{h}^{\varepsilon,\mathtt{m},\delta}(t,\bullet) - c_0^{\varepsilon,\delta} t$ | (8.32) | (8.10) | (8.5b), (8.5c), (8.15b) |
| $\bar{h}^{\varepsilon,\mathtt{m},\delta}(t,\bullet)$ | (8.14) | — | (8.15a), (8.15b) |
| $h^\delta(t,\bullet) - c_0^\delta t$ | (3.22) | (8.8) | (8.5a), (8.5b), (8.15a) |
| $h(t,\bullet) - c_0 t$ | (1.11) | Theorem 3.1 | (8.1), (8.2a), (8.5a) |

Table 4: Successive approximations used to prove the main convergence result. The last column refers to the limits in which the function is involved.

Let the function $\varrho : \mathbb{R}^2 \to \mathbb{R}$ and its rescaling $\varrho_\delta$ be as defined in (3.19), and let furthermore $\varphi : \mathbb{R} \to \mathbb{R}_+$ be a symmetric smooth function, supported on $[-1, 1]$ and satisfying $\int_\mathbb{R} \varphi(x)dx = 1$. For any $\delta > 0$ we define its rescaling $\varphi_\delta(x) := \delta^{-1}\varphi(\delta^{-1}x)$. Let $h^\delta$ to be the solution of the KPZ equation (3.22) with the smoothened noise $\zeta^\delta = \varrho_\delta \star \zeta$ and the initial condition $h_0^\delta = \varphi_\delta * h_0$. The function $h_0$ is the same as in the statement of Theorem 1.1. It is important to consider a smooth initial condition $h_0^\delta$ (rather than just a Hölder continuous $h_0$) for our following analysis. More precisely, we use the limit of $h_0^\delta$ in (8.35) with respect to a suitable modification of the $\mathcal{C}^1$ norm.

We need to modify the functions $\varrho_\delta$ and $\varphi_\delta$ in a way that their integrals/sums over the discretised domains equal 1. For this, we approximate these functions by their local averages as

$$\varrho_{\delta,\varepsilon}(t,x) := \varepsilon^{-1} \int_{y\in\mathbb{R}: |y-x|\leq\varepsilon/2} \varrho_\delta(t,y)dy, \qquad \varphi_{\delta,\varepsilon}(x) := \varepsilon^{-1} \int_{y\in\mathbb{R}: |y-x|\leq\varepsilon/2} \varphi_\delta(y)dy, \tag{8.3}$$

which satisfy $\int_{D_\varepsilon} \varrho_{\delta,\varepsilon}(z)dz = 1$ and $\varepsilon \sum_{x\in\Lambda_\varepsilon} \varphi_{\delta,\varepsilon}(x) = 1$. Then we define the processes

$$\tilde{\zeta}^{\varepsilon,\mathtt{m},\delta}(t,x) := \varepsilon \sum_{y\in\Lambda_\varepsilon} \int_0^\infty \varrho_{\delta,\varepsilon}(t-s, x-y) d\widetilde{M}^{\varepsilon,\mathtt{m}}(s,y), \tag{8.4}$$

where the integral is with respect to the time variable $s$. We note that these processes are not martingales, but they are smooth in the time variable. On the other hand, a convolution with this process can be interpreted as a stochastic integral. For example, a convolution with the kernel $G^\varepsilon$ may be written as

$$(G^\varepsilon \star_\varepsilon^+ \tilde{\zeta}^{\varepsilon,\mathtt{m},\delta})(t,x) = \varepsilon \sum_{y\in\Lambda_\varepsilon} \int_0^t G_{t-s}^{\varepsilon,\delta}(x-y) d\widetilde{M}^{\varepsilon,\mathtt{m}}(s,y)$$

for the new kernel

$$G^{\varepsilon,\delta} := G^\varepsilon \star_\varepsilon^+ \varrho_{\varepsilon,\delta}.$$

We also define $\tilde{h}^{\varepsilon,\mathtt{m},\delta}$ to be the solution of the semi-discrete equation (7.32) with the initial condition $\tilde{h}_0^{\varepsilon,\delta} := \varphi_{\delta,\varepsilon} *_\varepsilon \tilde{h}_0^\varepsilon$.

Recall furthermore that $h^\delta$ is the solution of the mollified and renormalised KPZ equation (3.22). Using these processes, (8.2a) follows if for some $\delta_0 > 0$ and for some constants $c_0^{\varepsilon,\delta}$ and $c_0^\delta$ we have

$$\lim_{\delta\to 0} \mathbb{E}|F(h^\delta - c_0^\delta \bullet) - F(h - c_0 \bullet)| = 0, \tag{8.5a}$$

$$\lim_{\varepsilon\to 0} \mathbb{E}[F(\mathcal{Y}_\varepsilon(\tilde{h}^{\varepsilon,\mathtt{m},\delta}) - c_0^{\varepsilon,\delta} \bullet)] = \mathbb{E}[F(h^\delta - c_0^\delta \bullet)], \tag{8.5b}$$

$$\lim_{\delta\to 0} \sup_{\varepsilon\in(0,\varepsilon_0)} \mathbb{E}|F(\mathcal{Y}_\varepsilon(\tilde{h}^{\varepsilon,\mathtt{m}}) - c_0^\varepsilon \bullet) - F(\mathcal{Y}_\varepsilon(\tilde{h}^{\varepsilon,\mathtt{m},\delta}) - c_0^{\varepsilon,\delta} \bullet)| = 0, \tag{8.5c}$$

where (8.5b) holds for every fixed $\delta \in (0, \delta_0)$. We used again that the pairs of processes in (8.5a) and (8.5c) are defined on the same probability spaces. We prove these three limits in the following sections.

### 8.1 Proof of (8.5a)

The limit (8.5a) follows from a much stronger convergence stated in Theorem 3.1, whenever

$$\lim_{\delta \to 0} c_0^\delta = c_0. \tag{8.6}$$

We recall from (3.23), (3.20) and (3.21) that the renormalisation constant $C^\delta$ in (3.22) may be written as

$$C^\delta = C^\delta_\star + C^\delta_0, \qquad C^\delta_\star := \nu \int_{\mathbb{R}^2} (\partial_x G^\delta)^2 dz, \tag{8.7}$$

where $G^\delta := G \star \varrho_\delta$ and $C^\delta_0$ has a finite limit $C_0$ as $\delta \to 0$. Then we make the following choice of the constant $c_0^\delta$:

$$c_0^\delta := \frac{1}{2}(C_0^\delta - C_0) + c_0, \tag{8.8}$$

which clearly satisfies (8.6). The choice of this constant will be clear from our proof of the limit (8.5b).

### 8.2 Proof of (8.5c)

Similarly, we can have the limit (8.5c) from a much stronger convergence stated in Proposition 7.7, but for this we need to make a particular choice of the constants $c_0^\varepsilon$ and $c_0^{\varepsilon,\delta}$. For this, we recall the definition of the renormalisation constant (6.41) and use $G^\varepsilon = K^\varepsilon + R^\varepsilon$ to write

$$\widetilde{C}_{\varepsilon,\delta}({}^{+}\!\vee^{-}) = \widetilde{C}^{\varepsilon,\delta}_\star + \widetilde{C}^{\varepsilon,\delta}_0, \qquad \widetilde{C}^{\varepsilon,\delta}_\star := \nu^\varepsilon \int_{D_\varepsilon} \nabla^+_\varepsilon G^{\varepsilon,\delta}(z) \nabla^-_\varepsilon G^{\varepsilon,\delta}(z)\, dz, \tag{8.9}$$

where $\widetilde{C}^{\varepsilon,\delta}_0$ has a finite limit $\widetilde{C}^\varepsilon_0$ as $\delta \to 0$ uniformly in $\varepsilon$. Using the constants (7.26) and (7.31) we then define

$$c_0^{\varepsilon,\delta} := \tilde{c}_0^{\varepsilon,\delta} + \frac{1}{2}(\widetilde{C}_0^{\varepsilon,\delta} - C_0) + c_0, \qquad c_0^\varepsilon := \tilde{c}_0^\varepsilon + \frac{1}{2}(\widetilde{C}_0^\varepsilon - C_0) + c_0, \tag{8.10}$$

which are bounded uniformly in $\varepsilon$. Then we clearly have

$$\lim_{\delta \to 0} (c_0^{\varepsilon,\delta} - \tilde{c}_0^{\varepsilon,\delta}) = c_0^\varepsilon - \tilde{c}_0^\varepsilon, \tag{8.11}$$

and the expression in (8.5c) can be written as

$$\lim_{\delta \to 0} \sup_{\varepsilon \in (0, \varepsilon_0)} \mathbb{E}|F(\mathcal{Y}_\varepsilon(\tilde{h}^{\varepsilon,\mathtt{m}} - \tilde{c}_0^\varepsilon \bullet) - (c_0^\varepsilon - \tilde{c}_0^\varepsilon)\bullet) - F(\mathcal{Y}_\varepsilon(\tilde{h}^{\varepsilon,\mathtt{m},\delta} - \tilde{c}_0^{\varepsilon,\delta}\bullet) - (c_0^{\varepsilon,\delta} - \tilde{c}_0^{\varepsilon,\delta})\bullet)|.$$

This limit vanishes, as follows from Proposition 7.7, and it is left to prove (8.5b).

### 8.3 Proof of (8.5b)

We are going to prove the desired limit (8.5b) on the level of the mild forms of the equations (7.32) and (3.22). Using the constant (8.8), the equation for $h^\delta - c_0^\delta \bullet$ in the mild form is

$$(h^\delta - c_0^\delta \bullet)(t,x) = (G_t * h_0^\delta)(x) - \frac{1}{2} G \star^{+} (\partial_x (h^\delta - c_0^\delta \bullet))^2 (t,x) + \frac{1}{2}(C^\delta_\star + C_0 - 2c_0)t + (G \star^{+} \zeta^\delta)(t,x), \tag{8.12}$$

where $\star^{+}$ is a convolution on $\mathbb{R}_+ \times \mathbb{R}$. Furthermore, we use the constant (8.10) and the equation (7.32) to get

$$\begin{aligned}(\tilde{h}^{\varepsilon,\mathtt{m},\delta} - c_0^{\varepsilon,\delta}\bullet)(t,x) = (G^\varepsilon_t *_\varepsilon \tilde{h}_0^{\varepsilon,\delta})(x) - \frac{1}{2}\widetilde{G}^\varepsilon \star^{+}_\varepsilon (\nabla^-_\varepsilon (\tilde{h}^{\varepsilon,\mathtt{m},\delta} - c_0^{\varepsilon,\delta}\bullet) \diamond \nabla^+_\varepsilon (\tilde{h}^{\varepsilon,\mathtt{m},\delta} - c_0^{\varepsilon,\delta}\bullet))(t,x) \\ + \frac{1}{2}(\widetilde{C}^{\varepsilon,\delta}_\star + C_0 - 2c_0)t + (G^\varepsilon \star^{+}_\varepsilon \tilde{\zeta}^{\varepsilon,\mathtt{m},\delta})(t,x).\end{aligned} \tag{8.13}$$

It will be convenient to introduce one more intermediate process $\bar{h}^{\varepsilon,\mathfrak{m},\delta}$ on $\mathbb{R}_+ \times \mathbb{R}$. We note that the formula (8.4) makes sense for all $(t,x) \in \mathbb{R}_+ \times \mathbb{R}$, which defines the noise $\tilde{\zeta}^{\varepsilon,\mathfrak{m},\delta}(t,x)$ on $\mathbb{R}_+ \times \mathbb{R}$. To avoid ambiguity, let us denote by $\bar{\zeta}^{\varepsilon,\mathfrak{m},\delta}$ such an extension of the noise $\tilde{\zeta}^{\varepsilon,\mathfrak{m},\delta}$ off the grid. Let then $\bar{h}^{\varepsilon,\mathfrak{m},\delta}$ be the solution of the following equation on $\mathbb{R}_+ \times \mathbb{R}$:

$$\bar{h}^{\varepsilon,\mathfrak{m},\delta}(t,x) = (G_t * h_0^\delta)(x) - \frac{1}{2} G \star^+ (\partial_x \bar{h}^{\varepsilon,\mathfrak{m},\delta})^2(t,x) + \frac{1}{2}(C_\star^\delta + C_0 - 2c_0)t + (G \star^+ \bar{\zeta}^{\varepsilon,\mathfrak{m},\delta})(t,x). \tag{8.14}$$

Then the limit (8.5b) follows from

$$\lim_{\varepsilon\to 0} \mathbb{E}[F(\bar{h}^{\varepsilon,\mathfrak{m},\delta})] = \mathbb{E}[F(h^\delta - c_0^\delta \bullet)], \tag{8.15a}$$

$$\lim_{\varepsilon\to 0} \mathbb{E}|F(\mathcal{Y}_\varepsilon(\tilde{h}^{\varepsilon,\mathfrak{m},\delta}) - c_0^{\varepsilon,\delta}\bullet) - F(\bar{h}^{\varepsilon,\mathfrak{m},\delta})| = 0, \tag{8.15b}$$

for every fixed $\delta \in (0, \delta_0)$. We are going to prove these limits in the following two sections.

### 8.3.1 Proof of (8.15a)

For a fixed $\delta > 0$, the processes involved in (8.15a) are almost surely smooth, and the limit can be proved by a standard fixed point argument in a space of sufficiently regular functions. More precisely, for a fixed $t > 0$, let us denote the function spaces

$$L^\infty_{t,x} := L^\infty([0,t] \times \mathbb{R}), \qquad L^\infty_t \mathcal{C}^1_x := L^\infty([0,t]; \mathcal{C}^1(\mathbb{R})).$$

We will show the continuous dependence of the solution of equation (8.14) on the driving noise. To state the problem slightly more generally, for $f_0 \in \mathcal{C}^1(\mathbb{R})$, $T > 0$ and $\zeta \in L^\infty_{T,x}$ we consider the PDE

$$\partial_t f = \frac{1}{2}\Delta f - \frac{1}{2}(\partial_x f)^2 + A + \zeta,$$

with the initial condition $f(0,\bullet) = f_0(\bullet)$ and a constant $A$. We can write it in mild form

$$f(t,x) = (G_t * f_0)(x) - \frac{1}{2} G \star^+ (\partial_x f)^2(t,x) + At + (G \star^+ \zeta)(t,x), \tag{8.16}$$

where $G$ is the heat kernel. Assume there exists $L > 0$ such that $\|f_0\|_{\mathcal{C}^1} \le L$ and $\|\zeta\|_{L^\infty_{T,x}} \le L$. We are going to prove that there is a unique solution $f \in L^\infty_T \mathcal{C}^1_x$, and the solution map $f = \mathcal{S}_T(\zeta, f_0)$ is locally continuous from $L^\infty_{T,x} \times \mathcal{C}^1(\mathbb{R})$ to $L^\infty_T \mathcal{C}^1_x$.

We denote by $\mathcal{M}_t(f)(x)$ the right-hand side of (8.16), and we are going to prove that $f \mapsto \mathcal{M}_t(f)$ is a contraction map on $\mathcal{B}_{L,t} := \{f \in L^\infty_T \mathcal{C}^1_x : \|f\|_{L^\infty_t \mathcal{C}^1_x} \le L+1\}$ for sufficiently small $t \in (0,T]$. Taking $f \in \mathcal{B}_{L,t}$, by (8.16) and $\|G_t\|_{L^1} = 1$ we have

$$\|\mathcal{M}_t(f)\|_{L^\infty} \le \|f_0\|_{L^\infty} + \frac{t}{2}(L+1)^2 + |A|t + Lt. \tag{8.17}$$

Furthermore, taking $\partial_x$ on the right-hand side of (8.16), we get

$$\partial_x \mathcal{M}_t(f) = (G_t * \partial_x f_0)(x) - \frac{1}{2}\partial_x G \star^+ (\partial_x f)^2(t,x) + (\partial_x G \star^+ \zeta)(t,x).$$

For any $t > 0$ we have $\|\partial_x G_t\|_{L^1} = \sqrt{\frac{1}{\pi t}}$, and hence $\int_0^t \|\partial_x G_s\|_{L^1} ds = \sqrt{\frac{4t}{\pi}}$ and

$$\|\partial_x \mathcal{M}_t(f)\|_{L^\infty} \le \|\partial_x f_0\|_{L^\infty} + \sqrt{\frac{t}{\pi}}(L+1)^2 + \sqrt{\frac{4t}{\pi}} L. \tag{8.18}$$

By (8.17)-(8.18), and $\|f_0\|_{\mathcal{C}^1} \le L$ we get $\|\mathcal{M}_t(f)\|_{L^\infty_t \mathcal{C}^1_x} \le L+1$ for $t > 0$ small enough, depending on $L$ and $\delta$. It means that $\mathcal{M}_t$ maps $\mathcal{B}_{L,t}$ to itself.

Let us now take $f, \bar f \in \mathcal{B}_{L,t}$ with the initial states $f_0 = \bar f_0$. Then, we have

$$(\mathcal{M}_t(f) - \mathcal{M}_t(\bar f))(x) = -\frac{1}{2} G \star^+ ((\partial_x f)^2 - (\partial_x \bar f)^2)(t,x),$$
$$(\partial_x \mathcal{M}_t(f) - \partial_x \mathcal{M}_t(\bar f))(x) = -\frac{1}{2} \partial_x G \star^+ ((\partial_x f)^2 - (\partial_x \bar f)^2)(t,x),$$

and hence

$$\|\mathcal{M}_t(f) - \mathcal{M}_t(\bar f)\|_{L_t^\infty \mathcal{C}_x^1} \leq L\Big(t + \sqrt{\frac{4t}{\pi}}\Big) \|f - \bar f\|_{L_t^\infty \mathcal{C}_x^1}.$$

Taking $t > 0$ small enough, we conclude that $\mathcal{M}_t$ is a contraction on $\mathcal{B}_{L,t}$. Then by the Banach fixed point theorem there exists a unique solution $f \in L_t^\infty \mathcal{C}_x^1$ of equation (8.16).

Let us now denote by $f = \mathcal{S}_t(\zeta, f_0)$ the solution map of (8.16) on $\mathcal{B}_{L,t}$. By similar estimates, we conclude that $S_t$ is Lipschitz continuous in both arguments. The extension of the solution to longer time intervals $[0, T]$ is the standard procedure, and is done by patching local solutions.

In order to apply the result which we have just proved to equation (8.14), we need to show that the driving noises converge.

**Lemma 8.1** *The martingales $(\widetilde{M}^{\varepsilon,\mathfrak{m}}(t,x))_{t\geq 0, x\in\Lambda_\varepsilon}$ weakly converge as $\varepsilon \to 0$ in the Skorokhod topology $\mathcal{D}(\mathbb{R}_+, \mathscr{D}'(\mathbb{R}))$ to a 2-periodic cylindrical Wiener process on $L^2(\mathbb{R})$ with variance $\nu$. The constant $\nu$ is defined in* (2.9).

*Proof.* We have verified in Lemma 6.10 that the martingales $\widetilde{M}^{\varepsilon,\mathfrak{m}}$ satisfy [GMW24, Assm. 1]. Hence, by [GMW24, Lem. 2.6], it suffices to show that for every twice continuously differentiable $\varphi : \mathbb{R} \to \mathbb{R}$, whose support radius does not exceed 1, and for every fixed $t > 0$ the following limit holds in distribution

$$\lim_{\varepsilon\to 0} \int_0^t (\iota_\varepsilon \widetilde{\mathbf{C}}_{\varepsilon,\mathfrak{m}})(s,\varphi) ds = t\nu \int_{\mathbb{R}} \varphi(x) dx, \tag{8.19}$$

where $\widetilde{\mathbf{C}}_{\varepsilon,\mathfrak{m}}$ is defined in (2.31) and where we use the map (1.22). Using the constant (2.8) and Lemma B.3, we get

$$\Big|\int_0^t (\iota_\varepsilon \widetilde{\mathbf{C}}_{\varepsilon,\mathfrak{m}})(s,\varphi) ds - t\nu^\varepsilon \varepsilon \sum_{x\in\Lambda_\varepsilon} \varphi(x)\Big| = \Big|\int_0^t \iota_\varepsilon (\widetilde{\mathbf{C}}_{\varepsilon,\mathfrak{m}} - \nu^\varepsilon)(s,\varphi) ds\Big| \leq C\varepsilon^\theta,$$

almost surely for some $\theta > 0$ and some constant $C > 0$ depending on $t$ and $\mathfrak{m}$. Hence, the limit on the left-hand side of (8.19) equals

$$\lim_{\varepsilon\to 0} t\nu^\varepsilon \varepsilon \sum_{x\in\Lambda_\varepsilon} \varphi(x) = t\nu \int_{\mathbb{R}} \varphi(x) dx,$$

where we used $\lim_{\varepsilon\to 0} \nu^\varepsilon = \nu$, as follows from (2.9). □

By Lemma 8.1 we readily conclude that as $\varepsilon \to 0$, the processes $\bar\zeta^{\varepsilon,\mathfrak{m},\delta}$ converge in law in the topology of $L^\infty_{T,x}$ to the smooth Gaussian noise $\zeta^\delta$, driving equation (3.22). Then, using equations (8.12) and (8.14), from the continuity of the solution map $\mathcal{S}_T$ we conclude that $\bar h^{\varepsilon,\mathfrak{m},\delta}$ converges in law to $h^\delta - c_0^\delta \bullet$, as $\varepsilon \to 0$, in the topology of $L_T^\infty \mathcal{C}_x^1$. Combining it with continuity of the solution map $\mathcal{S}_t$, we get

$$\lim_{\varepsilon\to 0} \mathbb{E}[F(\bar h^{\varepsilon,\mathfrak{m},\delta})] = \mathbb{E}[F(h^\delta - c_0^\delta \bullet)],$$

which is (8.15a) as desired.

### 8.3.2 Proof of (8.15b)

We observe that the two processes $\tilde{h}^{\varepsilon,\mathrm{m},\delta}$ and $\bar{h}^{\varepsilon,\mathrm{m},\delta}$ live on the same probability space. To prove their convergence we will use several $\varepsilon$-dependent norms and metrics. Namely, for a discrete $f^\varepsilon : \Lambda_\varepsilon \to \mathbb{R}$ we set

$$\|f^\varepsilon\|^{(\varepsilon)}_{L^\infty} := \sup_{x\in\Lambda_\varepsilon} |f^\varepsilon(x)|, \qquad \|f^\varepsilon\|^{(\varepsilon)}_{\mathcal{C}^1} := \|f^\varepsilon\|^{(\varepsilon)}_{L^\infty} + \|\nabla^+_\varepsilon f^\varepsilon\|^{(\varepsilon)}_{L^\infty}$$

and

$$\|f^\varepsilon\|^{(\varepsilon)}_{L^1} := \varepsilon \sum_{x\in\Lambda_\varepsilon} |f^\varepsilon(x)|.$$

For another function $f : \mathbb{R} \to \mathbb{R}$ we define

$$\|f^\varepsilon; f\|^{(\varepsilon)}_{L^\infty} := \sup_{\substack{x\in\Lambda_\varepsilon, \bar{x}\in\mathbb{R}\\ |\bar{x}-x|\le\varepsilon/2}} |f^\varepsilon(x) - f(\bar{x})|,$$

$$\|f^\varepsilon; f\|^{(\varepsilon)}_{\mathcal{C}^1} := \|f^\varepsilon; f\|^{(\varepsilon)}_{L^\infty} + \|\nabla^+_\varepsilon f^\varepsilon; f'\|^{(\varepsilon)}_{L^\infty} + \|\nabla^-_\varepsilon f^\varepsilon; f'\|^{(\varepsilon)}_{L^\infty}$$

and

$$\|f^\varepsilon; f\|^{(\varepsilon)}_{L^1} := \sum_{x\in\Lambda_\varepsilon} \int_{|\bar{x}-x|\le\varepsilon/2} |f^\varepsilon(x) - f(\bar{x})| d\bar{x}.$$

For time-dependent functions $f^\varepsilon : \mathbb{R}_+ \times \Lambda_\varepsilon \to \mathbb{R}$ and $f : \mathbb{R}_+ \times \mathbb{R} \to \mathbb{R}$ we define respectively

$$\|f^\varepsilon\|_{[0,t]} := \sup_{s\in[0,t]} \|f^\varepsilon_s\|, \qquad \|f\|_{[0,t]} := \sup_{s\in[0,t]} \|f_s\|, \qquad \|f^\varepsilon; f\|_{[0,t]} := \sup_{s\in[0,t]} \|f^\varepsilon_s; f_s\|,$$

where $\|\bullet\|$ and $\|\bullet;\bullet\|$ are any norm and metric.

Using the piece-wise linear extension $\mathcal{Y}_\varepsilon(f^\varepsilon)$, as in (8.1), one can see that the convergence $\lim_{\varepsilon\to0} \|f^\varepsilon; f\|^{(\varepsilon)}_{L^\infty,[0,T]} = 0$ implies the convergence $\lim_{\varepsilon\to0} \sup_{t\in[0,T]} \|\mathcal{Y}_\varepsilon(f^\varepsilon_t) - f_t\|_{L^\infty} = 0$. Then the limit (8.15b) will hold if we show the much stronger result

$$\lim_{\varepsilon\to0} \mathbb{E}\|\tilde{h}^{\varepsilon,\mathrm{m},\delta} - c^{\varepsilon,\delta}_0 \bullet; \bar{h}^{\varepsilon,\mathrm{m},\delta}\|^{(\varepsilon)}_{\mathcal{C}^1,[0,T]} = 0. \tag{8.20}$$

This is what we are going to do in the rest of this section.

In order to proceed, we need to compare the discrete and continuous heat kernels.

**Lemma 8.2** *For any $T > 0$, we have*

$$\lim_{\varepsilon\to0} \sup_{t\in[0,T]} \int_0^t \|G^\varepsilon_s; G_s\|^{(\varepsilon)}_{L^1} ds = 0, \qquad \lim_{\varepsilon\to0} \sup_{t\in[0,T]} \int_0^t \|\nabla^\pm_\varepsilon G^\varepsilon_s; \partial_x G_s\|^{(\varepsilon)}_{L^1} ds = 0. \tag{8.21}$$

*Proof.* For $x \in \Lambda_\varepsilon$ and $\bar{x} \in \mathbb{R}$ such that $|\bar{x} - x| \le \frac{\varepsilon}{2}$, we write

$$G_t(\bar{x}) - G^\varepsilon_t(x) = \big(G_t(\bar{x}) - G_t(x)\big) + \big(G_t(x) - G^\varepsilon_t(x)\big) \tag{8.22}$$

and bound the two differences separately. For the first term in (8.22) we have

$$G_t(\bar{x}) - G_t(x) = \int_0^1 \partial_\eta G_t(\eta\bar{x} + (1-\eta)x) d\eta = (\bar{x} - x) \int_0^1 G'_t(\eta\bar{x} + (1-\eta)x) d\eta.$$

Since $|\bar{x}| + \frac{\varepsilon}{2} \ge |\eta\bar{x} + (1-\eta)x| \ge \max\{|\bar{x}| - \frac{\varepsilon}{2}, 0\} =: (|\bar{x}| - \frac{\varepsilon}{2})_+$ for any $\eta \in [0,1]$, we get

$$|G_t(\bar{x}) - G_t(x)| \le \frac{\varepsilon}{2} \frac{|\bar{x}| + \varepsilon/2}{4\sqrt{\pi t^3}} e^{-\frac{(|\bar{x}|-\varepsilon/2)^2_+}{4t}}.$$

If $x \in \Lambda_\varepsilon \backslash \{0, \pm\varepsilon\}$, then we have $|\bar{x}| \geq \frac{3}{2}\varepsilon$, $|\bar{x}| + \frac{\varepsilon}{2} \leq 2|\bar{x}|$ and $|\bar{x}| - \frac{\varepsilon}{2} \geq \frac{2}{3}|\bar{x}|$, and hence

$$\sum_{x \in \Lambda_\varepsilon \backslash \{0, \pm\varepsilon\}} \int_{|\bar{x}-x| \leq \varepsilon/2} |G_t(\bar{x}) - G_t(x)| d\bar{x} \leq \int_{\mathbb{R}} \frac{\varepsilon}{2} \frac{2|\bar{x}|}{4\sqrt{\pi t^3}} e^{-\frac{(2|\bar{x}|/3)^2}{4t}} d\bar{x} \lesssim \frac{\varepsilon}{\sqrt{t}}.$$

If $x \in \{0, \pm\varepsilon\}$, then we simply bound $|G_t(\bar{x}) - G_t(x)| \leq G_t(\bar{x}) + G_t(x) \leq \frac{2}{\sqrt{4\pi t}}$. These two cases yield

$$\sum_{x \in \Lambda_\varepsilon} \int_{|\bar{x}-x| \leq \varepsilon/2} |G_t(\bar{x}) - G_t(x)| d\bar{x} \lesssim \frac{\varepsilon}{\sqrt{t}}. \tag{8.23}$$

For the second term in (8.22) we have, by (4.8),

$$\begin{aligned} G^\varepsilon_t(x) - G_t(x) &= \int_{-\frac{1}{2\varepsilon}}^{\frac{1}{2\varepsilon}} e^{-f_\varepsilon(\omega)t} e^{2\pi i \omega x} d\omega - \int_{\mathbb{R}} e^{-2\pi^2\omega^2 t} e^{2\pi i \omega x} d\omega \\ &= \int_{-\frac{1}{2\varepsilon}}^{\frac{1}{2\varepsilon}} (e^{-f_\varepsilon(\omega)t} - e^{-2\pi^2\omega^2 t}) e^{2\pi i \omega x} d\omega + \int_{|\omega| \geq \frac{1}{2\varepsilon}} e^{-2\pi^2\omega^2 t} e^{2\pi i \omega x} d\omega. \end{aligned} \tag{8.24}$$

To bound the first term in (8.24), we use the inequality $|e^{-a} - e^{-b}| \leq |a - b| \int_0^1 e^{-\eta a - (1-\eta)b} d\eta$ for any $a, b > 0$, and that for $|\varepsilon w| \leq \frac{1}{2}$, we have $|f_\varepsilon(\omega) - 2\pi^2\omega^2| \leq C\varepsilon^2\omega^4$ for some constant $C$. Then we absolutely bound the first term in (8.24) by

$$\begin{aligned} &\int_{-\frac{1}{2\varepsilon}}^{\frac{1}{2\varepsilon}} |f_\varepsilon(\omega) - 2\pi^2\omega^2| t \int_0^1 e^{-\eta f_\varepsilon(\omega)t - (1-\eta)2\pi^2\omega^2 t} d\eta d\omega \\ &\qquad \lesssim \int_{-\frac{1}{2\varepsilon}}^{\frac{1}{2\varepsilon}} \varepsilon^2\omega^4 t e^{-C\omega^2 t} d\omega \lesssim \frac{\varepsilon^2}{t^{3/2}} \int_{\mathbb{R}} \omega^4 e^{-C\omega^2} d\omega \lesssim \frac{\varepsilon^2}{t^{3/2}}, \end{aligned} \tag{8.25}$$

where we changed the variable $\omega$ to $\omega\sqrt{t}$. We can absolutely bound the second term in (8.24) by

$$\int_{|\omega| \geq \frac{1}{2\varepsilon}} e^{-2\pi^2\omega^2 t} d\omega = \frac{1}{\sqrt{t}} \int_{|\omega| \geq \frac{\sqrt{t}}{2\varepsilon}} e^{-2\pi^2\omega^2} dw \lesssim \frac{1}{\sqrt{t}} e^{-\frac{Ct}{\varepsilon^2}},$$

where in the first identity we changed $\omega$ to $\omega\sqrt{t}$. Combining, we get $|G^\varepsilon_t(x) - G_t(x)| \lesssim \frac{\varepsilon^2}{t^{3/2}} + \frac{1}{\sqrt{t}} e^{-C\varepsilon^{-2}t}$. This bound is not good enough since it is not integrable at $t = 0$. But, by (A.35) with $n = 4$, we also have

$$|G^\varepsilon_t(x) - G_t(x)| \leq |G^\varepsilon_t(x)| + |G_t(x)| \lesssim \frac{1}{\sqrt{4\pi t}} e^{-\frac{x^2}{4t}} + \frac{(\sqrt{t} + \varepsilon)^3}{(\sqrt{t} + |x| + \varepsilon)^4}.$$

Interpolating between these two bounds with weights $\frac{1}{2}$ and $\frac{1}{2}$, and using $\sqrt{a+b} \leq \sqrt{a} + \sqrt{b}$, we get

$$|G^\varepsilon_t(x) - G_t(x)| \lesssim \Big(\frac{\varepsilon}{t^{3/4}} + \frac{1}{t^{1/4}} e^{-C\varepsilon^{-2}t}\Big)\Big(\frac{1}{t^{1/4}} e^{-\frac{x^2}{8t}} + \frac{(\sqrt{t} + \varepsilon)^{3/2}}{(\sqrt{t} + |x| + \varepsilon)^2}\Big).$$

This yields

$$\sum_{x \in \Lambda_\varepsilon} \int_{|\bar{x}-x| \leq \varepsilon/2} |G^\varepsilon_t(x) - G_t(x)| d\bar{x} \lesssim (\sqrt{t} + 1)^{1/2} \Big(\frac{\varepsilon}{t^{3/4}} + \frac{1}{t^{1/4}} e^{-C\varepsilon^{-2}t}\Big), \tag{8.26}$$

where we estimated a Riemann sum by the respective integral. From (8.23) and (8.26) we get $\int_0^t \|G^\varepsilon_s; G_s\|^{(\varepsilon)}_{L^1} ds \lesssim \varepsilon$ with a proportionality constant bounded uniformly in $t \in [0, T]$.

Now we turn to the proof of the second limit in (8.21) using a similar approach. Clearly, it is sufficient to prove the bound only for $\nabla_\varepsilon^+$. On the one hand, we have

$$G_t'(\bar{x}) - G_t'(x) = \int_0^1 \partial_\eta G_t'(\eta\bar{x} + (1-\eta)x)d\eta = (\bar{x} - x)\int_0^1 G_t''(\eta\bar{x} + (1-\eta)x)d\eta.$$

Since $|\bar{x}| + \frac{\varepsilon}{2} \geq |\eta\bar{x} + (1-\eta)x| \geq (|\bar{x}| - \frac{\varepsilon}{2})_+$ for any $\eta \in [0,1]$, we have

$$|G_t'(\bar{x}) - G_t'(x)| \leq \frac{\varepsilon}{2}\frac{(|\bar{x}| + \varepsilon/2)^2}{8\sqrt{\pi t^5}} e^{-\frac{(|\bar{x}|-\varepsilon/2)_+^2}{4t}}.$$

If $x \in \Lambda_\varepsilon \backslash \{0, \pm\varepsilon\}$, then we have $|\bar{x}| \geq \frac{3}{2}\varepsilon$, $|\bar{x}| + \frac{\varepsilon}{2} \leq 2|\bar{x}|$ and $|\bar{x}| - \frac{\varepsilon}{2} \geq \frac{2}{3}|\bar{x}|$, and hence

$$\sum_{x\in\Lambda_\varepsilon\backslash\{0,\pm\varepsilon\}} \int_{|\bar{x}-x|\leq\varepsilon/2} |G_t'(\bar{x}) - G_t'(x)|d\bar{x} \leq \int_{\mathbb{R}} \frac{\varepsilon}{2}\frac{(2|\bar{x}|)^2}{8\sqrt{\pi t^5}} e^{-\frac{(\frac{2}{3}|\bar{x}|)^2}{4t}} d\bar{x} \lesssim \frac{\varepsilon}{t}.$$

If $x \in \{0, \pm\varepsilon\}$, we have $|\bar{x}| \leq \frac{3}{2}\varepsilon$ and we bound $|G_t'(\bar{x}) - G_t'(x)| \leq |G_t'(\bar{x})| + |G_t'(x)| \leq \frac{3\varepsilon}{4\sqrt{\pi t^3}}$. The two cases yield

$$\sum_{x\in\Lambda_\varepsilon} \int_{|\bar{x}-x|\leq\varepsilon/2} |G_t'(\bar{x}) - G_t'(x)|d\bar{x} \lesssim \varepsilon\Big(\frac{1}{t} + \frac{\varepsilon}{t^{3/2}}\Big). \tag{8.27}$$

This bound is not good enough since it is not integrable at $t = 0$. But we also have the bound

$$\sum_{x\in\Lambda_\varepsilon} \int_{|\bar{x}-x|\leq\varepsilon/2} |G_t'(\bar{x}) - G_t'(x)|d\bar{x} \leq \sum_{x\in\Lambda_\varepsilon} \int_{|\bar{x}-x|\leq\varepsilon/2} (|G_t'(\bar{x})| + |G_t'(x)|)d\bar{x} \lesssim \frac{1}{\sqrt{t}}. \tag{8.28}$$

Now we can interpolate between (8.27) and (8.28) with weights $\frac{1}{4}$ and $\frac{3}{4}$, and obtain

$$\sum_{x\in\Lambda_\varepsilon} \int_{|\bar{x}-x|\leq\varepsilon/2} |G_t'(\bar{x}) - G_t'(x)|d\bar{x} \lesssim \varepsilon^{1/4}\Big(\frac{1}{t^{3/8}} + \frac{\varepsilon^{1/4}}{t^{3/4}}\Big). \tag{8.29}$$

Now we turn to $\nabla_\varepsilon^+ G_t^\varepsilon(x) - \partial_x G_t(x)$ for $x \in \Lambda_\varepsilon$. By (4.8) we have

$$\begin{aligned}\nabla_\varepsilon^+ G_t^\varepsilon(x) - \partial_x G_t(x) &= \int_{-\frac{1}{2\varepsilon}}^{\frac{1}{2\varepsilon}} e^{-f_\varepsilon(\omega)t} e^{2\pi i\omega x}(m^+(\varepsilon\omega)i\omega)d\omega - \int_{\mathbb{R}} e^{-2\pi^2\omega^2 t} e^{2\pi i\omega x}(2\pi i\omega)d\omega \\ &= \int_{-\frac{1}{2\varepsilon}}^{\frac{1}{2\varepsilon}} e^{-f_\varepsilon(\omega)t}(m^+(\varepsilon\omega)i\omega - 2\pi i\omega)e^{2\pi i\omega x}d\omega + \int_{-\frac{1}{2\varepsilon}}^{\frac{1}{2\varepsilon}} (e^{-f_\varepsilon(\omega)t} - e^{-2\pi^2\omega^2 t})(2\pi i\omega)e^{2\pi i\omega x}d\omega \\ &\quad + \int_{|\omega|\geq\frac{1}{2\varepsilon}} e^{-2\pi^2\omega^2 t} e^{2\pi i\omega x}(2\pi i\omega)d\omega.\end{aligned} \tag{8.30}$$

To bound the first term in (8.30), we use $|m^+(\varepsilon\omega)i\omega - 2\pi i\omega| \leq C\varepsilon\omega^2$ for some constant $C$ and any $|\varepsilon\omega| \leq \frac{1}{2}$, and get

$$\Big|\int_{-\frac{1}{2\varepsilon}}^{\frac{1}{2\varepsilon}} e^{-f_\varepsilon(\omega)t}(m^+(\varepsilon\omega)i\omega - 2\pi i\omega)e^{2\pi i\omega x}d\omega\Big| \lesssim \int_{\mathbb{R}} \varepsilon\omega^2 e^{-C\omega^2 t}d\omega \lesssim \frac{\varepsilon}{t^{3/2}}.$$

The second term of (8.30) can be bounded similarly to (8.25):

$$\Big|\int_{-\frac{1}{2\varepsilon}}^{\frac{1}{2\varepsilon}} (e^{-f_\varepsilon(\omega)t} - e^{-2\pi^2\omega^2 t})(2\pi i\omega)e^{2\pi i\omega x}d\omega\Big| \lesssim \int_{-\frac{1}{2\varepsilon}}^{\frac{1}{2\varepsilon}} \varepsilon^2|\omega|^5 t e^{-C\omega^2 t}d\omega \lesssim \frac{\varepsilon^2}{t^2}.$$

The third term of (8.30) can be bounded absolutely by $\int_{|\omega|\geq\frac{1}{2\varepsilon}} 2\pi|\omega|e^{-2\pi^2\omega^2 t}d\omega \lesssim \frac{1}{t}e^{-\frac{Ct}{\varepsilon^2}}$. Combining, we get $|\nabla_\varepsilon^+ G_t^\varepsilon(x) - \partial_x G_t(x)| \lesssim \frac{\varepsilon}{t^{3/2}} + \frac{\varepsilon^2}{t^2} + \frac{1}{t}e^{-C\varepsilon^{-2}t}$. But, by (A.35) with $n = 4$, we also have the bound

$$|\nabla_\varepsilon^+ G_t^\varepsilon(x) - \partial_x G_t(x)| \leq |\nabla_\varepsilon^+ G_t^\varepsilon(x)| + |\partial_x G_t(x)| \lesssim \frac{|x|}{4\sqrt{\pi t^3}}e^{-\frac{x^2}{4t}} + \frac{(\sqrt{t}+\varepsilon)^2}{(\sqrt{t}+|x|+\varepsilon)^4}.$$

Interpolating between these two bounds with weights $\frac{1}{4}$ and $\frac{3}{4}$, and using $(a+b)^\beta \leq a^\beta + b^\beta$ for any $a, b > 0$ and $\beta \in (0, 1)$, we get

$$|\nabla_\varepsilon^+ G_t^\varepsilon(x) - \partial_x G_t(x)| \lesssim \Big(\frac{\varepsilon^{1/4}}{t^{3/8}} + \frac{\varepsilon^{1/2}}{t^{1/2}} + \frac{1}{t^{1/4}}e^{-C\varepsilon^{-2}t}\Big)\Big(\frac{|x|^{3/4}}{t^{9/8}}e^{-\frac{3x^2}{16t}} + \frac{(\sqrt{t}+\varepsilon)^{3/2}}{(\sqrt{t}+|x|+\varepsilon)^3}\Big).$$

This yields

$$\sum_{x\in\Lambda_\varepsilon}\int_{|\bar{x}-x|\leq\varepsilon/2} |\nabla_\varepsilon^+ G_t^\varepsilon(x) - \partial_x G_t(x)|d\bar{x} \lesssim \frac{1}{t^{1/4}}\Big(\frac{\varepsilon^{1/4}}{t^{3/8}} + \frac{\varepsilon^{1/2}}{t^{1/2}} + \frac{1}{t^{1/4}}e^{-C\varepsilon^{-2}t}\Big), \qquad (8.31)$$

where we estimated a Riemann sum by the respective integral. We get from (8.29) and (8.31) the bound $\int_0^t \sum_{x\in\Lambda_\varepsilon}\int_{|\bar{x}-x|\leq\varepsilon/2} |\nabla_\varepsilon^+ G_s^\varepsilon(x) - \partial_x G_s(x)|d\bar{x}ds \lesssim \varepsilon^{1/4}$ with a proportionality constant bounded uniformly in $t \in [0, T]$. □

The following lemma will be useful when we compare two equations.

**Lemma 8.3** *For $f^\varepsilon, g^\varepsilon : \Lambda_\varepsilon \to \mathbb{R}$ and $f, g : \mathbb{R} \to \mathbb{R}$, we have*

$$\|f^\varepsilon *_\varepsilon g^\varepsilon; f * g\|_{L^\infty}^{(\varepsilon)} \leq \|f^\varepsilon\|_{L^1(\Lambda_\varepsilon)}\|g^\varepsilon; g\|_{L^\infty}^{(\varepsilon)} + \|g\|_{L^\infty(\mathbb{R})}\|f^\varepsilon; f\|_{L^1}^{(\varepsilon)}.$$

*If moreover functions depend on the time variable, i.e. $f^\varepsilon, g^\varepsilon : \mathbb{R}_+ \times \Lambda_\varepsilon \to \mathbb{R}$ and $f, g : \mathbb{R}_+ \times \mathbb{R} \to \mathbb{R}$, then for any $t > 0$ we have*

$$\|f^\varepsilon \star_\varepsilon^+ g^\varepsilon; f \star^+ g\|_{L^\infty,[0,t]}^{(\varepsilon)} \leq \|g^\varepsilon; g\|_{L^\infty,[0,t]}^{(\varepsilon)} \int_0^t \|f_s^\varepsilon\|_{L^1(\Lambda_\varepsilon)}ds + \|g\|_{L^\infty_{t,x}} \int_0^t \|f_s^\varepsilon; f_s\|_{L^1}^{(\varepsilon)}ds.$$

*Proof.* The second claim follows from the first claim easily, so let us only prove the first one. Let $x \in \Lambda_\varepsilon$ and $\bar{x} \in \mathbb{R}$ be such that $|\bar{x} - x| \leq \frac{\varepsilon}{2}$. We have

$$(f^\varepsilon *_\varepsilon g^\varepsilon)(x) - (f * g)(\bar{x}) = \varepsilon \sum_{y\in\Lambda_\varepsilon} f^\varepsilon(x-y)g^\varepsilon(y) - \int_{\mathbb{R}} f(\bar{x}-\bar{y})g(\bar{y})d\bar{y}.$$

We add and subtract $\sum_{y\in\Lambda_\varepsilon}\int_{|\bar{y}-y|\leq\varepsilon/2} f^\varepsilon(x-y)g(\bar{y})d\bar{y}$ to get

$$\begin{aligned}&\Big|\varepsilon \sum_{y\in\Lambda_\varepsilon} f^\varepsilon(x-y)g^\varepsilon(y) - \int_{|\bar{y}-y|\leq\varepsilon/2} f^\varepsilon(x-y) \sum_{y\in\Lambda_\varepsilon} g(\bar{y})d\bar{y}\Big| \\ &\qquad = \Big|\sum_{y\in\Lambda_\varepsilon} f^\varepsilon(x-y)\int_{|\bar{y}-y|\leq\varepsilon/2} (g^\varepsilon(y) - g(\bar{y}))d\bar{y}\Big| \leq \|f^\varepsilon\|_{L^1(\Lambda_\varepsilon)}\|g^\varepsilon; g\|_{L^\infty}^{(\varepsilon)},\end{aligned}$$

and

$$\Big|\sum_{y\in\Lambda_\varepsilon}\int_{|\bar{y}-y|\leq\varepsilon/2} f^\varepsilon(x-y)g(\bar{y})d\bar{y} - \int_{\mathbb{R}} f(\bar{x}-\bar{y})g(\bar{y})d\bar{y}\Big|$$

$$= \left| \sum_{y \in \Lambda_\varepsilon} \int_{|\bar{y}-y| \leq \varepsilon/2} (f^\varepsilon(x-y) - f(\bar{x}-\bar{y}))g(\bar{y})d\bar{y} \right| \leq \|g\|_{L^\infty(\mathbb{R})} \|f^\varepsilon; f\|^{(\varepsilon)}_{L^1},$$

where in the last integral we necessarily have $|(x-y)-(\bar{x}-\bar{y})| \leq |x-\bar{x}| + |y-\bar{y}| \leq \frac{\varepsilon}{2} + \frac{\varepsilon}{2} = \varepsilon$. We thus obtained the required bound on $|(f^\varepsilon *_\varepsilon g^\varepsilon)(x) - (f * g)(\bar{x})|$, which after taking a supremum over $x \in \Lambda_\varepsilon$ and $\bar{x} \in \mathbb{R}$ such that $|\bar{x}-x| \leq \frac{\varepsilon}{2}$ yields our claim. □

Now, we will prove the limit (8.20). We will consider the difference of the solutions to (8.13) and (8.14), as well as their spatial derivatives. Recalling the definition of the renormalised product (7.33), we can write (8.13) as

$$\begin{aligned}(\tilde{h}^{\varepsilon,\mathtt{m},\delta} - c_0^{\varepsilon,\delta}\bullet)(t,x) &= (G_t^\varepsilon *_\varepsilon \tilde{h}_0^{\varepsilon,\delta})(x) - \frac{1}{2}\widetilde{G}^\varepsilon \star_\varepsilon^+ (\nabla_\varepsilon^-(\tilde{h}^{\varepsilon,\mathtt{m},\delta} - c_0^{\varepsilon,\delta}\bullet)\nabla_\varepsilon^+(\tilde{h}^{\varepsilon,\mathtt{m},\delta} - c_0^{\varepsilon,\delta}\bullet))(t,x) \\ &+ \frac{1}{2}(\widetilde{C}_\star^{\varepsilon,\delta} + C_0 - 2c_0)t + (G^\varepsilon \star_\varepsilon^+ \zeta^{\varepsilon,\mathtt{m},\delta})(t,x) + \frac{1}{2}(\widetilde{G}^\varepsilon \star_\varepsilon^+ \widetilde{C}_{\varepsilon,\mathtt{m},\delta})(t,x). \qquad (8.32)\end{aligned}$$

We need to show that each term on the right-hand side of (8.32) either vanishes or converges as $\varepsilon \to 0$ to the respective term on the right-hand side of (8.14). For this we are going to use the norms and metrics defined at the beginning of this section.

Using Lemma 8.3 and the estimate $\|G_t^\varepsilon\|_{L^1(\Lambda_\varepsilon)} = 1$, we get

$$\|G_t^\varepsilon *_\varepsilon \tilde{h}_0^{\varepsilon,\delta}; G_t * h_0^\delta\|^{(\varepsilon)}_{\mathcal{C}^1} \leq \|\tilde{h}_0^{\varepsilon,\delta}; h_0^\delta\|^{(\varepsilon)}_{\mathcal{C}^1} + \|G_t^\varepsilon; G_t\|^{(\varepsilon)}_{L^1} \|h_0^\delta\|_{\mathcal{C}^1},$$

where the right-hand sides vanish as $\varepsilon \to 0$ due to Lemma 8.2 and our assumption on the initial conditions in Theorem 1.1. To bound the next term in (8.32), we write

$$\begin{aligned}&\widetilde{G}^\varepsilon \star_\varepsilon^+ (\nabla_\varepsilon^-(\tilde{h}^{\varepsilon,\mathtt{m},\delta} - c_0^{\varepsilon,\delta}\bullet)\nabla_\varepsilon^+(\tilde{h}^{\varepsilon,\mathtt{m},\delta} - c_0^{\varepsilon,\delta}\bullet)) \\ &\qquad = G^\varepsilon \star_\varepsilon^+ (\nabla_\varepsilon^-(\tilde{h}^{\varepsilon,\mathtt{m},\delta} - c_0^{\varepsilon,\delta}\bullet)\nabla_\varepsilon^+(\tilde{h}^{\varepsilon,\mathtt{m},\delta} - c_0^{\varepsilon,\delta}\bullet)) \qquad (8.33) \\ &\qquad\quad + (\widetilde{G}^\varepsilon - G) \star_\varepsilon^+ (\nabla_\varepsilon^-(\tilde{h}^{\varepsilon,\mathtt{m},\delta} - c_0^{\varepsilon,\delta}\bullet)\nabla_\varepsilon^+(\tilde{h}^{\varepsilon,\mathtt{m},\delta} - c_0^{\varepsilon,\delta}\bullet)).\end{aligned}$$

We can compare the first term on the right-hand side of (8.33) with the second term in (8.14). Using Lemma 8.3 we get

$$\begin{aligned}&\|G^\varepsilon \star_\varepsilon^+ (\nabla_\varepsilon^-(\tilde{h}^{\varepsilon,\mathtt{m},\delta} - c_0^{\varepsilon,\delta}\bullet)\nabla_\varepsilon^+(\tilde{h}^{\varepsilon,\mathtt{m},\delta} - c_0^{\varepsilon,\delta}\bullet)); G \star (\partial_x \bar{h}^{\varepsilon,\mathtt{m},\delta})^2\|^{(\varepsilon)}_{\mathcal{C}^1,[0,t]} \\ &\lesssim (t+\sqrt{t})\|\tilde{h}^{\varepsilon,\mathtt{m},\delta} - c_0^{\varepsilon,\delta}\bullet; \bar{h}^{\varepsilon,\mathtt{m},\delta}\|^{(\varepsilon)}_{\mathcal{C}^1,[0,t]}\|\tilde{h}^{\varepsilon,\mathtt{m},\delta} - c_0^{\varepsilon,\delta}\bullet\|^{(\varepsilon)}_{\mathcal{C}^1,[0,t]} \\ &\quad + (t+\sqrt{t})\|\tilde{h}^{\varepsilon,\mathtt{m},\delta} - c_0^{\varepsilon,\delta}\bullet; \bar{h}^{\varepsilon,\mathtt{m},\delta}\|^{(\varepsilon)}_{\mathcal{C}^1,[0,t]}\|\bar{h}^{\varepsilon,\mathtt{m},\delta}\|_{\mathcal{C}^1,[0,t]} + \|\bar{h}^{\varepsilon,\mathtt{m},\delta}\|^2_{\mathcal{C}^1,[0,t]} \int_0^t \|G_s^\varepsilon; G_s\|^{(\varepsilon)}_{\mathcal{C}^1} ds,\end{aligned}$$

where we made use of $\|G_t^\varepsilon\|_{L^1(\Lambda_\varepsilon)} = 1$ and $\int_0^t \|\nabla_\varepsilon^\pm G_s^\varepsilon\|_{L^1(\Lambda_\varepsilon)} ds \lesssim \sqrt{t}$, which can be proved using (A.36). The last term vanishes as $\varepsilon \to 0$ due to Lemma 8.2. Similarly, for the second term in (8.33) we have

$$\begin{aligned}&\|(\widetilde{G}^\varepsilon - G^\varepsilon) \star_\varepsilon^+ (\nabla_\varepsilon^-(\tilde{h}^{\varepsilon,\mathtt{m},\delta} - c_0^{\varepsilon,\delta}\bullet)\nabla_\varepsilon^+(\tilde{h}^{\varepsilon,\mathtt{m},\delta} - c_0^{\varepsilon,\delta}\bullet))\|^{(\varepsilon)}_{\mathcal{C}^1,[0,t]} \\ &\qquad\qquad \leq (\|\tilde{h}^{\varepsilon,\mathtt{m},\delta} - c_0^{\varepsilon,\delta}\bullet\|^{(\varepsilon)}_{\mathcal{C}^1,[0,t]})^2 \int_0^t \|\widetilde{G}_s^\varepsilon - G_s^\varepsilon\|^{(\varepsilon)}_{\mathcal{C}^1} ds.\end{aligned}$$

By Lemma 4.10 and (A.36) we have for any $\vartheta \in (0,1)$

$$\int_0^t \|\widetilde{G}_s^\varepsilon - G_s^\varepsilon\|^{(\varepsilon)}_{\mathcal{C}^1} ds \lesssim \varepsilon^\vartheta \int_0^t (s^{-\frac{\vartheta}{2}} + s^{-\frac{1+\vartheta}{2}}) ds \lesssim \varepsilon^\vartheta (t^{\frac{2-\vartheta}{2}} + t^{\frac{1-\vartheta}{2}}),$$

which vanishes as $\varepsilon \to 0$.

We conclude from the definitions of the constants (8.7) and (8.9) and Vitali's convergence theorem that

$$C_\star^\delta = \lim_{\varepsilon\to 0} \widetilde{C}_\star^{\varepsilon,\delta}. \tag{8.34}$$

Hence, the term $\frac{1}{2}(\widetilde{C}_\star^{\varepsilon,\delta} + C_0 - 2c_0)t$ in (8.32) converges to the term $\frac{1}{2}(C_\star^\delta + C_0 - 2c_0)t$ in (8.14) as $\varepsilon \to 0$.

Let us now compare the next terms in (8.32) and (8.14). Using Lemma 8.3 we get

$$\|G^\varepsilon \star_\varepsilon^+ \tilde{\zeta}^{\varepsilon,\mathfrak{m},\delta}; G\star\bar{\zeta}^{\varepsilon,\mathfrak{m},\delta}\|^{(\varepsilon)}_{\mathcal{C}^1,[0,t]} \lesssim (t+\sqrt{t})\|\tilde{\zeta}^{\varepsilon,\mathfrak{m},\delta}; \bar{\zeta}^{\varepsilon,\mathfrak{m},\delta}\|^{(\varepsilon)}_{\mathcal{C}^1,[0,t]} + \|\bar{\zeta}^{\varepsilon,\mathfrak{m},\delta}\|_{L^\infty_{t,x}} \int_0^t \|G^\varepsilon_s; G_s\|^{(\varepsilon)}_{\mathcal{C}^1} ds,$$

where the last term vanishes as $\varepsilon \to 0$ due to Lemma 8.2.

Finally, we need to study the limit of the last term in (8.32). Definition (6.40) allows to write this term as $\frac{1}{2}(\widetilde{G}^\varepsilon \star_\varepsilon^+ Q^{\varepsilon,\delta} \star_\varepsilon^+ (\widetilde{\mathbf{C}}_{\varepsilon,\mathfrak{m}} - \nu^\varepsilon))$. Hence, we get

$$\|\widetilde{G}^\varepsilon \star_\varepsilon^+ \widetilde{C}_{\varepsilon,\mathfrak{m},\delta}\|^{(\varepsilon)}_{\mathcal{C}^1,[0,t]} \lesssim (t+\sqrt{t})\|Q^{\varepsilon,\delta} \star_\varepsilon^+ (\widetilde{\mathbf{C}}_{\varepsilon,\mathfrak{m}} - \nu^\varepsilon)\|^{(\varepsilon)}_{L^\infty,[0,t]}.$$

We note that $Q^{\varepsilon,\delta}$ is smooth and by Lemma B.3 the expectation of this expression is bounded by a multiple of $\varepsilon^\vartheta$ for some $\vartheta > 0$.

Combining the previous estimates we get

$$\begin{aligned}
&\mathbb{E}\|\tilde{h}^{\varepsilon,\mathfrak{m},\delta} - c_0^{\varepsilon,\delta}\bullet; \bar{h}^{\varepsilon,\mathfrak{m},\delta}\|^{(\varepsilon)}_{\mathcal{C}^1,[0,T]} \\
&\qquad \le \|\tilde{h}_0^{\varepsilon,\delta}; h_0^\delta\|^{(\varepsilon)}_{\mathcal{C}^1} + C(T+\sqrt{T})\mathbb{E}\|\tilde{h}^{\varepsilon,\mathfrak{m},\delta} - c_0^{\varepsilon,\delta}\bullet; \bar{h}^{\varepsilon,\mathfrak{m},\delta}\|^{(\varepsilon)}_{\mathcal{C}^1,[0,T]} \\
&\qquad\quad + C(T+\sqrt{T})\mathbb{E}\|\tilde{\zeta}^{\varepsilon,\mathfrak{m},\delta}; \bar{\zeta}^{\varepsilon,\mathfrak{m},\delta}\|^{(\varepsilon)}_{L^\infty,[0,T]} + C\varepsilon^\vartheta(1 + T^{\frac{2-\vartheta}{2}} + T^{\frac{1-\vartheta}{2}}),
\end{aligned}$$

where the constant $C$ is independent of $\varepsilon$ and $T$. For $T > 0$ small enough, we can absorb the term proportional to $\mathbb{E}\|\tilde{h}^{\varepsilon,\mathfrak{m},\delta} - c_0^{\varepsilon,\delta}\bullet; \bar{h}^{\varepsilon,\mathfrak{m},\delta}\|^{(\varepsilon)}_{\mathcal{C}^1,[0,T]}$ to the left-hand side. From our assumptions in Theorem 1.1 on the initial states we conclude that

$$\lim_{\varepsilon\to 0} \|\tilde{h}_0^{\varepsilon,\delta}; h_0^\delta\|^{(\varepsilon)}_{\mathcal{C}^1} = 0. \tag{8.35}$$

Furthermore, by (8.3), (8.4) and Burkholder-Davis-Gundy inequality [GMW24, Prop. 2.1] we have $\lim_{\varepsilon\to 0} \mathbb{E}\|\tilde{\zeta}^{\varepsilon,\mathfrak{m},\delta}; \bar{\zeta}^{\varepsilon,\mathfrak{m},\delta}\|^{(\varepsilon)}_{L^\infty,[0,t]} = 0$. Then we obtain

$$\lim_{\varepsilon\to 0} \mathbb{E}\|\tilde{h}^{\varepsilon,\mathfrak{m},\delta} - c_0^{\varepsilon,\delta}\bullet; \bar{h}^{\varepsilon,\mathfrak{m},\delta}\|^{(\varepsilon)}_{\mathcal{C}^1,[0,T]} = 0.$$

This argument can be iterated to prove (8.20) on longer time intervals $[0, T]$ and we get the desired limit (8.20).

## Appendix A Extension of discrete kernels off the lattice

Our goal is to extend the function $G^\varepsilon$, solving (4.7), smoothly off the lattice. A natural extension is by trigonometric polynomials in the spatial variable. While it is smoothly extended off the lattice, this extension has a slow decay as $|x| \to \infty$. Since (4.7) is a discretized heat equation, one would like to define an extension of the solution off the lattice which has a very fast decay at infinity, similar to the super-exponential decay of the heat kernel. We construct such an extension in this section for a more general discrete kernel and then apply the result to $G^\varepsilon$. Due to the scaling property (4.10), it will be convenient to work with a kernel on $\mathbb{R}_+ \times \mathbb{Z}$.

### A.1 Extension procedure

For a function $u : \mathbb{R}_+ \times \mathbb{Z} \to \mathbb{R}$, we would like to define an extension $U : \mathbb{R}_+ \times \mathbb{R} \to \mathbb{R}$ such that $U(t,k) = u(t,k)$ for all $k \in \mathbb{Z}$ and $U$ inherits the decay of $u$. Throughout this section we will use $k$ for an integer variable and $x$ for a continuous one. To define such an extension we follow the idea of [HM18, Sec. 5.1].

We lift the function $u$ to a distribution on $\mathbb{R}$ by

$$(\iota u)(t,x) = \sum_{k\in\mathbb{Z}} u(t,k)\delta(x-k),$$

where $\delta$ is the Dirac delta. Then the spatial Fourier transform of this lift is a 1-periodic function

$$\mathcal{F}(\iota u)(t,\omega) = \sum_{k\in\mathbb{Z}} u(t,k)e^{-2\pi i k\omega}, \qquad \omega \in \mathbb{R}.$$

We can recover $u$ by applying the inverse spatial Fourier transform:

$$u(t,k) = \int_{-\frac{1}{2}}^{\frac{1}{2}} \mathcal{F}(\iota u)(t,\omega)e^{2\pi i\omega k}d\omega, \qquad k \in \mathbb{Z}.$$

Let $\psi : \mathbb{R} \to \mathbb{R}_+$ be a smooth, even function, compactly supported on $[-\delta,\delta]$ for some $0 < \delta < \frac{1}{4}$ and such that $\int_{\mathbb{R}} \psi(x)dx = 1$. We define $\varphi : \mathbb{R} \to \mathbb{R}$ via its Fourier transform

$$(\mathcal{F}\varphi)(\omega) := \Big(\mathbf{1}_{[-\frac{1}{2},\frac{1}{2}]} * \psi\Big)(\omega), \qquad \omega \in \mathbb{R}. \tag{A.1}$$

Note that $\mathcal{F}^{-1}(\mathbf{1}_{[-\frac{1}{2},\frac{1}{2}]})(x) = \frac{\sin(\pi x)}{\pi x} =: D(x)$ is the Dirichlet kernel, which is defined at the origin by L'Hopital's rule as $D(0) = 1$ and satisfies $D(k) = 0$ for $k \in \mathbb{Z}\backslash\{0\}$ and $|D(x)| \lesssim |x|^{-1}$ for large $x$. On the other hand, $\mathcal{F}^{-1}\psi(0) = 1$ and $\mathcal{F}^{-1}\psi$ is Schwartz as follows from the properties of $\psi$. As a result, $\varphi$ is Schwartz and satisfies $\varphi(0) = 1$ and $\varphi(k) = 0$ for $k \in \mathbb{Z}\backslash\{0\}$.

We define the extension $U(t,x)$ of $u(t,k)$ to be

$$U(t,x) := \int_{-\infty}^{\infty} \mathcal{F}(\iota u)(t,\omega)\mathcal{F}\varphi(\omega)e^{2\pi i\omega x}d\omega = \sum_{k\in\mathbb{Z}} \varphi(x-k)u(t,k) \tag{A.2}$$

for $(t,x) \in \mathbb{R}_+ \times \mathbb{R}$. This choice clearly satisfies $U(t,k) = u(t,k)$ for all $k \in \mathbb{Z}$. It will be convenient to use the notation

$$U = \mathsf{Ext}(u).$$

We can also show that decay of $u$ transfers to decay of the same order for the extension $U$.

**Proposition A.1** *For any $j = (j_0, j_1) \in \mathbf{N}^2$ and some $\beta \geq 0$, suppose that there is constant $C = C(j,\beta) > 0$ such that*

$$|\partial_t^{j_0}(\nabla^-)^{j_1}u(t,x)| \leq C(\|(t,x)\|_{\mathfrak{s}} + 1)^{-|j|_{\mathfrak{s}}-\beta}, \tag{A.3}$$

*for any $(t,x) \in \mathbb{R}_+ \times \mathbb{Z}$, where the discrete operator $\nabla^-$ acts on the variable $x$. Then there is a constant $C' = C'(j,\beta) > 0$ such that*

$$|D^jU(z)| \leq C'(\|z\|_{\mathfrak{s}} + 1)^{-|j|_{\mathfrak{s}}-\beta} \tag{A.4}$$

*for any $z \in \mathbb{R}_+ \times \mathbb{R}$*

We are going to use this proposition for rescaled functions $u^\varepsilon : \mathbb{R}_+ \times \Lambda_\varepsilon \to \mathbb{R}$ defined as

$$u^\varepsilon(t,x) := \varepsilon^{-1} u(\varepsilon^{-2} t, \varepsilon^{-1} x).$$

Then the respective extension off the lattice

$$U^\varepsilon := \mathsf{Ext}(u^\varepsilon) \tag{A.5}$$

satisfies the same scaling identity. Let us denote the operator $D^j_\varepsilon := \partial_t^{j_0} (\nabla_\varepsilon^-)^{j_1}$ for $j = (j_0, j_1) \in \mathbf{N}^2$. Then we immediately get the following result.

**Corollary A.2** *If for any $j \in \mathbf{N}^2$ and some $\beta \geq 0$ there exists a constant $C = C(j,\beta) > 0$ such that*

$$|D^j_\varepsilon u^\varepsilon(z)| \leq C(\|z\|_\mathfrak{s} + \varepsilon)^{-|j|_\mathfrak{s} - \beta}, \tag{A.6}$$

*for any $z \in \mathbb{R}_+ \times \Lambda_\varepsilon$, then there is $C' = C'(j,\beta) > 0$ such that*

$$|D^j U^\varepsilon(z)| \leq C'(\|z\|_\mathfrak{s} + \varepsilon)^{-|j|_\mathfrak{s} - \beta} \tag{A.7}$$

*for any $z \in \mathbb{R}_+ \times \mathbb{R}$.*

We will also need the following result, which allows to compare two extensions.

**Corollary A.3** *Let $U^\varepsilon$ and $\widetilde{U}^\varepsilon$ be two extensions defined by* (A.5) *via two respective functions $u^\varepsilon$ and $\tilde{u}^\varepsilon$. If for any $j \in \mathbf{N}^2$ and some $\alpha, \beta \geq 0$ there exists a constant $C = C(j,\alpha,\beta) > 0$ such that*

$$|D^j_\varepsilon (\tilde{u}^\varepsilon - u^\varepsilon)(z)| \leq C\varepsilon^\alpha (\|z\|_\mathfrak{s} + \varepsilon)^{-|j|_\mathfrak{s} - \beta}, \tag{A.8}$$

*for any $z \in \mathbb{R}_+ \times \Lambda_\varepsilon$, then there is $C' = C'(j,\alpha,\beta) > 0$ such that*

$$|D^j (\widetilde{U}^\varepsilon - U^\varepsilon)(z)| \leq C'\varepsilon^\alpha (\|z\|_\mathfrak{s} + \varepsilon)^{-|j|_\mathfrak{s} - \beta} \tag{A.9}$$

*for any $z \in \mathbb{R}_+ \times \mathbb{R}$.*

*Proof.* The result follows from Corollary A.2, applied to the function $\varepsilon^{-\alpha}(\tilde{u}^\varepsilon - u^\varepsilon)(z)$. The definition (A.5) of the extension clearly satisfies $\varepsilon^{-\alpha}(\widetilde{U}^\varepsilon - U^\varepsilon) = \mathsf{Ext}(\varepsilon^{-\alpha}(\tilde{u}^\varepsilon - u^\varepsilon))$ by linearity. □

Before proving this result we need to prove in the following three lemmas several properties of the function $\varphi$. The first lemma shows that the extension (A.2) is given by a weighted average of $u(k)$, $k \in \mathbb{Z}$, where the weights sum to 1 and are allowed to be negative.

**Lemma A.4** *For any $x \in \mathbb{R}$ we have $\sum_{k\in\mathbb{Z}} \varphi(x+k) = 1$.*

*Proof.* By the Poisson summation formula

$$\sum_{k\in\mathbb{Z}} \varphi(x+k) = \sum_{\ell\in\mathbb{Z}} (\mathcal{F}\varphi)(\ell) e^{2\pi i \ell x}.$$

Since $\mathcal{F}\varphi = \mathbf{1}_{[-\frac{1}{2},\frac{1}{2}]} * \psi$ and the support of $\psi$ is in $[-\frac{1}{4},\frac{1}{4}]$, we conclude that the support of $\mathcal{F}\varphi$ is within $[-\frac{3}{4},\frac{3}{4}]$. Thus, $(\mathcal{F}\varphi)(\ell) = 0$ for all $\ell \in \mathbb{Z}\backslash\{0\}$. For $\ell = 0$ we have $(\mathcal{F}\varphi)(0) = \int_{-1/2}^{1/2} \psi(y)dy = 1$. We thus conclude that $\sum_{k\in\mathbb{Z}} \varphi(x+k) = 1$. □

**Lemma A.5** *For any $x \in \mathbb{R}$ and $p, q \in \mathbf{N}$ such that $p \neq q$ we have*

$$\sum_{\ell \in \mathbb{Z}} (x-\ell)^p \varphi^{(q)}(x-\ell) = 0. \tag{A.10}$$

*Furthermore, for any $\alpha \in \mathbb{R}$ and $k \in \mathbf{N}$, satisfying $\alpha > p+1$ and $k < x$, there exists a constant $C = C(\alpha, p, q) > 0$ such that*

$$\Big| \sum_{\ell \geq k} (x-\ell)^p \varphi^{(q)}(x-\ell) \Big| \leq C(x-k)^{p+1-\alpha}. \tag{A.11}$$

*Proof.* We note that the left-hand side of (A.10) is a 1-periodisation of the function $g(x) = x^p \varphi^{(q)}(x)$, whose Fourier transform is

$$(\mathcal{F}g)(\omega) = \frac{1}{(-2\pi i)^p} \frac{d^p}{d\omega^p} ((-2\pi i \omega)^q (\mathcal{F}\varphi)(\omega)).$$

Then by the Poisson summation formula we have

$$\begin{aligned} \sum_{\ell \in \mathbb{Z}} (x-\ell)^p \varphi^{(q)}(x-\ell) &= \sum_{\omega \in \mathbb{Z}} (\mathcal{F}g)(\omega) e^{2\pi i \omega x} = \sum_{\omega \in \mathbb{Z}} \frac{1}{(-2\pi i)^p} \frac{d^p}{d\omega^p} ((-2\pi i \omega)^q (\mathcal{F}\varphi)(\omega)) e^{2\pi i \omega x} \\ &= \sum_{\omega \in \mathbb{Z}} \sum_{j=0}^{p \wedge q} \frac{1}{(-2\pi i)^{p-j}} \binom{p}{j} \binom{q}{j} j! (-2\pi i \omega)^{q-j} (\mathcal{F}\varphi)^{(p-j)}(\omega) e^{2\pi i \omega x}, \end{aligned} \tag{A.12}$$

where we applied the Leibniz formula in the last identity. Furthermore, for any integer $m \geq 1$,

$$(\mathcal{F}\varphi)^{(m)}(\omega) = (\mathbf{1}_{[-\frac{1}{2}, \frac{1}{2}]} * \psi)^{(m)}(\omega) = \int_{\omega - \frac{1}{2}}^{\omega + \frac{1}{2}} \psi^{(m)}(y) dy = \psi^{(m-1)}(\omega + \tfrac{1}{2}) - \psi^{(m-1)}(\omega - \tfrac{1}{2}).$$

For any $\omega \in \mathbb{Z}$, $\omega \pm \frac{1}{2} \in \mathbb{Z} + \frac{1}{2}$, but $\psi$ is supported in $[-\frac{1}{4}, \frac{1}{4}]$, hence $\psi$ is identically 0 in any small enough neighborhood of $\mathbb{Z} + \frac{1}{2}$, whereby $\psi^{(m-1)}(\omega \pm \frac{1}{2}) = 0$ and $(\mathcal{F}\varphi)^{(m)}(\omega) = 0$.

Furthermore, from the definition (A.1) and properties of the function $\psi$ we have $(\mathcal{F}\varphi)(0) = 1$ and $(\mathcal{F}\varphi)(\omega) = 0$ for all $\omega \in \mathbb{Z} \backslash \{0\}$.

Going back to (A.12), we assumed $p \neq q$. If $p \geq q+1$, then $p-j \geq 1$ and we get $(\mathcal{F}\varphi)^{(p-j)}(\omega) = 0$ for any $\omega \in \mathbb{Z}$. If $q \geq p+1$, then for $j = p$ we have $(\mathcal{F}\varphi)^{(p-j)}(\omega) = (\mathcal{F}\varphi)(\omega) = 0$ for $\omega \in \mathbb{Z} \backslash \{0\}$ and in the case $\omega = 0$, the factor $\omega^{q-j} = \omega^{q-p}$ in (A.12) vanishes. We thus conclude that the expression in (A.12) vanishes.

To show (A.11), we note that $\varphi^{(q)}$ is Schwartz and hence $|\varphi^{(q)}(x)| \lesssim |x|^{-\alpha} \wedge 1$ for any $\alpha \geq 1$, where the proportionality constant depends on $q$ and $\alpha$. Using (A.10), for $k < x$ we can switch the index of the sum to get

$$\Big| \sum_{\ell \geq k} (x-\ell)^p \varphi^{(q)}(x-\ell) \Big| = \Big| \sum_{\ell < k} (x-\ell)^p \varphi^{(q)}(x-\ell) \Big| \lesssim \sum_{\ell < k} (x-\ell)^p ((x-\ell)^{-\alpha} \wedge 1). \tag{A.13}$$

We have $x - \ell > 0$ in this sum, and for $\alpha > p+1$ it is estimated by $C_1(x-k)^{p+1-\alpha}$. □

**Lemma A.6** *For any $x \in \mathbb{R}$ and integers $q \geq p \geq 1$, we have*

$$\sum_{\ell_1 \in \mathbb{Z}} \sum_{\ell_2 \geq \ell_1} \cdots \sum_{\ell_p \geq \ell_{p-1}} \varphi^{(q)}(x - \ell_p) = 0. \tag{A.14}$$

*Proof.* The left-hand side of (A.14) is the limit of $\sum_{\ell_1 \geq k} \sum_{\ell_2 \geq \ell_1} \cdots \sum_{\ell_p \geq \ell_{p-1}} \varphi^{(q)}(x - \ell_p)$ as $k \to -\infty$, if such a limit exists. We fix $x \in \mathbb{R}$ and consider $k \in \mathbb{Z}$ such that $k < x$. Changing the order of summations, we can write

$$\begin{aligned}\sum_{\ell_1 \geq k} \sum_{\ell_2 \geq \ell_1} \cdots \sum_{\ell_p \geq \ell_{p-1}} \varphi^{(q)}(x - \ell_p) &= \sum_{\ell_p \geq k} \varphi^{(q)}(x - \ell_p) \sum_{\substack{\ell_1, \ell_2, \dots, \ell_{p-1} \in \mathbb{Z}: \\ k \leq \ell_1 \leq \cdots \leq \ell_{p-1} \leq \ell_p}} 1 \\ &= \sum_{\ell_p \geq k} \varphi^{(q)}(x - \ell_p) \frac{1}{(p-1)!} \prod_{j=1}^{p-1} (\ell_p - k + j). \end{aligned} \tag{A.15}$$

Note that $\prod_{j=1}^{p-1}(\ell_p - k + j)$ is a polynomial of degree $p - 1$ in the variable $\ell_p - k$, where the coefficients of the polynomial depend on $p$. That is, (A.15) can be written as

$$\sum_{\ell \geq k} \varphi^{(q)}(x - \ell) \mathfrak{P}(\ell - k), \tag{A.16}$$

where $\mathfrak{P}$ is a polynomial of degree $p - 1$ (we use $\ell$ instead of $\ell_p$ to lighten notation). Writing $\ell - k = (\ell - x) + (x - k)$, we can write $\mathfrak{P}(\ell - k)$ as a linear combination of terms of the form $(x - \ell)^a (x - k)^b$, where $a, b \in \mathbf{N}$ such that $a + b \leq p - 1$. Hence, (A.16) can be written as a linear combination of terms of the form

$$(x - k)^b \sum_{\ell \geq k} \varphi^{(q)}(x - \ell)(x - \ell)^a.$$

Since $p \geq 1$ and $a \leq p - 1 \leq q - 1$ (hence $a \neq q$), the assumptions of Lemma A.5 are satisfied, and by (A.11) we have

$$\Big| (x - k)^b \sum_{\ell \geq k} \varphi^{(q)}(x - \ell)(\ell - x)^a \Big| \leq C (x - k)^{b + a + 1 - \alpha}$$

for any $\alpha > a + 1$. Taking $\alpha > p$, all these terms vanish as $k \to -\infty$ which means that (A.16) vanishes in the limit as well. □

Having the proved properties of the function $\varphi$, we can prove Proposition A.1.

*Proof of Proposition A.1.* We fix $j = (j_0, j_1) \in \mathbf{N}^2$ and separate the cases $j = (0, 0)$ and $j \geq (0, 1)$.

**Step I.** Consider first the case $j = (0, 0)$. We write (A.2) as

$$U(t, x) = \sum_{k \in \mathbb{Z}:\, |k - x| \leq |x|/2} \varphi(x - k) u(t, k) + \sum_{k \in \mathbb{Z}:\, |k - x| > |x|/2} \varphi(x - k) u(t, k). \tag{A.17}$$

If $|k - x| \leq |x|/2$, then $\|(t, k)\|_{\mathfrak{s}} \geq \frac{1}{2} \|(t, x)\|_{\mathfrak{s}}$. Using the assumption (A.3) with $j = (0, 0)$ we get

$$\Big| \sum_{k \in \mathbb{Z}:\, |k - x| \leq |x|/2} \varphi(x - k) u(t, k) \Big| \lesssim (\|(t, x)\|_{\mathfrak{s}} + 1)^{-\beta} \sum_{k \in \mathbb{Z}:\, |k - x| \leq |x|/2} |\varphi(x - k)|.$$

Since $\varphi$ is Schwartz, we have $|\varphi(x)| \lesssim |x|^{-\alpha} \wedge 1$ for any $\alpha \geq 2$ and with a proportionality constant depending on $\alpha$. Then the preceding expression is bounded by a constant times

$$(\|(t, x)\|_{\mathfrak{s}} + 1)^{-\beta} \sum_{k \in \mathbb{Z}:\, |k - x| \leq |x|/2} (|x - k|^{-\alpha} \wedge 1),$$

where the sum can be estimated by a Riemann integral and it is bounded independently of $x$ if $\alpha > 1$. Hence, the first sum in (A.17) is absolutely bounded by $C_1 (\|(t, x)\|_{\mathfrak{s}} + 1)^{-\beta}$.

Now, we will bound the second sum in (A.17). We have $\|(t,k)\|_{\mathfrak{s}} \geq \sqrt{t}$ and assumption (A.3) yields

$$\Big|\sum_{k\in\mathbb{Z}:\,|k-x|>|x|/2} \varphi(x-k)u(t,k)\Big| \lesssim (\sqrt{t}+1)^{-\beta} \sum_{k\in\mathbb{Z}:\,|k-x|>|x|/2} |\varphi(x-k)|$$
$$\lesssim (\sqrt{t}+1)^{-\beta} \sum_{k\in\mathbb{Z}:\,|k-x|>|x|/2} (|x-k|^{-\alpha}\wedge 1).$$

The sum can be estimated by a Riemann integral and the preceding expression is bounded by a constant times

$$(\sqrt{t}+1)^{-\beta}(|x|+1)^{1-\alpha} \lesssim (\sqrt{t}+1)^{-\beta}(|x|+1)^{-\beta},$$

for $\alpha > \beta$. Since $(\sqrt{t}+1)(|x|+1) \geq \|(t,x)\|_{\mathfrak{s}}+1$, the latter is bounded by $(\|(t,x)\|_{\mathfrak{s}}+1)^{-\beta}$. This completes the proof of (A.4) in the case $j=(0,0)$.

**Step II.** Consider the case $j=(j_0,j_1)\geq(0,1)$. We call $p=j_1$ to lighten notation. Fixing $x\in\mathbb{R}$, we define the $p$-fold "discrete-antiderivative" of $\varphi^{(p)}(x-\bullet)$ as

$$\Phi_x(k) := \sum_{\ell_1\geq k}\sum_{\ell_2\geq \ell_1}\cdots\sum_{\ell_p\geq \ell_{p-1}} \varphi^{(p)}(x-\ell_p), \tag{A.18}$$

which we expect to yield a function close to $\varphi^{(p)}(x-\bullet)$. It is straightforward to verify that

$$(\nabla^+)^p\Phi_x(k) = (-1)^p\varphi^{(p)}(x-k). \tag{A.19}$$

Furthermore, by (A.2), (A.19) and summation by parts

$$\begin{aligned} D^jU(t,x) &= \sum_{k\in\mathbb{Z}}\varphi^{(p)}(x-k)\partial_t^{j_0}u(t,k) \\ &= (-1)^p\sum_{k\in\mathbb{Z}}(\nabla^+)^p\Phi_x(k)\partial_t^{j_0}u(t,k) = \sum_{k\in\mathbb{Z}}\Phi_x(k)(\nabla^-)^p\partial_t^{j_0}u(t,k), \end{aligned} \tag{A.20}$$

where there are no boundary terms since the involved functions decay at $|k|\to\infty$.

Now, let us show that $k\mapsto\Phi_x(k)$ is highly concentrated around $x$, in the sense that

$$|\Phi_x(k)| \lesssim |x-k|^{p-\alpha}\wedge 1, \tag{A.21}$$

for any $\alpha$ large enough with a proportionality constant depending on $p$ and $\alpha$. Indeed, since $\varphi^{(p)}$ is Schwartz, we have $|\varphi^{(p)}(x)| \lesssim |x|^{-\alpha}\wedge 1$ for any $\alpha>0$. We fix some $\alpha \geq 2+|j|_{\mathfrak{s}}$ for the rest of the proof. If $k>x$ then in the expression (A.18) we have $\ell_p\geq\ell_{p-1}\geq\cdots\geq\ell_1\geq k>x$, and hence

$$|\Phi_x(k)| \lesssim \sum_{\ell_1\geq k}\sum_{\ell_2\geq \ell_1}\cdots\sum_{\ell_p\geq \ell_{p-1}} ((\ell_p-x)^{-\alpha}\wedge 1) \lesssim (k-x)^{p-\alpha}\wedge 1;$$

whereas if $k<x$ then by (A.14) we can switch iteratively the direction of the indices of all the sums and get

$$\begin{aligned} \sum_{\ell_1\geq k}\sum_{\ell_2\geq \ell_1}\cdots\sum_{\ell_p\geq \ell_{p-1}} \varphi^{(p)}(x-\ell_p) &= -\sum_{\ell_1\leq k-1}\sum_{\ell_2\geq \ell_1}\cdots\sum_{\ell_p\geq \ell_{p-1}} \varphi^{(p)}(x-\ell_p) \\ &= \sum_{\ell_1\leq k-1}\sum_{\ell_2\leq \ell_1-1}\cdots\sum_{\ell_p\geq \ell_{p-1}} \varphi^{(p)}(x-\ell_p) \\ &= \cdots = (-1)^p\sum_{\ell_1\leq k-1}\sum_{\ell_2\leq \ell_1-1}\cdots\sum_{\ell_p\leq \ell_{p-1}-1} \varphi^{(p)}(x-\ell_p). \end{aligned}$$

In this expression we have $\ell_p \leq \ell_{p-1} \leq \cdots \leq \ell_1 \leq k < x$, and thus

$$|\Phi_x(k)| \lesssim \sum_{\ell_1 \leq k} \sum_{\ell_2 \leq \ell_1} \cdots \sum_{\ell_p \leq \ell_{p-1}} ((x-\ell_p)^{-\alpha} \wedge 1) \lesssim (x-k)^{p-\alpha} \wedge 1.$$

This verifies our claim (A.21).

Now, let us go back to the expression (A.20) for $D^j U(t,x)$. Similarly to (A.17) we can decompose

$$D^j U(t,x) = \sum_{k \in \mathbb{Z}\colon |k-x| \leq |x|/2} \Phi_x(k)(\nabla^-)^p \partial_t^{j_0} u(t,k) + \sum_{k \in \mathbb{Z}\colon |k-x| > |x|/2} \Phi_x(k)(\nabla^-)^p \partial_t^{j_0} u(t,k). \tag{A.22}$$

We bound the first sum in the same way as we bounded the first sum in (A.22), with the only difference that now we apply assumption (A.3) with $j = (j_0, p)$. We get

$$\Big| \sum_{k \in \mathbb{Z}\colon |k-x| \leq |x|/2} \Phi_x(k)(\nabla^-)^p \partial_t^{j_0} u(t,k)\Big| \lesssim (\|(t,x)\|_{\mathfrak{s}} + 1)^{-|j|_{\mathfrak{s}} - \beta} \sum_{k \in \mathbb{Z}\colon |k-x| \leq |x|/2} |\Phi_x(k)|.$$

Applying (A.21), we bound this expression by a constant times

$$(\|(t,x)\|_{\mathfrak{s}} + 1)^{-|j|_{\mathfrak{s}} - \beta} \sum_{k \in \mathbb{Z}\colon |k-x| \leq |x|/2} (|x-k|^{p-\alpha} \wedge 1) \lesssim (\|(t,x)\|_{\mathfrak{s}} + 1)^{-|j|_{\mathfrak{s}} - \beta},$$

where the last bound holds for $\alpha > p + \beta$.

We bound the second sum in (A.22) similarly to how we bounded the second sum in (A.17), but using assumption (A.3) with $j = (j_0, p)$. We get

$$\Big| \sum_{k \in \mathbb{Z}\colon |k-x| > |x|/2} \Phi_x(k)(\nabla^-)^p \partial_t^{j_0} u(t,k)\Big| \lesssim (\sqrt{t} + 1)^{-|j|_{\mathfrak{s}} - \beta} \sum_{k \in \mathbb{Z}\colon |k-x| > |x|/2} |\Phi_x(k)|.$$

The estimate (A.21) allows to bound this expression by a constant times

$$(\sqrt{t} + 1)^{-|j|_{\mathfrak{s}} - \beta} \sum_{k \in \mathbb{Z}\colon |k-x| > |x|/2} (|x-k|^{p-\alpha} \wedge 1) \lesssim (\sqrt{t} + 1)^{-|j|_{\mathfrak{s}} - \beta} (|x| + 1)^{p-\alpha+\beta},$$

where the last estimate holds for $\alpha > p + \beta$. By our choice of $\alpha \geq |j|_{\mathfrak{s}} + 1 + \beta$, the latter is smaller than $(\sqrt{t}+1)^{-|j|_{\mathfrak{s}}-\beta}(|x|+1)^{-|j|_{\mathfrak{s}}-\beta} \leq (\|(t,x)\|_{\mathfrak{s}} + 1)^{-|j|_{\mathfrak{s}} - \beta}$. □

## A.2 Extension of $G^\varepsilon$

We write $G^\varepsilon_{\text{ext}}$ for the extension (A.5) of the discrete heat kernel (4.10), i.e. $G^\varepsilon_{\text{ext}} = \mathsf{Ext}(G^\varepsilon)$. The goal of this section is to prove bounds on $G^\varepsilon_{\text{ext}}$, which will follow if we prove the bounds (A.6). Before deriving these bounds, we show the following bounds on the Fourier transform of $G^\varepsilon$.

**Lemma A.7** *Let $g_\varepsilon(t,\omega) := e^{-f_\varepsilon(\omega)t}$, where the function $f_\varepsilon$ is defined in (4.9). Then there exist $\varepsilon_\star > 0$ and $a > 0$ such that for any $j = (j_0, j_1) \in \mathbf{N}^2$ there is a constant $A = A(j) > 0$ such that the following bounds hold uniformly in $t \in \mathbb{R}_+$, $\omega \in [-\frac{1}{2\varepsilon}, \frac{1}{2\varepsilon}]$ and $\varepsilon \in (0, \varepsilon_\star)$:*

1. *for $|\omega|\sqrt{t} \geq 1$,*

$$e^{a\omega^2 t} |\partial_t^{j_0} \partial_\omega^{j_1} g_\varepsilon(t,\omega)| \leq \begin{cases} A|\omega|^{2j_0 - j_1} & \text{if } 2j_0 \geq j_1, \\ A(\sqrt{t})^{j_1 - 2j_0} & \text{if } 2j_0 < j_1, \end{cases} \tag{A.23}$$

2. *for $|\omega|\sqrt{t} < 1$ and $t \geq \varepsilon^2$,*

$$|\partial_t^{j_0} \partial_\omega^{j_1} g_\varepsilon(t,\omega)| \leq A(\sqrt{t})^{j_1 - 2j_0}, \tag{A.24}$$

3. *for* $|\omega|\sqrt{t} < 1$ *and* $0 < t < \varepsilon^2$,

$$|\partial_t^{j_0}\partial_\omega^{j_1} g_\varepsilon(t,\omega)| \le \begin{cases} A|\omega|^{2j_0-j_1} & \text{if } 2j_0 \ge j_1, \\ A\varepsilon^{j_1-2j_0} & \text{if } 2j_0 < j_1. \end{cases} \tag{A.25}$$

*Proof.* By the Leibniz rule we have

$$\begin{aligned} \partial_\omega^{j_1}\partial_t^{j_0} g_\varepsilon(t,\omega) &= (-1)^{j_0}\partial_\omega^{j_1}(f_\varepsilon(\omega)^{j_0}e^{-f_\varepsilon(\omega)t}) \\ &= (-1)^{j_0}\sum_{n_0+\dots+n_{j_0}=j_1}\binom{j_0}{n_1,\dots,n_{j_0}}\Big(\prod_{i=1}^{j_0} f_\varepsilon^{(n_i)}(\omega)\Big)\partial_\omega^{n_0}e^{-f_\varepsilon(\omega)t}. \end{aligned} \tag{A.26}$$

From (4.9), Leibniz rule, Lemma 4.1 and the bounds $|\omega| \lesssim \varepsilon^{-1}$ we get

$$|f_\varepsilon^{(n)}(\omega)| \lesssim \begin{cases} |\omega|^{2-n} & \text{if } n \in \{0,1\}, \\ \varepsilon^{n-2} & \text{if } n \ge 2. \end{cases} \tag{A.27}$$

This bound suggests that we need to distinguish the values $n_i$ in (A.26) satisfying $n_i \in \{0,1\}$ and $n_i \ge 2$. Then (A.26) can be bounded a

$$|\partial_\omega^{j_1}\partial_t^{j_0} g_\varepsilon(t,\omega)| \lesssim \sum_{m=0}^{j_0}\sum_{\substack{n_0+\dots+n_{j_0}=j_1 \\ n_i\in\{0,1\} \text{ for } 1\le i\le m \\ n_i\ge 2 \text{ for } i>m}}\Big(\prod_{i=1}^{m}|\omega|^{2-n_i}\Big)\Big(\prod_{i=m+1}^{j_0}\varepsilon^{n_i-2}\Big)|\partial_\omega^{n_0}e^{-f_\varepsilon(\omega)t}|. \tag{A.28}$$

We can consider only the first indices taking values in $\{0,1\}$ because the order of the indices does not play any role up to a combinatorial factor. Furthermore, Faà di Bruno's formula yields

$$\partial_\omega^{n_0}e^{-f_\varepsilon(\omega)t} = \sum_{1\cdot p_1+2\cdot p_2+\dots+n_0\cdot p_{n_0}=n_0}\frac{n_0!}{p_1!\cdots p_{n_0}!}e^{-f_\varepsilon(\omega)t}\prod_{j=1}^{n_0}\Big(\frac{f_\varepsilon^{(j)}(\omega)t}{j!}\Big)^{p_j}, \tag{A.29}$$

where $p_1,\dots,p_{n_0} \ge 0$ in the sum. By Lemma 4.1 we have $f_\varepsilon(\omega) \ge 8\omega^2$. Then (A.29) can be bounded as

$$|\partial_\omega^{n_0}e^{-f_\varepsilon(\omega)t}| \lesssim \sum_{1\cdot p_1+2\cdot p_2+\dots+n_0\cdot p_{n_0}=n_0} e^{-8\omega^2 t}(|\omega|t)^{p_1}\prod_{j=2}^{n_0}(\varepsilon^{j-2}t)^{p_j}. \tag{A.30}$$

In order to estimate these expressions, we need to consider different regions for $t$ and $\omega$.

We need to consider the two cases: $|\omega|\sqrt{t} \ge 1$ and $|\omega|\sqrt{t} < 1$. We start with the case $|\omega|\sqrt{t} \ge 1$ (which implies $t \gtrsim \varepsilon^2$ because $|\omega| \lesssim \varepsilon^{-1}$). We use $\varepsilon \lesssim |\omega|^{-1}$ to bound (A.30) by

$$\sum_{1\cdot p_1+2\cdot p_2+\dots+n_0\cdot p_{n_0}=n_0} e^{-8\omega^2 t}(|\omega|t)^{p_1}\prod_{j=2}^{n_0}(|\omega|^{2-j}t)^{p_j} = |\omega|^{-n_0}e^{-8\omega^2 t}P_{n_0}(\omega^2 t),$$

where $P_{n_0}$ is a polynomial of degree $n_0$. Then the simple bound $e^{-7\omega^2 t}P_{n_0}(\omega^2 t) \lesssim 1$ yields $|\partial_\omega^{n_0}e^{-f_\varepsilon(\omega)t}| \lesssim |\omega|^{-n_0}e^{-\omega^2 t}$ and hence (A.28) is bounded by a constant multiple of

$$\sum_{m=0}^{j_0}\sum_{\substack{n_0+\dots+n_{j_0}=j_1 \\ n_i\in\{0,1\} \text{ for } 1\le i\le m \\ n_i\ge 2 \text{ for } i>m}}\Big(\prod_{i=1}^{m}|\omega|^{2-n_i}\Big)\Big(\prod_{i=m+1}^{j_0}\varepsilon^{n_i-2}\Big)|\omega|^{-n_0}e^{-\omega^2 t}. \tag{A.31}$$

The power of $|\omega|$ inside the sum is $2m-(n_0+\cdots+n_m)$, and the power of $\varepsilon$ is $(n_{m+1}+\cdots+n_{j_0})-2(j_0-m)=j_1-2j_0+2m-(n_0+\cdots+n_m)$. Introducing a new variable $r=n_0+\cdots+n_m$, we bound the preceding expression by a constant times

$$\sum_{m=0}^{j_0}\sum_{r=0}^{j_1-2j_0+2m}|\omega|^{2m-r}\varepsilon^{j_1-2j_0+2m-r}e^{-\omega^2 t}\lesssim\sum_{p=2j_0-j_1}^{2j_0}|\omega|^p\varepsilon^{j_1-2j_0+p}e^{-\omega^2 t},\tag{A.32}$$

where we changed the variable to $p=2m-r$. Using $\varepsilon\lesssim|\omega|^{-1}$ we bound (A.32) by a constant times $|\omega|^{2j_0-j_1}e^{-\omega^2 t}$ which is the desired bound (A.23) in the case $2j_0\geq j_1$. In the case $2j_0<j_1$ we use $|\omega|\lesssim\varepsilon^{-1}$ to bound this expression by a constant times $\sum_{p=2j_0-j_1}^{0}|\omega|^p\varepsilon^{j_1-2j_0+p}e^{-\omega^2 t}$. Furthermore, we use $|\omega|\geq 1/\sqrt{t}$ and $\varepsilon\lesssim\sqrt{t}$ to bound the latter by a constant multiple of $(\sqrt{t})^{j_1-2j_0}e^{-\omega^2 t}$, which is the desired bound (A.23).

Now, we consider the case $|\omega|\sqrt{t}<1$. Then we estimate $e^{-8\omega^2 t}\leq 1$ and bound (A.30) by

$$|\partial_\omega^{n_0}e^{-f_\varepsilon(\omega)t}|\lesssim\sum_{1\cdot p_1+2\cdot p_2+\cdots+n_0\cdot p_{n_0}=n_0}t^{p\frac{1}{2}}\prod_{j=2}^{n_0}(\varepsilon^{j-2}t)^{p_j}.\tag{A.33}$$

We need to distinguish the two cases: $t>\varepsilon^2$ and $0<t\leq\varepsilon^2$.

Let us consider the case $t>\varepsilon^2$. Since $\varepsilon^{-1}\sqrt{t}>1$, we can bound (A.33) by

$$\sum_{p_1=0}^{n_0}\varepsilon^{n_0-p_1}t^{p\frac{1}{2}}\sum_{2\cdot p_2+\cdots+n_0\cdot p_{n_0}=n_0-p_1}\prod_{j=2}^{n_0}(\varepsilon^{-1}\sqrt{t})^{2p_j}\lesssim\sum_{p_1=0}^{n_0}\varepsilon^{n_0-p_1}t^{p\frac{1}{2}}(\varepsilon^{-1}\sqrt{t})^{n_0-p_1}\lesssim t^{n_0/2}.$$

In the first inequality we have estimated a polynomial with respect to $\varepsilon^{-1}\sqrt{t}$ by its leading monomial. Then we use $|\omega|\leq t^{-\frac{1}{2}}$ to bound (A.28) by a constant multiple of

$$\sum_{m=0}^{j_0}\sum_{\substack{n_0+\cdots+n_{j_0}=j_1\\ n_i\in\{0,1\}\text{ for }1\leq i\leq m\\ n_i\geq 2\text{ for }i>m}}\Big(\prod_{i=1}^{m}t^{n_i/2-1}\Big)\Big(\prod_{i=m+1}^{j_0}\varepsilon^{n_i-2}\Big)t^{n_0/2}.$$

This expression can be obtained from (A.31) if we remove the exponential and replace $|\omega|$ by $t^{-\frac{1}{2}}$. Hence, in the same way as we derived (A.32), we bound the preceding expression by a constant times

$$\sum_{p=2j_0-j_1}^{2j_0}t^{-p/2}\varepsilon^{j_1-2j_0+p}.\tag{A.34}$$

Using $t>\varepsilon^2$ we bound (A.32) by a constant times $t^{j\frac{1}{2}-j_0}$ which is the desired bound (A.24).

Finally, in the case $0<t\leq\varepsilon^2$ we estimate (A.33) by a constant times $\varepsilon^{n_0}$. Then (A.28) is bounded by a constant multiple of

$$\sum_{m=0}^{j_0}\sum_{\substack{n_0+\cdots+n_{j_0}=j_1\\ n_i\in\{0,1\}\text{ for }1\leq i\leq m\\ n_i\geq 2\text{ for }i>m}}\Big(\prod_{i=1}^{m}|\omega|^{2-n_i}\Big)\Big(\prod_{i=m+1}^{j_0}\varepsilon^{n_i-2}\Big)\varepsilon^{n_0}.$$

The power of $|\omega|$ inside the sum is $2m-(n_1+\cdots+n_m)$, and the power of $\varepsilon$ is $(n_0+n_{m+1}+\cdots+n_{j_0})-2(j_0-m)=j_1-2j_0+2m-(n_1+\cdots+n_m)$. Using the new variable $r=n_1+\cdots+n_m$, we bound the preceding expression by a constant times

$$\sum_{m=0}^{j_0}\sum_{r=0}^{(j_1-2j_0+2m)\wedge m}|\omega|^{2m-r}\varepsilon^{j_1-2j_0+2m-r}\lesssim\sum_{p=(2j_0-j_1)\vee 0}^{2j_0}|\omega|^p\varepsilon^{j_1-2j_0+p},$$

where we changed the variable to $p = 2m - r$. If $2j_0 \geq j_1$, then we use $\varepsilon \lesssim |\omega|^{-1}$ to bound this expression by a constant times $|\omega|^{2j_0 - j_1}$. If $2j_0 < j_1$, then we use $|\omega| \lesssim \varepsilon^{-1}$ to bound the preceding expression by a constant times $\varepsilon^{j_1 - 2j_0}$. These are the required bounds (A.25) in the case $|\omega|\sqrt{t} < 1$ and $0 < t \leq \varepsilon^2$. □

Now, we can prove the bounds (A.6) for the kernel $G^\varepsilon$.

**Lemma A.8** *For any $n \in \mathbf{N}$ and $j \in \mathbf{N}^2$ there is $C = C(n, j) > 0$, such that*

$$|D_\varepsilon^j G^\varepsilon(t,x)| \leq C(\sqrt{t} + \varepsilon)^{-|j|_\mathfrak{s} - 1} \left( \frac{\sqrt{t} + \varepsilon}{\|(t,x)\|_\mathfrak{s} + \varepsilon} \right)^n, \tag{A.35}$$

*uniformly on $(t,x) \in \mathbb{R}_+ \times \Lambda_\varepsilon$. In particular, taking $n = |j|_\mathfrak{s} + 1$ we get*

$$|D_\varepsilon^j G^\varepsilon(t,x)| \leq C(\|(t,x)\|_\mathfrak{s} + \varepsilon)^{-|j|_\mathfrak{s} - 1}. \tag{A.36}$$

*Proof.* The proof consists of two steps: in the first step we establish the bound

$$|D_\varepsilon^j G^\varepsilon(t,x)| \leq C(\sqrt{t} + \varepsilon)^{-|j|_\mathfrak{s} - 1}, \tag{A.37}$$

and in the second step we prove

$$|D_\varepsilon^j G^\varepsilon(t,x)| \leq C(\sqrt{t} + \varepsilon)^{-|j|_\mathfrak{s} - 1} \left( \frac{\sqrt{t} + \varepsilon}{|x|} \right)^n, \tag{A.38}$$

for $x \neq 0$. These two bounds imply (A.35), since if $|x| > \sqrt{t} + \varepsilon$ then we use (A.38) and $\frac{\sqrt{t}+\varepsilon}{2|x|} \leq \frac{\sqrt{t}+\varepsilon}{\sqrt{t}+|x|+\varepsilon}$, whereas if $|x| \leq \sqrt{t} + \varepsilon$ then we use (A.37) and $\frac{1}{2} \leq \frac{\sqrt{t}+\varepsilon}{\sqrt{t}+|x|+\varepsilon}$.

**Step I.** We will prove the bound (A.37). Using $g_\varepsilon(t,\omega) = e^{-f_\varepsilon(\omega)t}$, (4.8) yields

$$D_\varepsilon^j G^\varepsilon(t,x) = \int_{-\frac{1}{2\varepsilon}}^{\frac{1}{2\varepsilon}} \partial_t^{j_0} g_\varepsilon(t,\omega) (i\omega m^-(\varepsilon\omega))^{j_1} e^{2\pi i \omega x} d\omega, \tag{A.39}$$

where we use the function $m^-$ defined in (4.3). Applying Lemma 4.1, we get

$$|D_\varepsilon^j G^\varepsilon(t,x)| \leq \int_{-\frac{1}{2\varepsilon}}^{\frac{1}{2\varepsilon}} |\partial_t^{j_0} g_\varepsilon(t,\omega)| (2\pi|\omega|)^{j_1} d\omega. \tag{A.40}$$

In order to bound this integral, we will consider two integration regions: $|\omega| \geq 1/\sqrt{t}$ and $|\omega| < 1/\sqrt{t}$. In the case $|\omega| \geq 1/\sqrt{t}$ the property (A.23) allows to bound $|\partial_t^{j_0} g_\varepsilon(t,\omega)| \lesssim |\omega|^{2j_0} e^{-a\omega^2 t}$. Then

$$\int_{1/\sqrt{t} \leq |\omega| \leq 1/2\varepsilon} |\partial_t^{j_0} g_\varepsilon(t,\omega)| (2\pi|\omega|)^{j_1} d\omega \lesssim \int_{1/\sqrt{t} \leq |\omega| \leq 1/2\varepsilon} |\omega|^{|j|_\mathfrak{s}} e^{-a\omega^2 t} d\omega. \tag{A.41}$$

By the change of variable $u = \omega\sqrt{t}$, this integral is bounded above by a constant times

$$(\sqrt{t})^{-|j|_\mathfrak{s} - 1} \int_1^\infty u^{|j|_\mathfrak{s}} e^{-au^2} du \lesssim (\sqrt{t})^{-|j|_\mathfrak{s} - 1},$$

which is a bound of the desired order (A.37) because $t \gtrsim \varepsilon^2$.

Now, we will bound the part of the integral (A.40) for which $|\omega| < 1/\sqrt{t}$. In the case $t \geq \varepsilon^2$ we use (A.24) to bound

$$\int_{|\omega| < 1/\sqrt{t}} |\partial_t^{j_0} g_\varepsilon(t,\omega)| (2\pi|\omega|)^{j_1} d\omega \lesssim (\sqrt{t})^{-2j_0} \int_{|\omega| < 1/\sqrt{t}} |\omega|^{j_1} d\omega \lesssim (\sqrt{t})^{-|j|_\mathfrak{s} - 1}, \tag{A.42}$$

which is of the desired order (A.37). In the case $t < \varepsilon^2$ we apply (A.25) to bound the integral on the left-hand side (A.42) by a constant multiple of

$$\int_{|\omega|\le 1/2\varepsilon} |\omega|^{|j|_{\mathfrak{s}}} d\omega \lesssim \varepsilon^{-|j|_{\mathfrak{s}}-1}.$$

This is a bound of the desired order (A.37).

**Step II.** We now prove the bound (A.38). Throughout the proof we assume $x \neq 0$. Since $e^{2\pi i\omega x} = \partial_\omega^n\big(\frac{e^{2\pi i\omega x}}{(2\pi ix)^n}\big)$ for any integer $n \ge 1$, (A.39) yields

$$\begin{aligned} D_\varepsilon^j G^\varepsilon(t,x) &= \int_{-\frac{1}{2\varepsilon}}^{\frac{1}{2\varepsilon}} \partial_t^{j_0} g_\varepsilon(t,\omega)(i\omega m^-(\varepsilon\omega))^{j_1} \partial_\omega^n \Big(\frac{e^{2\pi i\omega x}}{(2\pi ix)^n}\Big) d\omega \\ &= (-1)^n \int_{-\frac{1}{2\varepsilon}}^{\frac{1}{2\varepsilon}} \partial_\omega^n \Big(\partial_t^{j_0} g_\varepsilon(t,\omega)(i\omega m^-(\varepsilon\omega))^{j_1}\Big) \frac{e^{2\pi i\omega x}}{(2\pi ix)^n} d\omega, \end{aligned} \tag{A.43}$$

where in the last line we used integration by parts with periodic boundary values. By a generalized Leibniz rule we write (A.43) as

$$\sum_{\substack{n_1+n_2=n,\\ n_1\ge 0, 0\le n_2\le j_1}} \binom{n}{n_1}(-1)^n \int_{-\frac{1}{2\varepsilon}}^{\frac{1}{2\varepsilon}} (\partial_\omega^{n_1}\partial_t^{j_0} g_\varepsilon(t,\omega))(\partial_\omega^{n_2}(i\omega m^-(\varepsilon\omega))^{j_1}) \frac{e^{2\pi i\omega x}}{(2\pi ix)^n} d\omega. \tag{A.44}$$

To bound the integrals, we will consider two cases: $n_1 \le 2j_0$ and $n_1 > 2j_0$.

**Step II-1.** Let us first consider the terms in (A.44) with $n_1 \le 2j_0$. We have $|\partial_\omega^{n_2}\omega^{j_1}| \lesssim |\omega|^{j_1-n_2}$ for any $n_2 \le j_1$. Then the bound on the integral (A.44) is derived in exactly the same way as in Step I, but with $2j_0$ replaced by $2j_0 - n_1$ and $j_1$ replaced by $j_1 - n_2$. More precisely, we get in the case $t \ge \varepsilon^2$ a bound of the order $\frac{1}{|x|^n}(\sqrt{t})^{n_1+n_2-|j|_{\mathfrak{s}}-1} \le \frac{1}{|x|^n}(\sqrt{t})^{n-|j|_{\mathfrak{s}}-1}$. In the case $t < \varepsilon^2$ we get a bound of the order $\frac{1}{|x|^n}\varepsilon^{n-|j|_{\mathfrak{s}}-1}$. Combining these two cases, we get a bound of order $(\sqrt{t}+\varepsilon)^{-1-|j|_{\mathfrak{s}}}\big(\frac{\sqrt{t}+\varepsilon}{|x|}\big)^n$, as desired in (A.38).

**Step II-2.** Now, we will bound the terms in (A.44) with $n_1 > 2j_0$. For this, we use $|\partial_\omega^{n_2}\omega^{j_1}| \lesssim |\omega|^{j_1-n_2}$ to bound the integral absolutely by a constant times

$$\frac{1}{|x|^n}\int_{-\frac{1}{2\varepsilon}}^{\frac{1}{2\varepsilon}} |\partial_\omega^{n_1}\partial_t^{j_0} g_\varepsilon(t,\omega)||\omega|^{j_1-n_3} d\omega, \tag{A.45}$$

and we consider the two domains of integration: $|\omega| \ge 1/\sqrt{t}$ and $|\omega| < 1/\sqrt{t}$. In the case $|\omega| \ge 1/\sqrt{t}$, we apply (A.23) to bound

$$\begin{aligned} &\frac{1}{|x|^n}\int_{1/\sqrt{t}\le|\omega|\le 1/2\varepsilon} |\partial_\omega^{n_1}\partial_t^{j_0} g_\varepsilon(t,\omega)||\omega|^{j_1-n_2} d\omega \\ &\qquad\qquad \lesssim \frac{1}{|x|^n}(\sqrt{t})^{n_1-2j_0}\int_{1/\sqrt{t}\le|\omega|\le 1/2\varepsilon} |\omega|^{j_1-n_2} e^{-a\omega^2 t} d\omega. \end{aligned} \tag{A.46}$$

We change the integration variable to $u = \omega\sqrt{t}$ and bound this expression by

$$2\frac{1}{|x|^n}(\sqrt{t})^{n_1+n_2-|j|_{\mathfrak{s}}-1}\int_1^\infty |u|^{j_1-n_2} e^{-au^2} du \lesssim \frac{1}{|x|^n}(\sqrt{t})^{n_1+n_2-|j|_{\mathfrak{s}}-1}.$$

Since the bounds $|\omega| \ge 1/\sqrt{t}$ and $|\omega| \lesssim \varepsilon^{-1}$ imply $\varepsilon \lesssim \sqrt{t}$, the preceding expression is bounded by a constant times $\frac{1}{|x|^n}(\sqrt{t})^{n-|j|_{\mathfrak{s}}-1}$. This is a bound of the desired order (A.38) because $\varepsilon \lesssim \sqrt{t}$.

In the case $|\omega| < 1/\sqrt{t}$ we need to bound the respective part of the integral (A.45) for $t \geq \varepsilon^2$ and $t < \varepsilon^2$. For $t \geq \varepsilon^2$, (A.24) yields

$$\frac{1}{|x|^n} \int_{|\omega|<(1/\sqrt{t})\wedge(1/2\varepsilon)} |\partial_\omega^{n_1} \partial_t^{j_0} g_\varepsilon(t,\omega)| |\omega|^{j_1-n_2} d\omega \tag{A.47}$$
$$\lesssim \frac{\varepsilon^{n_2}}{|x|^n} (\sqrt{t})^{n_1-2j_0} \int_{|\omega|<1/\sqrt{t}} |\omega|^{j_1-n_2} d\omega \lesssim \frac{1}{|x|^n} (\sqrt{t})^{n_1+n_2-|j|_\mathfrak{s}-1},$$

where we used the restriction $j_1 \geq n_2$ in the sum in (A.44). Using $\varepsilon \leq \sqrt{t}$, we bound this expression by $\frac{1}{|x|^n}(\sqrt{t})^{n-|j|_\mathfrak{s}-1}$, which is a bound of the desired order (A.38).

If $t < \varepsilon^2$, then (A.25) allows to bound the integral on the left-hand side of (A.47) by

$$\frac{\varepsilon^{n_1-2j_0}}{|x|^n} \int_{|\omega|\leq 1/2\varepsilon} |\omega|^{j_1-n_2} d\omega \lesssim \frac{\varepsilon^{n-|j|_\mathfrak{s}-1}}{|x|^n},$$

where we used again the restriction $j_1 \geq n_3$. This is a bound of the desired order (A.38).

Combining the preceding bounds, we conclude that the terms in the sum (A.44) with $n_1 > 2j_0$ are absolutely bounded by a constant multiple of $(\sqrt{t}+\varepsilon)^{-1-|j|_\mathfrak{s}}(\frac{\sqrt{t}+\varepsilon}{|x|})^n$. Hence, the whole sum in (A.44) is bounded by the same expression and we get (A.38). □

Combining Lemma A.8 and Corollary A.2 we get the following result:

**Proposition A.9** *For any $j \in \mathbf{N}^2$ there is a constant $C = C(j) > 0$ such that the extension $G^\varepsilon_{\text{ext}} = \mathsf{Ext}(G^\varepsilon)$, defined in* (A.5)*, satisfies*

$$|D^j G^\varepsilon_{\text{ext}}(t,x)| \leq C(\|(t,x)\|_\mathfrak{s} + \varepsilon)^{-|j|_\mathfrak{s}-1},$$

*uniformly over $(t,x) \in \mathbb{R}_+ \times \mathbb{R}$.*

## Appendix B Bounds on functions defined via the symmetric exclusion process

The stopping time $\tau_{\varepsilon,\mathfrak{m}}$ introduced in Section 2.2 guarantees bounds on the random functions $\widetilde{\mathbf{C}}_{\varepsilon,\mathfrak{m}}$ and $\widetilde{\mathsf{C}}_{\varepsilon,\mathfrak{m}}$, which we prove in Lemmas B.3 and B.4 below. Before that we need to derive some estimates on the rescaled height function $\tilde{h}^{\varepsilon,\mathfrak{m}}_{\text{sym}}$ of the symmetric exclusion process, defined in Section 2.2. Throughout this section we use the map (5.10) of a semidiscrete function into a function in time and a distribution in space. We write $(\iota_\varepsilon f_\varepsilon)(t,\varphi)$ for the duality pairing of this distribution with a test function $\varphi$.

**Lemma B.1** *For any $p \geq 1$ and $T > 0$ there is a constant $C > 0$ such that the bound*

$$\Big(\mathbb{E}\Big|\sup_{t\in[\tau_{\varepsilon,\mathfrak{m}},T]} (\lambda \vee \varepsilon \vee \sqrt{t-\tau_{\varepsilon,\mathfrak{m}}})^{2-\alpha} |\iota_\varepsilon(\nabla^+_\varepsilon \tilde{h}^{\varepsilon,\mathfrak{m}}_{\text{sym}})(t, \nabla^+_\varepsilon \varphi^\lambda_x)|\Big|^p\Big)^{\frac{1}{p}} \leq C\mathfrak{m} \tag{B.1}$$

*holds uniformly over $t \in [0,T]$, $\varphi \in \mathcal{B}^2$, $x \in \Lambda_\varepsilon$, $\lambda \in (0,1]$ and $\varepsilon \in (0,1)$. The constant $\alpha$ is the same as in the statement of Theorem 1.1.*

*Proof.* By the definition in Section 2.2 we have the equation

$$d\tilde{h}^{\varepsilon,\mathfrak{m}}_{\text{sym}}(t,x) = \frac{1}{2}(1+\sqrt{\varepsilon})(1-\varrho_\varepsilon^2)\Delta_\varepsilon \tilde{h}^{\varepsilon,\mathfrak{m}}_{\text{sym}}(t,x)dt + d\widetilde{M}^{\varepsilon,\mathfrak{m}}_{\text{sym}}(t,x), \tag{B.2}$$

for $t > \tau_{\varepsilon,\mathfrak{m}}$, with the initial condition $\tilde{h}^{\varepsilon,\mathfrak{m}}_{\text{sym}}(\tau_{\varepsilon,\mathfrak{m}}, \bullet) = \tilde{h}^{\varepsilon}(\tau_{\varepsilon,\mathfrak{m}}-, \bullet)$. Let us use the shorthand $t_{\varepsilon,\mathfrak{m}} = t - \tau_{\varepsilon,\mathfrak{m}}$. Then the mild form of the equation is

$$\tilde{h}^{\varepsilon,\mathfrak{m}}_{\text{sym}}(t,x) = (G^{\varepsilon}_{a_\varepsilon t_{\varepsilon,\mathfrak{m}}} *_\varepsilon \tilde{h}^{\varepsilon}(\tau_{\gamma,\mathfrak{m}}-))(x) + \varepsilon \sum_{y\in\Lambda_\varepsilon} \int_{\tau_{\varepsilon,\mathfrak{m}}}^{t} G^{\varepsilon}_{a_\varepsilon(t-s)}(x-y)\, d\widetilde{M}^{\varepsilon,\mathfrak{m}}_{\text{sym}}(s,y), \tag{B.3}$$

where we use $a_\varepsilon = 1 - \varrho_\varepsilon^2$ and the discrete heat kernel $G^\varepsilon$ as in (2.14). This yields

$$\begin{aligned}\iota_\varepsilon(\nabla^+_\varepsilon \tilde{h}^{\varepsilon,\mathfrak{m}}_{\text{sym}})(t, \nabla^+_\varepsilon \varphi^\lambda_x) &= \varepsilon \sum_{y\in\Lambda_\varepsilon} (\nabla^+_\varepsilon \varphi^\lambda_x *_\varepsilon \nabla^+_\varepsilon G^\varepsilon_{a_\varepsilon t_{\varepsilon,\mathfrak{m}}})(y) \tilde{h}^\varepsilon(\tau_{\gamma,\mathfrak{m}}-, y) \\ &\quad + \varepsilon \sum_{y\in\Lambda_\varepsilon} \int_{\tau_{\varepsilon,\mathfrak{m}}}^{t} (\nabla^+_\varepsilon \varphi^\lambda_x *_\varepsilon \nabla^+_\varepsilon G^\varepsilon_{a_\varepsilon(t-s)})(y)\, d\widetilde{M}^{\varepsilon,\mathfrak{m}}_{\text{sym}}(s,y).\end{aligned}$$

We denote the two terms on the right-hand side by $A_{\varepsilon,\lambda}(t)$ and $B_{\varepsilon,\lambda}(t)$. Minkowski's inequality allows to bound the expression on the left side of (B.1) by

$$\Big(\mathbb{E}\Big|\sup_{t\in[\tau_{\varepsilon,\mathfrak{m}},T]} (\lambda \vee \varepsilon \vee \sqrt{t_{\varepsilon,\mathfrak{m}}})^2 |A_{\varepsilon,\lambda}(t)|\Big|^p\Big)^{\frac{1}{p}} + \Big(\mathbb{E}\Big|\sup_{t\in[\tau_{\varepsilon,\mathfrak{m}},T]} (\lambda \vee \varepsilon \vee \sqrt{t_{\varepsilon,\mathfrak{m}}})^2 |B_{\varepsilon,\lambda}(t)|\Big|^p\Big)^{\frac{1}{p}}, \tag{B.4}$$

and we are going to estimate these two terms separately.

Let us start with bounding the first term in (B.4). We note that we can write

$$A_{\varepsilon,\lambda}(t) = \varepsilon \sum_{y\in\Lambda_\varepsilon} (\nabla^+_\varepsilon \varphi^\lambda_x *_\varepsilon \nabla^+_\varepsilon G^\varepsilon_{a_\varepsilon t_{\varepsilon,\mathfrak{m}}})(y)(\tilde{h}^\varepsilon(\tau_{\gamma,\mathfrak{m}}-, y) - \tilde{h}^\varepsilon(\tau_{\gamma,\mathfrak{m}}-, x)),$$

where the extra term $\tilde{h}^\varepsilon(\tau_{\gamma,\mathfrak{m}}-, x)$ is independent of $y$, what makes the sum over $y$ vanish. The definition of the stopping time (2.30) yields $\|\tilde{h}^\varepsilon(\tau_{\gamma,\mathfrak{m}}-)\|^{(\varepsilon)}_{\mathcal{C}^\alpha} \leq \mathfrak{m}$. Then we get

$$|A_{\varepsilon,\lambda}(t)| \leq \mathfrak{m}\varepsilon \sum_{y\in\Lambda_\varepsilon} |(\nabla^+_\varepsilon \varphi^\lambda_x *_\varepsilon \nabla^+_\varepsilon G^\varepsilon_{a_\varepsilon t_{\varepsilon,\mathfrak{m}}})(y)||y-x|^\alpha. \tag{B.5}$$

Lemma A.8 allows to write $G^\varepsilon_{a_\varepsilon t_{\varepsilon,\mathfrak{m}}}(x) = \frac{1}{t^{1/2}_{\varepsilon,\mathfrak{m}}+\varepsilon} \bar{G}^\varepsilon(\frac{x}{t^{1/2}_{\varepsilon,\mathfrak{m}}+\varepsilon})$, where $\bar{G}^\varepsilon$ is a Schwartz function. Then we have a convolution in (B.5) of two scaling functions, centered at $x$: one is rescaled by $\lambda$ and the other one is rescaled by $t^{1/2}_{\varepsilon,\mathfrak{m}} + \varepsilon$. Hence, we get

$$|A_{\varepsilon,\lambda}(t)| \lesssim \mathfrak{m}(\lambda \vee (t^{1/2}_{\varepsilon,\mathfrak{m}} + \varepsilon))^{\alpha-2} \lesssim \mathfrak{m}(\lambda \vee \varepsilon \vee t^{1/2}_{\varepsilon,\mathfrak{m}})^{\alpha-2}. \tag{B.6}$$

To define the second term in (B.4), we define

$$B_{\varepsilon,\lambda}(t,t') = \varepsilon \sum_{y\in\Lambda_\varepsilon} \int_{\tau_{\varepsilon,\mathfrak{m}}}^{t'} (\nabla^+_\varepsilon \varphi^\lambda_x *_\varepsilon \nabla^+_\varepsilon G^\varepsilon_{a_\varepsilon(t-s)})(y)\, d\widetilde{M}^{\varepsilon,\mathfrak{m}}_{\text{sym}}(s,y),$$

which is a martingale with respect to the variable $t' \in [\tau_{\varepsilon,\mathfrak{m}}, t]$. In order to apply the Burkholder-Davis-Gundy inequality [GMW24, Prop. 2.1] to this martingale, we need to bound its jumps and bracket process. The jump times coincide with those of the martingales $\widetilde{M}^{\varepsilon,\mathfrak{m}}_{\text{sym}}$. Let $\Delta_s f = f(s) - f(s-)$ denote the jump size of a function $f$ at time $s$. Then

$$\Delta_{t'} B_{\varepsilon,\lambda}(t,\bullet) = \varepsilon \sum_{y\in\Lambda_\varepsilon} (\nabla^+_\varepsilon \varphi^\lambda_x *_\varepsilon \nabla^+_\varepsilon G^\varepsilon_{a_\varepsilon(t-t')})(y) \Delta_{t'} \widetilde{M}^{\varepsilon,\mathfrak{m}}_{\text{sym}}(\bullet, y).$$

The jump size of $\widetilde{M}^{\varepsilon,\mathfrak{m}}_{\text{sym}}$ is proportional to $\sqrt{\varepsilon}$, and we get

$$|\Delta_{t'} B_{\varepsilon,\lambda}(t,\bullet)| \lesssim \varepsilon^{\frac{3}{2}} \sum_{y\in\Lambda_\varepsilon} |(\nabla^+_\varepsilon \varphi^\lambda_x *_\varepsilon \nabla^+_\varepsilon G^\varepsilon_{a_\varepsilon(t-t')})(y)| \lesssim \sqrt{\varepsilon}(\lambda \vee \varepsilon \vee (t-t')^{1/2})^{-2},$$

where we used the preceding properties of the involved functions, as in (B.6). Furthermore, we have

$$\begin{aligned}\langle B_{\varepsilon,\lambda}(t,\bullet)\rangle_{t'} &= \varepsilon \sum_{y\in\Lambda_\varepsilon} \int_{\tau_{\varepsilon,\mathfrak{m}}}^{t'} (\nabla_\varepsilon^+ \varphi_x^\lambda *_\varepsilon \nabla_\varepsilon^+ G^\varepsilon_{a_\varepsilon(t-s)})(y)^2 \widetilde{\mathbf{C}}_{\varepsilon,\mathfrak{m}}(s,y)\, ds \\ &\lesssim \varepsilon \sum_{y\in\Lambda_\varepsilon} \int_{\tau_{\varepsilon,\mathfrak{m}}}^{t'} (\nabla_\varepsilon^+ \varphi_x^\lambda *_\varepsilon \nabla_\varepsilon^+ G^\varepsilon_{a_\varepsilon(t-s)})(y)^2\, ds,\end{aligned}$$

where we use the function (2.32), which is bounded uniformly in $\varepsilon$ and its arguments. Using again the properties of the involved functions, as above, we get

$$\langle B_{\varepsilon,\lambda}(t,\bullet)\rangle_{t'} \lesssim \int_{\tau_{\varepsilon,\mathfrak{m}}}^{t'} (\lambda \vee \varepsilon \vee (t-s)^{1/2})^{-5} ds \lesssim (\lambda\vee\varepsilon\vee t_{\varepsilon,\mathfrak{m}}^{1/2})^{-3}.$$

Applying the Burkholder-Davis-Gundy inequality [GMW24, Prop. 2.1] and using the just derived bounds, we get

$$\Big(\mathbb{E}\Big[\Big|\sup_{t\in[\tau_{\varepsilon,\mathfrak{m}},T]} (\lambda\vee\varepsilon\vee\sqrt{t_{\varepsilon,\mathfrak{m}}})^{2-\alpha} |B_{\varepsilon,\lambda}(t)|\Big|^p\Big]\Big)^{\frac{1}{p}} \lesssim 1, \tag{B.7}$$

where the proportionality constant depends on $T$. Applying (B.6) and (B.7), we bound (B.4) by $C\mathfrak{m}$, as desired. □

We also need to control the product appearing in the function (2.32).

**Lemma B.2** *There is $\varepsilon_0 \in (0,1)$ such that for any $p\geq 1$ and $T>0$ there is a constant $C>0$ such that the bound*

$$\Big(\mathbb{E}\Big|\sup_{t\in[\tau_{\varepsilon,\mathfrak{m}},T]} \frac{1}{\mathbf{w}_\lambda^{\varepsilon,\mathfrak{m}}(t)} |\iota_\varepsilon(\nabla_\varepsilon^- \tilde h^{\varepsilon,\mathfrak{m}}_{\mathrm{sym}} \nabla_\varepsilon^+ \tilde h^{\varepsilon,\mathfrak{m}}_{\mathrm{sym}})(t,\varphi_x^\lambda)|\Big|^p\Big)^{\frac{1}{p}} \leq C(|\ln\varepsilon| + \mathfrak{m}^2) \tag{B.8}$$

*holds uniformly over $t\in[0,T]$, $\varphi\in\mathcal{B}^2$, $x\in\Lambda_\varepsilon$, $\lambda\in(0,1]$ and $\varepsilon\in(0,\varepsilon_0)$, where the random weight is*

$$\mathbf{w}_\lambda^{\varepsilon,\mathfrak{m}}(t) = (\sqrt{t-\tau_{\varepsilon,\mathfrak{m}}}+\varepsilon)^{\alpha-1}(\lambda\wedge\sqrt{t-\tau_{\varepsilon,\mathfrak{m}}}+\varepsilon)^{-1}. \tag{B.9}$$

*The constant $\alpha$ is the same as in the statement of Theorem 1.1.*

*Proof.* In the setting of (B.3) we define $\mathcal{N}^\varepsilon(t,x) = (G^\varepsilon_{a_\varepsilon t_{\varepsilon,\mathfrak{m}}} *_\varepsilon \tilde h^\varepsilon(\tau_{\gamma,\mathfrak{m}}-))(x)$ and

$$\mathcal{M}^\varepsilon_{t'}(t,x) = \varepsilon\sum_{y\in\Lambda_\varepsilon}\int_{\tau_{\varepsilon,\mathfrak{m}}}^{t'} G^\varepsilon_{a_\varepsilon(t-s)}(x-y)\, d\widetilde{M}^{\varepsilon,\mathfrak{m}}_{\mathrm{sym}}(s,y), \tag{B.10}$$

which is a martingale with respect to $t'\in[\tau_{\varepsilon,\mathfrak{m}},t]$. W recall that $t_{\varepsilon,\mathfrak{m}} = t-\tau_{\varepsilon,\mathfrak{m}}$ and $a_\varepsilon$ is a positive constant which is bounded uniformly in $\varepsilon$. Then (B.3) gives $\tilde h^{\varepsilon,\mathfrak{m}}_{\mathrm{sym}}(t,x) = \mathcal{N}^\varepsilon(t,x) + \mathcal{M}^\varepsilon_t(t,x)$ and

$$\begin{aligned}(\nabla_\varepsilon^- \tilde h^{\varepsilon,\mathfrak{m}}_{\mathrm{sym}} \nabla_\varepsilon^+ \tilde h^{\varepsilon,\mathfrak{m}}_{\mathrm{sym}})(t,x) &= (\nabla_\varepsilon^-\mathcal{N}^\varepsilon\nabla_\varepsilon^+\mathcal{N}^\varepsilon)(t,x) + (\nabla_\varepsilon^-\mathcal{N}^\varepsilon\nabla_\varepsilon^+\mathcal{M}^\varepsilon_t)(t,x) \\ &\quad + (\nabla_\varepsilon^-\mathcal{M}^\varepsilon_t\nabla_\varepsilon^+\mathcal{N}^\varepsilon)(t,x) + (\nabla_\varepsilon^-\mathcal{M}^\varepsilon_t\nabla_\varepsilon^+\mathcal{M}^\varepsilon_t)(t,x).\end{aligned} \tag{B.11}$$

By Minkowski's inequality we will prove (B.8) if we show that the same bounds hold for each of these four terms. We are going to bound these terms one-by-one.

As in the proof of the preceding lemma, we write $G^\varepsilon_{a_\varepsilon t_{\varepsilon,\mathfrak{m}}}(x) = \frac{1}{t_{\varepsilon,\mathfrak{m}}^{1/2}+\varepsilon}\bar G^\varepsilon(\frac{x}{t_{\varepsilon,\mathfrak{m}}^{1/2}+\varepsilon})$, where $\bar G^\varepsilon$ is a Schwartz function. Then we can write

$$\nabla_\varepsilon^\pm \mathcal{N}^\varepsilon(t,x) = \varepsilon\sum_{y\in\Lambda_\varepsilon} \nabla_\varepsilon^\pm G^\varepsilon_{a_\varepsilon t_{\varepsilon,\mathfrak{m}}}(x-y)(\tilde h^\varepsilon(\tau_{\gamma,\mathfrak{m}}-,y) - \tilde h^\varepsilon(\tau_{\gamma,\mathfrak{m}}-,x)),$$

where we can subtract $\tilde{h}^\varepsilon(\tau_{\gamma,\mathfrak{m}}-,x)$ because the sum over $y$ vanishes. The definition of the stopping time (2.30) yields $\|\tilde{h}^\varepsilon(\tau_{\gamma,\mathfrak{m}}-)\|^{(\varepsilon)}_{\mathcal{C}^\alpha} \le \mathfrak{m}$, and we get

$$|\nabla^\pm_\varepsilon \mathcal{N}^\varepsilon(t,x)| \le \mathfrak{m}\varepsilon \sum_{y\in\Lambda_\varepsilon} |\nabla^\pm_\varepsilon G^\varepsilon_{a_\varepsilon t_{\varepsilon,\mathfrak{m}}}(x-y)||y-x|^\alpha \lesssim \mathfrak{m}(t^{1/2}_{\varepsilon,\mathfrak{m}}+\varepsilon)^{\alpha-1}, \tag{B.12}$$

where we used the scaling property of the discrete heat kernel. One can see that the same argument gives the bound

$$|(\nabla^-_\varepsilon)^2 \mathcal{N}^\varepsilon(t,x)| \le \mathfrak{m}\varepsilon \sum_{y\in\Lambda_\varepsilon} |(\nabla^-_\varepsilon)^2 G^\varepsilon_{a_\varepsilon t_{\varepsilon,\mathfrak{m}}}(x-y)||y-x|^\alpha \lesssim \mathfrak{m}(t^{1/2}_{\varepsilon,\mathfrak{m}}+\varepsilon)^{\alpha-2}, \tag{B.13}$$

which we are going to use later. Then the first term in (B.11) satisfies

$$|(\nabla^-_\varepsilon \mathcal{N}^\varepsilon \nabla^+_\varepsilon \mathcal{N}^\varepsilon)(t,x)| \lesssim \mathfrak{m}^2(t^{1/2}_{\varepsilon,\mathfrak{m}}+\varepsilon)^{2(\alpha-1)} \lesssim \mathfrak{m}^2 \mathbf{w}^{\varepsilon,\mathfrak{m}}_\lambda(t),$$

where we use the weight (B.9), and we get

$$\Big(\mathbb{E}\Big|\sup_{t\in[\tau_{\varepsilon,\mathfrak{m}},T]} \frac{1}{\mathbf{w}^{\varepsilon,\mathfrak{m}}_\lambda(t)} |(\nabla^-_\varepsilon \mathcal{N}^\varepsilon \nabla^+_\varepsilon \mathcal{N}^\varepsilon)(t,\varphi^\lambda_x)|\Big|^p\Big)^{\frac{1}{p}} \lesssim \mathfrak{m}^2. \tag{B.14}$$

To bound the second term in (B.11) we note that the process $t' \mapsto \iota_\varepsilon(\nabla^-_\varepsilon \mathcal{N}^\varepsilon \nabla^+_\varepsilon \mathcal{M}^\varepsilon_{t'})(t,\varphi^\lambda_x)$ is a martingale on $t' \in [\tau_{\varepsilon,\mathfrak{m}}, t]$. This is because $\nabla^-_\varepsilon \mathcal{N}^\varepsilon(t,y)$ is measurable with respect to the $\sigma$-algebra $\mathscr{F}_{\tau_{\varepsilon,\mathfrak{m}}}$. To bound this martingale, we use Burkholder-Davis-Gundy inequality. For this, we need to bound the jumps and the bracket process of the martingale. The jump is (we use the operator $\Delta_{t'}$ introduced in the previous lemma)

$$\begin{aligned}\Delta_{t'}\iota_\varepsilon(\nabla^-_\varepsilon \mathcal{N}^\varepsilon \nabla^+_\varepsilon \mathcal{M}^\varepsilon_\bullet)(t,\varphi^\lambda_x) &= \iota_\varepsilon(\nabla^-_\varepsilon \mathcal{N}^\varepsilon \Delta_{t'} \nabla^+_\varepsilon \mathcal{M}^\varepsilon_\bullet)(t,\varphi^\lambda_x)\\ &= \varepsilon^2 \sum_{y,y'\in\Lambda_\varepsilon} \varphi^\lambda_x(y)\nabla^-_\varepsilon \mathcal{N}^\varepsilon(t,y)\nabla^+_\varepsilon G^\varepsilon_{a_\varepsilon(t-t')}(y-y')\, \Delta_{t'}\widetilde{M}^{\varepsilon,\mathfrak{m}}_{\text{sym}}(\bullet,y'),\end{aligned}$$

where we used (B.10). Applying summation by parts, we get

$$-\varepsilon^2 \sum_{y,y'\in\Lambda_\varepsilon} \nabla^-_\varepsilon(\varphi^\lambda_x(\bullet)\nabla^-_\varepsilon \mathcal{N}^\varepsilon(t,\bullet))(y) G^\varepsilon_{a_\varepsilon(t-t')}(y-y')\, \Delta_{t'}\widetilde{M}^{\varepsilon,\mathfrak{m}}_{\text{sym}}(\bullet,y').$$

The jump size of $\widetilde{M}^{\varepsilon,\mathfrak{m}}_{\text{sym}}$ is of order $\sqrt{\varepsilon}$ and we get

$$\begin{aligned}|\Delta_{t'}\iota_\varepsilon(\nabla^-_\varepsilon \mathcal{N}^\varepsilon \nabla^+_\varepsilon \mathcal{M}^\varepsilon_\bullet)(t,\varphi^\lambda_x)| &\lesssim \varepsilon^{\frac{5}{2}} \sum_{y,y'\in\Lambda_\varepsilon} |\nabla^-_\varepsilon(\varphi^\lambda_x(\bullet)\nabla^-_\varepsilon \mathcal{N}^\varepsilon(t,\bullet))(y) G^\varepsilon_{a_\varepsilon(t-t')}(y-y')|\\ &\lesssim \varepsilon^{\frac{3}{2}} \sum_{y\in\Lambda_\varepsilon} |\nabla^-_\varepsilon(\varphi^\lambda_x(\bullet)\nabla^-_\varepsilon \mathcal{N}^\varepsilon(t,\bullet))(y)|,\end{aligned} \tag{B.15}$$

where we used the fact that the discrete heat kernel sums up to $1$ in the spatial variable. Furthermore,

$$|\nabla^-_\varepsilon(\varphi^\lambda_x(\bullet)\nabla^-_\varepsilon \mathcal{N}^\varepsilon(t,\bullet))(y)| \le |\nabla^-_\varepsilon \varphi^\lambda_x(y)\nabla^-_\varepsilon \mathcal{N}^\varepsilon(t,y)| + |\varphi^\lambda_x(y-\varepsilon)(\nabla^-_\varepsilon)^2 \mathcal{N}^\varepsilon(t,y)|.$$

Applying (B.12) and (B.13) we get

$$|\nabla^-_\varepsilon(\varphi^\lambda_x(\bullet)\nabla^-_\varepsilon \mathcal{N}^\varepsilon(t,\bullet))(y)| \le \mathfrak{m}(t^{1/2}_{\varepsilon,\mathfrak{m}}+\varepsilon)^{\alpha-1}((\lambda\vee\varepsilon)^{-1} + (t^{1/2}_{\varepsilon,\mathfrak{m}}+\varepsilon)^{-1})\tilde{\varphi}^{\varepsilon,\lambda}_x(y), \tag{B.16}$$

where $\tilde{\varphi}^{\varepsilon,\lambda}_x$ is a compactly supported, positive, continuous function, centered at $x$ and rescaled by $\lambda$. This function is not differentiable and is uniformly bounded by $c\lambda^{-1}$ for some constant $c > 0$. Then (B.15) is bounded as

$$|\Delta_{t'}\iota_\varepsilon(\nabla^-_\varepsilon \mathcal{N}^\varepsilon \nabla^+_\varepsilon \mathcal{M}^\varepsilon_\bullet)(t,\varphi^\lambda_x)| \lesssim \varepsilon^{\frac{1}{2}}\mathfrak{m}(t^{1/2}_{\varepsilon,\mathfrak{m}}+\varepsilon)^{\alpha-1}((\lambda\vee\varepsilon)^{-1} + (t^{1/2}_{\varepsilon,\mathfrak{m}}+\varepsilon)^{-1}). \tag{B.17}$$

Now, we will derive a bound on the bracket process of the martingale $t' \mapsto \iota_\varepsilon(\nabla_\varepsilon^- \mathcal{N}^\varepsilon \nabla_\varepsilon^+ \mathcal{M}_{t'}^\varepsilon)(t, \varphi_x^\lambda)$. Using (B.10) and (2.32) we get

$$\begin{aligned}
&\langle \iota_\varepsilon(\nabla_\varepsilon^- \mathcal{N}^\varepsilon \nabla_\varepsilon^+ \mathcal{M}_\bullet^\varepsilon)(t, \varphi_x^\lambda)\rangle_t \\
&\quad = \varepsilon^2 \sum_{y' \in \Lambda_\varepsilon} \int_{\tau_{\varepsilon,\mathfrak{m}}}^t \Big(\varepsilon \sum_{y \in \Lambda_\varepsilon} \varphi_x^\lambda(y) \nabla_\varepsilon^- \mathcal{N}^\varepsilon(t, y) \nabla_\varepsilon^+ G_{a_\varepsilon(t-s)}^\varepsilon(y - y')\Big)^2 d\langle \widetilde{M}^{\varepsilon,\mathfrak{m}}(\bullet, y')\rangle_s \\
&\quad \lesssim \int_{\tau_{\varepsilon,\mathfrak{m}}}^t \varepsilon \sum_{y' \in \Lambda_\varepsilon} \Big(\varepsilon \sum_{y \in \Lambda_\varepsilon} \varphi_x^\lambda(y) \nabla_\varepsilon^- \mathcal{N}^\varepsilon(t, y) \nabla_\varepsilon^+ G_{a_\varepsilon(t-s)}^\varepsilon(y - y')\Big)^2 ds. \qquad \text{(B.18)}
\end{aligned}$$

We apply summation by parts and write the expression inside the integral as

$$\begin{aligned}
&\varepsilon \sum_{y' \in \Lambda_\varepsilon} \Big(\varepsilon \sum_{y \in \Lambda_\varepsilon} \nabla_\varepsilon^- (\varphi_x^\lambda(\bullet) \nabla_\varepsilon^- \mathcal{N}^\varepsilon(t, \bullet))(y) G_{a_\varepsilon(t-s)}^\varepsilon(y - y')\Big)^2 \qquad \text{(B.19)} \\
&\quad = \varepsilon^2 \sum_{y_1, y_2 \in \Lambda_\varepsilon} \nabla_\varepsilon^- (\varphi_x^\lambda(\bullet) \nabla_\varepsilon^- \mathcal{N}^\varepsilon(t, \bullet))(y_1) \nabla_\varepsilon^- (\varphi_x^\lambda(\bullet) \nabla_\varepsilon^- \mathcal{N}^\varepsilon(t, \bullet))(y_2) G_{2a_\varepsilon(t-s)}^\varepsilon(y_1 - y_2),
\end{aligned}$$

where we have used the identity $G_{a_\varepsilon(t-s)}^\varepsilon *_\varepsilon G_{a_\varepsilon(t-s)}^\varepsilon = G_{2a_\varepsilon(t-s)}^\varepsilon$. Using (B.16), we estimate (B.19) by a constant times

$$\mathfrak{m}^2 (t_{\varepsilon,\mathfrak{m}}^{1/2} + \varepsilon)^{2(\alpha-1)} ((\lambda \vee \varepsilon)^{-1} + (t_{\varepsilon,\mathfrak{m}}^{1/2} + \varepsilon)^{-1})^2 \varepsilon^2 \sum_{y_1, y_2 \in \Lambda_\varepsilon} \tilde{\varphi}_x^{\varepsilon,\lambda}(y_1) \tilde{\varphi}_x^{\varepsilon,\lambda}(y_2) G_{2a_\varepsilon(t-s)}^\varepsilon(y_1 - y_2).$$

The scaling property of the discrete heat kernel yields

$$\varepsilon^2 \sum_{y_1, y_2 \in \Lambda_\varepsilon} \tilde{\varphi}_x^{\varepsilon,\lambda}(y_1) \tilde{\varphi}_x^{\varepsilon,\lambda}(y_2) G_{2a_\varepsilon(t-s)}^\varepsilon(y_1 - y_2) \lesssim 1,$$

and hence (B.19) is bounded by a constant times $\mathfrak{m}^2 (t_{\varepsilon,\mathfrak{m}}^{1/2} + \varepsilon)^{2(\alpha-1)} ((\lambda \vee \varepsilon)^{-1} + (t_{\varepsilon,\mathfrak{m}}^{1/2} + \varepsilon)^{-1})^2$. Then the bracket process (B.18) is bounded as

$$\langle \iota_\varepsilon(\nabla_\varepsilon^- \mathcal{N}^\varepsilon \nabla_\varepsilon^+ \mathcal{M}_\bullet^\varepsilon)(t, \varphi_x^\lambda)\rangle_t \lesssim \mathfrak{m}^2 t_{\varepsilon,\mathfrak{m}} (t_{\varepsilon,\mathfrak{m}}^{1/2} + \varepsilon)^{2(\alpha-1)} ((\lambda \vee \varepsilon)^{-1} + (t_{\varepsilon,\mathfrak{m}}^{1/2} + \varepsilon)^{-1})^2. \qquad \text{(B.20)}$$

Applying the Burkholder-Davis-Gundy inequality [GMW24, Prop. 2.1] to the martingale $t' \mapsto \iota_\varepsilon(\nabla_\varepsilon^- \mathcal{N}^\varepsilon \nabla_\varepsilon^+ \mathcal{M}_{t'}^\varepsilon)(t, \varphi_x^\lambda)$, and using (B.17) and (B.20) we get

$$\begin{aligned}
&\Big(\mathbb{E}\Big|\sup_{t \in [\tau_{\varepsilon,\mathfrak{m}}, T]} \frac{1}{\mathbf{w}_\lambda^{\varepsilon,\mathfrak{m}}(t)} |\iota_\varepsilon(\nabla_\varepsilon^- \mathcal{N}^\varepsilon \nabla_\varepsilon^+ \mathcal{M}_t^\varepsilon)(t, \varphi_x^\lambda)|\Big|^p\Big)^{\frac{1}{p}} \\
&\quad \le \Big(\mathbb{E}\Big[\sup_{t \in [\tau_{\varepsilon,\mathfrak{m}}, T]} \frac{1}{\mathbf{w}_\lambda^{\varepsilon,\mathfrak{m}}(t)^p} \mathbb{E}\Big[\sup_{t' \in [\tau_{\varepsilon,\mathfrak{m}}, t]} |\iota_\varepsilon(\nabla_\varepsilon^- \mathcal{N}^\varepsilon \nabla_\varepsilon^+ \mathcal{M}_{t'}^\varepsilon)(t, \varphi_x^\lambda)|^p \Big| \mathscr{F}_{\tau_{\varepsilon,\mathfrak{m}}}\Big]\Big]\Big)^{\frac{1}{p}} \\
&\quad \lesssim \Big(\mathbb{E}\Big[\sup_{t \in [\tau_{\varepsilon,\mathfrak{m}}, T]} \frac{1}{\mathbf{w}_\lambda^{\varepsilon,\mathfrak{m}}(t)^p} \mathfrak{m}^p (t_{\varepsilon,\mathfrak{m}}^{1/2} + \varepsilon)^{(\alpha-1)p} ((\lambda \vee \varepsilon)^{-1} + (t_{\varepsilon,\mathfrak{m}}^{1/2} + \varepsilon)^{-1})^p\Big]\Big)^{\frac{1}{p}} \lesssim \mathfrak{m}. \qquad \text{(B.21)}
\end{aligned}$$

This bound is the reason to choose the weight (B.9) to be so complicated. The same bound holds on the third term in (B.11).

It is left to bound the last term in (B.11). Let us define the functions

$$F^\varepsilon((t, y); z_1, z_2) = \nabla_\varepsilon^- G_{a_\varepsilon(t-t_1)}^\varepsilon(y - x_1) \nabla_\varepsilon^+ G_{a_\varepsilon(t-t_2)}^\varepsilon(y - x_2),$$

where $z_i = (t_i, x_i)$, and $\bar{F}_{t-s}^\varepsilon(y - y_1) = F^\varepsilon((t, y); (s, y_1), (s, y_1))$. Recalling the definition (B.10) and applying [GMW24, Eq. 2.16], we write this term as

$$\iota_\varepsilon(\nabla_\varepsilon^- \mathcal{M}_t^\varepsilon \nabla_\varepsilon^+ \mathcal{M}_t^\varepsilon)(t, \varphi_x^\lambda) = \mathcal{I}_2^\varepsilon((\iota_\varepsilon F^\varepsilon)(t, \varphi_x^\lambda)) + \varepsilon^2 \sum_{y \in \Lambda_\varepsilon} \int_{\tau_{\varepsilon,\mathfrak{m}}}^t (\varphi_x^\lambda *_\varepsilon \bar{F}_{t-s}^\varepsilon)(y)\, d[\widetilde{M}_{\mathrm{sym}}^{\varepsilon,\mathfrak{m}}(y)]_s, \qquad \text{(B.22)}$$

where $\mathcal{I}_2^\varepsilon$ is a second-order stochastic integral with respect to the martingales $\widetilde{M}^{\varepsilon,\mathfrak{m}}_{\text{sym}}$ and $[\widetilde{M}^{\varepsilon,\mathfrak{m}}_{\text{sym}}(y)]_s$ is the quadratic variation of the martingale.

Let us bound the first term in (B.22). We note that $|\nabla_\varepsilon^\pm G^\varepsilon_{a_\varepsilon(t-s)}(y)| \lesssim (\sqrt{t-s}+\varepsilon)^{-2}$ which follows from the scaling property of the kernel. Using this bound to estimate the $L^\infty$ norms of the discrete heat kernel and applying [GMW24, Thm. 3.5] we get

$$
\begin{aligned}
&\Big(\mathbb{E}\Big[\sup_{t\in[\tau_{\varepsilon,\mathfrak{m}},T]} \frac{1}{\mathbf{w}_\lambda^{\varepsilon,\mathfrak{m}}(t)^p}|\mathcal{I}_2^\varepsilon((\iota_\varepsilon F^\varepsilon)(t,\varphi_x^\lambda))|^p\Big|\mathscr{F}_{\tau_{\varepsilon,\mathfrak{m}}}\Big]\Big)^{\frac{1}{p}}\\
&\quad\lesssim \left(\int_{\tau_{\varepsilon,\mathfrak{m}}}^T\int_{\tau_{\varepsilon,\mathfrak{m}}}^T \frac{1}{\mathbf{w}_\lambda^{\varepsilon,\mathfrak{m}}(T)^2}\varepsilon^2\sum_{y_1,y_2\in\Lambda_\varepsilon}(\iota_\varepsilon F^\varepsilon)(T,\varphi_x^\lambda;(s_1,y_1),(s_2,y_2))^2\,ds_1ds_2\right)^{\frac12}\\
&\qquad+\varepsilon^{\frac32}\left(\int_{\tau_{\varepsilon,\mathfrak{m}}}^T \frac{(\sqrt{T-s}+\varepsilon)^{-4}}{\mathbf{w}_\lambda^{\varepsilon,\mathfrak{m}}(T)^2}\varepsilon\sum_{y\in\Lambda_\varepsilon}(|\nabla_\varepsilon^-\varphi_x^\lambda|*_\varepsilon|G^\varepsilon_{a_\varepsilon(T-s)}|)(y)^2\,ds\right)^{\frac12} \qquad \text{(B.23)}\\
&\qquad+\varepsilon^{\frac32}\left(\int_{\tau_{\varepsilon,\mathfrak{m}}}^T \frac{(\sqrt{T-s}+\varepsilon)^{-4}}{\mathbf{w}_\lambda^{\varepsilon,\mathfrak{m}}(T)^2}\varepsilon\sum_{y\in\Lambda_\varepsilon}(|\nabla_\varepsilon^+\varphi_x^\lambda|*_\varepsilon|G^\varepsilon_{a_\varepsilon(T-s)}|)(y)^2\,ds\right)^{\frac12}\\
&\qquad+\varepsilon^3\sup_{s\in[\tau_{\varepsilon,\mathfrak{m}},T]}\sup_{y\in\Lambda_\varepsilon}\frac{(\sqrt{T-s}+\varepsilon)^{-2}}{\mathbf{w}_\lambda^{\varepsilon,\mathfrak{m}}(T)}(|\nabla_\varepsilon^-\varphi_x^\lambda|*_\varepsilon|G^\varepsilon_{a_\varepsilon(T-s)}|)(y).
\end{aligned}
$$

We can write

$$
\begin{aligned}
&\varepsilon^2\sum_{y_1,y_2\in\Lambda_\varepsilon}(\iota_\varepsilon F^\varepsilon)(t,\varphi_x^\lambda;(s_1,y_1),(s_2,y_2))^2\\
&\quad=\varepsilon^2\sum_{y_1,y_2\in\Lambda_\varepsilon}\Big(\varepsilon\sum_{y\in\Lambda_\varepsilon}\varphi_x^\lambda(y)\nabla_\varepsilon^-G^\varepsilon_{a_\varepsilon(t-s_1)}(y-y_1)\nabla_\varepsilon^+G^\varepsilon_{a_\varepsilon(t-s_2)}(y-y_2)\Big)^2\\
&\quad=\varepsilon^2\sum_{y,y'\in\Lambda_\varepsilon}\varphi_x^\lambda(y)\varphi_x^\lambda(y')\Big(\varepsilon\sum_{y_1\in\Lambda_\varepsilon}\nabla_\varepsilon^-G^\varepsilon_{a_\varepsilon(t-s_1)}(y-y_1)\nabla_\varepsilon^-G^\varepsilon_{a_\varepsilon(t-s_1)}(y'-y_1)\Big)\\
&\qquad\qquad\times\Big(\varepsilon\sum_{y_2\in\Lambda_\varepsilon}\nabla_\varepsilon^+G^\varepsilon_{a_\varepsilon(t-s_2)}(y-y_2)\nabla_\varepsilon^+G^\varepsilon_{a_\varepsilon(t-s_2)}(y'-y_2)\Big).
\end{aligned}
$$

The scaling property of the discrete heat kernel implies that the two expressions in the parentheses can be written as $(\frac{1}{(T-s_1)^{1/2}+\varepsilon})^3\psi^\varepsilon(\frac{y-y'}{(T-s_1)^{1/2}+\varepsilon})$ and $(\frac{1}{(T-s_2)^{1/2}+\varepsilon})^3\widetilde{\psi}^\varepsilon(\frac{y-y'}{(T-s_2)^{1/2}+\varepsilon})$ respectively, where the functions $\psi^\varepsilon$ and $\widetilde{\psi}^\varepsilon$ are Schwartz and bounded uniformly in $\varepsilon$. Then the whole preceding expression is bounded by a constant times $(\frac{1}{(T-s_1)^{1/2}+\varepsilon})^2(\frac{1}{(T-s_2)^{1/2}+\varepsilon})^2$. Hence, the first term on the right-hand side of (B.23) is bounded by a constant multiple of

$$
\left(\int_{\tau_{\varepsilon,\mathfrak{m}}}^T\int_{\tau_{\varepsilon,\mathfrak{m}}}^T\frac{1}{\mathbf{w}_\lambda^{\varepsilon,\mathfrak{m}}(T)^2}\Big(\frac{1}{(T-s_1)^{1/2}+\varepsilon}\Big)^2\Big(\frac{1}{(T-s_2)^{1/2}+\varepsilon}\Big)^2ds_1ds_2\right)^{\frac12}\lesssim|\ln\varepsilon|.
$$

Similarly, we have $\varepsilon\sum_{y\in\Lambda_\varepsilon}(|\nabla_\varepsilon^\pm\varphi_x^\lambda|*_\varepsilon|G^\varepsilon_{a_\varepsilon(T-s)}|)(y)^2\lesssim(\lambda\vee\varepsilon)^{-2}$. Using then $\varepsilon\le\sqrt{T-s}+\varepsilon$, we estimate the second and third terms on the right-hand side of (B.23) by a constant times

$$
\left(\int_{\tau_{\varepsilon,\mathfrak{m}}}^T\frac{(\sqrt{T-s}+\varepsilon)^{-1}}{\mathbf{w}_\lambda^{\varepsilon,\mathfrak{m}}(T)^2}(\lambda\vee\varepsilon)^{-2}\,ds\right)^{\frac12}\lesssim(\lambda\vee\varepsilon)^{-1}.
$$

Using the bounds $(|\nabla_\varepsilon^\pm\varphi_x^\lambda|*_\varepsilon|G^\varepsilon_{a_\varepsilon(T-s)}|)(y)\lesssim(\lambda\vee\varepsilon)^{-1}$ and $\varepsilon\le\sqrt{T-s}+\varepsilon$, we estimate the last term in (B.23) by a constant times $\varepsilon(\lambda\vee\varepsilon)^{-1}$.

We have just proved the bound

$$\Big(\mathbb{E}\Big[\sup_{t\in[\tau_{\varepsilon,\mathfrak{m}},T]}\frac{1}{\mathbf{w}_\lambda^{\varepsilon,\mathfrak{m}}(t)^p}|\mathcal{I}_2^\varepsilon((\iota_\varepsilon F^\varepsilon)(t,\varphi_x^\lambda))|^p\Big|\mathscr{F}_{\tau_{\varepsilon,\mathfrak{m}}}\Big]\Big)^{\frac{1}{p}}\lesssim|\ln\varepsilon|+(\lambda\vee\varepsilon)^{-1}$$

on the first term in (B.22), and we are going to bound the second term. We recall from Section 2 that the martingales $\bar{M}_{\text{sym}}^{\varepsilon,\mathfrak{m}}(t,x)=\varepsilon^{-\frac{1}{2}}([\widetilde{M}_{\text{sym}}^{\varepsilon,\mathfrak{m}}(x)]_t-\langle\widetilde{M}_{\text{sym}}^{\varepsilon,\mathfrak{m}}(x)\rangle_t)$ have essentially the same properties as the martingales $\widetilde{M}_{\text{sym}}^{\varepsilon,\mathfrak{m}}$. Then we write the second term in (B.22) as

$$\varepsilon^{\frac{5}{2}}\sum_{y\in\Lambda_\varepsilon}\int_{\tau_{\varepsilon,\mathfrak{m}}}^t(\varphi_x^\lambda*_\varepsilon\bar{F}_{t-s}^\varepsilon)(y)\,d\bar{M}_{\text{sym}}^{\varepsilon,\mathfrak{m}}(s,y)+\varepsilon^2\sum_{y\in\Lambda_\varepsilon}\int_{\tau_{\varepsilon,\mathfrak{m}}}^t(\varphi_x^\lambda*_\varepsilon\bar{F}_{t-s}^\varepsilon)(y)\,d\langle\widetilde{M}_{\text{sym}}^{\varepsilon,\mathfrak{m}}(y)\rangle_s,\quad\text{(B.24)}$$

and we denote the two terms by $A_\lambda^\varepsilon(t)$ and $B_\lambda^\varepsilon(t)$ respectively. The first term is a stochastic integral, and applying [GMW24, Thm. 3.5] we bound it as

$$\Big(\mathbb{E}\Big[\sup_{t\in[\tau_{\varepsilon,\mathfrak{m}},T]}\frac{1}{\mathbf{w}_\lambda^{\varepsilon,\mathfrak{m}}(t)^p}|A_\lambda^\varepsilon(t)|^p\Big|\mathscr{F}_{\tau_{\varepsilon,\mathfrak{m}}}\Big]\Big)^{\frac{1}{p}}\lesssim\varepsilon^{\frac{3}{2}}\left(\int_{\tau_{\varepsilon,\mathfrak{m}}}^T\frac{1}{\mathbf{w}_\lambda^{\varepsilon,\mathfrak{m}}(T)^2}\varepsilon\sum_{y\in\Lambda_\varepsilon}(\varphi_x^\lambda*_\varepsilon\bar{F}_{T-s}^\varepsilon)(y)^2\,ds\right)^{\frac{1}{2}}$$
$$+\varepsilon^3\sup_{s\in[\tau_{\varepsilon,\mathfrak{m}},T]}\sup_{y\in\Lambda_\varepsilon}\frac{1}{\mathbf{w}_\lambda^{\varepsilon,\mathfrak{m}}(T)}|(\varphi_x^\lambda*_\varepsilon\bar{F}_{T-s}^\varepsilon)(y)|.$$

As follows from the definition of the function $\bar{F}^\varepsilon$ and from the scaling property of the discrete heat kernel, we can write $\bar{F}_{T-s}^\varepsilon(y)=(\frac{1}{(T-s)^{1/2}+\varepsilon})^4\bar{\psi}^\varepsilon(\frac{y}{(T-s)^{1/2}+\varepsilon})$, where $\bar{\psi}^\varepsilon$ is Schwartz. Then we have $|(\varphi_x^\lambda*_\varepsilon\bar{F}_{T-s}^\varepsilon)(y)|\lesssim(\frac{1}{(T-s)^{1/2}+\varepsilon})^3$ and $\varepsilon\sum_{y\in\Lambda_\varepsilon}(\varphi_x^\lambda*_\varepsilon\bar{F}_{T-s}^\varepsilon)(y)^2\lesssim(\frac{1}{(T-s)^{1/2}+\varepsilon})^5$. Using $\varepsilon\le\sqrt{T-s}+\varepsilon$, we bound the preceding expression by a constant times

$$\varepsilon^{\frac{3}{2}}\left(\int_{\tau_{\varepsilon,\mathfrak{m}}}^T\frac{1}{\mathbf{w}_\lambda^{\varepsilon,\mathfrak{m}}(T)^2}\Big(\frac{1}{(T-s)^{1/2}+\varepsilon}\Big)^5ds\right)^{\frac{1}{2}}$$
$$+\varepsilon^3\sup_{s\in[\tau_{\varepsilon,\mathfrak{m}},T]}\sup_{y\in\Lambda_\varepsilon}\frac{1}{\mathbf{w}_\lambda^{\varepsilon,\mathfrak{m}}(T)}\Big(\frac{1}{(T-s)^{1/2}+\varepsilon}\Big)^3$$
$$\lesssim\left(\int_{\tau_{\varepsilon,\mathfrak{m}}}^T\frac{1}{\mathbf{w}_\lambda^{\varepsilon,\mathfrak{m}}(T)^2}\Big(\frac{1}{(T-s)^{1/2}+\varepsilon}\Big)^2ds\right)^{\frac{1}{2}}+1\lesssim|\ln\varepsilon|.$$

Hence, we have just proved

$$\Big(\mathbb{E}\Big[\sup_{t\in[\tau_{\varepsilon,\mathfrak{m}},T]}\frac{1}{\mathbf{w}_\lambda^{\varepsilon,\mathfrak{m}}(t)^p}|A_\lambda^\varepsilon(t)|^p\Big|\mathscr{F}_{\tau_{\varepsilon,\mathfrak{m}}}\Big]\Big)^{\frac{1}{p}}\le C_1|\ln\varepsilon|\quad\text{(B.25)}$$

for some constant $C_1>0$, and we will now analyse the second term in (B.24).

The formula (2.32) for the bracket process yields (recall the constant (2.8))

$$\begin{aligned}B_\lambda^\varepsilon(t)&=\nu^\varepsilon\int_{\tau_{\varepsilon,\mathfrak{m}}}^t\varepsilon\sum_{y\in\Lambda_\varepsilon}(\varphi_x^\lambda*_\varepsilon\bar{F}_{t-s}^\varepsilon)(y)(1-\varepsilon(\nabla_\varepsilon^-\tilde{h}_{\text{sym}}^{\varepsilon,\mathfrak{m}}\nabla_\varepsilon^+\tilde{h}_{\text{sym}}^{\varepsilon,\mathfrak{m}})(s,y))\,ds\\&=\nu^\varepsilon\int_{\tau_{\varepsilon,\mathfrak{m}}}^t\varepsilon\sum_{y\in\Lambda_\varepsilon}(\varphi_x^\lambda*_\varepsilon\bar{F}_{t-s}^\varepsilon)(y)\,ds\\&\quad-\nu^\varepsilon\int_{\tau_{\varepsilon,\mathfrak{m}}}^t\varepsilon^2\sum_{y\in\Lambda_\varepsilon}\bar{F}_{t-s}^\varepsilon(y)\iota_\varepsilon(\nabla_\varepsilon^-\tilde{h}_{\text{sym}}^{\varepsilon,\mathfrak{m}}\nabla_\varepsilon^+\tilde{h}_{\text{sym}}^{\varepsilon,\mathfrak{m}})(t,\varphi_{x+y}^\lambda)\,ds.\end{aligned}\quad\text{(B.26)}$$

We note that $\bar{F}^\varepsilon_t(x) = Q^\varepsilon_{a_\varepsilon t}(x)$, where the latter is defined in (2.16). Then Lemma 4.7 allows to bound the first term in (B.26) by a constant multiple of $(\lambda \vee \varepsilon)^{-1}$. Applying furthermore Lemma 4.5, we obtain

$$\Big(\mathbb{E}\Big[\sup_{t\in[\tau_{\varepsilon,\mathfrak{m}},T]} \frac{1}{\mathbf{w}^{\varepsilon,\mathfrak{m}}_\lambda(t)^p} |B^\varepsilon_\lambda(t)|^p \Big| \mathscr{F}_{\tau_{\varepsilon,\mathfrak{m}}}\Big]\Big)^{\frac{1}{p}} \leq C_2 \\ + (1+\sqrt{\varepsilon})\Theta \sup_{y\in\Lambda_\varepsilon} \Big(\mathbb{E}\Big[\sup_{t\in[\tau_{\varepsilon,\mathfrak{m}},T]} \frac{1}{\mathbf{w}^{\varepsilon,\mathfrak{m}}_\lambda(t)^p} |\iota_\varepsilon(\nabla^-_\varepsilon \tilde{h}^{\varepsilon,\mathfrak{m}}_{\text{sym}} \nabla^+_\varepsilon \tilde{h}^{\varepsilon,\mathfrak{m}}_{\text{sym}})(t,\varphi^\lambda_{x+y})|^p \Big| \mathscr{F}_{\tau_{\varepsilon,\mathfrak{m}}}\Big]\Big)^{\frac{1}{p}}, \tag{B.27}$$

for some constant $C_2 > 0$.

Let us now denote the quantity which we need to bound in (B.8) by

$$D^\varepsilon_\lambda = \sup_{x\in\Lambda_\varepsilon} \Big(\mathbb{E}\Big|\sup_{t\in[\tau_{\varepsilon,\mathfrak{m}},T]} \frac{1}{\mathbf{w}^{\varepsilon,\mathfrak{m}}_\lambda(t)} |\iota_\varepsilon(\nabla^-_\varepsilon \tilde{h}^{\varepsilon,\mathfrak{m}}_{\text{sym}} \nabla^+_\varepsilon \tilde{h}^{\varepsilon,\mathfrak{m}}_{\text{sym}})(t,\varphi^\lambda_x)|\Big|^p\Big)^{\frac{1}{p}}.$$

Then the bounds (B.25) and (B.27) which we obtained on the two terms in (B.22) and Minkowski inequality yield

$$\Big(\mathbb{E}\Big|\sup_{t\in[\tau_{\varepsilon,\mathfrak{m}},T]} \frac{1}{\mathbf{w}^{\varepsilon,\mathfrak{m}}_\lambda(t)} |\iota_\varepsilon(\nabla^-_\varepsilon \mathcal{M}^\varepsilon_t \nabla^+_\varepsilon \mathcal{M}^\varepsilon_t)(t,\varphi^\lambda_x)|\Big|^p\Big)^{\frac{1}{p}} \\ \leq C_1|\ln\varepsilon| + C_2 + (1+\sqrt{\varepsilon})\Theta D^\varepsilon_\lambda.$$

Combining this bound with the bounds (B.14) and (B.21) on the other three terms in (B.11), we get

$$D^\varepsilon_\lambda \leq C_1|\ln\varepsilon| + C_3\mathfrak{m}^2 + (1+\sqrt{\varepsilon})\Theta D^\varepsilon_\lambda,$$

for some constant $C_3 > 0$, where we used our assumption $\mathfrak{m} \geq 1$. Since $\Theta \in (0,1)$, we can take $\varepsilon$ small enough so that $(1+\sqrt{\varepsilon})\Theta < 1$. Then the preceding inequality yields $D^\varepsilon_\lambda \lesssim |\ln\varepsilon| + \mathfrak{m}^2$, which is the required bound (B.8). □

Now, we are ready to estimate the function $\widetilde{\mathbf{C}}_{\varepsilon,\mathfrak{m}}$, defined in Section 2.2.

**Lemma B.3** *There is $\theta_0 \in (0,1)$ and $\varepsilon_0 \in (0,1)$ such that for any $p \geq 1$, $T > 0$ and $\theta \in [0,\theta_0)$ there is a constant $C > 0$ for which the bounds*

$$\Big(\mathbb{E}\Big|\sup_{t\in[0,\tau_{\varepsilon,\mathfrak{m}}\wedge T]} |\iota_\varepsilon(\widetilde{\mathbf{C}}_{\varepsilon,\mathfrak{m}} - \nu^\varepsilon)(t,\varphi^\lambda_x)|\Big|^p\Big)^{\frac{1}{p}} \leq C\mathfrak{m}\varepsilon^\theta(\lambda\vee\varepsilon)^{2\alpha-1-\underline{\kappa}-\theta}, \tag{B.28}$$

$$\Big(\mathbb{E}\Big|\sup_{t\in[\tau_{\varepsilon,\mathfrak{m}},T]} \frac{1}{\mathbf{w}^{\varepsilon,\mathfrak{m}}_\lambda(t)} |\iota_\varepsilon(\widetilde{\mathbf{C}}_{\varepsilon,\mathfrak{m}} - \nu^\varepsilon)(t,\varphi^\lambda_x)|\Big|^p\Big)^{\frac{1}{p}} \leq C\varepsilon(|\ln\varepsilon| + \mathfrak{m}^2) \tag{B.29}$$

*hold uniformly over $\varphi \in \mathcal{B}^2$, $x \in \Lambda_\varepsilon$, $\lambda \in (0,1]$ and $\varepsilon \in (0,\varepsilon_0)$. The constants $\alpha$ and $\underline{\kappa}$ are the same as in* (2.29)*, the constant $\nu^\varepsilon$ is defined in* (2.8) *and the random weight $\mathbf{w}^{\varepsilon,\mathfrak{m}}_\lambda(t)$ is defined in* (B.9).

*Proof.* We have from (2.13)

$$\iota_\varepsilon(\widetilde{\mathbf{C}}_{\varepsilon,\mathfrak{m}} - \nu^\varepsilon)(t,\varphi^\lambda_x) = -\varepsilon(1+\sqrt{\varepsilon})\iota_\varepsilon(\nabla^-_\varepsilon \tilde{h}^\varepsilon \nabla^+_\varepsilon \tilde{h}^\varepsilon)(t,\varphi^\lambda_x) \\ - 2\varrho_\varepsilon\varepsilon^{\frac{1}{2}}(1+\sqrt{\varepsilon})\iota_\varepsilon(\nabla_\varepsilon \tilde{h}^\varepsilon)(t,\varphi^\lambda_x) - \varepsilon^2\iota_\varepsilon(\nabla^+_\varepsilon \tilde{h}^\varepsilon)(t,\nabla^+_\varepsilon\varphi^\lambda_x) \tag{B.30}$$

for $t \in [0,\tau_{\varepsilon,\mathfrak{m}})$, where we applied summation by parts in the last term.

The definition of the stopping times (2.28) (see the comment after this definitions) yields $|\nabla_\varepsilon \tilde{h}^\varepsilon(t,x)| \leq (\mathfrak{m}+2)\varepsilon^{\alpha-1} \leq 3\mathfrak{m}\varepsilon^{\alpha-1}$ (we used the assumption $\mathfrak{m} \geq 1$) and

$$|\iota_\varepsilon(\nabla_\varepsilon \tilde{h}^\varepsilon)(t,\varphi^\lambda_x)| = |(\iota_\varepsilon \tilde{h}^\varepsilon)(t,\nabla_\varepsilon\varphi^\lambda_x)| \leq 3\mathfrak{m}(\lambda\vee\varepsilon)^{-1}.$$

Then, by interpolation, for any $\gamma \in [0,1]$ we have

$$|\iota_\varepsilon(\nabla_\varepsilon \tilde{h}^\varepsilon)(t, \varphi_x^\lambda)| \leq 3\mathfrak{m}\varepsilon^{(\alpha-1)\gamma}(\lambda \vee \varepsilon)^{\gamma-1} \tag{B.31}$$

and similarly

$$|\iota_\varepsilon(\nabla_\varepsilon^+ \tilde{h}^\varepsilon)(t, \nabla_\varepsilon^+ \varphi_x^\lambda)| \leq 3\mathfrak{m}\varepsilon^{(\alpha-1)\tilde{\gamma}}(\lambda \vee \varepsilon)^{\tilde{\gamma}-2}, \tag{B.32}$$

for any $\tilde{\gamma} \in [0,1]$.

Let us bound the first term on the right-hand side of (B.30). Using the renormalised product (2.17) and the function (2.18) we write

$$\begin{aligned}\iota_\varepsilon(\nabla_\varepsilon^- \tilde{h}^\varepsilon \nabla_\varepsilon^+ \tilde{h}^\varepsilon)(t, \varphi_x^\lambda) &= \iota_\varepsilon(\nabla_\varepsilon^- \tilde{h}^\varepsilon \diamond \nabla_\varepsilon^+ \tilde{h}^\varepsilon)(t, \varphi_x^\lambda) + (\iota_\varepsilon \widetilde{C}_\varepsilon)(t, \varphi_x^\lambda) \\ &= \iota_\varepsilon(\nabla_\varepsilon^- \tilde{h}^\varepsilon \diamond \nabla_\varepsilon^+ \tilde{h}^\varepsilon)(t, \varphi_x^\lambda) + \nu^\varepsilon \int_0^\infty \varepsilon \sum_{x' \in \Lambda_\varepsilon} Q^\varepsilon(t-s, x')\varphi_x^\lambda(x') \\ &\quad - (1+\sqrt{\varepsilon})\varepsilon Q^\varepsilon \star_\varepsilon^+ \iota_\varepsilon(\nabla_\varepsilon^- \tilde{h}^\varepsilon \nabla_\varepsilon^+ \tilde{h}^\varepsilon)(t, \varphi_x^\lambda) - 2\varrho_\varepsilon(1+\sqrt{\varepsilon})\varepsilon^{\frac{1}{2}} Q^\varepsilon \star_\varepsilon^+ \iota_\varepsilon(\nabla_\varepsilon \tilde{h}^\varepsilon)(t, \varphi_x^\lambda) \\ &\quad - \varepsilon^2 Q^\varepsilon \star_\varepsilon^+ \iota_\varepsilon(\nabla_\varepsilon^+ \tilde{h}^\varepsilon)(t, \nabla_\varepsilon^+ \varphi_x^\lambda),\end{aligned} \tag{B.33}$$

where we used summation by parts in the last term. Let us bound the terms on the right-hand side one-by-one. The definitions of the stopping times (2.29) and (2.30) allow to absolutely bound the first term in (B.33) by $\mathfrak{m}\varepsilon^{-\underline{\kappa}}(\lambda \vee \varepsilon)^{2(\alpha-1)}$. Using Lemma 4.7, the second term in (B.33) is absolutely bounded by a constant multiple of $\varepsilon^{\delta-1}(\lambda \vee \varepsilon)^{-1-\delta}$ for any $\delta \in [0,1]$. Combining Lemma 4.5 with (B.31), we absolutely bound the fourth term in (B.33) by a constant times $\mathfrak{m}\varepsilon^{(\alpha-1)\gamma-\frac{1}{2}}(\lambda \vee \varepsilon)^{\gamma-1}$. Similarly, using (B.32), the last term in (B.33) is absolutely bounded by a constant times $\mathfrak{m}\varepsilon^{(\alpha-1)\tilde{\gamma}+1}(\lambda \vee \varepsilon)^{\tilde{\gamma}-2} \leq \mathfrak{m}\varepsilon^{(\alpha-1)\tilde{\gamma}}(\lambda \vee \varepsilon)^{\tilde{\gamma}-1}$. Hence, after taking $\tilde{\gamma} = \gamma$ we obtain the following bound on (B.33):

$$\begin{aligned}|\iota_\varepsilon(\nabla_\varepsilon^- \tilde{h}^\varepsilon \nabla_\varepsilon^+ \tilde{h}^\varepsilon)(t, \varphi_x^\lambda)| &\leq (1+\sqrt{\varepsilon})\varepsilon\Big|Q^\varepsilon \star_\varepsilon^+ \iota_\varepsilon(\nabla_\varepsilon^- \tilde{h}^\varepsilon \nabla_\varepsilon^+ \tilde{h}^\varepsilon)(t, \varphi_x^\lambda)\Big| \\ &\quad + C_1\mathfrak{m}\Big(\varepsilon^{-\underline{\kappa}}(\lambda \vee \varepsilon)^{2(\alpha-1)} + \varepsilon^{\delta-1}(\lambda \vee \varepsilon)^{-1-\delta} + \varepsilon^{(\alpha-1)\gamma-\frac{1}{2}}(\lambda \vee \varepsilon)^{\gamma-1}\Big).\end{aligned}$$

Let us take $\delta \in [0, \frac{1}{2})$ and $\gamma$ such that $(\alpha-1)\gamma - \frac{1}{2} = \delta - 1$. Then the second term in the parentheses dominates the last one. Let us denote furthermore $A_{\varepsilon,\lambda} = \sup_{t \in [0,T]} \sup_{x \in \Lambda_\varepsilon} |\iota_\varepsilon(\nabla_\varepsilon^- \tilde{h}^\varepsilon \nabla_\varepsilon^+ \tilde{h}^\varepsilon)(t, \varphi_x^\lambda)|$. Then Lemma 4.5 yields

$$A_{\varepsilon,\lambda} \leq (1+\sqrt{\varepsilon})\Theta A_{\varepsilon,\lambda} + C_1\mathfrak{m}\Big(\varepsilon^{-\underline{\kappa}}(\lambda \vee \varepsilon)^{2(\alpha-1)} + \varepsilon^{\delta-1}(\lambda \vee \varepsilon)^{-1-\delta}\Big).$$

Taking $\varepsilon$ small, so that $(1+\sqrt{\varepsilon})\Theta < 1$, we can move the first term from the right to the left side and get

$$A_{\varepsilon,\lambda} \leq C_2\mathfrak{m}(\varepsilon^{-\underline{\kappa}}(\lambda \vee \varepsilon)^{2(\alpha-1)} + \varepsilon^{\delta-1}(\lambda \vee \varepsilon)^{-1-\delta}) \tag{B.34}$$

for any $\delta \in [0, \frac{1}{2})$.

Combining the last bound with (B.31) and (B.32), we estimate (B.30) absolutely by a constant multiple of

$$\mathfrak{m}\Big(\varepsilon^{1-\underline{\kappa}}(\lambda \vee \varepsilon)^{2(\alpha-1)} + \varepsilon^{\delta}(\lambda \vee \varepsilon)^{-1-\delta} + \varepsilon^{(\alpha-1)\gamma+\frac{1}{2}}(\lambda \vee \varepsilon)^{\gamma-1} + \varepsilon^{(\alpha-1)\gamma+2}(\lambda \vee \varepsilon)^{\gamma-2}\Big).$$

We can clearly bound it by an expression of the order $\mathfrak{m}\varepsilon^{\theta}(\lambda \vee \varepsilon)^{2\alpha-1-\underline{\kappa}-\theta}$ for any $\theta > 0$ small enough. This is a bound of the desired order.

Now, we will consider the case $t \in [\tau_{\varepsilon,\mathfrak{m}}, T]$. The definition (2.32) yields

$$\iota_\varepsilon(\widetilde{\mathbf{C}}_{\varepsilon,\mathfrak{m}} - \nu^\varepsilon)(t, \varphi_x^\lambda) = -\varepsilon\nu^\varepsilon \iota_\varepsilon(\nabla_\varepsilon^- \tilde{h}_{\text{sym}}^\varepsilon \nabla_\varepsilon^+ \tilde{h}_{\text{sym}}^\varepsilon)(t, \varphi_x^\lambda). \tag{B.35}$$

Applying Lemma B.2, we bound the $p$-th moment of this expression by

$$\Big(\mathbb{E}\Big[|\iota_\varepsilon(\widetilde{\mathbf{C}}_{\varepsilon,\mathfrak{m}} - \nu^\varepsilon)(t, \varphi_x^\lambda)|^p \Big| \mathscr{F}_{\tau_{\varepsilon,\mathfrak{m}}}\Big)^{\frac{1}{p}} \lesssim \varepsilon(|\ln \varepsilon| + \mathfrak{m}^2)\mathbf{w}_\lambda^{\varepsilon,\mathfrak{m}}(t),$$

where the random weight $\mathbf{w}_\lambda^{\varepsilon,\mathfrak{m}}(t)$ is defined in (B.9). Then the desired bound (B.29) follows. □

Similarly to the previous lemma, we can estimate the function $\widetilde{\mathsf{C}}_{\varepsilon,\mathfrak{m}}$, defined in Section 2.2.

**Lemma B.4** *In the setting of Lemma B.3 one has the bound*

$$\Big(\mathbb{E}|(\iota_\varepsilon \widetilde{\mathsf{C}}_{\varepsilon,\mathfrak{m}})(t, \varphi_x^\lambda)|^p\Big)^{\frac{1}{p}} \leq C\mathfrak{m}\varepsilon^\theta(\lambda \vee \varepsilon)^{-\frac{1}{2}-\theta}, \tag{B.36}$$

*for any* $\theta \in [0, \theta_0)$.

*Proof.* In the case $t \in [0, \tau_{\varepsilon,\mathfrak{m}})$ we need to bound the function (6.35). We write

$$\begin{aligned}(\iota_\varepsilon \widetilde{\mathsf{C}}_{\varepsilon,\mathfrak{m}})(t, \varphi_x^\lambda) = -\frac{1}{2}\varepsilon^{\frac{3}{2}}\Big((1 + \sqrt{\varepsilon})\iota_\varepsilon(\nabla_\varepsilon^+ \tilde{h}^\varepsilon)(t, \nabla_\varepsilon^+ \varphi_x^\lambda) - \iota_\varepsilon(\nabla_\varepsilon^- \tilde{h}^\varepsilon \nabla_\varepsilon^+ \tilde{h}^\varepsilon)(t, \varphi_x^\lambda)\\ - 2\varrho_\varepsilon \varepsilon^{-\frac{1}{2}}\iota_\varepsilon(\nabla_\varepsilon \tilde{h}^\varepsilon)(t, \varphi_x^\lambda) - \varrho_\varepsilon^2 \sum_{y \in \Lambda_\varepsilon} \varphi_x^\lambda(y)\Big),\end{aligned}$$

where we used summation by parts in the first term. The terms involved in this expression are the same which we bounded in the preceding lemma. Then, applying (B.31), (B.31) and (B.34) we get

$$\begin{aligned}|(\iota_\varepsilon \widetilde{\mathsf{C}}_{\varepsilon,\mathfrak{m}})(t, \varphi_x^\lambda)| \lesssim \mathfrak{m}\varepsilon^{\frac{3}{2}}\Big(\varepsilon^{(\alpha-1)\tilde{\gamma}}(\lambda \vee \varepsilon)^{\tilde{\gamma}-2} + \varepsilon^{-\underline{\kappa}}(\lambda \vee \varepsilon)^{2(\alpha-1)} + \varepsilon^{\delta-1}(\lambda \vee \varepsilon)^{-1-\delta}\\ + \varepsilon^{(\alpha-1)\gamma-\frac{1}{2}}(\lambda \vee \varepsilon)^{\gamma-1} + \varepsilon^{-1}\Big),\end{aligned}$$

for any $\tilde{\gamma}, \gamma \in [0, 1]$ and $\delta \in [0, \frac{1}{2})$. Taking $\tilde{\gamma} = \gamma = 1$ and $\delta = 0$, and estimating positive powers of $\varepsilon$ by the powers of $\lambda \vee \varepsilon$, we get

$$|(\iota_\varepsilon \widetilde{\mathsf{C}}_{\varepsilon,\mathfrak{m}})(t, \varphi_x^\lambda)| \lesssim \mathfrak{m}\varepsilon^\theta(\lambda \vee \varepsilon)^{-\frac{1}{2}-\theta},$$

for any $\theta \geq 0$ small enough.

In the case $t \in [\tau_{\varepsilon,\mathfrak{m}}, T]$ we use (6.36) to write

$$(\iota_\varepsilon \widetilde{\mathsf{C}}_{\varepsilon,\mathfrak{m}})(t, \varphi_x^\lambda) = -\frac{1}{2}\varepsilon^{\frac{3}{2}}\nu^\varepsilon \iota_\varepsilon(\nabla_\varepsilon^+ \tilde{h}_{\text{sym}}^\varepsilon)(t, \nabla_\varepsilon^+ \varphi_x^\lambda),$$

where we used summation by parts. Applying Lemma B.1, we get

$$\Big(\mathbb{E}\Big|(\lambda \vee \varepsilon)^{2-\alpha}(\iota_\varepsilon \widetilde{\mathsf{C}}_{\varepsilon,\mathfrak{m}})(t, \varphi_x^\lambda)\Big|^p\Big)^{\frac{1}{p}} \leq \Big(\mathbb{E}\Big|(\lambda \vee \varepsilon \vee \sqrt{t - \tau_{\varepsilon,\mathfrak{m}}})^{2-\alpha}(\iota_\varepsilon \widetilde{\mathsf{C}}_{\varepsilon,\mathfrak{m}})(t, \varphi_x^\lambda)\Big|^p\Big)^{\frac{1}{p}} \lesssim \varepsilon^{\frac{3}{2}}\mathfrak{m},$$

which yields

$$\Big(\mathbb{E}|(\iota_\varepsilon \widetilde{\mathsf{C}}_{\varepsilon,\mathfrak{m}})(t, \varphi_x^\lambda)|^p\Big)^{\frac{1}{p}} \lesssim \varepsilon^{\frac{3}{2}}(\lambda \vee \varepsilon)^{\alpha-2}\mathfrak{m} \lesssim \varepsilon^\theta(\lambda \vee \varepsilon)^{\alpha-\frac{1}{2}-\theta}\mathfrak{m},$$

for any $\theta \geq 0$ small enough, which is a stronger bound than we need in (B.36). □

## Appendix C An improved bound on generalised convolutions

The bounds provided in [GMW24, Thm. 4.3 & Cor. 4.5] can be slightly improved. We prefer not to duplicate all complicated definitions but we rather refer to [GMW24].

**Theorem C.1** *In the setting of [GMW24, Thm. 4.3], for any $p \geq 2$, $\theta_1 > 0$ and $\theta_2 \geq 0$ sufficiently small there is a constant $C$ for which the following bound holds*

$$\left(\mathbb{E}\left[\sup_{t\in\mathbb{R}_+} |(\mathcal{I}_\gamma^{\varepsilon,\mathrm{L}} \mathcal{K}_{\mathbb{G}}^{\lambda,\varepsilon})_t|^p\right]\right)^{\frac{1}{p}} \leq C(\lambda \vee \mathfrak{e})^{\nu_\gamma} \sum_{\substack{\mathbf{p}\in\mathscr{P}_m:\, \mathbf{p}^{-1}(1)\cap\Gamma=\varnothing,\\ \mathrm{L}^{-1}(\triangledown)\subset\mathbf{p}^{-1}(1),\, \mathbf{p}^{-1}(\infty)=\varnothing}} \varepsilon^{\alpha_\gamma(\mathbf{p})}\mathfrak{e}^{-\delta_\gamma(\mathbf{p})} \tag{C.1}$$
$$+ C(\lambda \vee \mathfrak{e})^{\nu_\gamma-\theta_2} \sum_{\substack{\mathbf{p}\in\mathscr{P}_m:\, \mathbf{p}^{-1}(1)\cap\Gamma=\varnothing,\\ \mathrm{L}^{-1}(\triangledown)\subset\mathbf{p}^{-1}(1),\, \mathbf{p}^{-1}(\infty)\neq\varnothing}} \varepsilon^{\alpha_\gamma(\mathbf{p})-\theta_1}\mathfrak{e}^{-\delta_\gamma(\mathbf{p})+\theta_2}$$

*uniformly in $\lambda \in (0,1]$, $\mathfrak{e} \in [\varepsilon, 1]$ and $\varepsilon \in (0,1]$.*

The difference between this theorem and [GMW24, Thm. 4.3] is in a slightly different second term in (C.1). Decreasing the power of $\lambda$ by $\theta_2$ allows to increase the power of $\mathfrak{e}$. As we see in our application of this bound in Section 6.4 this increase of the power allows to compensate a divergence as $\varepsilon \to 0$. [GMW24, Cor. 4.4] is a restatement of this result for a shifted function, and we prefer not to state it here. Instead, we will always refer to Theorem C.1.

*Proof of Theorem C.1.* We can slightly change the bound in [GMW24, Lem. 4.7]. More precisely, we bounded the norm of the function in the proof of this lemma by a constant multiple of

$$\Big(\prod_{\nu\in T^\circ} 2^{-\ell_\nu\tilde{\eta}(\nu)}\Big)\Big(\prod_{\nu\in T^\circ_\Gamma} 2^{\ell_\nu|\mathfrak{s}|}\Big)^{\frac12}\Big(\prod_{\nu\in T^\circ_{\mathbf{p}^{-1}(\infty)}\setminus T^\circ_\Gamma} 2^{\ell_\nu|\mathfrak{s}|}\Big)\Big(\prod_{\nu\in T^\circ_{\mathbf{p}^{-1}(\infty)}\cap T^\circ_\Gamma} 2^{\ell_\nu|\mathfrak{s}|}\Big)^{\frac12}\Big(\prod_{\nu\in T^\circ_{\mathbf{p}^{-1}(2)}\setminus T^\circ_\Gamma} 2^{\ell_\nu|\mathfrak{s}|}\Big)^{\frac12}.$$

In the proof of [GMW24, Lem. 4.7] we used $2^{-\ell_\nu} \geq \mathfrak{e}$ to bound the last three factors by $\mathfrak{e}^{-\delta_\gamma(\mathbf{p})}$. This time we do the following trick: we choose any $\nu' \in T^\circ$ and set $\tilde{\tilde{\eta}}(\nu') = \tilde{\eta}(\nu') - \theta_2$ and $\tilde{\tilde{\eta}}(\nu) = \tilde{\eta}(\nu)$ for all $\nu \in T^\circ$ such that $\nu \neq \nu'$. Then the preceding expression equals

$$\Big(\prod_{\nu\in T^\circ} 2^{-\ell_\nu\tilde{\tilde{\eta}}(\nu)}\Big)\Big(\prod_{\nu\in T^\circ_\Gamma} 2^{\ell_\nu|\mathfrak{s}|}\Big)^{\frac12}$$
$$\times\, 2^{-\ell_{\nu'}\theta_2}\Big(\prod_{\nu\in T^\circ_{\mathbf{p}^{-1}(\infty)}\setminus T^\circ_\Gamma} 2^{\ell_\nu|\mathfrak{s}|}\Big)\Big(\prod_{\nu\in T^\circ_{\mathbf{p}^{-1}(\infty)}\cap T^\circ_\Gamma} 2^{\ell_\nu|\mathfrak{s}|}\Big)^{\frac12}\Big(\prod_{\nu\in T^\circ_{\mathbf{p}^{-1}(2)}\setminus T^\circ_\Gamma} 2^{\ell_\nu|\mathfrak{s}|}\Big)^{\frac12}.$$

As soon as $\theta_2 \geq 0$ is chosen sufficiently small and $\delta_\gamma(\mathbf{p}) > 0$, we use $2^{-\ell_\nu} \geq \mathfrak{e}$ to bound it by

$$\Big(\prod_{\nu\in T^\circ} 2^{-\ell_\nu\tilde{\tilde{\eta}}(\nu)}\Big)\Big(\prod_{\nu\in T^\circ_\Gamma} 2^{\ell_\nu|\mathfrak{s}|}\Big)^{\frac12}\mathfrak{e}^{-\delta_\gamma(\mathbf{p})+\theta_2}.$$

After that we can use this bound, instead of [GMW24, Lem. 4.7], in the proof of [GMW24, Thm. 4.3]. The power of $\lambda$ is computed in [GMW24, Lem. 4.8] but via the function $\tilde{\tilde{\eta}}$, which yields the desired decrease of the power by $\theta_2$. □

## References

[ACH24] A. AGGARWAL, I. CORWIN, and M. HEGDE. KPZ fixed point convergence of the ASEP and stochastic six-vertex models (2024). arXiv:2412.18117.

[ACQ11] G. Amir, I. Corwin, and J. Quastel. Probability distribution of the free energy of the continuum directed random polymer in $1+1$ dimensions. *Comm. Pure Appl. Math.* **64**, no. 4, (2011), 466–537. doi:10.1002/cpa.20347.

[AKQ14] T. Alberts, K. Khanin, and J. Quastel. The intermediate disorder regime for directed polymers in dimension $1+1$. *Ann. Probab.* **42**, no. 3, (2014), 1212–1256. doi:10.1214/13-AOP858.

[Bar17] M. T. Barlow. *Random walks and heat kernels on graphs*, vol. 438 of *London Mathematical Society Lecture Note Series*. Cambridge University Press, Cambridge, 2017. doi:10.1017/9781107415690.

[BFS21] C. Bernardin, T. Funaki, and S. Sethuraman. Derivation of coupled KPZ-Burgers equation from multi-species zero-range processes. *Ann. Appl. Probab.* **31**, no. 4, (2021), 1966–2017. doi:10.1214/20-aap1639.

[BG97] L. Bertini and G. Giacomin. Stochastic Burgers and KPZ equations from particle systems. *Comm. Math. Phys.* **183**, no. 3, (1997), 571–607. doi:10.1007/s002200050044.

[CGST20] I. Corwin, P. Ghosal, H. Shen, and L.-C. Tsai. Stochastic PDE limit of the six vertex model. *Comm. Math. Phys.* **375**, no. 3, (2020), 1945–2038. doi:10.1007/s00220-019-03678-z.

[Col51] J. D. Cole. On a quasi-linear parabolic equation occurring in aerodynamics. *Quart. Appl. Math.* **9**, (1951), 225–236. doi:10.1090/qam/42889.

[Cor12] I. Corwin. The Kardar-Parisi-Zhang equation and universality class. *Random Matrices Theory Appl.* **1**, no. 1, (2012), 1130001, 76. doi:10.1142/S2010326311300014.

[Cor18] I. Corwin. Exactly solving the KPZ equation. In *Random growth models*, vol. 75 of *Proc. Sympos. Appl. Math.*, 203–254. Amer. Math. Soc., Providence, RI, 2018.

[CS18] I. Corwin and H. Shen. Open ASEP in the weakly asymmetric regime. *Comm. Pure Appl. Math.* **71**, no. 10, (2018), 2065–2128. doi:10.1002/cpa.21744.

[CST18] I. Corwin, H. Shen, and L.-C. Tsai. ASEP$(q,j)$ converges to the KPZ equation. *Ann. Inst. Henri Poincaré Probab. Stat.* **54**, no. 2, (2018), 995–1012. doi:10.1214/17-AIHP829.

[CT17] I. Corwin and L.-C. Tsai. KPZ equation limit of higher-spin exclusion processes. *Ann. Probab.* **45**, no. 3, (2017), 1771–1798. doi:10.1214/16-AOP1101.

[CT20] I. Corwin and L.-C. Tsai. SPDE limit of weakly inhomogeneous ASEP. *Electron. J. Probab.* **25**, (2020), Paper No. 156, 55. doi:10.1214/20-ejp565.

[DGP17] J. Diehl, M. Gubinelli, and N. Perkowski. The Kardar-Parisi-Zhang equation as scaling limit of weakly asymmetric interacting Brownian motions. *Comm. Math. Phys.* **354**, no. 2, (2017), 549–589. doi:10.1007/s00220-017-2918-6.

[DPD03] G. Da Prato and A. Debussche. Strong solutions to the stochastic quantization equations. *Ann. Probab.* **31**, no. 4, (2003), 1900–1916. doi:10.1214/aop/1068646370.

[DT16] A. Dembo and L.-C. Tsai. Weakly asymmetric non-simple exclusion process and the Kardar-Parisi-Zhang equation. *Comm. Math. Phys.* **341**, no. 1, (2016), 219–261. doi:10.1007/s00220-015-2527-1.

[EH19] D. Erhard and M. Hairer. Discretisation of regularity structures. *Ann. Inst. Henri Poincaré Probab. Stat.* **55**, no. 4, (2019), 2209–2248. doi:10.1214/18-AIHP947.

[EL15] A. M. Etheridge and C. Labbé. Scaling limits of weakly asymmetric interfaces. *Comm. Math. Phys.* **336**, no. 1, (2015), 287–336. doi:10.1007/s00220-014-2243-2.

[FH20] P. K. Friz and M. Hairer. *A course on rough paths. With an introduction to regularity structures*. Universitext. Springer, Cham, second ed., [2020] ©2020. doi:10.1007/978-3-030-41556-3.

[Gär88] J. Gärtner. Convergence towards Burgers' equation and propagation of chaos for weakly asymmetric exclusion processes. *Stochastic Process. Appl.* **27**, no. 2, (1988), 233–260. doi:10.1016/0304-4149(87)90040-8.

[GIP15] M. Gubinelli, P. Imkeller, and N. Perkowski. Paracontrolled distributions and singular PDEs. *Forum Math. Pi* **3**, (2015), e6, 75. doi:10.1017/fmp.2015.2.

[GJ13] M. Gubinelli and M. Jara. Regularization by noise and stochastic Burgers equations. *Stoch. Partial Differ. Equ. Anal. Comput.* **1**, no. 2, (2013), 325–350. doi:10.1007/s40072-013-0011-5.

[GJ14] P. Gonçalves and M. Jara. Nonlinear fluctuations of weakly asymmetric interacting particle systems. *Arch. Ration. Mech. Anal.* **212**, no. 2, (2014), 597–644. doi:10.1007/s00205-013-0693-x.

[GJ17] P. Gonçalves and M. Jara. Stochastic Burgers equation from long range exclusion interactions. *Stochastic Process. Appl.* **127**, no. 12, (2017), 4029–4052. doi:10.1016/j.spa.2017.03.022.

[GJS15] P. Gonçalves, M. Jara, and S. Sethuraman. A stochastic Burgers equation from a class of microscopic interactions. *Ann. Probab.* **43**, no. 1, (2015), 286–338. doi:10.1214/13-AOP878.

[GMW24] P. Grazieschi, K. Matetski, and H. Weber. Martingale-driven integrals and singular SPDEs. *Probab. Theory Relat. Fields* **190**, no. 3-4, (2024), 1063–1120. doi:10.1007/s00440-024-01311-2.

[GMW25] P. Grazieschi, K. Matetski, and H. Weber. The dynamical Ising-Kac model in 3D converges to $\Phi^4_3$. *Probability Theory and Related Fields* **191**, (2025), 671–778. doi:10.1007/s00440-024-01316-x.

[GP17] M. Gubinelli and N. Perkowski. KPZ reloaded. *Comm. Math. Phys.* **349**, no. 1, (2017), 165–269. doi:10.1007/s00220-016-2788-3.

[GP18] M. Gubinelli and N. Perkowski. Energy solutions of KPZ are unique. *J. Amer. Math. Soc.* **31**, no. 2, (2018), 427–471. doi:10.1090/jams/889.

[GP20] M. Gubinelli and N. Perkowski. The infinitesimal generator of the stochastic Burgers equation. *Probab. Theory Related Fields* **178**, no. 3-4, (2020), 1067–1124. doi:10.1007/s00440-020-00996-5.

[GPS20] P. Gonçalves, N. Perkowski, and M. Simon. Derivation of the stochastic Burgers equation with Dirichlet boundary conditions from the WASEP. *Ann. H. Lebesgue* **3**, (2020), 87–167. doi:10.5802/ahl.28.

[Hai13] M. Hairer. Solving the KPZ equation. *Ann. of Math. (2)* **178**, no. 2, (2013), 559–664. doi:10.4007/annals.2013.178.2.4.

[Hai14] M. Hairer. A theory of regularity structures. *Invent. Math.* **198**, no. 2, (2014), 269–504. doi:10.1007/s00222-014-0505-4.

[HM18] M. Hairer and K. Matetski. Discretisations of rough stochastic PDEs. *Ann. Probab.* **46**, no. 3, (2018), 1651–1709. doi:10.1214/17-AOP1212.

[Hop50] E. Hopf. The partial differential equation $u_t + uu_x = \mu u_{xx}$. *Comm. Pure Appl. Math.* **3**, (1950), 201–230. doi:10.1002/cpa.3160030302.

[HQ18] M. Hairer and J. Quastel. A class of growth models rescaling to KPZ. *Forum Math. Pi* **6**, (2018), e3, 112. doi:10.1017/fmp.2018.2.

[HS17] M. Hairer and H. Shen. A central limit theorem for the KPZ equation. *Ann. Probab.* **45**, no. 6B, (2017), 4167–4221. doi:10.1214/16-AOP1162.

[JMF20] M. Jara and G. R. Moreno Flores. Stationary directed polymers and energy solutions of the Burgers equation. *Stochastic Process. Appl.* **130**, no. 10, (2020), 5973–5998. doi:10.1016/j.spa.2020.04.012.

[KM17] A. Kupiainen and M. Marcozzi. Renormalization of generalized KPZ equation. *J. Stat. Phys.* **166**, no. 3-4, (2017), 876–902. doi:10.1007/s10955-016-1636-3.

[KPZ86] M. Kardar, G. Parisi, and Y.-C. Zhang. Dynamical scaling of growing interfaces. *Phys. Rev. Lett.* **56**, no. 9, (1986), 889–892. doi:10.1103/PhysRevLett.56.889.

[Kup16] A. Kupiainen. Renormalization group and stochastic PDEs. *Ann. Henri Poincaré* **17**, no. 3, (2016), 497–535. doi:10.1007/s00023-015-0408-y.

[Lab17] C. Labbé. Weakly asymmetric bridges and the KPZ equation. *Comm. Math. Phys.* **353**, no. 3, (2017), 1261–1298. doi:10.1007/s00220-017-2875-0.

[Lig05] T. M. Liggett. *Interacting particle systems*. Classics in Mathematics. Springer-Verlag, Berlin, 2005. Reprint of the 1985 original.

[Lin20] Y. Lin. KPZ equation limit of stochastic higher spin six vertex model. *Math. Phys. Anal. Geom.* **23**, no. 1, (2020), Paper No. 1, 118. doi:10.1007/s11040-019-9325-5.

[MQR21] K. Matetski, J. Quastel, and D. Remenik. The KPZ fixed point. *Acta Math.* **227**, no. 1, (2021), 115–203. doi:10.4310/acta.2021.v227.n1.a3.

[Nak19] M. Nakashima. Free energy of directed polymers in random environment in $1+1$-dimension at high temperature. *Electron. J. Probab.* **24**, (2019), Paper No. 50, 43. doi:10.1214/19-EJP292.

[OQR16] J. Ortmann, J. Quastel, and D. Remenik. Exact formulas for random growth with half-flat initial data. *Ann. Appl. Probab.* **26**, no. 1, (2016), 507–548. doi:10.1214/15-AAP1099.

[Par19] S. Parekh. The KPZ limit of ASEP with boundary. *Comm. Math. Phys.* **365**, no. 2, (2019), 569–649. doi:10.1007/s00220-018-3258-x.

[Par23] S. Parekh. Convergence of ASEP to KPZ with basic coupling of the dynamics. *Stochastic Process. Appl.* **160**, (2023), 351–370. doi:10.1016/j.spa.2023.03.007.

[QS23] J. Quastel and S. Sarkar. Convergence of exclusion processes and the KPZ equation to the KPZ fixed point. *J. Amer. Math. Soc.* **36**, no. 1, (2023), 251–289. doi:10.1090/jams/999.

[Qua12] J. Quastel. Introduction to KPZ. In *Current developments in mathematics, 2011*, 125–194. Int. Press, Somerville, MA, 2012.

[Spo14] H. Spohn. KPZ scaling theory and the semidiscrete directed polymer model. In *Random matrix theory, interacting particle systems, and integrable systems*, vol. 65 of *Math. Sci. Res. Inst. Publ.*, 483–493. Cambridge Univ. Press, New York, 2014.

[SSSX21] H. Shen, J. Song, R. Sun, and L. Xu. Scaling limit of a directed polymer among a Poisson field of independent walks. *J. Funct. Anal.* **281**, no. 5, (2021), Paper No. 109066, 55. doi:10.1016/j.jfa.2021.109066.

[Wal86] J. B. Walsh. An introduction to stochastic partial differential equations. In *École d'été de probabilités de Saint-Flour, XIV—1984*, vol. 1180 of *Lecture Notes in Math.*, 265–439. Springer, Berlin, 1986. doi:10.1007/BFb0074920.

[Yan22] K. Yang. KPZ equation from non-simple variations on open ASEP. *Probab. Theory Related Fields* **183**, no. 1-2, (2022), 415–545. doi:10.1007/s00440-022-01133-0.

[Yan23a] K. Yang. Fluctuations for some nonstationary interacting particle systems via Boltzmann-Gibbs principle. *Forum Math. Sigma* **11**, (2023), Paper No. e32, 86. doi:10.1017/fms.2023.27.

[Yan23b] K. Yang. Kardar-Parisi-Zhang equation from long-range exclusion processes. *Comm. Math. Phys.* **400**, no. 3, (2023), 1535–1663. doi:10.1007/s00220-022-04628-y.

[Yan24] K. Yang. Non-stationary KPZ equation from ASEP with slow bonds. *Ann. Inst. Henri Poincaré Probab. Stat.* **60**, no. 2, (2024), 1246–1294. doi:10.1214/23-aihp1364.